\documentclass[envcountsame,envcountchap]{svmono}  

\usepackage{amsmath}
\usepackage{amsfonts}
\usepackage{latexsym}
\usepackage{pstricks}
\usepackage{graphics}
\usepackage{amscd}
\usepackage{a4}
\usepackage[mathcal]{euscript}
\usepackage{mathrsfs}
\usepackage{theorem}
\usepackage{amssymb}
\usepackage{times}
\usepackage{amsrefs} 
\usepackage{makeidx}         

\input{boxedeps}  
\SetEPSFDirectory{figures/}  
\HideDisplacementBoxes

\SetRokickiEPSFSpecial  

\DeclareMathAlphabet{\ams}{U}{msb}{m}{n}
\DeclareMathAlphabet{\goth}{U}{euf}{m}{n}
\DeclareMathAlphabet{\mathgothic}{U}{euf}{m}{n}
\DeclareMathAlphabet{\ams}{U}{msb}{m}{n}
 at 12truept   
\input cyracc.def

\newfont{\Ma}{msam10}

\def\id{\text{id}}

\def\aut{\text{Aut}}
\def\sym{\text{Sym}}

\def\gal{\text{Gal}}

\def\ov{\overline}

\def\wtl{\widetilde}
\def\qedsymbol{$\Box$}
\def\LL{\mathcal L}

\def\SSS{\mathcal S}
\def\SS{\mathgothic{S}}

\def\ve{\varepsilon}
\def\aa{\alpha}
\def\ww{\omega}

\def\ss{\sigma}

\def\ve{\varepsilon}
\def\Z{\ams{Z}}\def\E{\ams{E}}
\def\R{\ams{R}}

\def\F{\ams{F}}
\def\P{\ams{P}}

\def\v{\mbox{\boldmath $v$}}
\def\u{\mbox{\boldmath $u$}}

\def\quo{/\kern -.45em\sim}
\def\hcl{[\kern-1.5pt[}
\def\hcr{]\kern-1.5pt]}
\newcmykcolor{cyan}{1 0 0 0}
\newcmykcolor{myyellow}{0 0 1 0}
\renewcommand{\theparagraph}{}

\begin{document}

%
%
%
%
\author{Brent Everitt}
\title{The Combinatorial Topology of Groups\\
}
\maketitle  
%
%
%

\frontmatter

\tableofcontents

\mainmatter

\chapstarthook{By graph and map of graphs, I mean something purely combinatorial or
algebraic. Pictures can be drawn, but one has to understand that maps are rigid and not just
continuous, maps do not $\ldots$ wrap edges around several edges.--John Stallings
\cite{Stallings83}.}
\chapter{Combinatorial Complexes}\label{chapter1}

\section{$1$-Complexes (ie: Graphs)}\label{chapter1:graphs}

\subsection{The category of graphs}\label{chapter1:graphs:category}

\begin{definition}[$\mathbf{1}$-complex: first go]\label{chapter1:graphs:definition100}
A \emph{$1$-complex\/} or \emph{graph\/} is a non-empty set $X$ together with an involutary map 
$i:X\rightarrow X$ (ie: $i^2=\id_X$) and an idempotent map $s:X\rightarrow X^0$ (ie: $s^2=s$)
where $X^0$ is the set of fixed points of $i$.
\end{definition}

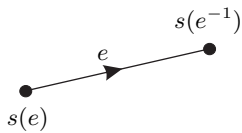
\begin{figure}
\begin{pspicture}(0,0)(13,1.75)
\rput(.5,0){
\rput(6,1){\BoxedEPSF{fig5.eps scaled 1000}}
\rput(4.8,.35){$s(e)$}
\rput(7.2,1.65){$s(e^{-1})$}
\rput(5.8,1.2){$e$}
}
\end{pspicture}
\caption{an edge of a graph with its start and 
terminal vertices.}\label{chapter1:graphs:figure100}
\end{figure}

Thus a graph has $0$-cells or {\em vertices\/} $X^0$ 
and $1$-cells or {\em edges\/} $X^1=X\setminus X^0$.
From now on we will write $x^{-1}$ for $i(x)$, and say
that the edge $e\in X^1$ has {\em start vertex\/} $s(e)$ and 
{\em terminal vertex\/}
$s(e^{-1})$. One thinks of the {\em inverse edge\/} $e^{-1}$ as just
$e$, but traversed in the reverse direction (or with the reverse
orientation). 
The edge $e$ is \emph{incident\/} with the vertex $v$ if $e\in s^{-1}(v)$.
We draw pictures like Figure \ref{chapter1:graphs:figure100},
although
they are purely for illustrative purposes. If the vertex set has 
cardinality that
of the power set of the continuum for instance, 
then there are not enough points on a piece of paper
for a picture to fit! Definition \ref{chapter1:graphs:definition100} is quite terse,
and it is sometimes useful to spell it out a little more:

\begin{definition}[$\mathbf{1}$-complex: second go]
\label{chapter1:graphs:definition200}
A \emph{$1$-complex\/} or \emph{graph\/} consists of two disjoint non-empty sets $X^0$ 
and $X^1$, together with two incidence maps and an inverse map,
$$
s,t:X^1\rightarrow X^0\text{ and }\,\,^{-1}:X^1\rightarrow X^1,
$$
such that,
(i). $e^{-1}\not= e=(e^{-1})^{-1}$ for all 
$e\in X^1$, and 
(ii). $t(e)=s(e^{-1})$ for all $e\in X$.
\end{definition}

More terminology: an {\em arc\/} is an edge/inverse edge pair,
and an {\em orientation\/} for $X$ is a set $\mathcal{O}$ consisting
of all the vertices and
exactly one edge from each arc. 
Write $\ov{e}$ for the arc
containing the edge $e$, so that $\ov{e^{-1}}=\ov{e}$.
The graph $X$ is {\em finite\/} when $X^0$ is finite and 
{\em locally finite\/}
when the set $s^{-1}(v)$ is finite for every $v\in X^0$. Thus a 
finite graph may have infinitely
many edges, a situation that possibly differs from that in combinatorics. 
The cardinality
of the set $s^{-1}(v)$ is the {\em valency\/} 
of the vertex $v$.
A \emph{pointed graph\/} is a pair
$X_v:=(X,v)$ for $v\in X$ a vertex.

\begin{exercise}\label{chapter1:graphs:exercise100}
Here is another definition of graph
more in the spirit of \S\ref{chapter1:2complexes}. 
A $0$-complex is a non-empty set $X$ and a map of $0$-complexes is a map 
$f:X\rightarrow Y$ of sets.
The $0$-sphere $S^0$ is the $0$-complex with \emph{two\/} elements. 

A graph $X$
is a graded set $X=X^0,X^1$ 
with $X^1\not=\varnothing$, such that
\begin{description}
\item[(C1).] $X^0$ is a $0$-complex;
\item[(C2).] there is an involutory map $^{-1}:X\rightarrow X$ with 
fixed point set $X^0$;
\item[(C3).] each $e\in X^1$ has \emph{boundary\/} $\partial e=(X^e,\aa_e)$
with $X^e$ the $0$-sphere $S^{0}$ and 
$\aa_e:X^e\rightarrow X^{(0)}$ a map of $0$-complexes.
\end{description}
Here is the exercise: show that all three definitions 
of graph are equivalent.
\end{exercise}

The \emph{trivial graph\/} has a single vertex and no edges.
Figure \ref{chapter1:graphs:figure300} shows some more examples of graphs,
including some with countably many edges, which will tend to
be more interesting than finite graphs.

\begin{figure}\label{chapter1:graphs:figure300}
\begin{pspicture}(0,0)(13,3)
\rput(1,1.5){\BoxedEPSF{chapter3.fig2000.eps scaled 800}}
\rput(4.5,1.5){\BoxedEPSF{chapter1.fig300.eps scaled 450}}
\rput(10,1.5){\BoxedEPSF{chapter1.fig400.eps scaled 450}}
\rput(10.25,1.9){$\ldots$}
\end{pspicture}
\caption{examples of graphs.}
\label{chapter1:graphs:figure300}
\end{figure}
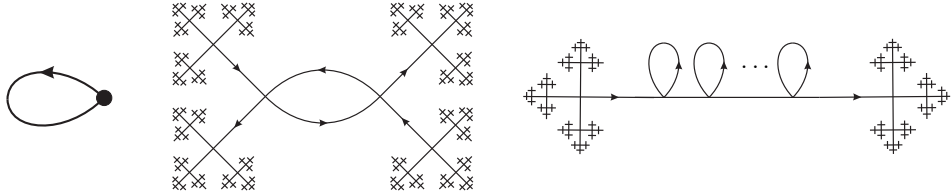

\begin{figure}
\begin{pspicture}(0,0)(13,2)
\rput(-3.75,0.3){
\rput(0.7,-0.55){
\rput(5,2){$X$}\rput(6.55,2){$Y$}
\rput(5,0.45){$X$}\rput(6.55,0.45){$Y$}
\psline[linewidth=.2mm]{->}(5.3,2)(6.3,2)
\psline[linewidth=.2mm]{->}(5.3,0.45)(6.3,0.45)
\psline[linewidth=.2mm]{->}(5,1.7)(5,.7)
\psline[linewidth=.2mm]{->}(6.55,1.7)(6.55,.7)
\rput(4.7,1.2){$\ss_X^{}$}\rput(6.85,1.2){$\ss_Y^{}$}
\rput(5.8,2.2){$f$}\rput(5.8,.65){$f$}
}}
\rput(6,0){
\rput(2,1){\BoxedEPSF{fig5.eps scaled 1000}}
\rput(1,.4){$s_X^{}(e)$}
\rput(3.1,1.65){$t_X^{}(e)$}
\rput(1.8,1.2){$e$}\rput(5.35,1){$u$}
\rput(5,1){\BoxedEPSF{fig79.eps scaled 1000}}
\psline[linewidth=.3mm,linestyle=dotted]{->}(2.4,1)(4.8,1)
\pscurve[linewidth=.3mm,linestyle=dotted]{->}(3.4,1.2)(3.6,1.2)(4,1.4)(4.5,1.5)(4.8,1.2)
\pscurve[linewidth=.3mm,linestyle=dotted]{->}(1.1,.7)(1.4,.7)(2.5,.4)(3.5,.4)(4.8,.8)
}
\end{pspicture}
\caption{graph map $f:X\rightarrow Y$}\label{chapter1:graphs:figure200}
\end{figure}
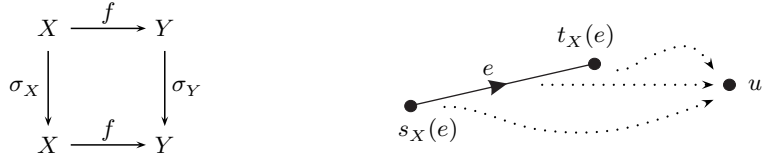

A {\em map\/} of graphs is a set map $f :X\rightarrow Y$
with $f(X^0)\subseteq Y^0$, such that the diagram on the left of
Figure \ref{chapter1:graphs:figure200} 
commutes, where $\ss_X^{}$ is one of the maps $s_X^{}$ or $^{-1}$ for $X$, 
and $\ss_Y^{}$ similarly, ie: $f s_X^{}(x)=s_Y^{} f(x)$ and 
$f(x^{-1})=f(x)^{-1}$.
Notice that a map can send edges to vertices, and so we
call $f$ \emph{dimension preserving\/} if 
we also have $f(X^1)\subseteq Y^1$.
A map $f:X_v\rightarrow Y_u$ of pointed graphs is a graph map
$f:X\rightarrow Y$ with $f(v)=u$.

The commuting of $f$ with $s$ and $^{-1}$ is a combinatorial version
of continuity: an edge incident with a vertex is either mapped to an edge incident
with the image of the vertex, or to the image vertex itself.
In the second case, 
if $e$ is an edge mapped by $f$ to a vertex $u$
as on the right in Figure \ref{chapter1:graphs:figure200},
then the commuting condition becomes $f s_X=f$, 
and in particular the start vertex $s_X^{}(e)$ must also be mapped to $u$
(the condition also ensures that $t_X^{}(e)$ is mapped to $u$, and not 
left ``hanging'').

\begin{exercise}\label{chapter1:graphs:exercise150}
Show that graphs and their mappings form a category.
\end{exercise}

For a fixed vertex $v\in X$, and $s_X^{-1}(v)\in X^1$ the edges starting at $v$,
the map $f:X\rightarrow Y$ induces
a map $s_X^{-1}(v)\rightarrow s_Y^{-1}(u)$, where $u=f(v)$.
We call this induced map the \emph{local continuity of $f$ at the vertex $u$\/}.

A graph map $f:X\rightarrow Y$ \emph{preserves orientation\/} whenever there are 
orientations $\mathcal{O}_X$ for $X$ and $\mathcal{O}_Y$ for $Y$ with 
$f(\mathcal{O}_X)\subset\mathcal{O}_Y$. We leave it as an exercise to show that
it is always possible to choose orientations for $X$ and $Y$ making a map
$f:X\rightarrow Y$ orientation preserving (although a map may not be orientation
preserving with respect to fixed orientations).

A map $f:X\rightarrow Y$ is an \emph{isomorphism\/} if it 
is dimension preserving and a bijection on the vertex and edge sets. 

\begin{exercise}\label{chapter1:graphs:exercise200}
Show that if $f:X\rightarrow Y$ is an isomorphism then the inverse map
$f^{-1}:Y\rightarrow X$ is 
also a graph isomorphism. Show that the set $\aut(X)$ of graph isomorphisms
$X\rightarrow X$ forms a group under composition.
\end{exercise}

A group $G$ acts on a graph $X$ if there is a 
homomorphism
$G\stackrel{\varphi}{\rightarrow}\aut(X)$. We abbreviate 
$\varphi(g)(x)$ to $g(x)$.
An action \emph{preserves orientation\/} if there is an orientation $\mathcal{O}$ for $X$ with 
$g(\mathcal{O})=\mathcal{O}$
for all $g\in G$. 

\begin{exercise}\label{chapter1:graphs:exercise300}
Let $X$ be a set and $^{-1}:X\rightarrow X$ a bijective map without fixed points. 
Let $G\rightarrow\text{Sym}(X)$ be a group action (here $\text{Sym}(X)$ is the
symmetric group on $X$) that commutes with $^{-1}$, ie: $g(x^{-1})=g(x)^{-1}$ for
all $x\in X$. An \emph{inversion\/} is a $g\in G$ such that $g(x)=x^{-1}$ for some
$x$, and $G$ is said to act without inversions if no $g\in G$ is an inversion
(equivalently, no $G$-orbit contains both some $x$ and its inverse $x^{-1}$).
Show that there exists an $\mathcal{O}\subset X$ with $X=\mathcal{O}\cup\mathcal{O}^{-1}$
a disjoint union and $g(\mathcal{O})=\mathcal{O}$ 
if and only if $G$ acts without inversions on $X$.
\end{exercise}

A group $G$ acts {\em freely\/} if and only if the action is free on the vertices, ie: 
if $g\in G$ and $v$ a vertex with $g(v)=v$ implies $g$
is the identity element. 

\begin{exercise}\label{chapter1:graphs:exercise400}
If $G$ acts freely and orientation preservingly on a graph, then
show that the action is free on the edges too.
\end{exercise}

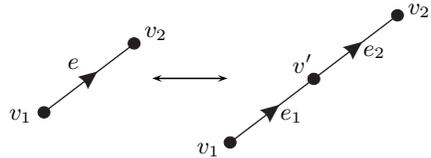
\begin{figure}
\begin{pspicture}(0,0)(13,2)
\rput(4,0){
\rput(-5,0){
\rput(6,1){\BoxedEPSF{fig10.eps scaled 1000}}
\rput(5.1,.5){$v_1$}
\rput(6.9,1.6){$v_2$}
\rput(5.8,1.2){$e$}
}
\psline[linewidth=.2mm]{<->}(1.85,1)(2.85,1)
\rput(-2,0){
\rput(6,1){\BoxedEPSF{fig78.eps scaled 1000}}
\rput(4.6,.05){$v_1$}
\rput(7.4,1.9){$v_2$}
\rput(5.85,1.2){$v'$}
\rput(5.7,.5){$e_1$}
\rput(6.8,1.35){$e_2$}
}}
\end{pspicture}
\caption{Subdividing an edge}
\label{chapter1:graphs:figure400}
\end{figure}

Graph isomorphisms are pretty rigid, and it is useful to have 
a relation with a bit more ``slack''.
Thus, a \emph{subdivision\/} of an edge replaces it by
two new edges and a new vertex as in Figure \ref{chapter1:graphs:figure400},
or is the reverse of this process.
Write $X\leftrightarrow X'$ when two graphs differ by the subdivision of a single edge.
Two graphs $X$ and $Y$ are then \emph{homeomorphic\/}, written $X\approx Y$, when there
is a finite sequence $X=X_0\leftrightarrow X_1\leftrightarrow
\cdots\leftrightarrow X_k=Y$ of subdivisions connecting them. It is easy to
see that homeomorphism is an equivalence relation for graphs.
A \emph{topological invariant\/} is a property of graphs that is invariant under homeomorphism
in the sense that if $X\approx Y$, then
$X$ has the property if and only if $Y$ does.

\subsection{Quotients and subgraphs}\label{chapter1:graphs:quotients}

The most useful construction in the category of graphs is the quotient:

\begin{definition}[quotient relation and quotient graph]\label{chapter1:graphs:definition300}
If $X$ is a graph, 
then a {\em quotient relation\/} is an equivalence relation $\sim$ on $X$
such that
$$
\text{(i)}.\, x\sim y \Rightarrow s(x)\sim s(y)\text{ and }
x^{-1}\sim y^{-1}\hspace{2em}
\text{(ii)}.\, x\sim x^{-1}\Rightarrow[x]\cap X^0\not=\varnothing,
$$
where $[x]$ is the equivalence class of $x$. 
If $\sim$ is a quotient relation on a graph $X$ then 
define $s$ and $^{-1}$ on the equivalence classes $X\quo$ by
$$
\text{(i)}.\,s[x]=[s(x)],
\hspace{2em}
\text{(ii)}.\,[x]^{-1}=[x^{-1}].
$$
\end{definition}

Notice that edges can be equivalent to vertices, but if an edge is equivalent to its
inverse then it must also be equivalent to a vertex. This ensures that
in the quotient we have $[e]\not=[e]^{-1}$.

\begin{proposition}\label{chapter1:graphs:result100}
If $\sim$ is a quotient relation then $X\quo$ with the maps
$s$ and $^{-1}$ defined above is a graph,
and the quotient map $q:X\rightarrow X\quo$ given by $q(x)=[x]$ is a map of graphs.
\end{proposition}

The proof is a straight forward exercise. In particular, the fixed points 
in $X\quo$ of the new inverse map $^{-1}$
are precisely those equivalence classes $[x]$ where $x\sim v$ for
some $v\in X^0$. Thus the quotient has vertices the $[v]$ for $v\in X^0$ (and these classes
may include some of the edges of the old graph $X$) and edges those $[e]$ with 
$[e]\cap X^0=\varnothing$.


The two main examples of graph quotients arise by factoring out the action of a group,
or by squashing a subgraph down to a vertex.
For the first we have the following,

\begin{proposition}\label{chapter1:graphs:result200}
Let $\sim$ be the equivalence
relation on $X$ given by the orbits of the action of a group $G$.
Then $\sim$ is a quotient relation if and only if the group action is orientation 
preserving, and we write $X/G:=X\quo$ for the quotient.
\end{proposition}

Again, the proof is left as an exercise (see Exercise \ref{chapter1:graphs:exercise300}).
A graph $X$ is a \emph{subgraph\/} of $Y$ if there is a 
mapping $X\hookrightarrow Y$ that is an isomorphism
onto its image. Equivaltently, it is a subset $X\subset Y$,
such that the maps $s$ and $^{-1}$ give a graph when restricted to 
$X$. 

Let $X\subset Y$ be a subgraph
and define a relation $\sim$ on $Y$ by $x\sim y$ if and only if $x=y$ or both 
$x$ and $y$ lie in $X$. Then this is a quotient relation
and we write $Y/X$ for $Y\quo$, the \emph{quotient of $Y$ by 
the subgraph $X$}. It is what results by squashing
$X$ to a vertex. 

Extending this a little, if $X_\aa, (\aa\in A)$ is a family of disjoint subgraphs in $Y$
then define $\sim$ by $x\sim y$ iff $x=y$ or $x$ and $y$ lie in the \emph{same\/} $X_\aa$,
and write $Y/X_\aa\,(\aa\in A)$, or just $Y/X_\aa$, 
for the corresponding quotient. Note the difference between this,
where each $X_\aa$ has been squashed to a distinct vertex $v_\aa$, and $Y/(\bigcup X_\aa)$,
where the whole union is squashed to just the one vertex.


\begin{exercise}\label{chapter1:2complexes:exercise400}
Recall that an equivalence relation on a set $X$ is 
a subset $S\subset X\times X$ such that, (i). $S$ contains the diagonal, $(x,x)\in S$ for
all $x\in X$, (ii). $(x,y)\in S\Rightarrow (y,x)\in S$, and 
(iii). $(x,y),(y,z)\in S\Rightarrow (x,z)\in S$. Show that if $S_\aa (\aa\in A)$ are
equivalence relations on $X$ then so is $\bigcap S_\aa$, and hence if $Y$ is 
any subset of $X$ we may define the \emph{equivalence relation generated by\/} $Y$
to be the intersection of all equivalence relations $S$ with $Y\subset S$.
\end{exercise}

\begin{exercise}\label{chapter1:graphs:exercise500}
Let $X_1,X_2$ and $Y$ be graphs and $f_i:Y\rightarrow X_i\,(i=1,2)$ dimension
preserving maps of graphs. 
Let $\sim$ on the disjoint union $X_1\bigcup X_2$ be the equivalence
relation generated by $x\sim y$ if and only if there is a $z\in Y$ with
$x=f_1(z)$ and $y=f_2(z)$. 
Show that $\sim$ is a quotient relation if there are orientations
$\mathcal{O}$ for $Y$, and $\mathcal{O}_i$ 
for $X_i\,(i=1,2)$ with $f_i(\mathcal{O})\subseteq \mathcal{O}_i$.
\end{exercise}

\subsection{Balls, spheres, paths and homotopies}\label{chapter1:graphs:paths}

\begin{figure}
\begin{pspicture}(0,0)(13,1)
\rput(4,.5){\BoxedEPSF{fig10.eps scaled 1000}}
\rput(8,.5){\BoxedEPSF{chapter3.fig2000.eps scaled 800}}
\rput(3.2,.5){$I^1$}\rput(9.2,.5){$S^1$}
\end{pspicture}
\caption{the $I^1$ and $S^1$-graphs.}\label{chapter1:graphs:figure500}
\end{figure}
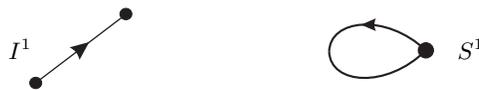

The $I^1$-graph and $S^1$-graph are shown in Figure \ref{chapter1:graphs:figure500}.
A graph $X$ is a \emph{$1$-ball\/} if it 
is homeomorphic to $I^1$ (or is trivial), and a \emph{$1$-sphere\/} if it is homeomorphic
to $S^1$ (or is trivial). 

\begin{figure}
\begin{pspicture}(0,0)(13,4)
\rput(3,2){\BoxedEPSF{chapter1.fig6000.eps scaled 1000}}
\rput(0.6,1.4){$v_0$}\rput(1.9,1.8){$v_1$}\rput(5.4,1.4){$v_n$}
\rput(1.3,2.1){$e_0$}\rput(2.5,2.25){$e_1$}\rput(5,2.3){$e_{n-1}$}
\rput(9,2){\BoxedEPSF{chapter1.fig6100.eps scaled 1000}}
\rput(9,3.6){$v_0$}\rput(10.2,3){$v_1$}\rput(10.6,2.1){$v_2$}
\rput(7.6,2.1){$v_{n-2}$}\rput(7.95,3){$v_{n-1}$}
\rput(9.8,3.9){$e_0$}\rput(11,2.8){$e_1$}
\rput(8.1,3.9){$e_{n-1}$}\rput(6.8,2.8){$e_{n-2}$}
\end{pspicture}
\caption{a $1$-ball and a $1$-sphere.}\label{chapter1:graphs:figure550}
\end{figure}
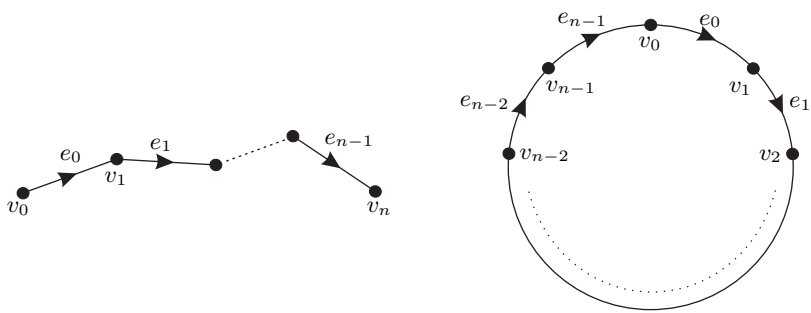

It is easy to see that a $1$-ball and $1$-sphere have the form shown in Figure 
\ref{chapter1:graphs:figure550}, and that 
the vertices of a $1$-sphere can be labelled $v_0,\ldots,v_{n-1}$
and the edges $e_1^\pm\ldots,e_{n-1}^\pm$ with $s(e_i)=v_{i}$ and $t(e_i)=v_{i+1}\,(i<n-1)$
and $t(e_{n-1})=v_0$. Similarly for 
a $1$-ball, so that it has \emph{end vertices\/} $v_0,v_n$ in an obvious sense.
The following is easily proved by induction:

\begin{lemma}\label{complexes:graphs:result300}
A graph $X$ is a $1$-sphere if and only if either $X$ is trivial or $X=S^1$, or
there are non-trivial $1$-balls $B_i\,(i=1,2)$ with end vertices $v_{i1},v_{i2}$,
such that $X=B_1\bigcup B_2\quo$, 
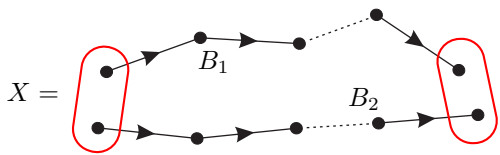
$$
\begin{pspicture}(0,0)(12,2)
\rput(6,1){\BoxedEPSF{chapter1.fig200.eps scaled 1000}}
\rput(.1,-.1){\pspolygon[linearc=.3,linecolor=red,
showpoints=false](3,.2)(3.6,0)(3.8,1.4)(3.2,1.6)}
\rput(5.45,-1.1){\rput{20}{\pspolygon[linearc=.3,linecolor=red,
showpoints=false](3,.2)(3.6,0)(3.8,1.4)(3.2,1.6)}}
\rput(5,1.2){$B_1$}\rput(7,.7){$B_2$}
\rput(2.6,.8){$X=$}
\end{pspicture}
$$
where the equivalence classes
of $\sim$ are $\{v_{11},v_{21}\}, \{v_{12},v_{22}\}$ and the $\{x\}$ for all other
cells $x\in X_1\bigcup X_2$.
\end{lemma}

The \emph{standard orientation\/} $\mathcal{O}$ for a $1$-sphere consists of all the vertices and
$\{e_0,e_1,\ldots,$ $e_{n-1}\}$, ie: the edges taken in a clockwise direction in
Figure \ref{chapter1:graphs:figure550}. From now on, an orientation preserving map
between $1$-spheres \emph{preserves the standard orientations\/} on each.

\begin{exercise}\label{chapter1:graphs:exercise550}
Let $X,Y$ be $1$-spheres with their standard orientations $\mathcal{O}_X=\{e_{1i}\}$,
$\mathcal{O}_Y=\{e_{2i}\}$ and $f:X\rightarrow Y$ an orientation preserving map. If
$f(e_{1i})=e_{2j}$, then show that $f(e_{1,i+1})=e_{2,j+1}$, or is the vertex $t(e_{2j})$.
Deduce that if $f$ is an orientation preserving
isomorphism then it is a rotation (in the obvious sense).
\end{exercise}

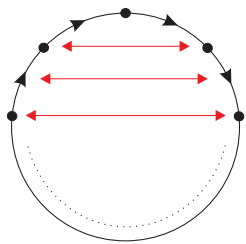
\begin{figure}
\begin{pspicture}(0,0)(13,3)
\rput(6.5,1.5){\BoxedEPSF{chapter1.fig6100a.eps scaled 800}}
\end{pspicture}
\caption{the map $\iota:X\rightarrow X$ for $X\approx S^1$.}
\label{chapter1:graphs:figure575}
\end{figure}

There is one particular map of $1$-spheres that is not orientation preserving but 
will be nevertheless useful later on. If $X$ is a $1$-sphere, 
let $\iota:X\rightarrow X$ be the
map interchanging the edges $e_i$ and $e_{n-i+1}^{-1}$ as in Figure 
\ref{chapter1:graphs:figure575}.

A \emph{path\/} in $X$ is a graph mapping $\gamma:B\rightarrow X$ with $B$ a $1$-ball.
It is convenient not to insist that the map preserve dimension, but by 
Exercise \ref{chapter1:graphs:exercise600} 
below, we can always replace the $1$-ball by another 
so that the map is dimension preserving.
In any case, the image in $X$ is a sequence of edges $e_1\ldots e_k$
(which we will also call $\gamma$),
that are consecutively incident in the obvious way: $s(e_i^{-1})=s(e_{i+1})$, and there is no
harm in thinking about paths in terms of their images. 
A path 
joins the vertices $s(e_1),s(e_k^{-1})$ 
that are the images of the end vertices of the $1$-ball, and is \emph{closed\/} if these
end vertices have the same image. If $\gamma:B\rightarrow X$ is the path
$e_1\ldots e_k$ then the \emph{inverse path\/} $\gamma^{-1}:B\rightarrow X$
has edges $e_k^{-1}\ldots e_1^{-1}$.

\begin{exercise}\label{chapter1:graphs:exercise600}
let $\gamma:B\rightarrow X$ be a path with the edges of $B$ labelled 
$e_1^\pm\ldots,e_n^\pm$ as in the comments before Lemma \ref{complexes:graphs:result300},
and image edges $e'_1\ldots e'_k$ in $X$. Show there are 
$1\leq i_1\leq i_1\leq\cdots\leq i_\ell\leq n$
with $f(e_{i_j})=e'_j$ and all other edges mapped to vertices. Thus, $B$ can be replaced
by a $1$-ball $B'$ and dimension preserving map $\gamma':B'\rightarrow X$ having the same
image path.
\end{exercise}

\begin{exercise}\label{chapter1:graphs:exercise650}
Show that a closed path $\gamma:B\rightarrow X$ gives a mapping $S\rightarrow X$
with $S$ a sphere.
\end{exercise}

If $f:X\rightarrow Y$ is a graph map and $\gamma:B\rightarrow X$ a path, then there
is an induced path in $Y$ given by the composition $f\gamma:B\rightarrow X\rightarrow Y$.
Thus, a graph mapping sends paths to paths. An example is when we have a quotient
relation $\sim$ on a graph $X$ with quotient map $q:X\rightarrow X\quo$. Then
the quotient relation can be easily extended to paths in $X$: if $\gamma,\mu$
are two such, then $\gamma\sim\mu$ precisely when $q(\gamma)=q(\mu)$ give the 
same path in the quotient $X\quo$.

The graph $X$ is \emph{connected\/} if any two vertices can be joined by a path. 
The \emph{connected component\/} of $X$ containing the vertex $v$ consists of those 
vertices that can be joined to $v$, together with all their incident edges.
A connected graph has finitely many edges if and only if it is finite and locally finite.

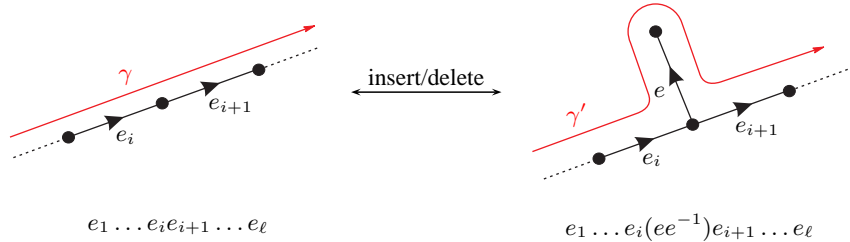
\begin{figure}
\begin{pspicture}(0,0)(13,3)
\rput(3,2){\BoxedEPSF{fig21.eps scaled 1000}}
\rput(10,2){\BoxedEPSF{fig22.eps scaled 1000}}
\psline[linewidth=.2mm]{<->}(5.5,2)(7.5,2)
\rput(6.5,2.2){insert/delete}
\rput(3.2,.2){$e_1\ldots e_ie_{i+1}
\ldots e_\ell$}
\rput(10,.2){$e_1\ldots e_i
(ee^{-1})e_{i+1}\ldots e_\ell$}
\rput(2.5,1.4){$e_i$}
\rput(3.9,1.8){$e_{i+1}$}
\rput(9.5,1.1){$e_i$}
\rput(10.9,1.5){$e_{i+1}$}
\rput(9.6,2){$e$}
\rput(2.5,2.2){${\red\gamma}$}\rput(8.5,1.7){${\red\gamma'}$}
\end{pspicture}
\caption{elementary homotopy}\label{chapter1:graphs:figure600}
\end{figure}

A \emph{spur\/} is a path of the form $ee^{-1}$,
ie: a path that consecutively traverses an edge and its inverse. 
An \emph{elementary homotopy\/} of a path, 
$e_1\ldots e_ie_{i+1}\ldots e_\ell\leftrightarrow 
e_1\ldots e_i(ee^{-1})e_{i+1}\ldots e_\ell$
inserts or deletes a {\em spur\/} so that incidence is preserved as in Figure \ref{chapter1:graphs:figure600}.
Two paths are \emph{freely homotopic\/} 
if and only if there is a finite sequence of elementary homotopies
taking one to the other. A path is {\em homotopically trivial\/} 
if it is homotopic to the trivial path. 
We leave it as an exercise
to formulate these notions for a path $B\rightarrow X$ thought of as a mapping.

\begin{exercise}\label{chapter1:graphs:exercise700}
Show that homotopic paths have the same start and end vertices, and thus homotopically 
trivial paths are necessarily closed. Show that homotopy is an equivalence relation on the
set of paths joining two fixed vertices.
\end{exercise}

\subsection{Forests and Trees}

A path in a graph is \emph{reduced\/} when it contains no spurs. 
By removing spurs we can ensure that if there exists a 
path between two vertices, then there must exist a reduced path; indeed 
for any two vertices, there is a reduced path between them if and only if 
they lie in the same component.

\begin{proposition}
\label{chapter1:graphs:result400}
The following are equivalent for a graph $X$:
\begin{enumerate}
\item There is at most one reduced path joining any two vertices;
\item any closed path is homotopically trivial;
\item any non-trivial closed path contains a spur.
\end{enumerate}
\end{proposition}

A graph satisfying any of the conditions of Proposition \ref{chapter1:graphs:result400} 
is called a \emph{forest\/}, and a connected forest is a {\em tree\/}. 

\begin{proof}
$(1\Rightarrow 2)$: a closed path $\gamma$ is necessarily contained in a component of the graph, 
thus by assumption there is a unique reduced path connecting any two of it's vertices. 
The path cannot contain just a single
vertex $u$ and edge $e$ (which it circumnavigates some number of times), for if so, 
then the edge $e$ and the trivial path at $u$ are distinct reduced paths from $u$ to $u$. 
Thus any closed path contains at least two distinct vertices. 
We show the path $\gamma$ is homotopic to the trivial path based at one of them, say $u$. 
Let $v$ be another vertex of $\gamma$ and $w$ 
the reduced path running from $u$ to $v$ as in Figure \ref{chapter1:graphs:figure700} (left). 
Then $\gamma$ decomposes into two parts, $w_1$ running 
from $u$ to $v$ and $w_2$ running from $v$ to $u$. If $w_1\not= w$ then it cannot be reduced, 
hence must contain a spur. Removing it and continuing, we have a series of homotopies that 
reduces $w_1$ to $w$ as in Figure \ref{chapter1:graphs:figure700} (middle).
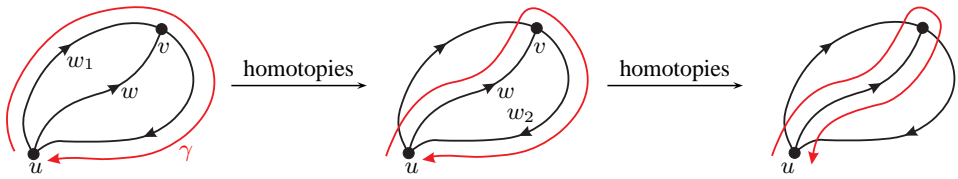
\begin{figure}
\begin{pspicture}(0,0)(13,2.5)
\rput(-.5,0){
\rput(0,0){
\rput(2,1.5){\BoxedEPSF{fig147a.eps scaled 350}}
\rput(1,.4){$u$}\rput(2.7,2){$v$}\rput(1.6,1.8){$w_1$}
\rput(2.25,1.4){$w$}\rput(3,.6){${\red\gamma}$}
}
\psline[linewidth=.2mm]{->}(3.6,1.5)(5.4,1.5)
\rput(4.5,1.7){homotopies}
\rput(5,0){
\rput(2,1.5){\BoxedEPSF{fig147b.eps scaled 350}}
\rput(1,.4){$u$}\rput(2.7,2){$v$}\rput(2.45,1.1){$w_2$}
\rput(2.25,1.4){$w$}
}
\psline[linewidth=.2mm]{->}(8.6,1.5)(10.4,1.5)
\rput(9.5,1.7){homotopies}
\rput(10,0){
\rput(2,1.5){\BoxedEPSF{fig147c.eps scaled 350}}
\rput(1,.4){$u$}
}}
\end{pspicture}
\caption{homotoping a closed path to the trivial path}\label{chapter1:graphs:figure700}
\end{figure}
Similarly for $w_2$ and $w^{-1}$ (which is the unique reduced path from $v$ to $u$)
as on the right of Figure \ref{chapter1:graphs:figure700}. 
Thus our path is homotopic to $ww^{-1}$, which in turn can be reduced by the removal of 
spurs to the trivial path based at $u$ as required.

$(2\Rightarrow 1)$: if $u$ and $v$ lie in different components then there is
no reduced path connecting them. Otherwise, if $w_1$ and $w_2$ are reduced paths 
running from $u$ to $v$ then the path $\gamma=w_1w_2^{-1}$
is by assumption homotopically trivial, hence contains a spur. As the $w_i$ are reduced, the 
spur must be at the beginning or the end of $\gamma$, ie: involve the first edges 
of $w_1$ and $w_2$, or the last edges as in Figure \ref{chapter1:graphs:figure800}.
\begin{figure}
\begin{pspicture}(0,0)(13,2.5)
\rput(4.5,1.5){\BoxedEPSF{fig148a.eps scaled 250}}
\rput(8.5,1.5){\BoxedEPSF{fig148b.eps scaled 250}}
\rput(6.5,1.5){or}
\rput(3.4,.25){$u$}\rput(5.55,2.55){$v$}\rput(7.35,.3){$u$}\rput(9.25,2.25){$v$}
\end{pspicture}
\caption{uniqueness of reduced paths}\label{chapter1:graphs:figure800}
\end{figure}
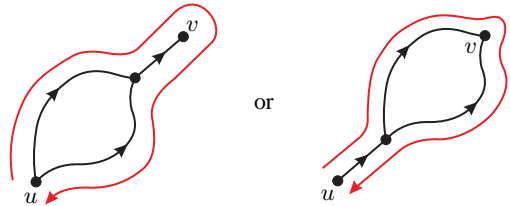
Shifting attention to the reduced subpaths not involving this spur and continuing, 
we get $w_1=w_2$. We leave the equivalence of parts 2 and 3 as an Exercise.
\qed
\end{proof}

\begin{exercise}\label{chapter1:graphs:exercise800}
Let $X$ be a finite graph, remembering that this means that the 
vertex set $X^0$ is finite. If each vertex has valency at 
least two, show that $X$ contains a homotopically non-trivial closed path.
Deduce that if $T$ is a finite tree, then $|T^1|=2(|T^0|-1)$.  
\end{exercise}

We'll have more to say about trees later. We finish this section by considering
how to approximate a graph by a tree:
if $X$ is a connected graph, then a \emph{spanning tree\/} is a subgraph 
$T\subset X$ that is a tree and contains all the vertices of $X$ (ie: $T^0=X^0$). 
The following exercise shows that under some mild set-theoretic
assumptions, spanning trees always exist.

\begin{exercise}\label{chapter1:graphs:exercise900}
Recall the well-ordering principle from set theory: any set $X$ can be 
partially ordered $\leq$ (see \S\ref{chapter3:lattices:aside})
so that for any $x,y\in X$, either $x\leq y$ or $y\leq x$,
and for any $S\subset X$ there is a $s\in S$ with $s\leq x$ for all 
$x\in X$.
In particular, choose a well ordering of the edges set of a graph.
Choose a basepoint vertex $v_0$, and consider those vertices at distance one
from $v_0$, ie: the $v\not= v_0$ with $s(e)=v_0$, $s(e^{-1})=v$ for some
edge $e$. For each such, choose
an edge $e_v$ that is minimal in the well-ordering amongst the edges 
joining $v_0$ to $v$. Let $T_1$ be the subgraph consisting of $v_0$, its distance 
$1$ neighbours  and the edges so chosen.
\begin{enumerate}
\item Show that $T_1$ is a tree. Continue the construction inductively:
at step $k$, take the tree $T_{k-1}$ constructed at step
$k-1$, and for each vertex $v$ of $X$ a distance $1$ from a vertex of 
$T_{k-1}$, choose a minimal edge as above. 
Let $T_k$ be the subgraph consisting of $T_{k-1}$ together 
with the distance $1$ vertices and minimal edges. Show that $T_k$ is a tree.
\item Show that $T=\bigcup T_k$, is the required spanning tree.
\end{enumerate}
\end{exercise}

In a graph the edge 
set can have a wildly different cardinality from the vertex set, causing difficulties 
with some arguments. This shortcoming is avoided by spanning trees which 
have a number of edges that is ``roughly'' the same as the number of vertices
of the graph they span:

\begin{proposition}\label{chapter1:graphs:result500}
Let $X$ be a connected graph and $T\subset X$ a spanning tree. Then
$$
|T^1|=\left\{
\begin{array}{ll}
2|X^0|-1,&\text{if $X$ is finite},\\
|X^0|,&\text{if $X$ is infinite}.
\end{array}\right.
$$
\end{proposition}

\begin{proof}
The result for finite graphs is the content of Exercise \ref{chapter1:graphs:exercise800}. 
If $X$ is an infinite graph with spanning tree $T$, then the edge set of $T$ must be infinite, 
as a finite edge set only spans
$|T^1|+1$ vertices.
Then $X^0=T^0=\bigcup_{e\in T^1}\{s(e),s(e^{-1})\}$ has the same cardinality as
$T^1$.
\qed
\end{proof}

\begin{exercise}\label{chapter1:graphs:exercise1000}
Let $T_\aa\subset X\,(\aa\in A)$ be a family of mutually disjoint trees
in a connected graph $X$. Show there is a spanning tree $T\subset X$
containing the $T_\aa$ as subgraphs, and such that 
$q(T)$ is a spanning tree for $X/T_\aa$, where $q:X\rightarrow X/T_\aa$ is the quotient map.

A \emph{spanning forest\/} is a subgraph $\Phi\subset X$
that is a forest and contains all the vertices of $X$. 
By considering $q^{-1}(T')$ for some spanning tree $T'$ of 
the (connected) graph $X/\Phi$, show that any spanning forest
can be extended to a spanning tree.
\end{exercise}

\section{The category of $2$-complexes}\label{chapter1:2complexes}

\subsection{$2$-complexes}\label{chapter1:2complexes:2complexes}

\begin{definition}[$\mathbf{2}$-complex]
\label{chapter1:2complexes:definition100}
A {\em combinatorial $2$-complex\/} $X$
is a graded set $X=X^0,X^1,$ $X^2$ 
with $X^i\not=\varnothing$, such that
\begin{description}
\item[(C1).] if $X^{(1)}:=X^0\bigcup X^1$, there are maps $^{-1}:X\rightarrow X$
and $s:X^{(1)}\rightarrow X^0$ making $X^{(1)}$ a  graph; 
\item[(C2).] each $\ss\in X^2$ has \emph{boundary\/} $\partial\ss=(X^\ss,\aa_\ss)$
with $X^\ss$ a $1$-sphere and 
$\aa_\ss:X^\ss\rightarrow X^{(1)}$ a dimension preserving
map of graphs. The map $^{-1}$ extends to all of $X$, with $\ss^{-1}\not=\ss=(\ss^{-1})^{-1}$
and $\partial\ss^{-1}=(X^\ss,\aa_\ss\iota)$, with $\iota$ the map given in
\S\ref{chapter1:graphs:paths}. 
\end{description}
\end{definition}

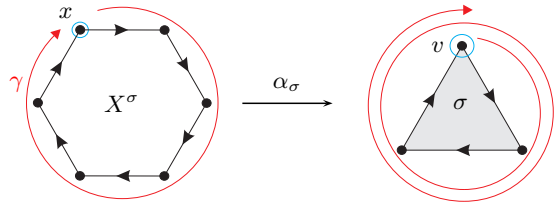
\begin{figure}
\begin{pspicture}(0,0)(12.5,2.5)
\rput(4,1.25){\BoxedEPSF{chapter1.fig1100.eps scaled 750}}
\rput(8.5,1.25){\BoxedEPSF{chapter1.fig1800a.eps scaled 750}}
\psline[linewidth=.2mm]{->}(5.6,1.25)(6.8,1.25)
\rput(6.2,1.5){$\aa_\ss$}\rput(4,1.25){$X^\ss$}\rput(8.5,1.25){$\ss$}
\rput(8.2,2){$v$}\rput(3.25,2.45){$x$}
\rput(2.6,1.5){${\red\gamma}$}
\end{pspicture}
\caption{boundary of a face $\ss$ consisting of a $1$-sphere and an attaching map, which
may wrap the sphere around the face several times. The vertex $v$ appears in the boundary
of $\ss$ and the red path wrapping twice around is a boundary path of $\ss$ starting at $v$.}
\label{chapter1:2complexes:figure100}
\end{figure}

The elements of $X^2$ are the $2$-cells or \emph{faces\/} and $\aa_\ss$
is the \emph{attaching map\/} of the face $\ss$
(see Figure \ref{chapter1:2complexes:figure100}). The underlying graph $X^{(1)}$
is called the \emph{$1$-skeleton\/}. One thinks of a face as a disc sewn by its boundary onto the $1$-skeleton 
as in Figure \ref{chapter1:2complexes:figure100}.

We say that the vertex $v$ \emph{appears in the boundary\/} of the face $\ss$ whenever
$\aa_\ss^{-1}(v)\not=\varnothing$. Thus there is a vertex $x$ of $X^\ss$ mapping to 
$v$ via the attaching map $\aa_\ss$ as in Figure \ref{chapter1:2complexes:figure100}, 
and indeed there may be several of them. 
We will call the vertices in $\aa_\ss^{-1}(v)\subset X^\ss$ 
\emph{the appearances of $v$ in the boundary of 
the face $\ss$\/}.
Similiarly for edges.
For an appearance $x$ of $v$, if we take
a path $\gamma$ consisting of the standard orientation of $X^\ss$,
then we call its image a \emph{boundary path of
$\ss$ starting at $v$\/}. 
We will often write $\gamma$ for both the path in $X^\ss$ and its image under the 
attaching map.
The vertex $v$ appears a total of $|\aa_\ss^{-1}(v)|$ times in the boundary
of $\ss$, and each appearance gives rise to a pair of boundary paths starting at $v$. 

Many of the concepts and adjectives pertaining to graphs can be applied 
to $2$-complexes by considering the $1$-skeleton. Thus we have
finite $2$-complexes, connected $2$-complexes, paths in
$2$-complexes, etc. An \emph{orientation\/} $\mathcal{O}$ for a
$2$-complex is a set consisting of all the vertices
and exactly one edge
from each arc and exactly one face from each face/inverse face pair.


\begin{exercise}
\label{chapter1:2complexes:definition200}
There is a more elementary notion of $2$-complex where the boundary $(X^\ss,\aa_\ss)$
is identified with its image under the attaching map:
we have a graded set $X=X^0,X^1,$ $X^2$ 
with $X^i\not=\varnothing$, such that
\begin{description}
\item[(C1).] if $X^{(1)}:=X^0\bigcup X^1$, there are maps $^{-1}:X\rightarrow X$
and $s:X^{(1)}\rightarrow X^0$ making $X^{(1)}$ a  graph; 
\item[(C2).] each $\ss\in X^2$ has \emph{boundary\/} $\partial\ss$
all cyclic permutations of some fixed closed path $\gamma_\ss$;
\item[(C3).] the map $^{-1}$ has no fixed points in $X^2$ and
$\partial\ss^{-1}$ consists of all cyclic permutations of the inverse
path $\gamma_\ss^{-1}$.
\end{description}
Show that a $2$-complex in our sense gives rise to such a $2$-complex, and
vice-versa (although we will prefer to keep track of the attaching maps
of the faces).
\end{exercise}

\begin{figure}
\begin{pspicture}(0,0)(13,4.5)
\rput(0,.1){
\rput(0,2){
\rput(1,1.5){\BoxedEPSF{chapter1.fig4000.eps scaled 1000}}
\rput(4,1.5){\BoxedEPSF{chapter1.fig4000.eps scaled 1000}}
\rput(7,1.5){\BoxedEPSF{chapter1.fig4000.eps scaled 1000}}
\psline[linewidth=.2mm]{->}(2.1,1.5)(2.9,1.5)
\psline[linewidth=.2mm]{<-}(5.1,1.5)(5.9,1.5)
\rput(1,1.5){$X^{\ss_1}$}\rput(4,.3){$X^{(1)}$}\rput(7,1.5){$X^{\ss_2}$}
\rput(2.5,1.7){$\aa_{\ss_{1}}$}\rput(5.5,1.7){$\aa_{{\ss_2}}$}
\rput(1,2.1){$e_1$}\rput(1,.9){$e_2$}
\rput(7.1,2.1){$e_1$}\rput(7.1,.95){$e_2$}
\rput(3.45,1.5){$v_1$}\rput(4.6,1.5){$v_2$}\rput(4,2.1){$e_1$}\rput(4,.9){$e_2$}
}
\rput(1,-.5){
\rput(3.3,1.5){\BoxedEPSF{chapter1.fig4100.eps scaled 1000}}
\rput(.5,1.5){$v_1$}\rput(3.75,1.5){$v_1$}
\rput(2.9,1.5){$v_2$}\rput(6.1,1.5){$v_2$}
\rput(1.7,.4){$e_2$}\rput(1.7,2.6){$e_1$}\rput(1.7,1.5){$\ss_1$}
\rput(4.9,.4){$e_2$}\rput(4.9,2.6){$e_1$}\rput(4.9,1.5){$\ss_2$}
}
\rput(-2,0.8){
\rput(12.5,1.5){\BoxedEPSF{chapter1.fig4200.eps scaled 500}}
\rput(2.5,-4.5){
\rput(11.6,7.4){$\ss_1$}
\rput(11.6,4.6){$\ss_2$}
\rput(12.45,6.4){${\red\partial\ss_1}$}
\rput(12.45,5.8){${\red\partial\ss_2}$}
\rput(9.7,5.45 ){${v}_1$}
\rput(10.2,6.55){${v}_2$}
\rput(10.6,5.45){${e}_2$}
\rput(9.4,6.6){${e}_1$}
}}
}
\end{pspicture}
\caption{a $2$-complex: pictorial version 
of definition \ref{chapter1:2complexes:definition100} (top left);
face-centric version with faces sewn on (bottom left) 
and topologically suggestive version (right).}
\label{chapter1:2complexes:figure200}
\end{figure}
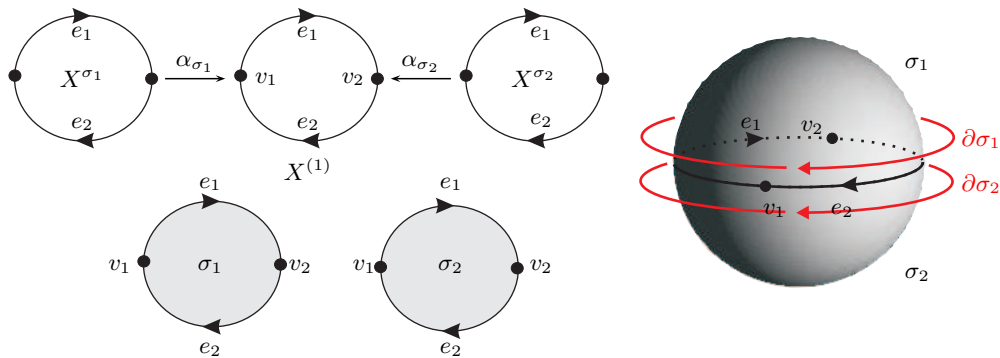

\subsection{Examples}\label{chapter1:2complexes:examples}

Figure \ref{chapter1:2complexes:figure200} gives three different versions of a 
$2$-complex that we will call the $2$-sphere. The first version (top left) is 
a straight pictorial version of definition \ref{chapter1:2complexes:definition100}:
the $1$-skeleton in the middle is a graph with two vertices and two
edges; there are two faces and the attaching maps are described by
labeling the edges of the $X^{\ss_i}$ with their images in the $1$-skeleton. This version
is both the most accurate and the most cumbersome.

In the second version (bottom left) we have adopted the
convention that parts of the complex with the
same label give the same cell and have drawn the complex ``face-centrically'' with the faces
thought of as discs sewn onto the $1$-skeleton. 
The third version is more along the lines of Exercise
\ref{chapter1:2complexes:definition200}, with the face boundaries
given by closed paths and their cyclic permutations.

\begin{figure}
\begin{pspicture}(0,0)(13,4.5)
\rput(0,.1){
\rput(0,2){
\rput(1,1.5){\BoxedEPSF{chapter1.fig4000.eps scaled 1000}}
\rput(4,1.5){\BoxedEPSF{chapter1.fig4000.eps scaled 1000}}
\rput(7,1.5){\BoxedEPSF{chapter1.fig4000.eps scaled 1000}}
\psline[linewidth=.2mm]{->}(2.1,1.5)(2.9,1.5)
\psline[linewidth=.2mm]{<-}(5.1,1.5)(5.9,1.5)
\rput(1,1.5){$X^{\ss_1}$}\rput(4,.3){$X^{(1)}$}\rput(7,1.5){$X^{\ss_2}$}
\rput(2.5,1.7){$\aa_{\ss_{1}}$}\rput(5.5,1.7){$\aa_{{\ss_2}}$}
\rput(1,2.1){$e_1$}\rput(1,.9){$e_2$}
\rput(7.1,2.1){$e_2^{-1}$}\rput(7.1,.95){$e_1^{-1}$}
\rput(3.45,1.5){$v_1$}\rput(4.6,1.5){$v_2$}\rput(4,2.1){$e_1$}\rput(4,.9){$e_2$}
}
\rput(1,-.5){
\rput(3.3,1.5){\BoxedEPSF{chapter1.fig4100.eps scaled 1000}}
\rput(.5,1.5){$v_1$}\rput(3.75,1.5){$v_1$}
\rput(2.9,1.5){$v_2$}\rput(6.1,1.5){$v_2$}
\rput(1.7,.4){$e_2$}\rput(1.7,2.65){$e_1$}\rput(1.7,1.5){$\ss_1$}
\rput(4.9,.4){$e_1^{-1}$}\rput(4.9,2.65){$e_2^{-1}$}\rput(4.9,1.5){$\ss_2$}
}
\rput(-2,0.8){
\rput(12.5,1.5){\BoxedEPSF{chapter1.fig4200a.eps scaled 500}}
\rput(2.5,-4.5){
\rput(11.6,7.4){$\ss_1$}
\rput(11.6,4.6){$\ss_2$}
\rput(12.45,6.4){${\red\partial\ss_1}$}
\rput(12.45,5.8){${\red\partial\ss_2}$}
\rput(9.7,5.45 ){${v}_1$}
\rput(10.2,6.55){${v}_2$}
\rput(10.6,5.45){${e}_2$}
\rput(9.4,6.6){${e}_1$}
}}
}
\end{pspicture}
\caption{another $2$-complex}
\label{chapter1:2complexes:figure150}
\end{figure}
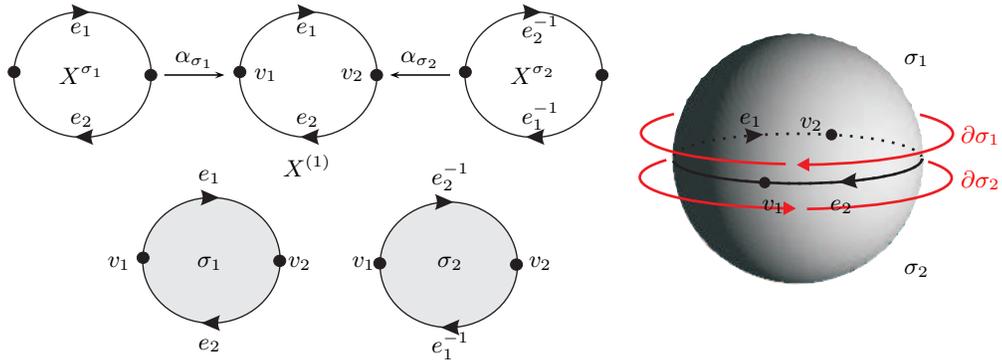

Figure \ref{chapter1:2complexes:figure150} gives a very similar
example, except that one of the face boundaries goes the other way
around the $1$-skeleton. 

\begin{figure}
\begin{pspicture}(0,0)(12.5,2.5)
\rput(8.2,0){
\rput(1.25,1.25){\BoxedEPSF{chapter1.fig4300.eps scaled 500}}
\rput(2.2,.2){$\ss$}\rput(.4,1.55){$v$}\rput(2.2,2.4){$v$}\rput(1.2,1.05){$e$}
\rput(1.2,2.5){$e$}
}
\rput(-1.8,-0.05){
\rput(4.1,1.5){\BoxedEPSF{chapter3.fig2000.eps scaled 650}}
\rput(7,1.5){\BoxedEPSF{chapter1.fig4000.eps scaled 1000}}
\psline[linewidth=.2mm]{<-}(5.1,1.5)(5.9,1.5)
\rput(4.1,1){$X^{(1)}$}\rput(7,1.5){$X^{\ss}$}
\rput(5.5,1.7){$\aa_{\ss}$}
\rput(4.9,1.5){$v$}
\rput(3.35,1.5){$e$}
\rput(7,2.1){$e$}\rput(7,.9){$e$}
}
\end{pspicture}
  \caption{projective plane complex}
  \label{chapter1:2complexes:figure400}
\end{figure}
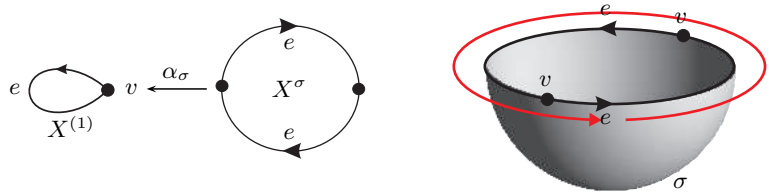

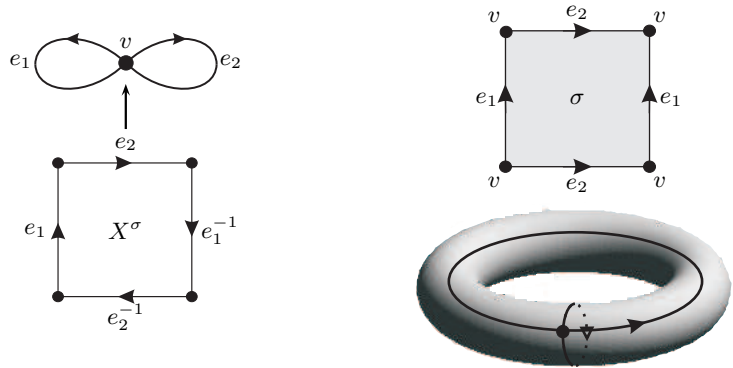
\begin{figure}
\begin{pspicture}(0,0)(12.5,5)
\rput(0,.5){
\rput(2,0.25){
\rput(1,1){\BoxedEPSF{chapter3.fig3000a.eps scaled 900}}
\rput(1,1){$X^\ss$}
\rput(-.2,1){$e_1$}
\rput(2.25,1){$e_1^{-1}$}
\rput(1,2.15){$e_2$}
\rput(1,-.15){$e_2^{-1}$}
}
\psline[linewidth=.3mm]{->}(3,2.6)(3,3.2)
\rput(0,0){
\rput(3,3.75){$v$}
\rput(1.6,3.5){$e_1$}
\rput(4.4,3.5){$e_2$}
\rput(3,3.5){\BoxedEPSF{chapter3.fig2100.eps scaled 750}}
}
}
\rput(9,1){\BoxedEPSF{fig121.eps scaled 500}}
\rput(8,2.5){
\rput(1,1){\BoxedEPSF{fig9.eps scaled 1000}}
\rput(1,1){$\ss$}
\rput(2.1,2.1){$v$}
\rput(2.1,-0.1){$v$}
\rput(-0.1,2.1){$v$}
\rput(-0.1,-0.1){$v$}
\rput(-.2,1){$e_1$}
\rput(2.25,1){$e_1$}
\rput(1,2.15){$e_2$}
\rput(1,-.15){$e_2$}
}
\end{pspicture}
\caption{various versions of the torus.}
 \label{chapter1:2complexes:figure500}
\end{figure}

Figure \ref{chapter1:2complexes:figure400} is the (real) projective plane $\R\P^2$, a
combinatorial model for the disc with antipodal points on the boundary
identified.
We have
again drawn the complex both face-centrically 
and \emph{ala\/} Definition \ref{chapter1:2complexes:definition100}.
Similarly, Figure \ref{chapter1:2complexes:figure500} shows various versions
of the torus complex.

\subsection{Maps of $2$-complexes}\label{chapter1:2complexes:maps}

To complete the definition of the category of $2$-complexes we need mappings.
The principle is the same as for graph maps
in \S\ref{chapter1:graphs:category}: continuity is captured 
by making maps commute with the various attaching maps of the 
cells. The definition is complicated
slightly as we allow a map to squash a face down to a path.

\begin{definition}[maps of $\mathbf{2}$-complexes]
\label{chapter1:2complexes:definition250}
A {\em map\/} $f:X\rightarrow Y$ of $2$-complexes is a map
$f:X^{(1)}\rightarrow Y^{(1)}$ of the underlying graphs such that for
each face $\ss\in X^2$ we have either $f(\ss)$ is a face $\tau\in Y^2$ or
$f(\ss)$ is a closed path in the graph $Y^{(1)}$. There are the conditions:
\begin{description}
\item[(M1).] Let $f(\ss)=\tau$, a face in $Y$, with 
$\partial\ss=(X^\ss,\aa_\ss)$ and 
$\partial\tau=(Y^\tau,\aa_\tau)$. Then
there is an \emph{orientation preserving\/} map
$\ve=\ve(f,\ss):X^\ss\rightarrow Y^{\tau}$ making the diagram
below left commute;
$$
\begin{pspicture}(0,0)(14,2)
\rput(-.5,-.8){
\rput(-1.3,.5){
\rput(5,2){$X^\ss$}\rput(6.65,2){$Y^{\tau}$}
\rput(4.95,0.45){$X^{(1)}$}\rput(6.65,0.45){$Y^{(1)}$}
\psline[linewidth=.2mm]{->}(5.3,2)(6.3,2)
\psline[linewidth=.2mm]{->}(5.3,0.45)(6.3,0.45)
\psline[linewidth=.2mm]{->}(5,1.7)(5,.7)
\psline[linewidth=.2mm]{->}(6.55,1.7)(6.55,.7)
\rput(4.7,1.2){$\aa_\ss$}\rput(6.9,1.175){$\aa_{\tau}$}
\rput(5.8,2.2){$\ve$}
\rput(5.8,.7){$f$}
}
\rput(3.7,.5){
\rput(5,2){$X^\ss$}\rput(6.55,2){$S$}
\rput(4.95,0.45){$X^{(1)}$}\rput(6.65,0.45){$Y^{(1)}$}
\psline[linewidth=.2mm]{->}(5.3,2)(6.3,2)
\psline[linewidth=.2mm]{->}(5.3,0.45)(6.3,0.45)
\psline[linewidth=.2mm]{->}(5,1.7)(5,.7)
\psline[linewidth=.2mm]{->}(6.55,1.7)(6.55,.7)
\rput(4.7,1.2){$\aa_\ss$}\rput(6.85,1.2){$\gamma$}
\rput(5.8,2.2){$\ve$}
\rput(5.8,.7){$f$}
}
}
\end{pspicture}
$$
Moreover $f(\ss^{-1})=\tau^{-1}$ and $\ve(f,\ss^{-1})=\ve(f,\ss)$;
\item[(M2).] Let $f(\ss)$ be
the closed path $\gamma:S\rightarrow Y$ where $S$ is a $1$-sphere. Then there
is an orientation preserving map $\ve=\ve(f,\ss):X^\ss\rightarrow S$ 
making the diagram above right commute.
Moreover, $f(\ss^{-1})$ is the inverse path and the path 
$\gamma:S\rightarrow X$ is \emph{homotopically trivial\/}.
\end{description}
\end{definition}

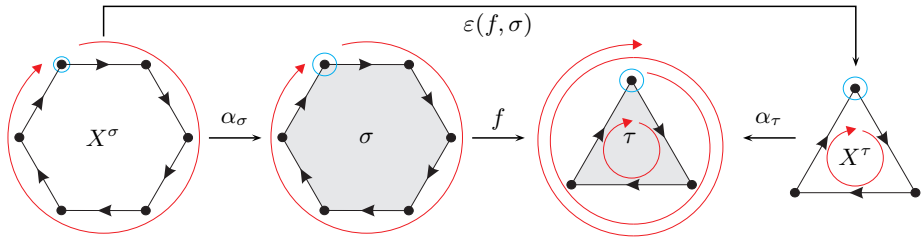
\begin{figure}
  \centering
\begin{pspicture}(0,0)(12.5,2.75)
\rput(1.5,1.25){\BoxedEPSF{chapter1.fig1100.eps scaled 750}}
\rput(5,1.25){\BoxedEPSF{chapter1.fig1700.eps scaled 750}}
\rput(8.5,1.25){\BoxedEPSF{chapter1.fig1800.eps scaled 750}}
\rput(11.5,1.25){\BoxedEPSF{chapter1.fig1800b.eps scaled 750}}
\psline[linewidth=.2mm]{->}(2.9,1.25)(3.6,1.25)
\psline[linewidth=.2mm]{->}(6.4,1.25)(7.1,1.25)
\psline[linewidth=.2mm]{->}(10.7,1.25)(10,1.25)
\psline[linewidth=.2mm]{->}(1.5,2.6)(1.5,3)(11.5,3)(11.5,2.2)
\rput(1.5,1.25){$X^\ss$}\rput(5,1.25){$\ss$}
\rput(11.5,1){$X^\tau$}\rput(8.5,1.25){$\tau$}
\rput(3.25,1.5){$\aa_\ss$}\rput(10.35,1.5){$\aa_\tau$}
\rput(6.75,1.5){$f$}\rput(6.75,2.75){$\ve(f,\ss)$}
\end{pspicture}  
  \caption{a boundary path for $\ss$ is mapped via $f$ to a boundary path of $\tau$, 
possibly repeated.}
\label{chapter1:2complexes:figure550}
\end{figure}

Thus, if $f(\ss)=\tau$, then a boundary path for $\ss$ is mapped to a path
that circumnavigates a boundary path for $\tau$, possibly a number of times.
See Figure \ref{chapter1:2complexes:figure550}.

If a face is mapped to a path then this path 
must be the image (possibly repeated) of the boundary path of the face, hence closed. 
Moreover it is
homotopically trivial. The motivating example is squashing a face flat: the boundary gets
squashed too, into a path of the form $e_1\ldots e_\ell e_\ell^{-1}\ldots e_1^{-1}$.
In this example the result is clearly homotopically trivial, but we will really
require this condition in Chapter \ref{chapter2} for the fundamental group 
$\pi_1$ to be a functor.

\begin{definition}[map conventions]
Let $X,Y$ and $Z$ be $2$-complexes.
\begin{enumerate}
\item \emph{When are two maps the same?} Let 
$f,g:X\rightarrow Y$ be
maps of $2$-complexes. Then $f=g$ when $f(x)=g(x)$ for all
cells $x\in X$ \emph{and\/} if $f(\ss)=g(\ss)=\tau$ for faces $\ss\in X$ and $\tau\in Y$,
then $\ve(f,\ss)=\ve(g,\ss)$.
\item \emph{Forming compositions}: Suppose that
$X\stackrel{f}{\longrightarrow}Y\stackrel{g}{\longrightarrow}Z,$
are maps of $2$-complexes. Then the composition $gf$ is formed in the usual way, with the
additional proviso that if $\ss\in X$ is a face with $f(\ss)$ a face
of $Y$ and $gf(\ss)$ a face of $Z$, then
$\ve(gf,\ss)=\ve(g,f(\ss))\ve(f,\ss)$.
\end{enumerate}
\end{definition}

These conventions have ramifications for commuting diagrams of complexes
and maps, which are after all, just statements about maps being the same. For example,
when we say that the diagram of maps and complexes on the left commutes,
$$
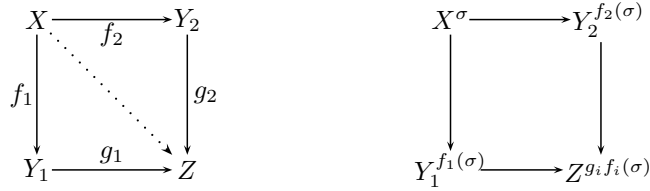
\begin{pspicture}(0,0)(12.5,2.25)
\rput(2.5,0){
\rput(0,2){$X$}
\rput(0,0){$Y_1$}
\rput(2,2){$Y_2$}
\rput(2,0){$Z$}
\psline[linewidth=.2mm]{->}(0,1.8)(0,.2)
\psline[linewidth=.2mm]{->}(.2,2)(1.8,2)
\psline[linewidth=.2mm]{->}(.2,0)(1.8,0)
\psline[linewidth=.2mm]{->}(2,1.8)(2,.2)
\psline[linewidth=.3mm,linestyle=dotted]{->}(.1,1.9)(1.8,.2)
\rput(-0.2,1){$f_1$}
\rput(1,1.8){$f_2$}
\rput(1,.2){$g_1$}
\rput(2.25,1){$g_2$}
}
\rput(8,0){
\rput(0,2){$X^\ss$}
\rput(0,0){$Y_1^{f_1(\ss)}$}
\rput(2.1,2){$Y_2^{f_2(\ss)}$}
\rput(2.1,0){$Z^{g_if_i(\ss)}$}
\psline[linewidth=.2mm]{->}(0,1.8)(0,.25)
\psline[linewidth=.2mm]{->}(.25,2)(1.6,2)
\psline[linewidth=.2mm]{->}(.4,0)(1.5,0)
\psline[linewidth=.2mm]{->}(2,1.7)(2,.2)
}
\end{pspicture}
$$
then the compositions $g_if_i\,(i=1,2)$ are the same as the dotted map
across the middle. If $\ss$ is a face of $X$ that maps to a face $g_if_i(\ss)$ of $Z$, 
then the diagram of maps on the right must also commute.

Let $f:X\rightarrow Y$ be a map of $2$-complexes,
$u,v$ vertices with $f(u)=v$, and $\tau$ a face in $Y$. 
If $\ss$ is a face of $X$ that maps to $\tau$,
let $\ve(f,\ss):X^\ss\rightarrow X^\tau$ be the orientation preserving map
from (M1) of Definition \ref{chapter1:2complexes:definition250}. This map
induces a map
$\ve(f,\ss):\aa_\ss^{-1}(u)\subset X^\ss\rightarrow\aa_\tau^{-1}(v)\subset Y^\tau$
from the appearances of $u$ in the boundary of $\ss$ to the appearances of $v$ in 
the boundary of $\tau$.
As this is true for all the $\ss$ mapping to $\tau$, we have,

\begin{definition}[local continuity]
\label{chapter1:2complexes:definition250}
Let $f:X\rightarrow Y$ be a map of $2$-complexes, $u,v$ vertices with $f(u)=v$, 
and $\tau$ a face in $Y$. The \emph{local continuity\/} of $f$ at $v$ is
$$\
\amalg\,\ve(f,\ss):\bigcup_{f(\ss)=\tau}\,\aa_\ss^{-1}(u)\rightarrow\aa_\tau^{-1}(v),
$$
where $\amalg\,\ve(f,\ss)$ is the (disjoint) union of the maps $\ve(f,\ss)$ over the 
faces $\ss$ mapping to $\tau$. 
\end{definition}

\begin{figure}
\begin{pspicture}(0,0)(12.5,4)
\rput(4,2){\BoxedEPSF{chapter1.fig700.eps scaled 800}}
\rput(3.4,2.6){$\ss_1^{}$}\rput(1.1,3.2){$X^{\ss_1^{}}$}
\rput(4.6,2.6){$\ss_2^{}$}\rput(7.1,3.2){$X^{\ss_2^{}}$}
\rput(4.6,1.4){$\ss_3^{}$}\rput(7.1,.8){$X^{\ss_3^{}}$}
\rput(3.4,1.4){$\ss_4^{}$}\rput(1.1,.8){$X^{\ss_4^{}}$}
\rput(4.2,2.2){$u$}
\rput(1,3.2){\BoxedEPSF{chapter1.fig600c.eps scaled 750}}
\psline[linewidth=.2mm]{->}(2,3)(2.6,2.8)
\rput(7,3.2){\BoxedEPSF{chapter1.fig600d.eps scaled 750}}
\psline[linewidth=.2mm]{->}(2,1)(2.6,1.2)
\rput(7,.8){\BoxedEPSF{chapter1.fig600a.eps scaled 750}}
\psline[linewidth=.2mm]{->}(6,3)(5.4,2.8)
\rput(1,.8){\BoxedEPSF{chapter1.fig600b.eps scaled 750}}
\psline[linewidth=.2mm]{->}(6,1)(5.4,1.2)
\rput(11.5,.8){\BoxedEPSF{fig121.eps scaled 300}} 
\psline[linewidth=.2mm]{->}(11.5,2.2)(11.5,1.6)
\rput(11.5,3.2){\BoxedEPSF{chapter1.fig600e.eps scaled 750}}
\rput(11.6,3.2){$Y^{\tau^{}}$}
\psline[linewidth=.2mm]{->}(8,2)(10,2)
\rput(9,2.2){$f$}
\end{pspicture}
\caption{local continuity for a map of $2$-complexes: the inifinite plane complex
maps via $f$ to the torus; $\tau$ is the single face of the torus and for 
the single vertex $v$, the set $\aa_\tau^{-1}(v)\in Y^{\tau}$ consists of the four vertices
marked with little red and blue circles and squares.
A fixed vertex $u$ mapping to $v$ appears in the boundary of four
faces $\ss_i^{}$ with $\aa_{\ss_i}^{-1}(u)$ consisting of a single vertex
in each case. For all other faces $\ss$ we have $\aa_{\ss}^{-1}(u)=\varnothing$. The red
and blue circles and squares in the $X^{\ss_i}$ map via the local continuity
of Definition \ref{chapter1:2complexes:definition250} to the corresponding ones in $Y^{\tau}$.}
 \label{chapter1:2complexes:figure600}
\end{figure}
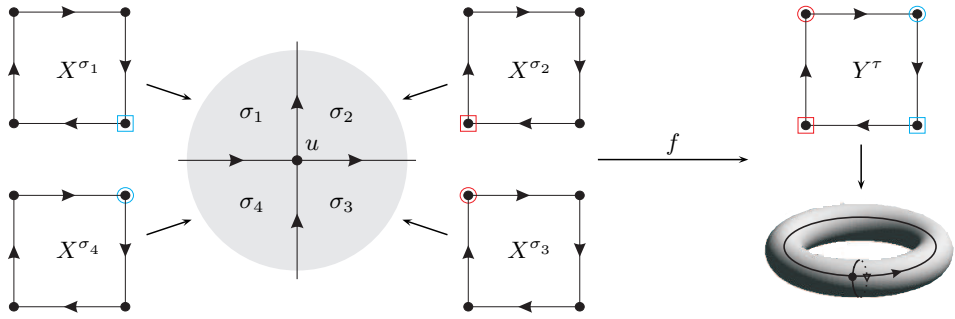

An example is given in Figure \ref{chapter1:2complexes:figure600} with the mapping of the 
``plane'' complex to the torus.

\begin{exercise}\label{chapter1:2complexes:exercise150}
Show that if the right hand side of the set in Definition \ref{chapter1:2complexes:definition250} 
is empty, then so is the left hand side.
\end{exercise}

\begin{definition}[dimension preserving maps]
\label{chapter1:2complexes:definition275}
A map $f:X\rightarrow Y$ is {\em dimension preserving\/}  if and only if
\begin{enumerate}
\item the graph map $f:X^{(1)}\rightarrow Y^{(1)}$ is dimension preserving;
\item $f(X^2)\subset Y^2$;
\item if $\ss\in X$ and $f(\ss)\in Y$ are faces then $\ve(f,\ss):X^\ss\rightarrow Y^\tau$
is an orientation preserving \emph{isomorphism\/}.
\end{enumerate}
\end{definition}
 
A map is an \emph{isomorphism\/} if it  is dimension preserving, and a
bijection on the vertex, edge and face sets. 
In this case one easily sees that, as for graphs, the inverse map $f^{-1}$
is also an isomorphism $f^{-1}:Y\rightarrow X$ (just reverse the horizontal
arrows in the left commuting diagram of Definition
\ref{chapter1:2complexes:definition250}(M1)!) so
that the set of automorphisms $f:X\rightarrow X$ forms a group $\aut(X)$ under composition.

\begin{exercise}\label{chapter1:2complexes:exercise175}
Show that the $2$-complexes of Figures \ref{chapter1:2complexes:figure200} and 
\ref{chapter1:2complexes:figure150} are isomorphic.
\end{exercise}

A group action $G\stackrel{\varphi}{\rightarrow}\aut(X)$  
\emph{preserves orientation\/} if there is an orientation
$\mathcal{O}$ for $X$ with $g(\mathcal{O})=\mathcal{O}$ for all $g\in
G$. One can then show, just as in Exercise
\ref{chapter1:graphs:exercise300}, that an action of a group preserves
orientation if and only if it acts without inversions and no $g\in G$
sends a face to its inverse.
A group acts \emph{freely\/} on a $2$-complex precisely when the
action on the underlying graph is free.

A $2$-complex $X$ is a 
\emph{subcomplex\/} of $Y$ if there
is a mapping $X\hookrightarrow Y$ of $2$-complexes that is an isomorphism
onto its image.

\subsection{Homotopies and homeomorphisms}\label{chapter1:2complexes:homotopies}

As with graphs we can deform paths, simulating in
a combinatorial manner the homotopies of paths in
topology.

Let $\gamma=e_1\ldots e_\ell$ be a path in the $2$-complex $X$. An \emph{elementary homotopy\/}
either inserts or deletes a spur as in \S\ref{chapter1:graphs:paths} or inserts or
deletes the boundary of a face in the following sense. If $\ss\in X^2$ is a face
with boundary $\partial\ss=(X^\ss,\aa_\ss)$, then by Exercise \ref{chapter1:graphs:exercise650},
$\aa_\ss(X^\ss)$ is a closed path, say $e_1'\ldots e_k'$. 
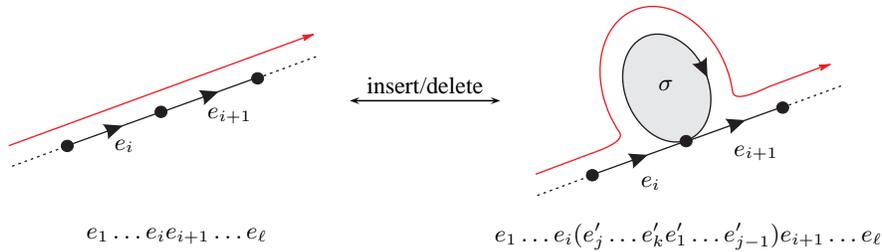
\begin{figure}
\begin{pspicture}(0,0)(13,3)
\rput(3,2){\BoxedEPSF{fig21.eps scaled 1000}}
\rput(10,2){\BoxedEPSF{fig23.eps scaled 1000}}
\psline[linewidth=.2mm]{<->}(5.5,2)(7.5,2)
\rput(6.5,2.2){insert/delete}
\rput(3.2,.2){$e_1\ldots e_ie_{i+1}\ldots e_\ell$}
\rput(10,.2){$e_1\ldots e_i(e_j'\ldots e_k'e_1'\ldots e_{j-1}')e_{i+1}\ldots e_\ell$}
\rput(2.5,1.4){$e_i$}
\rput(3.9,1.8){$e_{i+1}$}
\rput(9.5,.9){$e_i$}
\rput(10.9,1.3){$e_{i+1}$}
\rput(9.7,2.2){$\ss$}
\end{pspicture}  
  \caption{elementary homotopy--inserting/deleting the boundary of a face.}
  \label{chapter1:2complexes:figure650}
\end{figure}
The homotopy inserts into (or deletes from) $\gamma$ the result
of completely traversing this closed path, starting at one of its vertices, so that
all the incidences match up in the obvious way, ie: so that $s(e_j')=t(e_{j-1}')$
is the vertex $t(e_i)=s(e_{i+1})$.

Pictures such as the right hand side of Figure \ref{chapter1:2complexes:figure650} should
be approached with care. The \emph{entire\/} boundary path of $\ss$ must be traversed, 
including any repetitions. Later we will have faces with boundary a closed path that travels
a number of times around an edge loop. Any homotopy involving this boundary must then
travel the full number of times around the loop.

\begin{figure}
\begin{pspicture}(0,0)(13,3)
\rput(2.5,1.5){\BoxedEPSF{fig24.eps scaled 1000}}
\rput(6.6,1.5){\BoxedEPSF{fig25.eps scaled 1000}}
\rput(10.5,1.5){\BoxedEPSF{fig26.eps scaled 1000}}
\psline[linewidth=.2mm]{->}(4,1.5)(5,1.5)
\psline[linewidth=.2mm]{->}(8,1.5)(9,1.5)
\rput(2.5,1.4){$\ss$}
\rput(6.6,1.5){$\ss$}
\rput(10.4,1.6){$\ss$}
\end{pspicture}  
  \caption{homotoping a path across a face. To get from the first picture to the second, 
insert the boundary of the
face $\ss$; to get from the second to the third, remove the obvious spurs.}
  \label{chapter1:2complexes:figure700}
\end{figure}
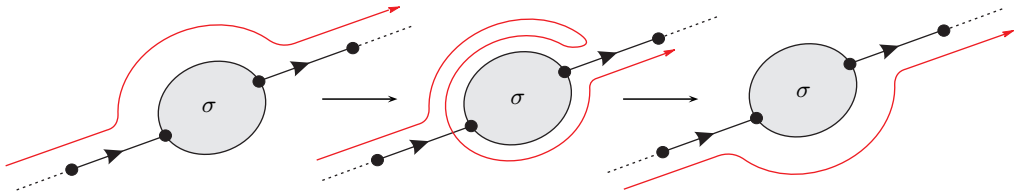

Two paths are \emph{homotopic\/} precisely when
there is a finite sequence of elementary homotopies, 
taking one to the other. 
A path homotopic to the trivial path 
is {\em homotopically trivial\/}. 
For example, two paths running
different ways around a face are homotopic as shown in Figure 
\ref{chapter1:2complexes:figure700}.

\begin{exercise}\label{chapter1:2complexes:exercise200}
Show that homotopic paths have the same start and end vertices, and thus homotopically 
trivial paths are necessarily closed. Show that homotopy is an equivalence relation on the
paths with common fixed endpoints.
\end{exercise}

We can also subdivide $2$-complexes to get homeomorphic ones, although this will play
less of a role than it does with graphs, where graphs homeomorphic to $S^1$ were essential to
the definition of $2$-complex.
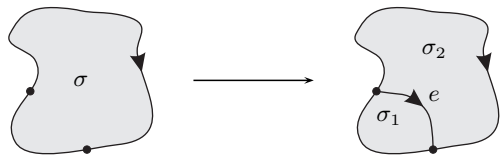
\begin{figure}
\begin{pspicture}(0,0)(12,2)
\rput(-4,0){
\rput(-1.3,0){
\rput(9,1){\BoxedEPSF{fig128a.eps scaled 1000}}
}
\rput(7.7,1){$\ss$}
\psline[linewidth=.2mm]{->}(9.2,1)(10.8,1)
\rput(3.3,0){
\rput(9,1){\BoxedEPSF{fig128b.eps scaled 1000}}
}
\rput(11.8,.5){$\ss_1^{}$}
\rput(12.4,1.4){$\ss_2^{}$}
\rput(12.4,.8){$e$}
}
\end{pspicture}
\caption{subdividing a face}\label{chapter1:2complexes:figure800}
\end{figure}
What we want to do is summarized by Figure \ref{chapter1:2complexes:figure800}: 
replace an existing face $\ss$ by two new faces $\ss_1,\ss_2$ by placing a new edge
running between vertices of $\ss$, or the reverse of this process.

\begin{exercise}\label{chapter1:2complexes:exercise300}
Formulate the definition of subdividing a face in the style of Definition 
\ref{chapter1:2complexes:definition100} by using the description of $1$-spheres given
by Lemma \ref{complexes:graphs:result300}.
\end{exercise}

\begin{figure}
\begin{pspicture}(0,0)(13,3)
\rput(-11,0){
\rput(12.5,1.5){\BoxedEPSF{fig118.eps scaled 500}}
\rput(2.5,-4.5){
\rput(11.4,7.4){$\ss_1$}
\rput(11.4,4.6){$\ss_2$}
\rput(9.7,5.45 ){${v}_1$}
\rput(10.2,6.55){${v}_2$}
\rput(10.6,5.45){${e}_1$}
\rput(9.4,6.55){${e}_2$}
}}
\rput(-6,0){
\rput(12.5,1.5){\BoxedEPSF{fig119.eps scaled 500}}
\rput(2.5,-4.5){
\rput(11.4,7.4){$\ss_1$}
\rput(11.4,4.6){$\ss_2$}
\rput(9.7,5.425 ){${v}_1$}
\rput(10.2,6.55){${v}_2$}
\rput(10.6,5.45){${e}_1$}
\rput(10.1,6.1){${e}_2$}
\rput(10.6,7){${e}_3$}
}}
\rput(-1,0){
\rput(12.5,1.5){\BoxedEPSF{fig120.eps scaled 500}}
\rput(2.5,-4.5){
\rput(11.4,7.4){$\ss_1$}
\rput(11.4,4.6){$\ss_2$}
\rput(9.8,5.45 ){${v}_1$}
\rput(10.2,6.55){${v}_2$}
\rput(10.65,5.45){${e}_1$}
\rput(10.1,6.1){${e}_2$}
\rput(10.6,7){${e}_3$}
}}
\end{pspicture}
\caption{subdividing a sphere}\label{chapter1:2complexes:figure900}
\end{figure}
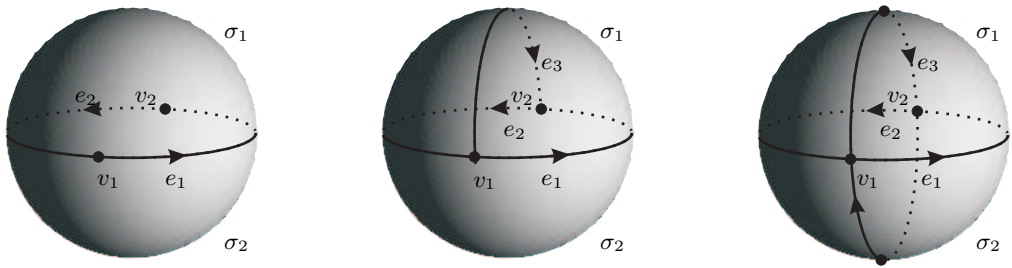

Write $X\leftrightarrow X'$ when the two complexes differ by the subdivision of an 
edge or face, so that
$X$ and $Y$ are then \emph{homeomorphic\/}, written $X\approx Y$, when there
is a finite sequence $X=X_0\leftrightarrow X_1\leftrightarrow
\cdots\leftrightarrow X_k=Y$ of subdivisions (of either type) connecting them. It is easy to
see that homeomorphism is an equivalence relation and 
a topological invariant is a well defined property of the equivalence classes.
Figure \ref{chapter1:2complexes:figure900} shows a series of subdivisions of the $2$-sphere.

\section{Quotients of $2$-complexes}\label{chapter1:quotients}

We often want to squash parts of a complex away, glue complexes 
together, factor out the action of a group, and so on.
In other words, we want to take quotients. We do this much as with graphs by 
defining an equivalence relation
on the cells of the complex and then defining a new complex whose cells are the 
equivalence classes of the
relation. It turns out that there are a number of subtleties complicating the 
exposition, arising when we want to
identify cells of different dimensions.

\subsection{Quotients in general}\label{chapter1:quotients:general}

All quotients start with an equivalence relation:

\begin{definition}[quotient relation on a $\mathbf{2}$-complex]
\label{chapter1:2complexes:definition300}
If $X$ is a $2$-complex, then a \emph{quotient relation\/} is an equivalence relation
on the vertices, edges and faces of $X$ 
such that
\begin{description}
\item[(Q1).] $\sim$ restricted to the $1$-skeleton $X^{(1)}$ is a graph quotient relation 
as in Definition \ref{chapter1:graphs:definition300}, with
quotient map $q:X^{(1)}\rightarrow X^{(1)}\quo$, and hence the induced relation
on the paths as in \S\ref{chapter1:graphs:paths}.
\item[(Q2).] if $\ss,\tau$ are faces with $\ss\sim\tau$ then
  $\ss^{-1}\sim\tau^{-1}$, and
if $\ss$ is a face with $\ss\sim\ss^{-1}$, then $[\ss]$ contains a vertex, edge 
or path, where $[\ss]$ is the equivalence class of $\ss$;
\item[(Q3).] if $\ss$ is a face with $[\ss]\subset X^2$, then there
is a $(X^{[\ss]},\aa_{[\ss]})$ with $X^{[\ss]}\approx S^1$ 
and $\aa_{[\ss]}:X^{[\ss]}\rightarrow X^{(1)}\quo$
dimension preserving,
such that for all $\tau\in[\ss]$ with $\partial\tau=(X^\tau,\aa_\tau)$,
there is an orientation preserving map $\ve=\ve(f,\tau):X^\tau\rightarrow X^{[\ss]}$
making the diagram below left commute:
$$
\begin{pspicture}(0,0)(14,2)
\rput(-.5,-.8){
\rput(-1.3,.5){
\rput(5,2){$X^\tau$}\rput(6.7,2){$X^{[\ss]}$}
\rput(4.95,0.45){$X^{(1)}$}\rput(6.7,0.45){$X^{(1)}\quo$}
\psline[linewidth=.2mm]{->}(5.3,2)(6.3,2)
\psline[linewidth=.2mm]{->}(5.3,0.45)(6.15,0.45)
\psline[linewidth=.2mm]{->}(5,1.7)(5,.7)
\psline[linewidth=.2mm]{->}(6.55,1.7)(6.55,.7)
\rput(4.7,1.2){$\aa_\tau$}\rput(6.9,1.175){$\aa_{[\ss]}$}
\rput(5.8,2.2){$\ve$}
\rput(5.8,.7){$q$}
}
\rput(3.7,.5){
\rput(5,2){$X^\tau$}\rput(6.55,2){$S$}
\rput(4.95,0.45){$X^{(1)}$}\rput(6.7,0.45){$X^{(1)}\quo$}
\psline[linewidth=.2mm]{->}(5.3,2)(6.3,2)
\psline[linewidth=.2mm]{->}(5.3,0.45)(6.15,0.45)
\psline[linewidth=.2mm]{->}(5,1.7)(5,.7)
\psline[linewidth=.2mm]{->}(6.55,1.7)(6.55,.7)
\rput(4.7,1.2){$\aa_\tau$}\rput(6.8,1.2){$\gamma$}
\rput(5.8,.7){$q$}
}
}
\end{pspicture}
$$
\item[(Q4).] if $\ss$ is a face with $[\ss]\not\subset X^2$, 
then there is a homotopically trivial path $\gamma:S\rightarrow X^{(1)}\quo$
such that for all faces $\tau\in[\ss]$ there is an
orientation preserving map $X^\tau\rightarrow S$
making the diagram above right commute.
\end{description}
\end{definition}

(Q2) ensures that when we form a quotient we have $[\ss]\not=[\ss]^{-1}$.
(Q3) says that equivalent faces have boundaries that fold up in the quotient
to give the same thing. Similarly, if a face is to be identified with a closed path
then (Q4) forces the boundary of the face to be identified with it as well. The
homotopically trivial condition is a little obscure at the moment. It's role will become
clearer in Chapter \ref{chapter2}.

\begin{definition}[quotient $\mathbf{2}$-complex]
\label{chapter1:2complexes:definition400}
If $\sim$ is a quotient relation on the $2$-complex $X$ then define the quotient
$X\quo$ as follows:
\begin{enumerate}
\item the $1$-skeleton $(X\quo)^{(1)}$ is the quotient graph
$X^{(1)}\quo$ with quotient map $q:X^{(1)}\rightarrow X^{(1)}\quo$;
\item the faces $(X\quo)^2$ are the $[\ss]\subset X^2$ for $\ss\in X^2$.
Such a face has boundary 
$\partial[\ss]=(X^{[\ss]},\aa_{[\ss]})$ as given by Definition 
\ref{chapter1:2complexes:definition300} (Q3).
\end{enumerate}
\end{definition}

\begin{proposition}
\label{chapter1:2complexes:result100}
If $\sim$ is a quotient relation then $X\quo$ 
is a $2$-complex
and the quotient map $q:X\rightarrow X\quo$ given by $q(x)=[x]$ is a map of $2$-complexes.
\end{proposition}

\begin{proof}
That the quotient is a $2$-complex is immediate from Definition
\ref{chapter1:2complexes:definition300}, and the commuting diagrams given there
are precisely what is needed for $q$ to be a map of $2$-complexes.
\qed
\end{proof}

\subsection{Quotients by a  group action}\label{chapter1:quotients:groupactions}

When a group $G$ acts on a $2$-complex $X$ we can replace $X$ by a complex on which the action
of $G$ is trivial, ie: every element of $G$ acts as the identity. We do this by factoring out
the group action, and we do this by forming a quotient. 
Recall from \S\ref{chapter1:2complexes:maps} the group $\aut(X)$ of automorphisms of the
$2$-complex $X$, and let 
a group action $G\stackrel{\varphi}{\rightarrow}\aut(X)$ be given.

Let $\sim$ be the equivalence relation on $X$ given by the orbits of the action, so that
$x\sim y$ if and only if 
$y=g(x)$ for some $g\in G$, where $x,y$ are cells of $X$, necessarily of the same dimension
(as automorphisms are dimension preserving).

\begin{proposition}[quotient by a group action]\label{chapter1:2complexes:result200}
Let $\sim$ be the equivalence
relation on the $2$-complex $X$ given by the orbits of a group action.
Then $\sim$ is a quotient relation if and only if the group action is orientation 
preserving. 
\end{proposition}

Although we can in principle consider group actions that don't preserve orientation, 
as the primary
purpose of such actions is to form quotients, we will only consider orientation preserving actions.
Compare this with Proposition \ref{chapter1:graphs:result200}, noting how the $2$-complex
structure imposes no new conditions for $\sim$ to be a quotient relation.
We write $X/G$ for the quotient complex $X\quo$.

\begin{proof}
By Exercise \ref{chapter1:graphs:exercise300}, Proposition \ref{chapter1:graphs:result200}
and basic properties of maps,
parts (Q1) and (Q2) of Definition \ref{chapter1:2complexes:definition300} 
are satisfied if and only if the $G$-action preserves orientation. 
As the elements of $G$ are dimension preserving, part (Q4) never arises.
If $[\ss]\subset X^2$ with $\ss$ a face,
let $(X^{[\ss]},\aa_{[\ss]})=(X^\ss,q\aa_\ss)$, where $q:X\rightarrow X\quo$ is the graph quotient map.
If $\tau\sim \ss$ then $\ss=g(\tau)$ for some $g\in G$, so we have an orientation
preserving isomorphism $X^\tau\rightarrow X^\ss$ with
the diagram below left commuting:
$$
\begin{pspicture}(0,0)(14,2)
\rput(-.5,-.8){
\rput(-1.3,.5){
\rput(5,2){$X^\tau$}\rput(6.7,2){$X^\ss$}
\rput(4.95,0.45){$X^{(1)}$}\rput(6.65,0.45){$X^{(1)}$}
\psline[linewidth=.2mm]{->}(5.3,2)(6.3,2)
\psline[linewidth=.2mm]{->}(5.3,0.45)(6.3,0.45)
\psline[linewidth=.2mm]{->}(5,1.7)(5,.7)
\psline[linewidth=.2mm]{->}(6.55,1.7)(6.55,.7)
\rput(4.7,1.2){$a^\tau$}\rput(6.9,1.175){$\aa_\ss$}
\rput(5.8,2.2){$\cong$}
\rput(5.8,.7){$g$}
}
\rput(3.7,.5){
\rput(5,2){$X^{(1)}$}\rput(6.6,2){$X^{(1)}$}
\psline[linewidth=.2mm]{->}(5.35,2)(6.25,2)
\psline[linewidth=.2mm]{->}(5,1.75)(5.65,.7)
\psline[linewidth=.2mm]{->}(6.55,1.75)(5.9,.7)
\rput(5.8,0.45){$X^{(1)}\quo$}
\rput(5.8,2.2){$g$}
\rput(5,1.3){$q$}\rput(6.5,1.3){$q$}
}
}
\end{pspicture}
$$
The triangular diagram on the right commutes by the nature of the quotient map: if $y=g(x)$
then $q(x)=q(y)$. Now glue the triangular diagram to the bottom of the square.
\qed
\end{proof}

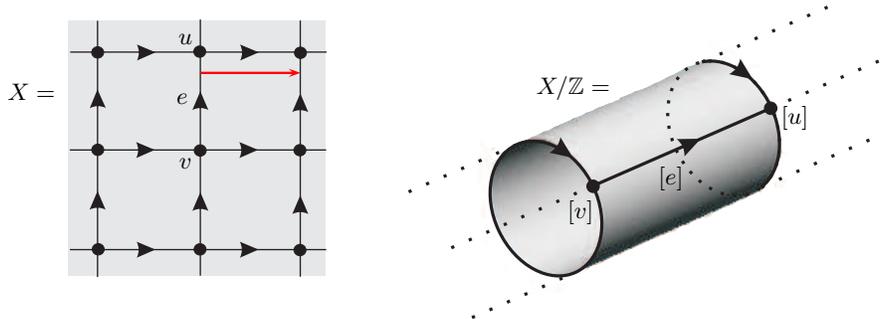
\begin{figure}
  \centering
\begin{pspicture}(0,0)(13,3.75)
\rput(0,0.25){
\rput(3,1.75){\BoxedEPSF{fig8.eps scaled 1000}}
\rput(.8,2.5){$X=$}
\rput(2.85,1.55){$v$}\rput(2.85,3.2){$u$}\rput(2.8,2.4){$e$}
}
\psline[linewidth=.3mm,linecolor=red]{->}(3.05,3)(4.375,3)
\rput(9,1.8){\BoxedEPSF{chapter1.fig900.eps scaled 600}}
\rput(8.1,1.15){$[v]$}\rput(10.95,2.4){$[u]$}\rput(9.3,1.6){$[e]$}
\rput(8,2.8){$X/\Z=$}
\end{pspicture}
\caption{let $\Z$ act on $X$ (left) with $1\in\Z$ 
acting as the translation of $X$ one step to the right, as shown by the red arrow:
The quotient $X/\Z$ is an infinite rolled up tube. If $m\in\Z$, then its effect on $X$ is to translate
$m$ steps to the right, whereas its effect on $X/\Z$ is to rotate a cell $m$ times around the tube,
bringing it back to itself. The induced $\Z$-action on $X/\Z$ is thus trivial.}
\label{chapter1:2complexes:figure1000}
\end{figure}

Consider as an example Figure \ref{chapter1:2complexes:figure1000}, where the Euclidean
plane complex is rolled into an infinite tube by the action of the integers $\Z$.

\begin{exercise}\label{chapter1:2complexes:exercise350}
In the proof of Proposition \ref{chapter1:2complexes:result200} 
we took $(X^{[\ss]},\aa_{[\ss]}s)=(X^\ss,q\aa_\ss)$.
Show that we are free to choose instead a different face from $[\ss]$: if 
$\tau\sim\ss$ and we take $(X^{[\ss]},\aa_{[\ss]})=(X^\tau,q\aa_\tau)$ instead, then this new version
of $X/G$ is isomorphic to the old one.
\end{exercise}

\subsection{Quotients by a subcomplex}\label{chapter1:quotients:subcomplex}

Now for a quotient that involves some serious squashing: if $Y\subset X$ is a
subcomplex we define a new complex where $Y$ has been compacted down to a single vertex,
extending the construction of \S\ref{chapter1:graphs:quotients} from graphs to $2$-complexes.

Define $\sim$ on $X$ to be the equivalence relation with the following equivalence classes:
(i). all the cells in $Y$ (of whatever dimension) form one class; (ii). every other class
has the form $[x]=\{x\}$. Thus, we have $x\sim y$ if and only if
either $x=y$ or both $x$ and $y$ lie in $Y$.

\begin{exercise}
\label{chapter1:2complexes:exercise360}
Let $X$ be a $1$-sphere and $Y\subset X$ a connected subcomplex. Show that the relation
just defined is a (graph) quotient relation with the quotient $X/Y$ another
$1$-sphere and the quotient map $q:X\rightarrow X/Y$ an orientation preserving map.
\end{exercise}

\begin{proposition}\label{chapter1:2complexes:result300}
The relation $\sim$ is a quotient relation.
\end{proposition}

Write $X/Y$ for the corresponding quotient, the \emph{quotient of $X$ by the subcomplex
$Y$\/}: it is what results from collapsing $Y$ to a vertex and propagating the effects of this
on the incidence of cells throughout $X$, but otherwise leaving the cells of $X$ unaffected.

\begin{proof}
The only part requiring more than a moments thought is the verification of the 
face conditions (Q3) and (Q4) in Definition 
\ref{chapter1:2complexes:definition300}.
A face $\ss\in Y$ is squashed to a vertex in $X/Y$, so
taking $S$ to be the trivial graph gives (Q4).
If $\ss\not\in Y$ then $[\ss]=\{\ss\}$ and $\aa_\ss^{-1}(Y)\subset X^\ss$ 
has connected components $T_1,\ldots,T_k$. By 
Exercise \ref{chapter1:2complexes:exercise360} we can form the successive quotients
$X^\ss/T_1,X^\ss/T_1/T_2,\ldots$, obtaining $1$-spheres at each stage.
Let $X^{[\ss]}$ be the end result of taking the $k$ quotients and 
$\ve(q,\ss):X^\ss\rightarrow X^{[\ss]}$ the composition of the quotient maps.
\qed
\end{proof}

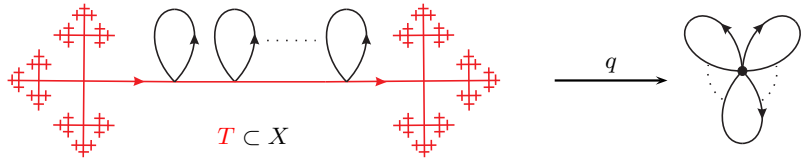
\begin{figure}
  \centering
\begin{pspicture}(0,0)(13,3)
\rput(-6,0){
\rput(10.25,1.5){\BoxedEPSF{chapter1.fig3000.eps scaled 600}}
\rput(10.25,.75){${\red T}\subset X$}
\rput(16.75,1.5){\BoxedEPSF{chapter1.fig3100.eps scaled 600}}
}
\psline[linewidth=.3mm]{->}(8.25,1.5)(9.75,1.5)
\rput(9,1.7){$q$}
\end{pspicture}  
  \caption{squashing a spanning tree down to a vertex}
\label{chapter1:2complexes:figure1100}
\end{figure}

A typical quotient by a subcomplex arises when $T\subset X$ is a spanning tree
for the $1$-skeleton and we form $X/T$ as in Figure
\ref{chapter1:2complexes:figure1100}.

As in \S\ref{chapter1:graphs:quotients} we write $X/Y_i\,(i\in I)$ for the result
of squashing each $Y_i\subset X$ to a vertex $v_i$, and $X/(\bigcup Y_i)$ for the
result of squashing $\bigcup Y_i$ to a single vertex $v$.

\subsection{Pushouts}\label{chapter1:quotients:pushouts}

The pushout is a pretty general construction which arises whenever a pair of 
complexes are glued together across a common subcomplex. 

\begin{definition}[pushout]
\label{chapter1:2complexes:definition500}
Let $X_1,X_2$ and $Y$ be $2$-complexes and 
$f_i:Y\rightarrow X_i\,(i=1,2)$ maps of $2$-complexes. 
Let $\sim$ be the equivalence
relation on $X_1\bigcup X_2$ \emph{generated\/} by 
$x\sim x'$ if and only if there is a $y\in Y$ with
$x=f_1(y)$ and $x'=f_2(y)$. 
If $\sim$ is a quotient
relation then call the quotient $X_1\bigcup X_2\quo$ the 
\emph{pushout\/} of the maps $f_i:Y\rightarrow X_i$,
and denote it by $X_1\coprod_Y Y_2$.
\end{definition}

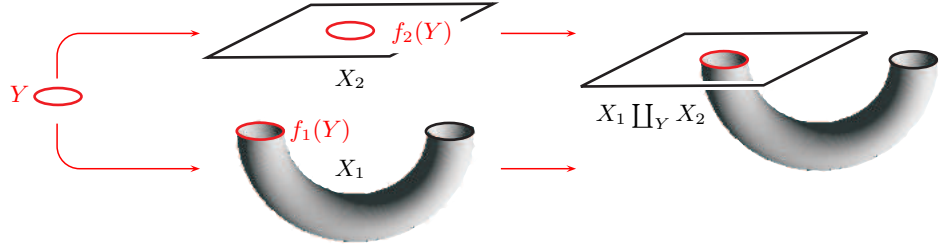
\begin{figure}
  \centering
\begin{pspicture}(0,0)(13,4)
\rput(1,2){\BoxedEPSF{chapter1.fig2300a.eps scaled 500}}
\rput(.6,2.2){${\red Y}$}
\rput(1,0){\rput(4,1){\BoxedEPSF{chapter1.fig2200.eps scaled 500}}
\rput(4,3){\BoxedEPSF{chapter1.fig2300.eps scaled 500}}
\rput(3.6,1.7){${\red f_1(Y)}$}\rput*(4.95,3){${\red f_2(Y)}$}
\rput(4,1.2){$X_1$}\rput(4,2.4){$X_2$}}
\rput(.5,0){\rput(10,2){\BoxedEPSF{chapter1.fig2100.eps scaled 500}}
\rput(8.5,1.9){$X_1\coprod_Y X_2$}}
\psline[linewidth=.2mm,linearc=.3,linecolor=red]{->}(1.1,2.4)(1.1,3)(3,3)
\psline[linewidth=.2mm,linearc=.3,linecolor=red]{->}(1.1,1.8)(1.1,1.2)(3,1.2)
\psline[linewidth=.2mm,linecolor=red]{->}(7,3)(8,3)
\psline[linewidth=.2mm,linecolor=red]{->}(7,1.2)(8,1.2)
\end{pspicture}  
\caption{typical pushout}
\label{chapter1:2complexes:figure1200}
\end{figure}

Figure \ref{chapter1:2complexes:figure1200} illustrates a typical pushout.

\begin{exercise}\label{chapter1:2complexes:exercise500}
If $\sim$ is the relation described in Definition \ref{chapter1:2complexes:definition500},
show that $x\sim x'$ iff
there are $x_0,x_1,\ldots,x_k\in X_1\bigcup X_2$ 
with $x_0=x$ and $x_k=x'$, and 
$y_1,\ldots,y_k\in Y$, such that 
$f_1(y_1)=x_0,f_2(y_1)=x_1,f_2(y_2)=x_1,f_1(y_2)=x_2,\ldots$ and so on.
\end{exercise}

Define $t_i:X_i\rightarrow X_1\coprod_YX_2\,(i=1,2)$ 
to be the composition $X_i\hookrightarrow 
X_1\bigcup X_2\rightarrow X_1\coprod X_2\quo$ of the 
inclusion of $X_i$ in the union and the quotient map.

\begin{theorem}[pushouts exist and are colimits]
\label{chapter1:2complexes:result400}
Let $Y,X_1,X_2$ be $2$-complexes, $\mathcal{O}\subset Y$,
$\mathcal{O}_i\subset X_i$ orientations and 
$f_i:Y\rightarrow X_i$ orientation and dimension preserving maps.
Then the quotient $q:X_1\bigcup X_2\rightarrow X_1\bigcup X_2\quo$,
and hence the $t_i$, are dimension preserving, and the pushout exists, with the 
diagram below left commuting.
$$
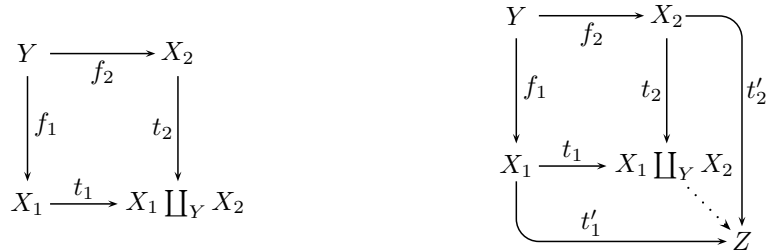
\begin{pspicture}(6,3)
\rput(0,-.5){
\rput(-1.5,1){
\rput(0,2){$Y$}
\rput(0,0){$X_1$}
\rput(2,2){$X_2$}
\rput(2.1,0){$X_1\coprod_Y X_2$}
\psline[linewidth=.2mm]{->}(0,1.7)(0,.3)
\psline[linewidth=.2mm]{->}(.3,2)(1.7,2)
\psline[linewidth=.2mm]{->}(.3,0)(1.2,0)
\psline[linewidth=.2mm]{->}(2,1.7)(2,.3)
\rput(.25,1.05){$f_1$}
\rput(1,1.75){$f_2$}
\rput(.75,.2){$t_1$}
\rput(1.8,1){$t_2$}
}
\rput(0,.5){
\rput(5,1){
\rput(0,2){$Y$}
\rput(0,0){$X_1$}
\rput(2,2){$X_2$}
\rput(2.1,0){$X_1\coprod_Y X_2$}
\psline[linewidth=.2mm]{->}(0,1.7)(0,.3)
\psline[linewidth=.2mm]{->}(.3,2)(1.7,2)
\psline[linewidth=.2mm]{->}(.3,0)(1.2,0)
\psline[linewidth=.2mm]{->}(2,1.7)(2,.3)
\rput(.25,1.05){$f_1$}
\rput(1,1.75){$f_2$}
\rput(.75,.2){$t_1$}
\rput(1.8,1){$t_2$}
}
\rput(8,0){$Z$}
\psline[linewidth=.2mm,linearc=.3]{->}(5,.8)(5,0)(7.8,0)
\psline[linewidth=.2mm,linearc=.3]{->}(7.25,3)(8,3)(8,.2)
\psline[linewidth=.3mm,linestyle=dotted]{->}(7.2,.8)(7.85,.15)
\rput(6,.25){$t'_1$}\rput(8.2,2){$t'_2$}
}
}
\end{pspicture}
$$
Moreover the pushout is universal in the sense that if $Z$,
$t'_1$, $t'_2$ are a $2$-complex and maps making 
such a square commute, then there
is a map $h:X_1\coprod_Y X_2\rightarrow Z$ making the 
diagram above right commute.
\end{theorem}

Thus the data $f_i:Y\rightarrow X_i$
forming the input to the pushout gives two sides of a
commutative square, and $X_1\coprod_Y X_2$ ``pushes out'' 
along the other two sides. Pushouts are thus examples of colimits in the 
category of $2$-complexes.

\begin{proof}
We show that under the assumptions $\sim$ is a quotient relation, with 
the result on the
$1$-skeleton given by Exercise \ref{chapter1:graphs:exercise500},
and a similar argument shows that we never have
$\ss\sim\ss^{-1}$ for a face $\ss$. Exercise \ref{chapter1:2complexes:exercise500}
and the definition of map gives that $\ss\sim\tau$ implies $\ss^{-1}\sim\tau^{-1}$.
As the $f_i$ 
are dimension preserving, all the cells in an equivalence class have the same 
dimension, leaving us with
(Q3) of Definition \ref{chapter1:2complexes:definition300} to do. 
Let
$(X^{[\ss]},\aa_{[\ss]})=(X^\ss,q\aa_\ss)$, with $q\aa_\ss$ dimension preserving.
Let $\tau\in[\ss]$ and suppose we are in the special
case 
$\ss=f_1(\rho), \tau=f_2(\rho)$ for some face $\rho$ in $Y$. Then we get a diagram
$$
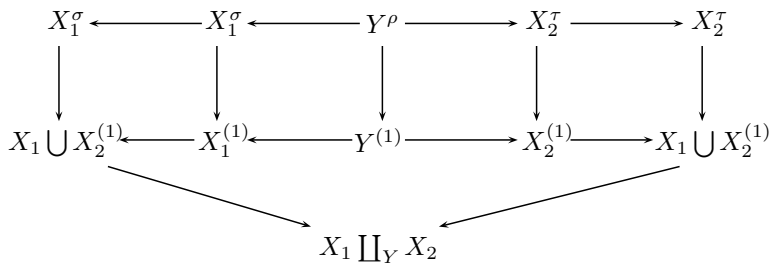
\begin{pspicture}(0,0)(13,3.5)
\rput(0,-.75){
\rput(1.25,1.8){
\rput(-2.1,0){
\rput(2.9,2){$X_1^\ss$}\rput(2.9,0.45){$X_1\bigcup X_2^{(1)}$}
\psline[linewidth=.2mm]{->}(2.8,1.7)(2.8,.7)
\psline[linewidth=.2mm]{->}(4.7,2)(3.2,2)
\psline[linewidth=.2mm]{->}(4.6,0.45)(3.6,0.45)
}
\rput(2.9,2){$X_1^\ss$}\rput(2.9,0.45){$X_1^{(1)}$}
\psline[linewidth=.2mm]{->}(2.8,1.7)(2.8,.7)
\psline[linewidth=.2mm]{->}(4.7,2)(3.2,2)
\psline[linewidth=.2mm]{->}(4.6,0.45)(3.2,0.45)
\rput(5,2){$Y^\rho$}\rput(4.95,0.45){$Y^{(1)}$}
\psline[linewidth=.2mm]{->}(5,1.7)(5,.7)
\rput(7.15,2){$X_2^\tau$}\rput(7.2,0.45){$X_2^{(1)}$}
\psline[linewidth=.2mm]{->}(7.05,1.7)(7.05,.7)
\psline[linewidth=.2mm]{->}(5.3,2)(6.8,2)
\psline[linewidth=.2mm]{->}(5.3,0.45)(6.8,0.45)
\rput(2.2,0){
\rput(7.15,2){$X_2^\tau$}\rput(7.2,0.45){$X_1\bigcup X_2^{(1)}$}
\psline[linewidth=.2mm]{->}(7.05,1.7)(7.05,.7)
\psline[linewidth=.2mm]{->}(5.3,2)(6.8,2)
\psline[linewidth=.2mm]{->}(5.3,0.45)(6.4,0.45)
}
}
\rput(6.2,.8){$X_1\coprod_YX_2$}
\psline[linewidth=.2mm]{->}(2.6,1.9)(5.3,1.1)
\psline[linewidth=.2mm]{->}(10.2,1.9)(7,1.1)
}
\end{pspicture}
$$
with the four squares commuting via the maps $f_i$ and the inclusions
$X_i\hookrightarrow X_1\bigcup X_2$. The $\ve(f_1,\ss)$, $\ve(f_2,\ss)$ 
are orientation preserving isomorphisms and the other two maps along the top 
are the identities. The diagram glued to the bottom commutes by the 
definition of $q$. The maps along the top (and their inverses) compose to give an isomorphism
$X_2^\tau\rightarrow X_1^\ss$, and this, together with the outside
circuit of the diagram, give condition (Q3).

When $\ss\sim\tau$ in general, we have $\ss=\ss_0=f_1(\rho_1), \ss_1=f_2(\rho_1),
\ldots,\ss_{k-1}=f_1(\rho_k),\tau=\ss_k=f_2(\rho_k)$, and the requirements for a quotient
relation can be verified by repeatedly applying the process of the previous paragraph.
In particular, for $\tau\in[\ss]$ the map $X^\tau\rightarrow X^{[\ss]}$ is an orientation
preserving isomorphism, so that the quotient map is dimension preserving.

If $y\in Y$ is a cell then its images under the
$t_if_i:Y\rightarrow X_i\hookrightarrow X_1\coprod_Y X_2\,(i=1,2)$ are equivalent by the definition of 
$\sim$, and so the square commutes. 

Suppose now that $Z,t'_1,t'_2$ are as in the statement of the Theorem.
After a moments thought it is clear what the map $h:X_1\coprod_Y X_2\rightarrow Z$
should be: every cell of the pushout has the form $[x]$ for some 
$x\in X_i$, so define $h[x]=t_i'(x)$. 
We leave it to the reader to show that this is well defined and gives a map of 
$2$-complexes. That (M1) is satisfied is very similar to the argument above.
\qed
\end{proof}

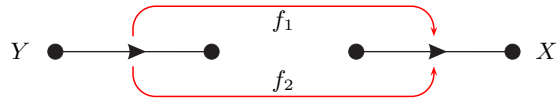
\begin{figure}
  \centering
\begin{pspicture}(0,0)(13,2)
\rput(4,1){\BoxedEPSF{chapter1.fig1500.eps scaled 500}}
\rput(8,1){\BoxedEPSF{chapter1.fig1500.eps scaled 500}}
\rput(2.5,1){$Y$}\rput(9.5,1){$X$}
\rput(6,1.4){$f_1$}\rput(6,.6){$f_2$}
\psline[linewidth=.2mm,linearc=.3,linecolor=red]{->}(4,1.2)(4,1.6)(8,1.6)(8,1.2)
\psline[linewidth=.2mm,linearc=.3,linecolor=red]{->}(4,.8)(4,.4)(8,.4)(8,.8)
\end{pspicture}
\caption{data $f_i:Y\rightarrow X_i$ for which the pushout does not exist.
Here $X_1=X_2=X$ and 
$Y$ are all graphs with a single edge joining two vertices and the map $f_1$
(respectively $f_2$)  sends the
edge of $Y$ to the edge $e$ (resp. $e^{-1}$). 
The resulting equivalence relation on $X$ has $e\sim e^{-1}$ but
$[e]\cap X^{0}=\varnothing$, and so is not a
quotient relation.}
\label{chapter1:2complexes:figure1300}
\end{figure}

\begin{exercise}
\label{chapter1:2complexes:exercise525}
With the conditions of Theorem \ref{chapter1:2complexes:result400}, 
show that if $X_1,X_2$ are connected then the
pushout is connected.
\end{exercise}

\begin{exercise}[pointed pushout]
\label{chapter1:2complexes:exercise550}
Formulate a pointed version of Theorem \ref{chapter1:2complexes:result400}
with all the complexes and 
maps in sight pointed.
\end{exercise}

Figure \ref{chapter1:2complexes:figure1300} is a simple example for which the pushout
doesn't exist, ie: the relation $\sim$ in the pushout construction is 
not a quotient relation.

Figure \ref{chapter1:2complexes:figure1350}
is the {\em Stallings fold\/}:
the graph $Y$ is a single edge joining two vertices.
Another example, shown in Figure \ref {chapter1:2complexes:figure1375},
is the \emph{wedge\/} of a pair of complexes: $Y$ is now the trivial complex 
consisting of just a single vertex.

\begin{figure}
  \centering
\begin{pspicture}(0,0)(13,3)
\rput(3,1.5){\BoxedEPSF{chapter1.fig1500.eps scaled 500}}
\rput(6,1.5){\BoxedEPSF{chapter1.fig6200.eps scaled 500}}
\rput(10,1.5){\BoxedEPSF{chapter1.fig6300.eps scaled 500}}
\psline[linewidth=.2mm,linearc=.3,linecolor=red]{->}(3,1.7)(3,2)(5.8,2)
\psline[linewidth=.2mm,linearc=.3,linecolor=red]{->}(3,1.3)(3,1)(5.8,1)
\psline[linewidth=.2mm,linecolor=red]{->}(7.4,1.5)(8.2,1.5)
\rput(1.5,1.5){$Y$}\rput(7.5,2.5){$X$}
\end{pspicture}  
  \caption{Stallings fold: here
$Y$ is a single edge $e$ joining two distinct vertices, and $X_1=X_2=X$,
$f_i:Y\rightarrow X, (i=1,2)$ 
with $f_1(s(e))=f_2(s(e))$ and $f_1(e)\not=f_2(e)^{-1}$}
  \label{chapter1:2complexes:figure1350}
\end{figure}
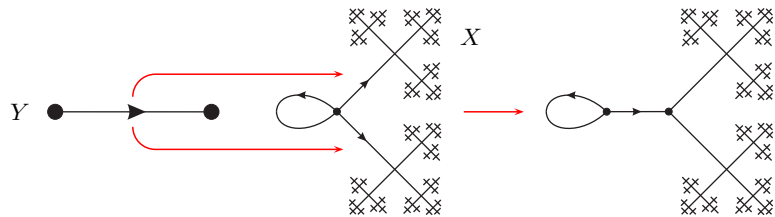

\begin{figure}
  \centering
\begin{pspicture}(0,0)(13,4)
\rput(2,2){${\red\bullet}$}
\rput(1.75,2){${\red Y}$}
\psline[linewidth=.2mm,linearc=.3,linecolor=red]{->}(2,2.2)(2,3)(3.8,3)
\psline[linewidth=.2mm,linearc=.3,linecolor=red]{->}(2,1.8)(2,1)(3.8,1)
\psline[linewidth=.2mm,linecolor=red]{->}(6.3,3)(7.7,3)
\psline[linewidth=.2mm,linecolor=red]{->}(6.3,1)(7.7,1)
\rput(0,0){
\rput(5,3){\BoxedEPSF{chapter1.fig1900.eps scaled 250}}
\rput(6,3.5){$X_1$}
}
\rput(0,-2){
\rput(5,3){\BoxedEPSF{chapter1.fig1900.eps scaled 250}}
\rput(6,2.5){$X_2$}
}
\rput(4.8,-1){
\rput(5,3){\BoxedEPSF{chapter1.fig2000.eps scaled 250}}
\rput(5,1.8){$X_1\coprod_Y X_2$}
}
\end{pspicture}
  \caption{wedge of spheres.}
  \label{chapter1:2complexes:figure1375}
\end{figure}
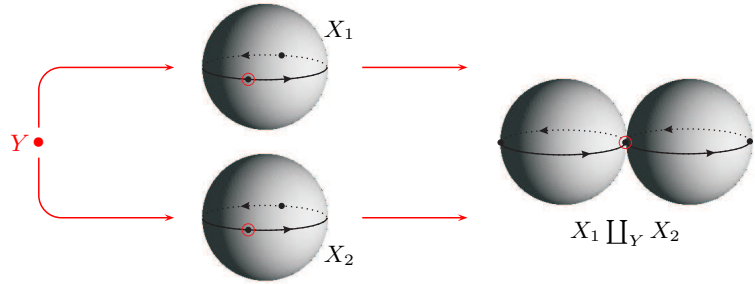

\begin{figure}
  \centering
\begin{pspicture}(0,0)(13,4)
\rput(10,2){\BoxedEPSF{chapter1.fig2100.eps scaled 500}}
\rput(4,1){\BoxedEPSF{chapter1.fig2200.eps scaled 500}}
\psline[linewidth=.2mm,linecolor=red]{->}(2.8,2)(2.8,3.1)
\rput(2.8,3){\BoxedEPSF{chapter1.fig2300.eps scaled 500}}
\rput(2,1.7){$f_1(Y)$}\rput(4,1.2){$X_1$}\rput(4,2.75){$X_2$}
\rput(8.5,1.9){$X_1\coprod_Y X_2$}
\rput(3,2.3){$f_2$}
\end{pspicture}
  \caption{glueing complexes via an attaching map is a pushout with $f_1$ an inclusion.}
  \label{chapter1:2complexes:figure1600}
\end{figure}
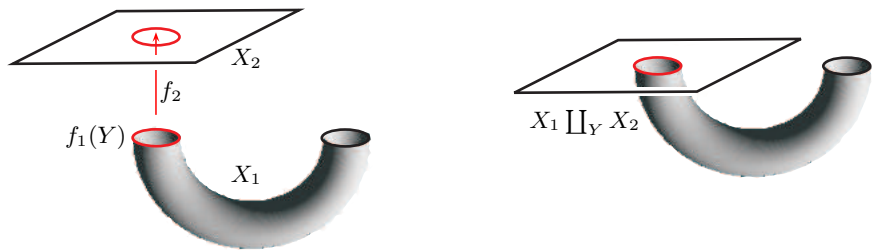

An important example arises when one of the maps $f_i$ is just an inclusion, 
so that the initial data consists of two complexes $X_1$ and $X_2$, and a map
$f$ from a \emph{subcomplex\/} of $X_1$ to $X_2$. The pushout (when it
exists) is the result of glueing $X_1$ and $X_2$ together via the attaching map $f$ 
as in Figure \ref{chapter1:2complexes:figure1600}.

\begin{figure}
  \centering
\begin{pspicture}(0,0)(12.5,4.5)
\rput(2,3){\BoxedEPSF{chapter1.fig1800a.eps scaled 750}}
\rput(2,3){$X_1$}
\rput(6.25,3){\BoxedEPSF{chapter1.fig1700.eps scaled 750}}
\rput(6.25,3){$Y$}
\rput(10.5,3){\BoxedEPSF{chapter1.fig1300a.eps scaled 750}}
\rput(10.5,3){$X_2$}
\rput(2,.5){\BoxedEPSF{chapter1.fig1000.eps scaled 350}}
\rput(6.25,.5){\BoxedEPSF{chapter1.fig1000.eps scaled 350}}
\rput(10.5,.5){\BoxedEPSF{chapter1.fig1000.eps scaled 350}}
\rput(1.1,.5){$X_1^\ss$}\rput(11.4,.5){$X_2^\ss$}
\psline[linewidth=.2mm]{->}(2,1.1)(2,1.6)
\psline[linewidth=.2mm]{->}(6.25,1.1)(6.25,1.6)
\psline[linewidth=.2mm]{->}(10.5,1.1)(10.5,1.6)
\psline[linewidth=.2mm,linecolor=blue]{->}(4.8,3)(3.4,3)
\psline[linewidth=.2mm,linecolor=blue]{->}(7.7,3)(9.1,3)
\rput(4.1,3.2){$g_1$}\rput(8.4,3.2){$g_2$}
\psline[linewidth=.2mm]{->}(5.5,.5)(2.7,.5)
\psline[linewidth=.2mm]{->}(7,.5)(9.8,.5)
\end{pspicture}  
\caption{initial data for a typical pushout of Chapter \ref{chapter3}: 
$Y$ has six vertices, edges and faces, $X_1$ has
three of everything and $X_2$ has two of everything. The boundaries of 
all faces are hexagonal, but we've only shown one in each case.
The attaching maps wrap around the $1$-skeletons as shown.}
\label{chapter1:2complexes:figure1400}
\end{figure}
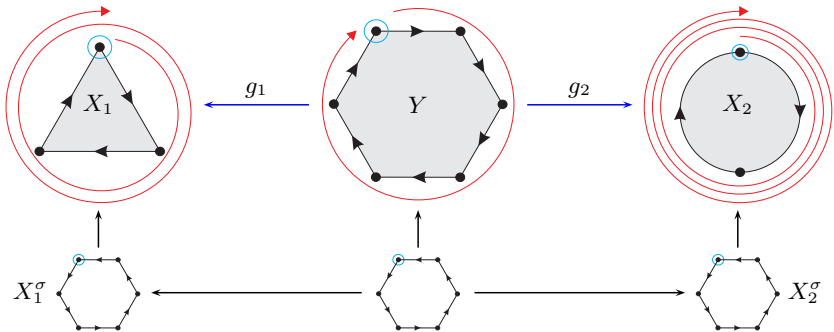

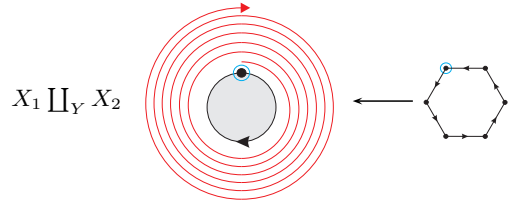
\begin{figure}
  \centering
\begin{pspicture}(0,0)(13,3)
\rput(5.5,1.5){\BoxedEPSF{chapter1.fig1400a.eps scaled 750}}
\rput(3.2,1.5){$X_1\coprod_Y X_2$}
\psline[linewidth=.2mm]{->}(7.8,1.5)(7,1.5)
\rput(8.5,1.5){\BoxedEPSF{chapter1.fig1000.eps scaled 350}}
\end{pspicture}  
\caption{pushout resulting from the set-up in Figure \ref{chapter1:2complexes:figure1400}. There
is a single vertex, edge and face in the quotient, and the face has a hexagonal boundary 
with the attaching map wrapping it around the $1$-skeleton six times as shown.}
\label{chapter1:2complexes:figure1500}
\end{figure}

Pushouts will really prove
their mettle in Chapters \ref{chapter3}-\ref{chapter4}, where the $f_i$ will be covering maps.
Figure \ref{chapter1:2complexes:figure1400} illustrates the kind
of initial set-up we will have, and Figure \ref{chapter1:2complexes:figure1500} 
the resulting pushout.

\section{Pullbacks and Higman composition}\label{chpt1:pullbacks}

We now come to a pair of constructions which both start with roughly the same kind of 
data: a complex $Y$, a (finite) family of complexes $X_i$, and a family of maps
$f_i:X_i\rightarrow Y$. The first of these, the pullback, is dual to the pushout:
it is a categorical limit.
It is what we get if we reverse
the directions of all the maps in the pushout. Pullbacks, like pushouts,
will play a crucial role in the theory of coverings in Chapters
\ref{chapter3}-\ref{chapter4}: they will act like a kind of ``union'' and
pushouts like a kind of ``intersection''. In \S\ref{chapter3:lattices:lattice}
we will be able to be much more precise
about what we mean by this.

The other construction, Higman composition, is less well known and can 
be performed only in very special circumstances. Nevertheless, when possible it will 
prove extremely powerful, and this makes its inclusion more than worthwhile.

\subsection{Pullbacks}\label{chapter1:pullhig:pullbacks}

\begin{figure}
  \centering
\begin{pspicture}(0,0)(13,4)
\psline[linewidth=.2mm,linecolor=red]{->}(5.5,1)(7,1)
\psline[linewidth=.2mm,linecolor=red]{->}(9,2.75)(9,1.25)
\psline[linewidth=.2mm,linecolor=red]{->}(3.5,2.75)(3.5,1.25)
\psline[linewidth=.2mm,linecolor=red]{->}(5.5,3)(7,3)
\rput(-2.5,2.25){
\rput(6,1){\BoxedEPSF{fig5.eps scaled 1000}}
\rput(4.8,.4){$u_1\times u_2$}
\rput(7.2,1.6){$v_1\times v_2$}
\rput(6,.75){$e_1\times e_2$}
\rput(4.8,1.4){$X_1\prod_Y X_2$}
}
\rput(-2.5,-.25){
\rput(6,1){\BoxedEPSF{fig5.eps scaled 1000}}
\rput(4.8,.4){$u_1$}
\rput(7.4,1.6){$v_1$}
\rput(5.8,1.2){$e_1$}
\rput(7.2,.75){$X_1$}
}
\rput(3.1,2.25){
\rput(6,1){\BoxedEPSF{fig5.eps scaled 1000}}
\rput(4.8,.4){$u_2$}
\rput(7.2,1.6){$v_2$}
\rput(5.8,1.2){$e_2$}
\rput(7.2,.6){$X_2$}
}
\rput(3.1,-.25){
\rput(6,1){\BoxedEPSF{fig5.eps scaled 1000}}
\rput(4.8,.4){$u$}
\rput(7.2,1.6){$v$}
\rput(5.8,1.2){$e$}
\rput(7.2,.6){$Y$}
}
\end{pspicture}  
  \caption{Construction of the pullback for graphs: the vertices $u_i$ and $v_i$ map via the
$f_i$ to $u$ and $v$, and the edges $e_i$ map via the $f_i$ to $e$. In the pullback we get 
vertices $u_1\times u_2$, $v_1\times v_2$ joined by an edge $e_1\times e_2$.}
    \label{chapter1:2complexes:figure1700}
\end{figure}
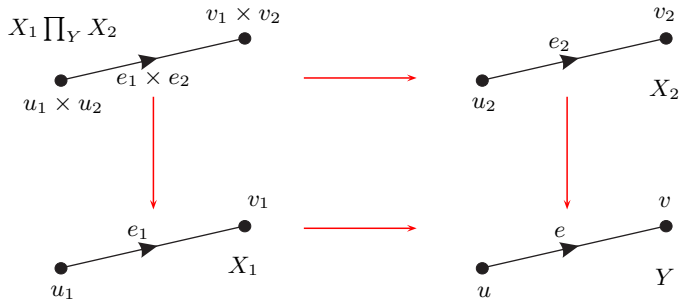

It is easier to do graphs and then extend to $2$-complexes:

\begin{definition}[pullback of graphs]\label{chapter1:2complexes:definition600}
Let $Y$ and $X_1,X_2$ be graphs and 
$f_i:X_i\rightarrow Y$ 
maps of graphs.
The {\em pullback\/} $X_1\prod_Y X_2$  has vertices
(respectively edges) the $x_1\times x_2$ for $x_i\in X_i^0$ (resp.
$x_i\in X_i^1$), such that $f_1(x_1)=f_2(x_2)$. The incidence maps are given by
$s(e_1\times e_2)=s(e_1)\times s(e_2)$ and $(e_1\times e_2)^{-1}
=e_1^{-1}\times e_2^{-1}$. See Figure \ref{chapter1:2complexes:figure1700}.
\end{definition}

We now set up the pullback of $2$-complexes by seeing how the boundaries of faces
in the $X_1,X_2$ behave when we ``pullback'' their boundaries. Suppose 
the $f_i$ are dimension preserving, and that
$\ss$ is a face of $Y$ and $\ss_i\,(i=1,2)$ faces of the $X_i$ mapping to $\ss$ via the $f_i$. We
get a by now familar commuting diagram:
$$
\begin{pspicture}(0,0)(13,2)
\rput(1.25,-.2){
\rput(2.9,2){$X_1^{\ss_1}$}\rput(5,2){$Y^\ss$}\rput(7.2,2){$X_2^{\ss_2}$}
\rput(2.8,0.45){$X_1^{(1)}$}
\rput(4.95,0.45){$Y^{(1)}$}\rput(7.2,0.45){$X_2^{(2)}$}
\psline[linewidth=.2mm]{<-}(4.7,2)(3.2,2)
\psline[linewidth=.2mm]{<-}(4.6,0.45)(3.2,0.45)
\psline[linewidth=.2mm]{<-}(5.3,2)(6.8,2)
\psline[linewidth=.2mm]{<-}(5.3,0.45)(6.8,0.45)
\psline[linewidth=.2mm]{->}(2.8,1.7)(2.8,.7)
\psline[linewidth=.2mm]{->}(5,1.7)(5,.7)
\psline[linewidth=.2mm]{->}(7.05,1.7)(7.05,.7)
\rput(4.7,1.2){$\aa_{\ss}$}
\rput(2.5,1.2){$\aa_{\ss_1}$}\rput(7.4,1.2){$\aa_{\ss_2}$}
}
\end{pspicture}
$$
The $\ve_i=\ve(f_i,\ss_i):X^{\ss_i}_i\rightarrow Y^\ss$ are orientation preserving 
isomorphisms as the $f_i$ preserve dimension. 
Let $\mathcal{O}=\{e_1,e_2,\ldots,e_n\}$ be the standard orientation of $Y^\ss$, so that
the edges of the $X_i^{\ss_i}$ can 
be labelled $e_{i1},e_{i2},\ldots,e_{ik}$, with $e_j=\ve_i(e_{ij})$, and
the closed path $\aa_{\ss_i}(e_{i1}),\aa_{\ss_i}(e_{i2}),\ldots,\aa_{\ss_i}(e_{ik})$ 
a boundary label for $\ss_i$ with $\aa_{\ss_i}(e_{ij})$ mapping via 
$f_i$ to $\aa_\ss(e_j)$ in $Y$.

The upshot is that the pullback contains a path of edges
\begin{equation}
  \label{chapter1:2complexes:equation200}
\aa_{\ss_1}(e_{11})\times\aa_{\ss_2}(e_{21}),\aa_{\ss_1}(e_{12})\times\aa_{\ss_2}(e_{22}),\ldots,
\aa_{\ss_1}(e_{1k})\times\aa_{\ss_2}(e_{2k}),  
\end{equation}
and since $s\aa_{\ss_1}(e_{11})=t\aa_{\ss_1}(e_{1k})$ and $s\aa_{\ss_2}(e_{21})=t\aa_{\ss_2}(e_{2k})$,
this path is closed. The idea is to ``sew a face'' into the $1$-skeleton
having boundary this closed path, by taking $Y^\ss$ and attaching
map $\aa_{\ss_1}\ve_1^{-1}\times\aa_{\ss_2}\ve_2^{-1}$ that sends
$e_j$ to $\aa_{\ss_1}(e_{1j})\times\aa_{\ss_2}(e_{2j})$.

\begin{definition}[pullback of $2$-complexes]
\label{chapter1:2complexes:definition700}
Let $X_1,X_2$ and $Y$ be $2$-complexes and 
$f_i:X_i\rightarrow Y$ dimension preserving
maps.
The {\em pullback\/} $X_1\prod_Y X_2$  has $k$-dimensional cells 
the $x_1\times x_2$, for $x_i\in X_i^k$, such that $f_1(x_1)=f_2(x_2)$. The incidence maps are
given by $s(e_1\times e_2)=s(e_1)\times s(e_2)$, $(e_1\times e_2)^{-1}
=e_1^{-1}\times e_2^{-1}$ and 
$$
\partial(\ss_1\times\ss_2)=(Y^\ss,\aa_{\ss_1}\ve_1^{-1}\times\aa_{\ss_2}\ve_2^{-1}),
$$
where $f_1(\ss_1)=\ss=f_2(\ss_2)$.
\end{definition}

\begin{figure}
  \centering
\begin{pspicture}(0,0)(13,2.2)
\rput(-.5,0){
\rput(3,1){\BoxedEPSF{fig128.eps scaled 1000}}
\rput(8,1){\BoxedEPSF{fig126.eps scaled 1000}}
\rput(11.75,1){\BoxedEPSF{fig127.eps scaled 1000}}
\rput(3.1,1){$\ss_1\times \ss_2$}
\rput(7.7,1.3){$\ss_1$}
\rput(11.85,.5){$\ss_2$}
\rput(6,1){$\Leftarrow$}\rput(10,1){and}
\rput(8.9,1.5){$\in X_1$}\rput(12.8,.5){$\in X_2$}
\rput(4.6,.65){$v_1\times v_2$}\rput(7.7,.45){$v_1$}\rput(11.15,1){$v_2$}
\rput(4.5,2.2){$\partial(\ss_1\times\ss_2)=\partial\ss_1\times\partial\ss_2$}
\rput(9.3,0){$\partial\ss_1$}\rput(12.75,1.2){$\partial\ss_2$}
}
\end{pspicture}
  \caption{pullback of faces whenever $f_1(\ss_1)=f_2(\ss_2)$.}
  \label{chapter1:2complexes:figure1800}
\end{figure}
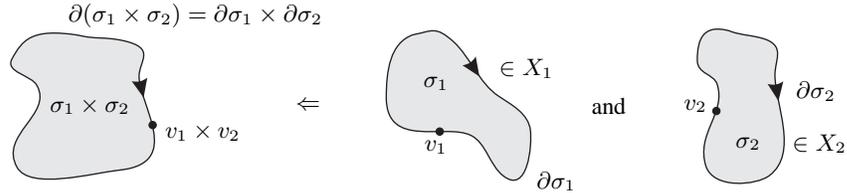

Notice that unlike the pushout, there is no question of whether 
the pullback exists or not.
Each cell of the pullback has the form $x_1\times x_2$ with the $x_i\in X_i$,
so define
$$
t_i:X_1\prod_YX_2\rightarrow X_i\,(i=1,2),
$$
by $t_i(x_1\times x_2)=x_i$. 
For $\ss_1\times\ss_2$ a face with $t_i(\ss_1\times\ss_2)=\ss_i$, we define
$\ve(t_i,\ss_1\times\ss_2)=\ve(f_i,\ss_i)^{-1}$.
We leave it as an exercise 
to see that the $t_i$ are dimension preserving maps of $2$-complexes.

\begin{theorem}[pullbacks are limits]
\label{chapter1:2complexes:result600}
The diagram below left
commutes,
$$
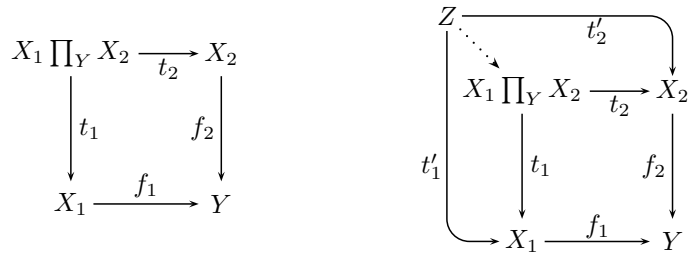
\begin{pspicture}(12.5,3)
\rput(0,-1){
\rput(2.5,1.5){
\rput(0,2){$X_1\prod_Y X_2$}
\rput(0,0){$X_1$}
\rput(2,2){$X_2$}
\rput(2,0){$Y$}
\psline[linewidth=.2mm]{->}(0,1.7)(0,.3)
\psline[linewidth=.2mm]{->}(.9,2)(1.7,2)
\psline[linewidth=.2mm]{->}(.3,0)(1.7,0)
\psline[linewidth=.2mm]{->}(2,1.7)(2,.3)
\rput(.25,1){$t_1$}
\rput(1.3,1.8){$t_2$}
\rput(1,.25){$f_1$}
\rput(1.75,1){$f_2$}
}
\rput(3.5,0){
\rput(5,1){
\rput(0,2){$X_1\prod_Y X_2$}
\rput(0,0){$X_1$}
\rput(2,2){$X_2$}
\rput(2,0){$Y$}
\psline[linewidth=.2mm]{->}(0,1.7)(0,.3)
\psline[linewidth=.2mm]{->}(.9,2)(1.7,2)
\psline[linewidth=.2mm]{->}(.3,0)(1.7,0)
\psline[linewidth=.2mm]{->}(2,1.7)(2,.3)
\rput(.25,1){$t_1$}
\rput(1.3,1.8){$t_2$}
\rput(1,.2){$f_1$}
\rput(1.75,1){$f_2$}
}
\rput(4,4){$Z$}
\psline[linewidth=.2mm,linearc=.3]{->}(4,3.8)(4,1)(4.7,1)
\psline[linewidth=.2mm,linearc=.3]{->}(4.2,4)(7,4)(7,3.2)
\rput(3.8,2){$t'_1$}
\rput(6,3.8){$t'_2$}
\psline[linewidth=.3mm,linestyle=dotted]{->}(4.1,3.9)(4.7,3.3)
}
}
\end{pspicture}
$$
Moreover, the pullback is universal in the sense that if $Z$,
$t'_1,t'_2$ are a $2$-complex and maps making such a square commute,
then
there is a map $h:Z\rightarrow X_1\prod_Y X_2$ making the
diagram above right commute.
\end{theorem}

\begin{proof}
We have $f_1t_1(x_1\times x_2)=f_2t_2(x_1\times x_2)$, and writing
$X=X_1\prod_YX_2$, the face isomorphisms $X^{\ss_1\times\ss_2}\rightarrow Y^\ss$ 
are identical for both $f_it_i$, and so the square commutes.
If $z$ is a cell of $Z$, then the commuting
of the large square of the righthand diagram that $f_1t_1'(z)=f_2t_2'(z)$,
so $t_1'(z)\times t_2'(z)$ is a cell of the pullback. 
Define $h:Z\rightarrow X_1\prod_Y X_2$ to be $z\mapsto t_1'(z)\times t_2'(z)$, and for the
isomorphism $Z^\ss\rightarrow X^{t'_1(\ss)\times t'_2(\ss)}=Y^{f_it'_i(\ss)}$, take the composition
$$
Z^\ss
\stackrel{}{\longrightarrow}X_i^{t_i'(\ss)}\stackrel{}{\longrightarrow}Y^{f_it'_i(\ss)}.
$$ 
We leave it to the reader to check
that this is a map of $2$-complexes with the required properties. 
\qed  
\end{proof}

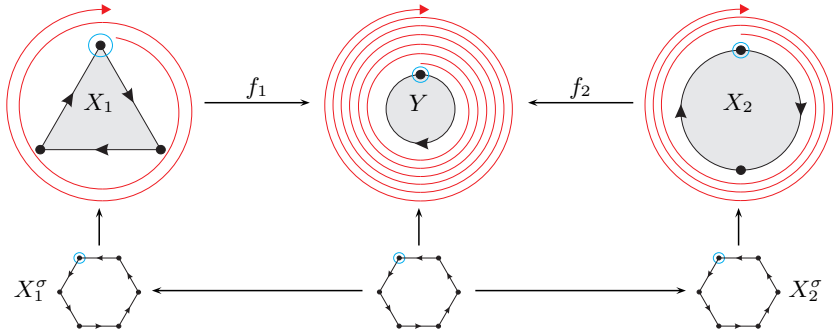
\begin{figure}
  \centering
\begin{pspicture}(0,0)(12.5,4)
\rput(2,3){\BoxedEPSF{chapter1.fig1800a.eps scaled 750}}
\rput(2,3){$X_1$}
\rput(6.25,3){\BoxedEPSF{chapter1.fig1400a.eps scaled 750}}
\rput(6.25,3){$Y$}
\rput(10.5,3){\BoxedEPSF{chapter1.fig1300a.eps scaled 750}}
\rput(10.5,3){$X_2$}
\rput(2,.5){\BoxedEPSF{chapter1.fig1000.eps scaled 350}}
\rput(6.25,.5){\BoxedEPSF{chapter1.fig1000.eps scaled 350}}
\rput(10.5,.5){\BoxedEPSF{chapter1.fig1000.eps scaled 350}}
\rput(1.1,.5){$X_1^\ss$}\rput(11.4,.5){$X_2^\ss$}
\psline[linewidth=.2mm]{->}(2,1.1)(2,1.6)
\psline[linewidth=.2mm]{->}(6.25,1.1)(6.25,1.6)
\psline[linewidth=.2mm]{->}(10.5,1.1)(10.5,1.6)
\psline[linewidth=.2mm]{<-}(4.8,3)(3.4,3)
\psline[linewidth=.2mm]{<-}(7.7,3)(9.1,3)
\rput(4.1,3.2){$f_1$}\rput(8.4,3.2){$f_2$}
\psline[linewidth=.2mm]{->}(5.5,.5)(2.7,.5)
\psline[linewidth=.2mm]{->}(7,.5)(9.8,.5)
\end{pspicture}  
\caption{initial data for a typical pullback of Chapters \ref{chapter3}-\ref{chapter4}: 
$Y$ has one vertex, edge and face, $X_1$ has
three of everything and $X_2$ has two of everything. The boundaries of 
all faces are hexagonal with just one shown. 
The attaching maps wrap around the $1$-skeletons as shown.}
\label{chapter1:2complexes:figure1900}
\end{figure}

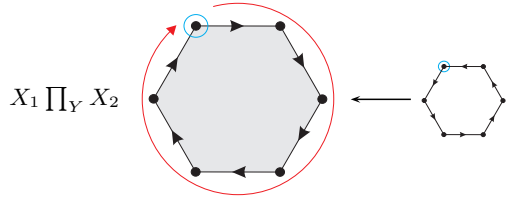
\begin{figure}
  \centering
\begin{pspicture}(0,0)(13,3)
\rput(5.5,1.5){\BoxedEPSF{chapter1.fig1700.eps scaled 750}}
\rput(3.2,1.5){$X_1\prod_Y X_2$}
\psline[linewidth=.2mm]{->}(7.8,1.5)(7,1.5)
\rput(8.5,1.5){\BoxedEPSF{chapter1.fig1000.eps scaled 350}}
\end{pspicture}  
\caption{Pullback resulting from the set-up in Figure \ref{chapter1:2complexes:figure1900}.}
\label{chapter1:2complexes:figure2000}
\end{figure}

Although pullbacks always exist, unlike pushouts, they are not necessarily connected,
unlike pushouts. This simple fact, even for graphs, has surprisingly far-reaching implications
as we will see when we study subgroups of free groups. There is a pointed version of all the
above which goes some way to fixing this:

\begin{exercise}[pointed pullbacks]
\label{chapter1:2complexes:exercise650}
Suppose the $f_i:(X_i)_{x_i}\rightarrow Y_y$ are maps of pointed
complexes. Then $x=x_1\times x_2$ is a vertex of the pullback.
Write $(X_1\prod_Y X_2)_x$ for the connected component
containing $x_1\times x_2$. Show that we have a Proposition
analogous to Proposition \ref{chapter1:2complexes:result600}
for the pointed pullback,
with all complexes and maps in sight pointed and connected.
\end{exercise}

The duality between pullbacks and pushouts can be seen by running
the example of \S\ref{chapter1:quotients:pushouts}
backwards, using the pullback to get back to where we started.
In Figure \ref{chapter1:2complexes:figure1900} we have the same two complexes
$X_1,X_2$ as in the pushout example, but this time $Y$ is the end result of that 
example and the maps $f_i$ are the 
$t_i$ from Proposition \ref{chapter1:2complexes:result400}. The resulting pullback, shown in 
Figure \ref{chapter1:2complexes:figure2000} is the starting point of the pushout example,
and the maps $t_i$ given by Proposition \ref{chapter1:2complexes:result600} are the starting
$f_i$ from there. The number of faces proliferates rather than 
declines as it does in the pushout: if $\ss_1,\ss_2,\ss_3$ are the faces of $X_1$ and 
$\tau_1,\tau_2$ of $X_2$, then we get six faces $\ss_i\times\tau_j\,(i=1,2,3;j=1,2)$ in the pullback.

\subsection{Higman composition}\label{chapter1:pullhig:higman}

The second of our two constructions starts with a (finite) family
of maps $f_i:Y_i\rightarrow X$ with the $Y_i$ disjoint
complexes. If there is a certain special configuration of edges in the
$Y_i$, then they can be threaded together into one large complex
$Y$ and a map $f:Y\rightarrow X$. Much as with pullbacks, the pay-off will not be evident
until Chapter \ref{chapter3}, when we show that if the $f_i$ are
covering maps then so is the new map $f$. 

\begin{definition}[handle configuration]
\label{chapter1:2complexes:definition800}
Let $X,Y_i\,(i=1,\ldots,n)$ be $2$-complexes and 
$f_i:Y_i\rightarrow X$ a collection of
dimension preserving maps of $2$-complexes. Let $e\in X$ be an edge
and $\{e_{j1},e_{j2}\}$ ($j=1,\ldots,m$), pairs of edges in $Y_0=\bigcup_i
Y_i$ such that for each $j$, the edges $\{e_{j1},e_{j2}\}$ lie in the
fiber $f_i^{-1}(e)$
for some $i$, and for each $i$, the complex $Y_i$ contains a pair
$\{e_{j1},e_{j2}\}$ for some $j$.

The pairs $\{e_{j1},e_{j2}\}$ form a \emph{handle configuration\/} if and only if
for every face $\ss\in X$ containing $e$ in its boundary we have,
\begin{description}
\item[(i).] if $\tau\in\bigcup_i f_i^{-1}(\ss)$ and for some $j$ we have
  $\aa_\tau(Y_0^\tau)=e_{j1}^k$ (respectively $e_{j2}^k$),
  then $k=m$ and there are faces 
  $\tau_{11},\ldots,\tau_{1m},\tau_{21},\ldots,\tau_{2m}\in\bigcup_i f_i^{-1}(\ss)$, with
  $\tau=\tau_{j1}$ (resp. $\tau_{j2}$), and $\aa_{\tau_{\ell i}}(Y^{\tau_{\ell i}})=e_{\ell i}^k$
  ($\ell=1,\ldots,m$ and $i=1,2$);
\item[(ii).] for any face $\tau\in\bigcup_i f_i^{-1}(\ss)$ not of the form (i) and containing 
$e_{j1}$ or $e_{j2}$ in its boundary we have
$\aa_\tau(Y_0^\tau)=(e_{j1}\gamma_{j1}e_{j2}\gamma_{j2})^k$ with $\gamma_{j1},\gamma_{j2}$ not
containing $e_{j1},e_{j2}$. Moreover, there are faces 
  $\tau_1,\ldots,\tau_m\in\bigcup_i f_i^{-1}(\ss)$, with
  $\tau=\tau_j$, the
$\aa_{\tau_\ell}(Y_0^{\tau_\ell})=(e_{\ell 1}\gamma_{\ell 1}e_{\ell 2}\gamma_{\ell 2})^k$ with 
$\gamma_{\ell 1},\gamma_{\ell 2}$ not containing $e_{\ell 1},e_{\ell 2}$, and the $\gamma_{\ell 1}$
(respectively the $\gamma_{\ell 2}$) all mapping to the same path in $X$.
\end{description}
\end{definition}

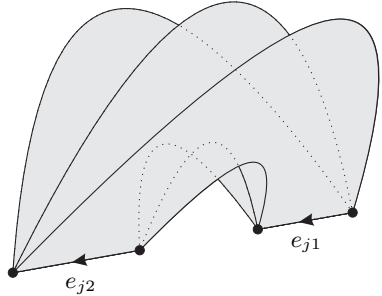
\begin{figure}
\begin{pspicture}(0,0)(12,4)
\rput(3.5,0){
\rput(2.5,2){\BoxedEPSF{fig133.eps scaled 400}}
\rput(1,0.05){$e_{j2}$}\rput(4,.6){$e_{j1}$}
}
\end{pspicture}
\caption{Faces around an edge pair in a handle configuration.}
\label{chapter1:2complexes:figure3000}
\end{figure}

Whenever Definition
\ref{chapter1:2complexes:definition800}(i) happens the edges
$e_{j1},e_{j2}$ in the handle configuration start and finish at the
same vertex. The faces in Definition
\ref{chapter1:2complexes:definition800}(ii) form a 
``taco'' arrangement as shown in Figure \ref{chapter1:2complexes:figure3000}:
whenever a face contains $e_{j1}$ in its boundary it also contains
$e_{j2}$ and vice-versa.

\begin{definition}[Higman composition]
\label{chapter1:2complexes:definition900}
Let $f_i:Y_i\rightarrow X$ ($i=1,\ldots,n$), be dimension preserving 
and $\{e_{j1},e_{j2}\}$ ($j=1,\ldots,m$), a handle configuration in $\bigcup_i
Y_i$. The \emph{Higman composition\/} produces a
$2$-complex $Y=\hcl Y_1,\ldots,Y_n\hcr$ and a map $f=\hcl
f_1,\ldots,f_n\hcr:Y\rightarrow X$ as follows:
\begin{description}
\item[(HC1).] Delete the edges $e_{j1},e_{j2}$ in the handle
  configurations and replace them by new edges $e'_{j1},e'_{j2}$, where
  for each $j$, the edge $e'_{j1}$ connects $s(e_{j1})$ to
  $t(e_{j+1,1})$, while $e'_{j2}$ connects $s(e_{j2})$ to
  $t(e_{j-1,2})$, subscripts modulo $m$ as in Figure
  \ref{chapter1:2complexes:figure3100}. 

\begin{figure}[h]
  \centering
\begin{pspicture}(0,0)(14,2.5)
\rput(-2,.2){
\rput(9,1.25){\BoxedEPSF{chapter1.fig5200.eps scaled 400}}
\rput(9.3,1.8){$e_{j1}$}\rput(6.7,1.3){$e_{j2}$}
\rput(10.2,.7){$e_{j+1,1}$}\rput(8.1,.3){$e_{j+1,2}$}
\rput(10.25,1.5){${\red e'_{j1}}$}\rput(6.6,.8){${\red e'_{j2}}$}
}
\end{pspicture}  
  \caption{Old (black) and new (red ) edges in a Higman composition.}
  \label{chapter1:2complexes:figure3100}
\end{figure}
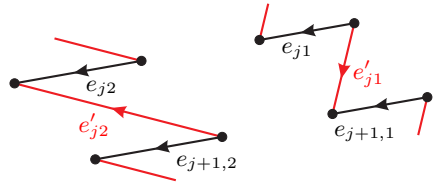

\item[(HC2).] For the faces $\tau$ in Definition
  \ref{chapter1:2complexes:definition800}(i), delete the
  $\tau_1,\ldots,\tau_m$ and replace then by new faces
  $\tau'_1,\ldots,\tau'_m$, with $\partial\tau'_\ell$ using the same
  $X^{\tau_\ell}$ but having attaching map sending it to $e'_{\ell 1}\ldots e'_{m1}
  e'_{11}\ldots e'_{\ell-1,1}$. 

\item[(HC3).] For the faces in Definition
  \ref{chapter1:2complexes:definition800}(ii), delete the
  $\tau_1,\ldots,\tau_m$ and replace then by new faces
  $\tau'_1,\ldots,\tau'_m$, with $\partial\tau'_\ell$ using the same
  $X^{\tau_\ell}$ but having attaching map sending it to the path
  shown in Figure \ref{chapter1:2complexes:figure3300}.

\begin{figure}[h]
  \centering
\begin{pspicture}(0,0)(14,5)
\rput(3.2,2.5){\BoxedEPSF{chapter1.fig5100.eps scaled 275}}
\rput(10.2,2.5){\BoxedEPSF{chapter1.fig5000a.eps scaled 275}}
\rput(2.2,3.5){$\tau_\ell$}\rput(5.1,2){$\tau_{\ell+1}$}
\rput(4,4){${\red\aa_{\tau_\ell}(X^{\tau_\ell})}$}
\rput(11,4){${\red\aa_{\tau'_\ell}(X^{\tau_\ell})}$}
\rput(10,2){$\tau'_\ell$}
\end{pspicture}  
  \caption{Old and new faces in a Higman composition (HC3).}
  \label{chapter1:2complexes:figure3300}
\end{figure}
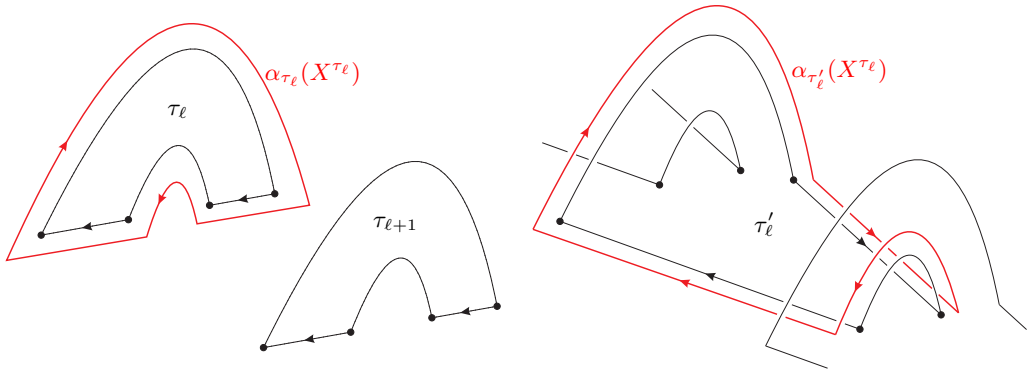

\item[(HC4).] A cell $x$ of $Y$ is either a cell $x$ of one of the $Y_i$ or replaces 
  a cell $x'$ of one of the $Y_i$. Define $f(x)$ to be $f_i(x)$ or $f_i(x')$ as 
  appropriate.
  The face isomorphisms
  remain unchanged.
\end{description}
\end{definition}

\begin{exercise}
\label{chapter1:2complexes:exercise700}
Give an example showing that if the $Y_i$ are all connected then the Higman composition
$\hcl Y_1,\ldots,Y_m\hcr$ is not necessarily connected (see also \S\ref{chapter3:operations:higman}).
\end{exercise}

We defer further exploration of Higman composition until Chapter 3,
where we will show that the Higman composition of a family of
coverings $f_i:Y_i\rightarrow X$ yields another covering 
$f:\hcl Y_1,\ldots,Y_m\hcr\rightarrow X$. 
\section{Notes on Chapter \thechapter}






\chapter{Topological Invariants}\label{chapter2}

\chapter{Coverings}\label{chapter3}

\section{Basics}\label{chapter3:basics}

\subsection{Coverings}\label{chapter3:basics:coverings}

\begin{definition}[covering]
\label{chapter3:basics:definition200}
A map $f:Y\rightarrow X$ of $2$-complexes is a \emph{covering\/} if and only if
\begin{description}
\item[(C1)] $f$ preserves dimension (see Definition \ref{chapter1:2complexes:definition275});
\item[(C2)] for every pair of vertices $u\in Y$ and $v\in X$, with $f(u)=v$, the local continuity 
of $f$ at $v$ (see \S\ref{chapter1:graphs:category}),
$$
s_Y^{-1}(u)\rightarrow s_X^{-1}(v),
$$
is a \emph{bijection\/}.
\item[(C3)] for every pair of vertices $u\in Y$ and $v\in X$ with $f(u)=v$, and 
every face $\tau$ of $X$,
the local continuity of $f$ at $v$ (see Definition \ref{chapter1:2complexes:definition250}),
$$
\amalg\,\ve(f,\ss):\bigcup_{f(\ss)=\tau}\kern-2mm\aa_\ss^{-1}(u)\rightarrow\aa_\tau^{-1}(v)  
$$
is a \emph{bijection\/}.
\end{description}
\end{definition}

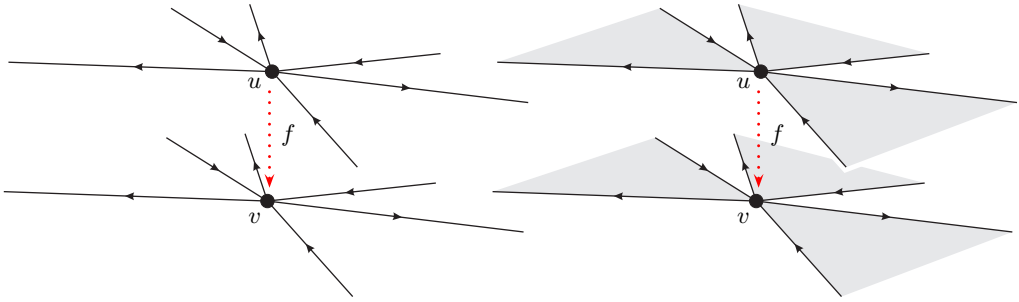
\begin{figure}
  \centering
\begin{pspicture}(0,0)(12.5,4)
\rput(-.2,-.6){
\rput(3.7,2.5){\BoxedEPSF{chapter3.fig50.eps scaled 250}}
\rput(3.4,3.4){$u$}\rput(3.4,1.6){$v$}\rput(3.85,2.7){$f$}
\psline[linewidth=.4mm,linestyle=dotted,linecolor=red]{->}(3.6,3.4)(3.6,2)
}
\rput(6.3,-.6){
\rput(3.7,2.5){\BoxedEPSF{chapter3.fig100.eps scaled 250}}
\rput(3.4,3.4){$u$}\rput(3.4,1.6){$v$}\rput(3.85,2.7){$f$}
\psline[linewidth=.4mm,linestyle=dotted,linecolor=red]{->}(3.6,3.4)(3.6,2)
}
\end{pspicture}  
  \caption{the local continuity maps are bijections for coverings}
  \label{chapter3:basics:figure100}
\end{figure}

A covering is a kind of ``local isomorphism'': if $f(u)=v$ then $Y$ looks
the same ``near'' $u$ as $X$ does ``near'' $v$. Thus (C2) ensures that the
configuration of edges around a vertex looks the
same both upstairs and downstairs (Figure \ref{chapter3:basics:figure100} left).
Similarly (C3) means that for a face $\tau$ downstairs containing the vertex $v$ in
its boundary and $f(u)=v$, this face
looks the same near $v$ as its pre-images do near $u$. Specifically, 
if $v$ appears $k$ times in the boundary of $\tau$, so
there are $k$ ``wedge-shaped'' pieces of $\tau$ fitting together around $v$,
then there are $k$ wedge-shaped pieces of face fitting together around
$u$, where these wedges belong to faces $\ss$ mapping to $\tau$
(see Figure \ref{chapter3:basics:figure100} right).

The terminology \emph{cover\/} and \emph{lift\/} is used for images
and pre-images of a covering map: if
$f(y)=x$, then one says that $y$ covers $x$,
or that $x$ lifts 
to $y$. The set of all lifts of $x$, or the set $f^{-1}(x)$ of all cells covering
$x$, is its \emph{fiber\/}.
Note that each $\ve(f,\ss)$ is an isomorphism as $f$ preserves dimension, 
so the local continuity map
is a disjoint union of a set of maps, each of which is the restriction of 
a bijection. In any case, each is individually an injection.

Part 3 of the definition gives in particular that
\begin{equation}
  \label{chapter3:equation100}
\sum_{f(\ss)=\tau} |\aa_\ss^{-1}(u)|=|\aa_\tau^{-1}(v)|,  
\end{equation}
so that $v$ appears the same number of times in the boundary of $\tau$
as $u$ does
in the boundaries of all the faces $\ss$ in the fiber of $\tau$ (although be sure
to take on board Example \ref{chapter3:basics:example400}).

It is easy to see that a covering $f:Y\rightarrow X$ can be restricted to a covering
$f:Y^{(1)}\rightarrow X^{(1)}$ of the $1$-skeletons.

One commonly sees the assumption that in a covering, both the covering complex 
$Y$ and the
covered complex $X$ are connected, but we won't assume this at the moment. Indeed,
we will find it useful in some situations to \emph{not\/} assume that a covering 
complex be connected.

\begin{exercise}\label{chapter3:basics:exercise100}
Let $f:Y\rightarrow X$ be a covering and $Y^\circ$ a connected component of $Y$. 
Show that restricting
$f$ to $Y^\circ$ gives a covering. Show that 
we may not restrict
a covering to an arbitrary subcomplex and still get a covering.
\end{exercise}

\begin{example}\label{chapter3:basics:example100}
Figure \ref{chapter3:basics:figure200} shows a simple graph covering, with the 
two vertices in $Y$ covering the single vertex of $X$, and the edges such that the
red path in $Y$ covers the red path in $X$. In particular, the two edges of $Y$ both cover the single edge of $X$.
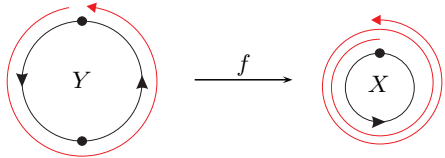
\begin{figure}[h]
  \centering
\begin{pspicture}(0,0)(12.5,2)
\rput(4,1){\BoxedEPSF{chapter3.fig200.eps scaled 750}}
\rput(8,1){\BoxedEPSF{chapter3.fig300.eps scaled 750}}
\rput(4,1){$Y$}\rput(7.95,.95){$X$}\rput(6.15,1.2){$f$}
\psline[linewidth=.2mm]{->}(5.5,1)(6.8,1)
\end{pspicture}
  \caption{a simple graph covering: the two vertices of $Y$ cover the 
single vertex of $X$ and the two arcs of $Y$ similarly.}
  \label{chapter3:basics:figure200}
\end{figure}
\end{example}

\begin{example}\label{chapter3:basics:example200}
Figure \ref{chapter3:basics:figure300} extends the graph covering of Example
\ref{chapter3:basics:example100} to a covering $f_1:Y_1\rightarrow X$ of $2$-complexes. 
\begin{figure}[h]
  \centering
\begin{pspicture}(0,0)(12.5,3)
\rput(-9.5,0){
\rput(12.5,1.5){\BoxedEPSF{chapter1.fig4200.eps scaled 500}}
\rput(2.5,-4.5){
\rput(11.6,7.4){$\ss_1$}
\rput(11.6,4.6){$\ss_2$}
\rput(12.45,6.4){${\red\partial\ss_1}$}
\rput(12.45,5.8){${\red\partial\ss_2}$}
\rput(9.7,5.45 ){${v}_1$}
\rput(10.2,6.55){${v}_2$}
\rput(10.6,5.45){${e}_2$}
\rput(9.4,6.6){${e}_1$}
}}
\rput(8.7,0){
\rput(1.25,1.25){\BoxedEPSF{chapter1.fig4300.eps scaled 500}}
\rput(2.2,.2){$\ss$}\rput(.4,1.55){$v$}\rput(2.2,2.4){$v$}\rput(1.2,1.05){$e$}
\rput(1.2,2.5){$e$}
}
\rput(-.25,0){
\psline[linewidth=.2mm]{->}(6,1.5)(8,1.5)
\rput(7,1.75){$f_1$}
}
\rput(0,1.5){$Y_1$}\rput(12.5,1.5){$X$}
\end{pspicture}
  \caption{the graph covering of Example \ref{chapter3:basics:example100} extended
to a covering of $2$-complexes: $Y$ is now the complex of Figure
\ref{chapter1:2complexes:figure200} and $X$ the complex of Figure
\ref{chapter1:2complexes:figure400}. The $\ss_i$ both cover the same face
$\ss\in X$, and the face isomorphisms are $\ve_{\ss_1}=\id$ and
$\ve_{\ss_2}$ a $1/2$-turn.
}
  \label{chapter3:basics:figure300}
\end{figure}
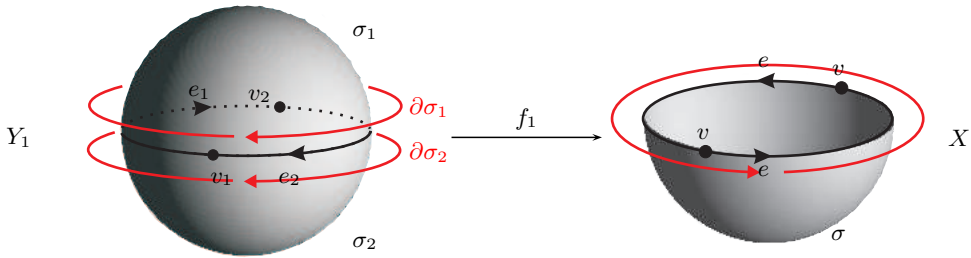
The two faces $\ss_1,\ss_2$ both cover the face $\ss$ of $X$, but the face isomorphisms
are different: $\ve_{\ss_1}$ is the identity map and $\ve_{\ss_2}$ is a clockwise $1/2$-turn. 
\end{example}

\begin{example}\label{chapter3:basics:example300}
Figure \ref{chapter3:basics:figure350} also extends the graph covering of Example
\ref{chapter3:basics:example100} to a covering $f_2:Y_2\rightarrow X$ of $2$-complexes. 
In this case $\ss_1$ covers the face $\ss$ of $X$ and $\ss_2$ covers $\ss^{-1}$. 
For the face isomorphisms, $\ve_{\ss_1}$ is the identity map and $\ve_{\ss_2}$ is a
$1/2$-turn.
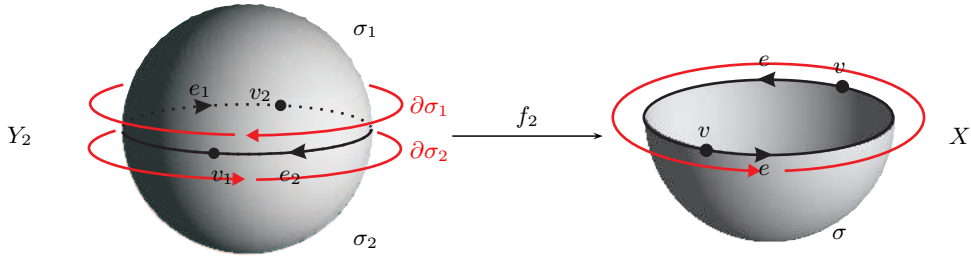
\begin{figure}[h]
  \centering
\begin{pspicture}(0,0)(12.5,3)
\rput(-9.5,0){
\rput(12.5,1.5){\BoxedEPSF{chapter1.fig4200a.eps scaled 500}}
\rput(2.5,-4.5){
\rput(11.6,7.4){$\ss_1$}
\rput(11.6,4.6){$\ss_2$}
\rput(12.45,6.4){${\red\partial\ss_1}$}
\rput(12.45,5.8){${\red\partial\ss_2}$}
\rput(9.7,5.45 ){${v}_1$}
\rput(10.2,6.55){${v}_2$}
\rput(10.6,5.45){${e}_2$}
\rput(9.4,6.6){${e}_1$}
}}
\rput(8.7,0){
\rput(1.25,1.25){\BoxedEPSF{chapter1.fig4300.eps scaled 500}}
\rput(2.2,.2){$\ss$}\rput(.4,1.55){$v$}\rput(2.2,2.4){$v$}\rput(1.2,1.05){$e$}
\rput(1.2,2.5){$e$}
}
\rput(-.25,0){
\psline[linewidth=.2mm]{->}(6,1.5)(8,1.5)
\rput(7,1.75){$f_2$}
}
\rput(0,1.5){$Y_2$}\rput(12.5,1.5){$X$}
\end{pspicture}
  \caption{the graph covering of Example \ref{chapter3:basics:example100} extended
to a different covering of $2$-complexes: $Y$ is now the complex of Figure
\ref{chapter1:2complexes:figure150} and $X$ the complex of Figure
\ref{chapter1:2complexes:figure400}. 
}
  \label{chapter3:basics:figure350}
\end{figure}
\end{example}

\begin{exercise}
If $\aa:Y_1\rightarrow Y_2$ is the isomorphism of Exercise 
\ref{chapter1:2complexes:exercise175}, show that it commutes with the 
coverings of the previous two examples: $f_1=f_2\aa$.
\end{exercise}

\begin{example}\label{chapter3:basics:example400}
Return to the complexes $Y_1$ and $X$ of Example \ref{chapter3:basics:example100} but tweak the 
covering slightly: define $f_1':Y_1\rightarrow X$ by $f_1'(x)=f_1(x)$ for all cells $x\in Y_1$.
Where the two coverings differ is in the face isomorphisms: define $\ve_{\ss_2}$ to
be the identity rather than a $1/2$-turn. One can then check that the 
(\ref{chapter3:equation100}) is satisfied, but the local continuity maps are not 
bijections, and so we do not have a covering.
\end{example}

\begin{exercise}[immersions]
\label{chapter3:basics:exercise150}
Call a map $f:Y\rightarrow X$ an \emph{immersion\/} when it preserves dimension
and the local continuity maps are \emph{injections\/}. Give examples of
immersions that are not coverings.
\end{exercise}

\begin{exercise}\label{chapter3:basics:exercise160}
Show that for any $X$ the identity map $\id: X\rightarrow X$ is a covering.
\end{exercise}

\subsection{Lifting}\label{chapter3:basics:lifting}

When we have a covering $f:Y\rightarrow X$ the complexes $Y$ and $X$
look the same so long as we restrict our attention to small pieces. 
If two complexes look the same as each other then we should be able to pull parts
of $X$ back through $f$ to find parts of $Y$ mapping to them. Putting these together, 
when we have a covering we should be able to
pull small pieces of $X$ back through $f$ and get 
identical small pieces of $Y$ covering them. The small pieces turn out to be 
paths and faces, and this process is called is called \emph{lifting\/}.

\begin{proposition}[path, spur and free homotopy lifting]
\label{chapter3:basics:result100}
Let $f:Y\rightarrow X$ be a covering with $f(u)=v$ vertices.
\begin{description}
\item[(i).] If $\gamma$ is a path in $X$ starting at $v$
then there is a path $\mu$ in $Y$ starting at $u$
and covering $\gamma$.
Moreover, if $\mu_1,\mu_2$ are paths in $Y$ starting at 
$u$ and covering the same path in $X$, 
then $\mu_1=\mu_2$.
\item[(ii).] A path in $Y$ covering a spur is itself a spur.
Consequently, two paths in $Y$ covering freely homotopic paths are themselves
freely homotopic.
\end{description}
\end{proposition}

Part (i) is called \emph{path lifting\/} and part (ii) is \emph{spur lifting\/}.
Call $\mu$ the \emph{lift\/} of $\gamma$ at $u$. Thus a path can be lifted to 
any vertex that covers its initial vertex to give a covering path,
and this lift is unique.
As with so many such results, it is
the uniqueness of the lift, rather than the existence, that turns out to be most 
useful.

\begin{proof}
The existence of $\mu$ is easily seen, as in Figure \ref{chapter3:basics:figure400},
\begin{figure}
  \centering
\begin{pspicture}(0,0)(12.5,3)
\rput(6.25,1.5){\BoxedEPSF{chapter3.fig400.eps scaled 1000}}
\rput(3.75,.7){$v$}
\rput(3.75,2){$u$}
\rput(4.8,.7){$e_1$}
\rput(4.8,2){$e'_1$}
\rput(5.8,.7){$e_2$}
\rput(5.8,2){$e'_2$}
\rput(8,.8){$e_n$}
\rput(8,2.1){$e'_n$}
\rput(2.5,.3){$X$}\rput(2.5,1.6){$Y$}
\rput(9,1){$\gamma$}\rput(9,2.4){$\mu$}
\rput(0,0){\psline[linewidth=.4mm,linestyle=dotted,linecolor=red]{->}(4,1.9)(4,.9)}
\rput(1.25,.4){\psline[linewidth=.4mm,linestyle=dotted,linecolor=red]{->}(4,1.9)(4,.9)}
\rput(2.5,0){\psline[linewidth=.4mm,linestyle=dotted,linecolor=red]{->}(4,1.9)(4,.9)}
\rput(3.4,.3){\psline[linewidth=.4mm,linestyle=dotted,linecolor=red]{->}(4,1.9)(4,.9)}
\rput(4.7,.4){\psline[linewidth=.4mm,linestyle=dotted,linecolor=red]{->}(4,1.9)(4,.9)}
\end{pspicture}
  \caption{lifting paths in a covering: individual edges can be lifted via the
bijection between the edges starting at a vertex upstairs and the edges starting at a 
vertex it covers downstairs. Paths are then lifted by repeated edge lifts.}
  \label{chapter3:basics:figure400}
\end{figure}
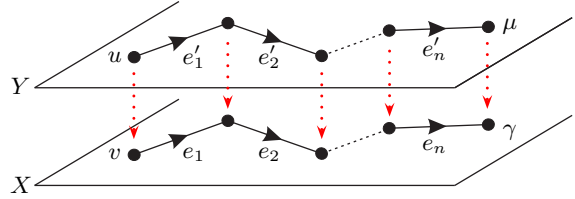
since if $\gamma=e_1\ldots e_n$, there is an edge $e'_1$ covering $e_1$ under the
bijection $s^{-1}_Y(u)\rightarrow s^{-1}_X(v)$.
This edge $e'_1$ must end at a vertex that covers the end vertex of $e_1$,
as coverings (being maps of complexes) preserve vertex-edge incidences. 
The process can be repeated starting at this new vertex to give $\mu$.
For the uniqueness, 
the first edges of the $\mu_i$ both have initial vertex $u$ and map to $e_1$, 
hence must be
the same edge. Continuing in this manner along the two paths gives their equality.
For part (ii), 
the path in $Y$ must have the form $e_1e_2$, where the middle vertex is the start of the 
edges $e_1^{-1}$ and $e_2$. Use the injectivity of $f$ on the edges starting
at this vertex to deduce that $e_1^{-1}=e_2$.
\qed
\end{proof}

\begin{exercise}\label{chapter3:basics:exercise150}
Let $f:Y\rightarrow X$ be a covering and $\gamma=\gamma_1\gamma_2$ a path in $X$, hence
$t(\gamma_1)=s(\gamma_2)$. Show that the lift at a vertex $u\in Y$ of $\gamma$ is the path
$\mu_1\mu_2$ consisting of the lift $\mu_1$ of $\gamma_1$ at $u$
followed by the lift $\mu_2$ of $\gamma_2$ at $t(\mu_1)$.
\end{exercise}

\begin{exercise}\label{chapter3:basics:exercise200}
Show that paths cannot necessarily be lifted by an immersion, but when they can, they are
unique. Show that spur lifting is a property enjoyed by immersions.  
\end{exercise}

\begin{proposition}[face lifting]\label{chapter3:basics:result300}
Let $f:Y\rightarrow X$ be a covering, $\tau\in X$ a face and $v$ a vertex
with $x\in X^\tau$ an appearance of $v$ in the boundary of $\tau$.
Let $\gamma$ be a boundary path of $\tau$ starting at $v$ and given by
$x$. 
Finally, let $u$ be a vertex of $Y$ covering $v$ and let
$\mu$ be the lift of $\gamma$ to $u$. Then there is a unique
appearance $y$ of $u$ in the boundary of some face $\ss$ covering
$\tau$, 
with the local continuity $\ve(f,\ss)(y)=x$
and $\mu$ a boundary path of $\ss$ starting at
$u$ and given by $y$.
\end{proposition}

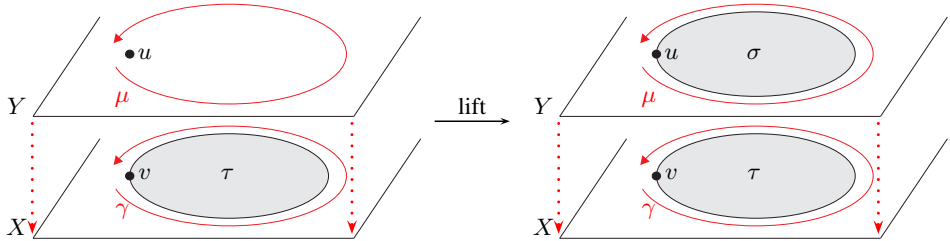
\begin{figure}
  \centering
\begin{pspicture}(0,0)(12.5,3)
\rput(-.25,0){
\rput(3,1.5){\BoxedEPSF{chapter3.fig500.eps scaled 700}}
\rput(0,0){\psline[linewidth=.4mm,linestyle=dotted,linecolor=red]{->}(.4,1.6)(.4,0)}
\rput(0,0){\psline[linewidth=.4mm,linestyle=dotted,linecolor=red]{->}(4.65,1.6)(4.65,0)}
\rput(3,.8){$\tau$}\rput(1.6,.3){${\red \gamma}$}\rput(1.6,1.8){${\red\mu}$}
\rput(1.9,.8){$v$}\rput(1.9,2.4){$u$}
\rput(.2,.1){$X$}\rput(.2,1.7){$Y$}
}
\rput(6.75,0){
\rput(3,1.5){\BoxedEPSF{chapter3.fig600.eps scaled 700}}
\rput(0,0){\psline[linewidth=.4mm,linestyle=dotted,linecolor=red]{->}(.4,1.6)(.4,0)}
\rput(0,0){\psline[linewidth=.4mm,linestyle=dotted,linecolor=red]{->}(4.65,1.6)(4.65,0)}
\rput(3,.8){$\tau$}\rput(3,2.4){$\ss$}
\rput(1.6,.3){${\red \gamma}$}\rput(1.6,1.8){${\red\mu}$}
\rput(1.9,.8){$v$}\rput(1.9,2.4){$u$}
\rput(.2,.1){$X$}\rput(.2,1.7){$Y$}
}
\psline[linewidth=.2mm]{->}(5.5,1.5)(6.5,1.5)
\rput(6,1.7){lift}
\end{pspicture}
  \caption{face lifting}
  \label{chapter3:basics:figure500}
\end{figure}

The boundaries of faces thus lift to the boundaries of faces, as 
in Figure \ref{chapter3:basics:figure500}. We will call the uniqueness
statement of Proposition \ref{chapter3:basics:result300}, ``uniqueness of face lifting''.

\begin{proof}
(C3) of Definition \ref{chapter3:basics:definition200} gives a unique vertex $y$
in $\bigcup_{f(\ss)=\tau}\,\aa_\ss^{-1}(u)$ mapping to $x$ via the local continuity of $f$,
ie: there is a face $\ss$ of $Y$ covering $\tau$,
containing $u$ in its boundary, and with the diagram of pointed maps
$$
\begin{pspicture}(0,0)(12.5,2)
\rput(1.75,-.8){
\rput(-1.3,.5){
\rput(5,2){$Y^\ss_y$}\rput(6.65,2){$X^{\tau}_x$}
\rput(4.95,0.45){$Y^{(1)}_u$}\rput(6.65,0.45){$X^{(1)}_v$}
\psline[linewidth=.2mm]{->}(5.3,2)(6.3,2)
\psline[linewidth=.2mm]{->}(5.3,0.45)(6.3,0.45)
\psline[linewidth=.2mm]{->}(5,1.7)(5,.7)
\psline[linewidth=.2mm]{->}(6.55,1.7)(6.55,.7)
\rput(4.7,1.2){$\aa_\ss$}\rput(6.9,1.175){$\aa_{\tau}$}
\rput(5.8,2.2){$\cong$}
\rput(5.8,.7){$f$}
}
}
\end{pspicture}
$$
commuting. In particular there is a boundary path of $\ss$ starting at $u$ that covers
$\gamma$, and by the uniqueness of lifts, this must be the path $\mu$.
\qed
\end{proof}

The first consequence of lifting justifies the usage of the word ``cover'',
and is not {\em a priori\/} obvious from the definition:

\begin{corollary}[surjectivity of coverings]
\label{chapter3:basics:result200}
If $f:Y\rightarrow X$ is a covering with $X$ connected then $f$ is a surjective map of
$2$-complexes, ie: every cell of $X$ is the image under $f$ of some cell of $Y$.
\end{corollary}

\begin{proof}
Path lifting gives the surjectivity on the vertices 
and edges, and face lifting on the faces:
fix a vertex $u$ of $Y$, and
by connectedness, we can join any vertex $v$ of $X$ to $f(u)$ by a path. Lift
this path to $u$, so that its terminal vertex in $Y$ maps via $f$ to $v$.
For an edge $e$ or face $\ss$ of $X$, let $v$ be a vertex in the boundary and lift
the edge or face to a vertex $u$ covering $v$.
\qed
\end{proof}

\begin{exercise}\label{chapter3:basics:exercise300}
Illustrate by an example why the connectedness of $X$ is necessary in Proposition
\ref{chapter3:basics:result200}.
\end{exercise}

Another result of being able to find pre-images of paths and faces is that homotopies
can be ``pulled back'' through a covering:

\begin{corollary}[homotopy lifting]\label{chapter3:basics:result400}
Let $Y\rightarrow X$ be a covering.
Then two paths that cover homotopic paths are themselves homotopic.
\end{corollary}

\begin{proof}
The homotopy between the covered paths is realised by a finite sequence of insertions or deletions of
spurs and face boundaries. By spur and face lifting, those sections
of the covering paths mapping to the spurs and face boundaries are themselves spurs and face boundaries, 
while
by uniqueness of path lifting, the remaining pieces are identical. Thus the same sequence of elementary
homotopies can be realised between the covering paths as between the covered ones.
\qed
\end{proof}

\begin{exercise}\label{chapter3:basics:exercise400}
Let $f:Y\rightarrow X$ be a covering and $\gamma_1,\gamma_2$ paths in $X$ related by
an elementary homotopy, ie: $\gamma_2$ is what results by inserting/deleting a spur or face
boundary into $\gamma_1$. Let $\mu_1$ be the lift of $\gamma_1$ to some vertex
$u$ of $Y$ and $\mu$ the result of lifting to the appropriate vertex
the elementary homotopy and performing
it on $\mu_1$. Show that $\mu$ is the lift $\mu_2$ of $\gamma_2$ at $u$.
\end{exercise}

Another spin-off of homotopy lifting 
is the following characterisation of the image of the induced homomorphism between 
fundamental groups:

\begin{corollary}\label{chapter3:basics:result500}  
Let $f:Y\rightarrow X$ be a covering with $f(u)=v$ and 
$$f_*:\pi_1(Y,u)\rightarrow\pi_1(X,v),$$
the induced homomorphism. Then
$f_*$ is injective, and
a closed path $\gamma$ at $v$ represents an element of 
$f_*\pi_1(Y,u)$ if and only if the lift $\mu$ of $\gamma$ to $u$ is closed.
\end{corollary}

The injectivity of the induced homomorphism is probably the single 
most important property of coverings: it means that the fundamental group of the
covering space can be identified with a subgroup of the fundamental group of the
covered space. The appropriate context in which to develop this idea properly
will be the Galois theory of coverings in Chapter \ref{chapter4}.

\begin{proof}
If two elements of $\pi_1(Y,u)$ map to the same element of 
$\pi_1(X,v)$ then they are represented by closed paths at $u$ covering homotopic paths
at $v$. Homotopy lifting gives that the 
paths in $Y$ are homotopic, and so the two elements of the fundamental group
coincide, thus establishing the injectivity of the homomorphism. For the second part, if
$\mu$ is closed then its homotopy class maps via $f_*$ to the homotopy class 
of $\gamma$. Conversely, if $\gamma$ represents an element in the image of the 
homomorphism then there is a closed path $\gamma_1$ at $v$, homotopic to $\gamma$,
with $f(\mu_1)=\gamma_1$ for $\mu_1$ closed at $u$ (and the lift of $\gamma_1)$. 
By homotopy lifting the lift $\mu$ of $\gamma$ is homotopic to $\mu_1$, 
hence has the same endpoints, ie: is closed.
\qed
\end{proof}

\begin{exercise}\label{chapter3:basics:exercise500}
Let $f:Y\rightarrow X$ be a covering with $f(u)=v$ and
$u'$ the terminal vertex of a path $\mu\in Y$ starting at $u$. Show that
$f_*\pi_1(Y,u)=hf_*\pi_1(Y,u')h^{-1}$, where $h$ is the 
homotopy class of $f(\mu)$.
\end{exercise}

This thread of ideas culminates in (and is subsumed by) the following general lifting result.

\begin{theorem}[map lifting]
\label{chapter3:basics:result600}  
If $f:Y\rightarrow X$ is a covering with $f(u)=v$ and $g:Z\rightarrow X$ a map
with $g(x)=v$ and $Z$ connected, then there is a map $\wtl{g}:Z\rightarrow Y$
making the diagram 
$$
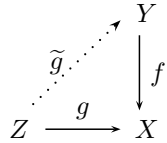
\begin{pspicture}(0,0)(12.5,2)
\rput(1.75,-.8){
\rput(-1.3,.5){
\rput(6.65,2){$Y$}
\rput(4.95,0.45){$Z$}\rput(6.65,0.45){$X$}
\psline[linewidth=.2mm]{->}(5.3,0.45)(6.3,0.45)
\psline[linewidth=.2mm]{->}(6.55,1.7)(6.55,.7)
\psline[linewidth=.3mm,linestyle=dotted]{->}(5.05,.7)(6.3,1.9)
\rput(6.8,1.175){$f$}
\rput(5.8,.7){$g$}\rput(5.45,1.35){$\wtl{g}$}
}
}
\end{pspicture}
$$
commute if and only $g_*\pi_1(Z,x)\subset f_*\pi_1(Y,u)$. If $\wtl{g}$ exists
then it is unique.
\end{theorem}

Think of this result is as a generalisation of path lifting: if $Z$ is a 
$1$-ball then the map $g:Z\rightarrow X$ is a path in $X$ starting at $v$. As the
fundamental group of a $1$-ball is trivial, the condition $g_*\pi_1(Z,x)\subset f_*\pi_1(Y,u)$
is trivially satisfied.
The resulting map $\wtl{g}:Z\rightarrow Y$ is a path in $Y$ starting at $u$,
and the commuting of the diagram just says that this new path is the lift of the old one.

\begin{proof}
The ``only if'' part can be dispensed with quickly as $\pi_1$ is a functor:
$f\wtl{g}=g$ gives $f_*\wtl{g}_*=g_*$
so that $g_*\pi_1(Z,x)=f_*\wtl{g}_*\pi_1(Z,x)\subset f_*\pi_1(Y,u)$.

Suppose we have the condition on the fundamental groups, which by 
Corollary \ref{chapter3:basics:result500}   
means that if $\gamma$ is a closed path at $x\in Z$
then the lift $\mu$ of $g(\gamma)$ to $u$ is also closed.
We proceed to define
a map $\wtl{g}$ having the required properties: if $z$ is a vertex of $Z$, then by connectedness
there is a path joining it to $x$. Take the image of this path by $g$ and then lift
the result via the covering $f$ to a path at $u$. Define $\wtl{g}(z)$ to be the 
end vertex of the resulting path in $Y$. Edges and faces are similar: 
choose a vertex $z$ in the boundary of the edge or face (for an edge $e$, 
choose $z=s(e)$) and then lift the image under $g$ of the edge/face via the covering 
$f$ to the vertex $\wtl{g}(z)$.

If this procedure is well defined, then it is easy to check that we have a map
which by definition makes the 
diagram commute (for the face isomorphism $Z^\ss\rightarrow Y^{\wtl{g}(\ss)}$, 
take the composition $\ve(\wtl{g},\ss)=\ve(f,\wtl{g}(\ss))^{-1}\ve(g,\ss)$). 
To show that $\wtl{g}$ is well defined on the vertices,
suppose that $\gamma_1,\gamma_2$ are paths in $Z$ from $x$ to the vertex $z$,
so that $\gamma_1\gamma_2^{-1}$ is a closed path at $x$. Thus, the lift of 
its $g$-image (which is $\mu_1\mu_2^{-1}$) is closed too, and so 
the image $\wtl{g}(z)$ does not depend on the choice of the path $\gamma$.
Edge images are well defined as there are no additional choices made. For a face
$\ss\in Z$ the construction involves a choice of 
vertex in its boundary, so suppose that $z_1,z_2$ are
two such. If $\gamma$ is a boundary path for $\ss$, then the lift of 
$g(\gamma)$ must pass through both $\wtl{g}(z_1)$ and $\wtl{g}(z_2)$ by
the well-definedness of $\wtl{g}$ on the vertices. Applying the uniqueness of face lifting 
to these two vertices gives what we want.

Finally, if $\wtl{g}_1,\wtl{g}_2$ are two maps making the diagram commute, then
path and face lifting gives $\wtl{g}_1(x)=\wtl{g}_2(x)$ for any cell $x\in Z$. The face 
isomorphisms of both are the compositions of face isomorphisms of $g$ and $f$, and so are
identical. Thus $\wtl{g}_1=\wtl{g}_2$.
\qed
\end{proof}

\subsection{Degree}\label{chapter3:basics:degree}

We now come to an important invariant that can be attached to a covering. Looking back
at Example \ref{chapter3:basics:example100},
we have a graph covering of $X$ where both the fiber of the vertex and the fiber of
the edge contain two cells. Extending this covering to one of $2$-complexes
in Example \ref{chapter3:basics:example200}, the fiber of the face also 
contains two cells. The fibers thus all have the same cardinality:

\begin{proposition}[covering degree]
\label{chapter3:basics:result700}  
If $f:Y\rightarrow X$ is a covering with $X$ connected, 
then any two fibers have the same cardinality.
\end{proposition}

This common cardinality of the fibers is called the 
\emph{degree\/} of the covering, written
$$
\deg(Y\rightarrow X).
$$
The connectedness of $X$ is easily seen to be essential, 
for if $X$ has components
$X_1$ and $X_2$, and $f_i:Y_i\rightarrow X_i\,(i=1,2)$ are coverings 
of different degree then we can cobble together a new
covering $f:Y=Y_1\bigcup Y_2\rightarrow X_1\bigcup X_2$ with $f|_{Y_i}=f_i$.
The cardinality of the fibers now depends on which component of $X_i$ they lie over.
Anyway, the connectedness of $X$ is used explicitly in the proof:

\begin{proof}
If $v,u$ are vertices of $X$ and $\gamma$ a path from $v$ to $u$, then lifting $\gamma$ 
to a path $\mu$ at any
vertex of the fiber of $v$ and taking its end vertex $t(\mu)$, gives a (set) mapping from 
the fiber of $v$ to the
fiber of $u$. Interchanging the roles of $v$ and $u$ and replacing $\gamma$ with 
$\gamma^{-1}$ gives the inverse of this set map, hence we have 
a bijection between the fibers of the two vertices.

If $e\in X$ is an edge, let $u$ be a vertex in the fiber of $s(e)$. Then there is a 
unique edge $e'$ in the fiber of $e$ with $s(e')=u$. It is easy to show that the map
$u\mapsto e'$ is a bijection $f^{-1}(s(e))\rightarrow f^{-1}(e)$.

Faces are similar: let $\tau\in X$ be a face, $v\in X$ a vertex in its boundary and
$x\in\aa_\tau^{-1}(v)$. If $u$ is a vertex in the fiber of $v$ there is a face 
$\ss\in f^{-1}(\tau)$ and a $y\in Y^\ss$ with $\ve(f,\ss)(y)=x$. Let 
$d:f^{-1}(v)\rightarrow f^{-1}(\tau)$ be the map defined by $d(u)=\ss$. 
Then $d$ is injective as the attaching map of the face $d(u)$ sends $y$ to
$u$. If $\ss'\in f^{-1}(\tau)$ then $\ss'=d(u')$ for 
$u'=\aa_{\ss'}\ve(f,\ss)^{-1}(x)$, so $d$ is a surjection.
\qed
\end{proof}

Degree plays a similar role for coverings as dimension does for vector spaces or
index does for groups.
For example, ``if $U$ is a subspace of $V$ and $\dim V/U=1$, then $U=V$'', or
``if $H$ is a subgroup of index one in a group $G$ then $H=G$'',
are arguments whose combinatorial topology version is,

\begin{corollary}\label{chapter3:basics:result800}  
A degree one covering of a connected complex is an isomorphism.
\end{corollary}

The Corollary follows immediately from the surjectivity of coverings Proposition
\ref{chapter3:basics:result200}, and the definition of degree. 
This simple little result will play a crucial role in the proof of the 
Galois correspondence of \S\ref{chapter4:correspondences:correspondence}.

\subsection{Lifting and excising simply connected subcomplexes}
\label{chapter3:basics:liftingsimplyconnected}

If $X$ is $2$-complex and $Z\subset X$ a simply connected subcomplex, then we saw
in Chapter \ref{chapter2} that the quotient map $q:X\rightarrow X/Z$ induces an isomorphism
$q_*:\pi_1(X,v)\rightarrow\pi_1(X/Z,q(v))$. If $f:Y\rightarrow X$ is a covering, then
$f^{-1}(Z)\subset Y$ is a collection of isomorphic copies of $Z$, each simply connected:

\begin{proposition}[lifting simply connected complexes]
\label{chapter3:basics:result650}
Let $f:Y\rightarrow X$ be a covering 
and $Z\subset X$ a connected, simply connected subcomplex. Then
$f^{-1}(Z)\subset Y$ is a disjoint union $f^{-1}(Z)=\bigcup_I
Z_i$ with the $Z_i$ connected, simply connected, and
$f$ maps each $Z_i$ isomorphically onto $Z$.
\end{proposition}

\begin{proof}
Let $v\in Z$ be a vertex with $\pi_1(Z,v)$ trivial and 
$u_i\in Z_i$ a vertex in the fiber $f^{-1}(v)$ with $\gamma$ a
closed path at $u_i$. Then $f(\gamma)$ is a closed path at $v$, hence
homotopically trivial. Homotopy lifting gives $\gamma$ is
homotopically trivial, and thus $Z_i$ is simply connected. Connectedness 
follows by path lifting. 

Observe that $f$ restricted to any one of the $Z_i$ is a covering
$f:Z_i\rightarrow Z$. 
Suppose $u_1,u_2\in Z_i$ are vertices with $f(u_1)=f(u_2)=v\in Z$. 
If $\gamma$ is a path in $Z_i$
from $u_1$ to $u_2$ then $f(\gamma)$ is a closed
path in $Z$ at $v$, hence homotopically trivial. Homotopy lifting
gives $\gamma$ is homotopically trivial. But homotopically
trivial paths are necessarily closed, so $u_1=u_2$, the covering
$f:Z_i\rightarrow Z$ has degree one and thus is an isomorphism
by Corollary \ref{chapter3:basics:result800}. 
\qed
\end{proof}

\begin{exercise}\label{chapter3:basics:exercise510}
In the situation of Proposition \ref{chapter3:basics:result650} let
$\ss$ be a face of $X$ and $T_1,\ldots,T_k\subset X^\ss$ balls such that
$\aa_\ss^{-1}(Z)=\bigcup T_\ell$. Let $\tau$ be a face of $Y$ in the fiber of $\ss$
and $S_1,\ldots,S_k\subset Y^\tau$ balls corresponding to the $T_\ell$ via the 
isomorphism $\ve(f,\ss):Y^\tau\rightarrow X^\ss$. Show that 
$\aa_\tau^{-1}f^{-1}(Z)=\bigcup S_\ell$ (ie: the $S_\ell$ are precisely those parts
of $X^\tau$ that attach into $f^{-1}(Z)$) and that each $S_\ell$ attaches into 
a $Z_i$ for some $i$.
\end{exercise}

\begin{theorem}[excising simply connected complexes]
\label{chapter3:basics:result660}
Let $f:Y\rightarrow X$ be a covering 
with $Z\subset X$ a connected, simply connected subcomplex,
$f^{-1}(Z)=\bigcup_I Z_i$ disjoint and with the $Z_i$ connected, and
$Y/Z_i$ and $X/Z$ the resulting quotients. Then there is an
induced covering $f':Y/Z_i\rightarrow X/Z$ making the
diagram,
$$
\begin{pspicture}(0,0)(12.5,2)
\rput(1.75,-.8){
\rput(-1.3,.5){
\rput(5,2){$Y$}\rput(6.675,2){$Y/Z_i$}
\rput(4.95,0.45){$X$}\rput(6.675,0.45){$X/Z$}
\psline[linewidth=.2mm]{->}(5.3,2)(6.3,2)
\psline[linewidth=.2mm]{->}(5.3,0.45)(6.3,0.45)
\psline[linewidth=.2mm]{->}(5,1.7)(5,.7)
\psline[linewidth=.2mm]{->}(6.575,1.7)(6.575,.7)
\rput(4.75,1.2){$f$}\rput(6.9,1.175){$f'$}
\rput(5.8,2.25){$q'$}
\rput(5.8,.65){$q$}
}
}
\end{pspicture}
$$
commute, where $q,q'$ are the quotient maps.%
\end{theorem}

\begin{proof}
Let $q:X\rightarrow X/Z$ and $q':Y\rightarrow Y/Z_i$ be the bottom and
top quotient maps, and let $x':=q'(x)$ be a cell of the quotient
$Y/Z_i$. Define $f':Y/Z_i\rightarrow X/Z$ by $f'(x')=qf(x)$. As each
$Z_i$ is mapped isomorphically onto $Z$, the map $f'$ is well defined. 

Suppose $v'$ is a vertex of $X/Z$ 
and let $u'\in Y/Z_i$ be a
vertex in the fiber of $v'$. We need to show that the local continuity
maps are bijections. We leave this as an exercise for the edges.
It is immediate for the faces if $v'$ is not the vertex
$q(Z)$, as the faces incident with $u'$ and $v'$ are
unaffected by passing to the quotient. Suppose then that $v'=q(Z)$,
$u'=q(Z_i)\in Y/Z_i$ for some $i$ and $q(\ss)$ is a face of the
quotient with $\aa_{q(\ss)}^{-1}(v')=\{x_1',\ldots,x_k'\}$.
Thus there are balls $T_1,\ldots,T_k\subset X^\ss$ with the 
$T_\ell$ those parts of $X^\ss$ attaching into $Z$ and $q(T_\ell)=x_\ell$.

We show the surjectivity of local continuity first. 
Fix $x_\ell'$ and let $x_\ell\in T_\ell$,
$v=\aa_\ss(x_\ell)$ and $u$ the unique vertex of $Z_i$ covering $v$ (unique as 
$Z_i$ covers $Z$ isomorphically). Applying the covering $f$ to the triple
$v,u,\ss$ yields a face $\tau$ in $f^{-1}(\ss)$ and a $y\in Y^\tau$ attaching
to $u$ and corresponding to $x_\ell$ via the isomorphism 
$\ve(f,\tau):Y^\tau\rightarrow X^\ss$. If $S_1,\ldots,S_k$ are the balls
in $Y^\tau$ of Exercise \ref{chapter3:basics:exercise510}, then we must have
$y\in S_\ell$. Thus $y'=q(S_\ell)\in Y^{q(\tau)}$ attaches to $u'$ and maps
via $\ve(f',q(\tau))$ to $x_\ell'$.

Now injectivity: let $y_1',y_2'\in Y^{q(\tau_1)},Y^{q(\tau_2)}$ attach to $u'$ and
map to $x_\ell\in X^{q(\ss)}$ via the $\ve(f',q(\tau_i))$. Thus there are 
balls $S_{\ell 1}, S_{\ell 2}\subset Y^{\tau_1},Y^{\tau_2}$ attaching into 
$Z_i$ and corresponding to $T_\ell\subset X^\ss$ via the isomorphisms
$\ve(f,\tau_i)$. Let $x_\ell\in T_\ell$, $v=\aa_\ss(x_\ell)\in Z$ and 
$y_1,y_2\in S_{\ell i}$ correspond to $x_\ell$. Then the $y_1,y_2$ attach to
a vertex in $Z_i$ that covers $v$, and as $f$ has degree one when 
restricted to $Z_i$, they must attach to the same vertex. Applying the
covering $f$ then gives $y_1=y_2$, hence $y_1'=y_2'$ as required.
\qed
\end{proof}

\begin{example}
\label{chapter3:basics:example410}
Figure \ref{chapter3:basics:figure600} shows a degree two covering 
$f:Y\rightarrow X$ and 
$Z$ (in red) a spanning tree for $X$ (at the bottom right). A vertex 
$v\in X$ is ringed (in blue) and the two vertices of $X^\ss$ attaching to it
are also ringed. The two balls $T_1,T_2\subset X^\ss$ with 
$\aa_\ss^{-1}(Z)=T_1\bigcup T_2$ are outlined in blue. The $Z_1,Z_2$
are the lifts of $Z$ to $Y$, and the balls $S_{ij}\subset Y^{\tau_i}$ attaching
to the $Z_j$ are outlined in blue. Finally, a vertex $u\in Y$ in the 
fiber of $v$ is ringed in blue and the vertices of the $Y^{\tau_i}$ attaching
to it also (one in each face). The quotients of $X$ by $Z$ and $Y$ by the 
$Z_i$ are on the left.
\begin{figure}[h]
  \centering
\begin{pspicture}(0,0)(12,5)
\rput(0.4,0){
\rput(11,4){\BoxedEPSF{chapter3.fig5100b.eps scaled 750}}
\rput(8.9,4){\BoxedEPSF{chapter3.fig5100.eps scaled 750}}
\rput(6.8,4){\BoxedEPSF{chapter3.fig5100a.eps scaled 750}}
\rput(8.9,4){$Y$}\rput(6.8,4){$Y^{\tau_1}$}\rput(11,4){$Y^{\tau_2}$}
\rput(8.1,2.8){${\red Z_2}$}\rput(9.6,2.8){${\red Z_1}$}
}
\rput(-.4,0){
\rput(4.8,4){\BoxedEPSF{chapter3.fig5400b.eps scaled 750}}
\rput(2.9,4){\BoxedEPSF{chapter3.fig5400.eps scaled 750}}
\rput(1,4){\BoxedEPSF{chapter3.fig5400a.eps scaled 750}}
\rput(2.9,4){$Y/Z_i$}\rput(1,3){$Y/Z_i^{q(\tau_1)}$}\rput(4.85,3){$Y/Z_i^{q(\tau_2)}$}
}
\rput(1.4,-3){
\rput(9,4){\BoxedEPSF{chapter3.fig5200.eps scaled 750}}
\rput(7,4){\BoxedEPSF{chapter3.fig5100c.eps scaled 750}}
\rput(9,4){$X$}\rput(7,4){$X^\ss$}\rput(9,3.1){${\red Z}$}
}
\rput(0.4,-3){
\rput(2.8,4){\BoxedEPSF{chapter3.fig5300.eps scaled 750}}
\rput(1.1,4){\BoxedEPSF{chapter3.fig5400c.eps scaled 750}}
\rput(2.8,4){$X/Z$}\rput(1.1,5){$X/Z^{q(\ss)}$}
}
\psline[linewidth=.1mm]{->}(7,1)(4.3,1)
\psline[linewidth=.1mm]{<-}(5.4,4)(6,4)
\psline[linewidth=.1mm]{->}(2.5,3)(2.5,2)
\psline[linewidth=.1mm]{->}(9.3,3)(9.3,2)
\rput(9.1,2.5){$f$}\rput(2.7,2.5){$f'$}\rput(5.65,1.2){$q$}\rput(5.7,4.2){$q'$}
\end{pspicture}
  \caption{lifting a spanning tree.}
\label{chapter3:basics:figure600}  
\end{figure}
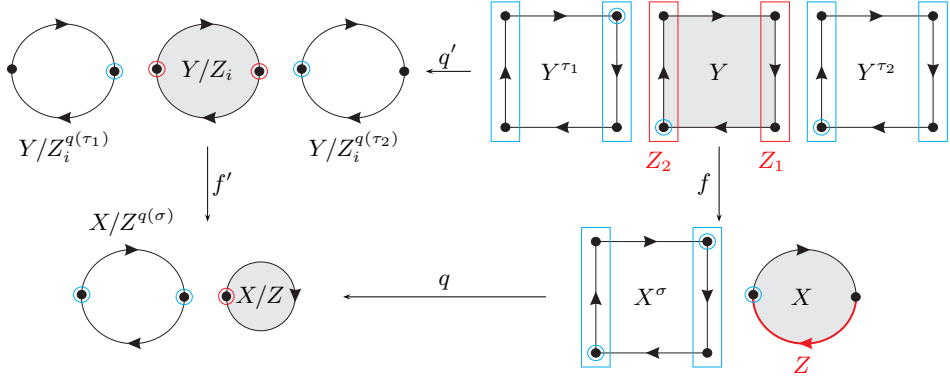
\end{example}

\begin{exercise}
\label{chapter3:basics:exercise700}
Let $f:Y\rightarrow X$ be a covering with $X$ connected and $Z\subset X$ a connected,
simply connected subcomplex with $f':Y/Z_i\rightarrow X/Z$ the induced covering
of Theorem \ref{chapter3:basics:result660}. Show that 
$$
\deg(Y\rightarrow X)=\deg(Y/Z_i\rightarrow X/Z).
$$
\end{exercise}

\section{Actions, intermediate and universal covers}\label{chapter3:actionsuniversal}

\subsection{Group actions}\label{chapter3:actionsuniversal:actions}

Recall from \S\ref{chapter1:2complexes:maps}
that a group acts freely on a $2$-complex precisely when it acts freely
on the vertices, ie: for any $g\in G$ and vertex $v$, if $g(v)=v$ then $g$ is the identity.
Such group actions give coverings:

\begin{proposition}[free actions give covers]
\label{chapter3:actionsuniversal:result100}
If a group $G$ acts orientation preservingly and freely on a $2$-complex $X$
then the quotient map $q:X\rightarrow X/G$ is a covering.
\end{proposition}

Our main supply of free group actions will come from the Galois group 
of a cover $Y\rightarrow X$ in Chapter \ref{chapter4}.

\begin{proof}
Let $[v]$ be a vertex in $X/G$ and $u\in X$ with $q(u)=[v]$, hence $u=g(v)$
for some $g\in G$. An edge of $X/G$ starting at $[v]$ has the form $[e]$
with $s(e)=v$. In particular $g(e)$ has start $u$ and $qg(e)=[e]$, giving the
surjectivity of (edge) local continuity. If $e_1,e_2$ start at $u$ with $q(e_1)=q(e_2)$
then $e_2=g'(e_1)$ for $g'\not=1$ and $g'(u)=u$, contradicting the freeness of the 
$G$-action. Thus we have injectivity of (edge) local continuity and 
$q:X^{(1)}\rightarrow X/G^{(1)}$ a graph covering.

Now let $[\ss]$ be a face of the quotient containing $[v]$ in its boundary. We have
$\partial[\ss]=(X^\ss,\aa_{[\ss]}=q\aa_\ss)$ with $\aa_{[\ss]}^{-1}[v]$ those vertices of 
$X^\ss$ sent by $\aa_\ss$ into the equivalence class $[v]$. For $u\in X$ with $q(u)=[v]$, 
the set $\bigcup_{\tau\in[\ss]}\aa_\tau^{-1}(u)$ consists of theose vertices of
$\bigcup_{\tau\in[\ss]}X^\tau$ that attach to $u$, and we require
\begin{equation}
  \label{chapter3:equation200}
\amalg\ve(g_\tau,\ss)^{-1}:\bigcup_{\tau\in[\ss]}\aa_\tau^{-1}(u)
\rightarrow\aa_{[\ss]}^{-1}[v]
\end{equation}
to be a bijection, where $g_\tau(\ss)=\tau$. Suppose that $\tau,\ww\in[\ss]$ and
$y\in X^\tau,z\in X^\ww$ map via (\ref{chapter3:equation200}) to $x\in X^\ss$. As
$\aa_\tau(y)=\aa_\ww(z)=u$, the elements $g_\tau,g_\ww\in G$ both map $\aa_\ss(x)$
to $u$. As the $\ve(g_\tau,\ss)$ are injective, we must have $\tau\not=\ww$,
hence $g_\tau\not=g_\ww$, so $g_\tau g_\ww^{-1}$ is a non-identity element fixing $u$,
contradicting the freeness of the $G$-action. The map (\ref{chapter3:equation200})
is thus an injection. 

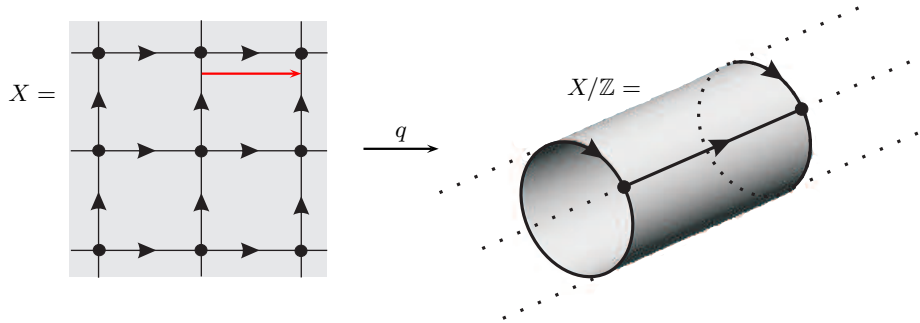
\begin{figure}
  \centering
\begin{pspicture}(0,0)(13,3.75)
\rput(0,0.25){
\rput(3,1.75){\BoxedEPSF{fig8.eps scaled 1000}}
\rput(.8,2.5){$X=$}
}
\psline[linewidth=.3mm,linecolor=red]{->}(3.05,3)(4.375,3)
\rput(.4,0){
\rput(9,1.8){\BoxedEPSF{chapter1.fig900.eps scaled 600}}
\rput(8,2.8){$X/\Z=$}
}
\psline[linewidth=.3mm]{->}(5.2,2)(6.2,2)
\rput(5.7,2.2){$q$}
\end{pspicture}
\caption{covering $X\rightarrow X/\Z$ given by the free $\Z$-action of Figure 
\ref{chapter1:2complexes:figure1000}.}
\label{chapter3:basics:figure650}  
\end{figure}

Now let $x\in X^\ss$ attach to $g(v)\in [v]$ and let $y=\ve(g,\ss)(x)\in X^{g(\ss)}$.
The attaching map of $X^{g(\ss)}$ is the composition $g^{-1}\aa_\ss\ve(g,\ss)^{-1}$,
sending $y$ to $v$, hence $y$ lies in the left hand side of (\ref{chapter3:equation200}).
\qed
\end{proof}

\subsection{Intermediate covers}
\label{chapter3:actionsuniversal:intermediate}

Suppose we have a commuting triangle of complexes and maps,
$$
\begin{pspicture}(0,0)(4,1.5)
\rput(0,-.25){
\rput(1,1.5){$Y$}
\rput(3,1.5){$Z$}
\rput(2,0.4){$X$}
\rput(2,1.7){$g$}
\rput(1.125,.825){$f$}
\rput(2.78,.825){$h$}
\psline[linewidth=.1mm]{->}(1.3,1.5)(2.7,1.5)
\psline[linewidth=.1mm]{->}(1.2,1.25)(1.8,.55)
\psline[linewidth=.1mm]{->}(2.8,1.25)(2.25,.55)
}
\end{pspicture}
$$
with all three maps coverings. We say that the covers 
$Y\stackrel{g}{\rightarrow}Z\stackrel{h}{\rightarrow}X$ are \emph{intermediate\/}
to $f:Y\rightarrow X$. We will see in \S\ref{chapter3:lattices} that the set of coverings 
intermediate to a fixed covering $f:Y\rightarrow X$ has a very nice structure. 

\begin{proposition}
\label{chapter3:actionsuniversal:result150}
Let $Y\stackrel{g}{\rightarrow}Z\stackrel{h}{\rightarrow}X$ be dimension preserving
maps of $2$-complexes with $f=hg:Y\rightarrow X$. If any two of $f,g,h$ are coverings, then
so is the third.
\end{proposition}

\begin{proof}
We do one of the three cases and leave the other two as an exercise. Suppose then that
$f$ and $h$ are coverings. We need to show that the local continuity maps for $g$ in 
(C2) and (C3) of Definition \ref{chapter3:basics:definition200} are bijections. Let $w\in Z$,
$v\in Y$ be vertices with $g(v)=w$ and $e\in Z$ an edge with $s(e)=w$. Injectivity 
is easiest: if $e_1,e_2\in Y$ with $s(e_i)=v$ and $g(e_i)=e$ then $f(e_i)=hg(e_i)=h(e)$.
As $f$ is a covering we get $e_1=e_2$. To find an edge at $v$ covering $e$, lift $h(e)$ 
to an edge $e''$ at $v$ and let $e'=g(e'')$. Then both $e,e'$ start at $w$ and cover $h(e)$,
so $e=e'$ as $h$ is a covering, and $g(e'')=e$ as required. Local continuity at a face $\ss$
with $w$ lying in its boundary is completely analogous.
\qed
\end{proof}

\begin{exercise}
\label{chapter3:actionsuniversal:exercise100}
Let $Y$ be a graph and $Y_1,Y_2\subset Y$
subgraphs of the form,
$$
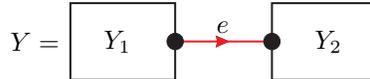
\begin{pspicture}(14,1)
\rput(4,0){
\rput(3,.5){\BoxedEPSF{free35.eps scaled 750}}
\rput(.5,.5){$Y=$}
\rput(1.6,.5){$Y_1$}\rput(4.4,.5){$Y_2$}
\rput(3,.7){$e$}
}
\end{pspicture}
$$
\begin{description}
\item[(i).] If $Y_1$ is a tree, $f:Y\rightarrow X$,
$h:Z\rightarrow X$ coverings with $X$ having a single vertex, and
$p:Y_2\hookrightarrow Z$ a subgraph,
then
there is an intermediate
covering $Y\stackrel{g}{\rightarrow}Z\stackrel{h}{\rightarrow}X$.
\item[(ii).] If $W\rightarrow Y$ is a covering and $Y_1$ a tree,
then $W$ also has the form shown above for some subgraphs 
$Y'_1,Y'_2\subset W$, and with $Y'_1$ a tree.
\end{description}
\end{exercise}

\subsection{Covers from the ``bottom up''}
\label{chapter3:actionsuniversal:bottom_up}

Much of the discussion of coverings so far as been in the abstract: we haven't seen many
actual covers! In this book we will construct specific examples in two ways that can be
broadly
described as ``bottom-up'' and ``top-down''. The first of these, 
which we describe in this section,
starts with a complex $X$ and builds upwards to give a covering of it. The other, which is
described in \S\ref{chapter4:galoisgroups:subgroupstocovers}, 
starts with a covering of $X$ and folds it down
into a smaller covering. In both cases how far to build up, or how far to fold down,
is governed by a subgroup of a certain group, although as it turns out,
a different group in the two cases.

For the bottom-up cover we imitate a standard construction in topology:

\begin{definition}[$\mathbf{1}$-skeleton of ``bottom-up'' cover]
\label{chapter3:actionsuniversal:definition50}
Let $X$ be a $2$-complex and $H\subset\pi_1(X,v)$ a subgroup. Define
$X\kern-2pt\uparrow\kern-2pt H$ to be the following graph:
\begin{description}
\item[(1).] The vertices of $X\kern-2pt\uparrow\kern-2pt H$
are the equivalence classes of paths starting at
$v$ under the following relation:
$\gamma_1\sim\gamma_2$ if and only if
$\gamma_1\gamma_2^{-1}$ represents the homotopy class of an element of $H$ (and
 so in particular is a closed path). Write $u_\gamma$ for the vertex with representative path
$\gamma$.
\item[(2).] Let $e$ be an edge of $X$ and $u_\gamma,u_\mu$ vertices of 
$X\kern-2pt\uparrow\kern-2pt H$.
Then there is an edge $e'$ of $X\kern-2pt\uparrow\kern-2pt H$ with start vertex
$u_\gamma$ and finish vertex $u_\mu$ if and only if $\gamma e\mu^{-1}$ represents
the homotopy class of an element of $H$.
\end{description}
\end{definition}

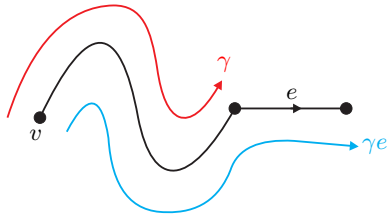
\begin{figure}
  \centering
\begin{pspicture}(12.5,3)
\rput(6.25,1.4){\BoxedEPSF{chapter3.fig2700.eps scaled 350}}
\rput(7.7,1.6){$e$}\rput(4.3,1.05){$v$}
\rput(6.8,2){${\red\gamma}$}\rput(8.8,.9){${\cyan\gamma e}$}
\end{pspicture}
  \caption{how an edge of $X$ gives rise to an edge of $X\kern-2pt\uparrow\kern-2pt H$: 
there is an edge $e'$
with start vertex the equivalence class of $\gamma$ and finish vertex the equivalence class
of $\gamma e$.}
  \label{chapter3:actionsuniversal:figure200}
\end{figure}

In particular, edges of $X\kern-2pt\uparrow\kern-2pt H$ arise via the scheme illustrated in
Figure \ref{chapter3:actionsuniversal:figure200}: in this case $\gamma e (\gamma e)^{-1}$
represents the homotopy class of the identity element of $H$.

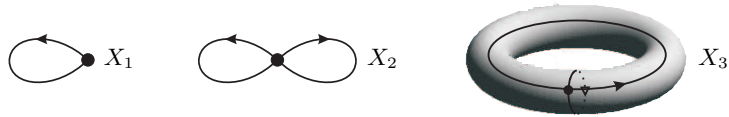
\begin{figure}
  \centering
\begin{pspicture}(12.5,2)
\rput(-.5,0){
\rput(3,1){\BoxedEPSF{chapter3.fig2000.eps scaled 650}}
\rput(6,1){\BoxedEPSF{chapter3.fig2100.eps scaled 650}}
\rput(10,1){\BoxedEPSF{fig121.eps scaled 350}}
\rput(3.9,1){$X_1$}\rput(7.4,1){$X_2$}\rput(11.8,1){$X_3$}
}
\end{pspicture}
  \caption{three complexes: $X_2$ and $X_3$ have the same $1$-skeleton.}
  \label{chapter3:actionsuniversal:figure300}
\end{figure}

\begin{figure}
  \centering
\begin{pspicture}(12.5,4.5)
\rput(0,.25){
\rput(-1.5,0){
\rput(2.5,2){\BoxedEPSF{chapter3.fig2600.eps scaled 800}}
\rput(6.5,2){\BoxedEPSF{fig18a.eps scaled 500}}
\rput(11.5,2){\BoxedEPSF{chapter3.fig2500.eps scaled 1000}}
}
}
\end{pspicture}
  \caption{the $1$-skeleton of $X_i\kern-2pt\uparrow\kern-2pt H$ for the $X_i$ of Figure
\ref{chapter3:actionsuniversal:figure300} and $H$ the identity subgroup.}
  \label{chapter3:actionsuniversal:figure400}
\end{figure}

Figure \ref{chapter3:actionsuniversal:figure300} gives complexes 
$X_i\,(i=1,2,3)$ and the
graphs $X_i\kern-2pt\uparrow\kern-2pt H$ for $H$ the identity subgroup
are given in Figure \ref{chapter3:actionsuniversal:figure400}.
Note that for $H$ trivial, paths $\gamma_1\sim\gamma_2$ if
and only if they are homotopic.
If $X_1$ is the $S^1$-graph on the left
of Figure \ref{chapter3:actionsuniversal:figure300}, 
then there
is a 1-1 correspondence between the homotopy classes of paths and paths of the form
$e\ldots e$ ($k$ times) or $e^{-1}\ldots e^{-1}$ ($k$ times). 
Thus $X\kern-2pt\uparrow\kern-2pt H$ has
vertices $u_k$ for $k\in\Z$. There is an edge connecting the vertex of the path
$e\ldots e$ ($k$ times) to the vertex of the path $e\ldots e$ ($k+1$ times), to give 
the infinite $2$-valent tree at left in Figure \ref{chapter3:actionsuniversal:figure400}.
Similarly the single-vertexed graph with two edges 
has $X\kern-2pt\uparrow\kern-2pt H$ the $4$-valent infinite
tree. 

The last of the three complexes is the torus, which has exactly the same $1$-skeleton as $X_2$,
but the presence of a face drastically changes the graph $X_i\kern-2pt\uparrow\kern-2pt H$.
Paths of the form $\gamma$ and $\gamma e_1e_2e_1^{-1}e_2^{-1}$ give distinct vertices
in $X_2\kern-2pt\uparrow\kern-2pt H$
but the \emph{same\/} vertex in $X_3\kern-2pt\uparrow\kern-2pt H$,
forcing the $4$-valent tree to bend into the grid shape shown.

\begin{exercise}
\label{chapter3:actionsuniversal:exercise300}
Let $u_\gamma$ and $u_\mu$ be vertices of $X\kern-2pt\uparrow\kern-2pt H$.
If there is a path 
$e_{1}'\ldots e_{k}'$ in $X\kern-2pt\uparrow\kern-2pt H$ from $u_\gamma$ to $u_\mu$,
then there is a path $e_1\ldots e_k$ in $X$ with 
the edge $e_{i}'$ arising from the edge $e_i$ as in Definition 
\ref{chapter3:actionsuniversal:definition50}, and 
$\gamma_1e_1\ldots e_k\gamma_2^{-1}$ an element of $H$.
\end{exercise}

\begin{proposition}
\label{chapter3:actionsuniversal:result400}
The graph $X\kern-2pt\uparrow\kern-2pt H$ is connected. 
Define $f:X\kern-2pt\uparrow\kern-2pt H\rightarrow X^{(1)}$ by
sending a vertex $u_\gamma$ to 
the terminal vertex of $\gamma$, 
and $f(e')=e$, where the edge $e'$ arises from $e$ as in Definition
\ref{chapter3:actionsuniversal:definition50}(2). Then
$f$ is a graph covering.
\end{proposition}

\begin{proof}
If $u_\gamma$ is a vertex of $X\kern-2pt\uparrow\kern-2pt H$ with
$\gamma=e_1\ldots e_k\in X$, then $e_{1}'\ldots e_{k}'\in X\kern-2pt\uparrow\kern-2pt H$ 
is a path 
connecting $u_\gamma$ to $u_\varnothing=v$, and
so $X\kern-2pt\uparrow\kern-2pt H$ is connected. 
To see that $f$ is a covering we need it to be
a dimension preserving map (which we leave to the reader) and
for every pair of vertices $u,x$ with $f(u)=x$, it induces
a bijection from the edges starting at $u$ to the edges starting at $x$.
Suppose then that $u=u_\gamma$ and $e_1$, $e_2$ are edges 
of $X\kern-2pt\uparrow\kern-2pt H$ connecting 
$u_\gamma$ to vertices $u_1$ and $u_2$. Let $x$ be the terminal
vertex of $\gamma$, so that $f(u_\gamma)=x$, and $e$ an edge starting at $x$
with $f(e_1)=f(e_2)=e$. 
Thus there are paths $\mu_1,\mu_2$ from $v$ to the terminal vertex of
$e$ with $\gamma e\mu_i^{-1} (i=1,2)$ representing an element of $H$.
In particular, $(\gamma e\mu_1^{-1})^{-1}\gamma e\mu_2^{-1}$, which 
is homotopic to $\mu_1\mu_2^{-1}$, represents an element of $H$, and so
$u_1=u_2$. The edges $e_i$ both arise by applying the
construction of Definition \ref{chapter3:actionsuniversal:definition50}(2)
to the pair of vertices $u$ and $u_1=u_2$,
and as only one edge can arise this way we have $e_1=e_2$. 
The local continuity maps are thus injective. For an edge $e$ starting at $v$,
and a vertex $u_\gamma$ in the fiber of $v$ with the path $\gamma$ from
$v$ to $x$, there is by definition an edge $e'$ connecting
$u_\gamma$ and $u_{\gamma e}$.
The local continuity maps are thus surjective.
\qed
\end{proof}

\begin{lemma}\label{chapter3:actionsuniversal:result450}
In the graph covering $f:X\kern-2pt\uparrow\kern-2pt H\rightarrow X^{(1)}$, 
the boundaries of faces of $X$
lift to closed paths in $X\kern-2pt\uparrow\kern-2pt H$. 
More precisely, let $x$ be a vertex
of $X$, $\ss$ a face containing $x$ in its boundary and $\gamma_\ss$ a boundary path of $\ss$
starting at $x$. If $u$ is a vertex covering $x$ via the graph covering $f$ and 
$\mu_\ss$ is the lift of $\gamma_\ss$ at $u$, then $\mu_\ss$ is a closed path
in $X\kern-2pt\uparrow\kern-2pt H$.
\end{lemma}

The proof is left as an excercise.
Thus the boundaries of faces in $X$ give rise to closed paths in 
$X\kern-2pt\uparrow\kern-2pt H$, and to
construct the $2$-skeleton of $X\kern-2pt\uparrow\kern-2pt H$
we ``sew'' faces into these closed paths:

\begin{definition}[``bottom-up'' cover]
\label{chapter3:actionsuniversal:definition75}
Let $X$ be a $2$-complex, $H\subset\pi_1(X,v)$ a subgroup and 
$X\kern-2pt\uparrow\kern-2pt H$ the graph of 
Definition \ref{chapter3:actionsuniversal:definition50}.
We add faces in the following way: 
for each face $\ss\in X$, fix a boundary label $\gamma_\ss$ with start vertex $x$ and let
$f^{-1}(x)=\{u_i\,|\,i\in I\}$ be the fiber of $x$ via the graph covering
$f:X\kern-2pt\uparrow\kern-2pt H\rightarrow X^{(1)}$. For the lift of $\gamma_\ss$ to 
each $u_i$, define a face $\ss_i$ with boundary this closed path: 
$\partial\ss_i=(X^\ss,\aa_{\ss_i})$ where $\aa_{\ss_i}=f^{-1}\aa_\ss$.
\end{definition}

There are quite a few choices made in this construction.
We will see that the complex is independent (upto 
isomorphism) of these choices in \S\ref{chapter4:correspondences:universal}.

\begin{proposition}
\label{chapter3:actionsuniversal:result500}
Let $X\kern-2pt\uparrow\kern-2pt H$ be the $2$-complex of
Definition \ref{chapter3:actionsuniversal:definition75} and
define $f:X\kern-2pt\uparrow\kern-2pt H\rightarrow X$ on the $1$-skeleton as in Lemma
\ref{chapter3:actionsuniversal:result400}, 
and for each face $\ss_i$ arising from the face $\ss$ of $X$ as in 
Definition \ref{chapter3:actionsuniversal:definition75}, define
$f(\ss_i)=\ss$.
Then $f$ is a covering of $2$-complexes with
$$
f_*\pi_1(X\kern-2pt\uparrow\kern-2pt H,u_\varnothing)=H\subset\pi_1(X,v).
$$
\end{proposition}

\begin{proof}
We have a path $\gamma'$ in $X\kern-2pt\uparrow\kern-2pt H$ 
from the vertex $u_\mu$ to the vertex $u_\nu$ if and only if there is a path 
$\gamma$ in $X$ from the terminal vertex of $\mu$ to the terminal vertex of $\nu$ with
$\mu\gamma\nu^{-1}$ homotopic to an element of $H$. 
In particular, if $\gamma'$ is a closed path at $u_\varnothing$, then $\gamma$ is a closed
path at $v$ homomtopic to an element of $H$, ie: we have
$f_*\pi_1(X\kern-2pt\uparrow\kern-2pt H,u)=H$ as claimed.
For $f$ to be a covering we need bijective local continuity on the faces. 
Thus, let $z\in X$ be a vertex, $\ss$ a face containing $z$ in its boundary and
$y\in X\kern-2pt\uparrow\kern-2pt H$ a vertex in the fiber of $z$. 
Suppose also that in the 
construction of $X\kern-2pt\uparrow\kern-2pt H$ we chose as boundary path for
$\ss$ the image under $\aa_\ss$ of the path $\gamma_\ss$ circumnavigating $X^\ss$, and 
suppose that $\gamma_\ss=\gamma_0\gamma_1\ldots\gamma_k$, where $z$ appears $k$ times in the
boundary of $\ss$ and $\gamma_0\ldots\gamma_i$ terminates at the $i$-th of these
appearences. Lift $\aa_\ss(\gamma_0\ldots\gamma_i)^{-1}$ to $y$ for each $i$. By 
definition, there is a face of $X\kern-2pt\uparrow\kern-2pt H$ arising by lifting $\gamma_\ss$
to the terminal vertices of each of these lifts. We leave it to the reader to show that
these are precisely the faces in the fiber of $\ss$ that $y$ appears in the boundary of, 
and that it appears exactly $k$ times.
\qed  
\end{proof}

There is an alternative construction of the complex
$X\kern-2pt\uparrow\kern-2pt H$ that is more group
theoretic, at least at the level of the $1$-skeleton.

\begin{figure}[h]
  \centering
  \begin{pspicture}(12.5,4)
\rput(2.5,2){\BoxedEPSF{chapter3.fig4000.eps scaled 500}}
\rput(9,2){\BoxedEPSF{chapter3.fig4000a.eps scaled 500}}
\rput(1.5,.5){$X$}\rput(3.4,3){$T$}\rput(3,.7){${\red e}$}
\rput(8,3){$T_i$}\rput(11.6,3){$T_j$}\rput(9.5,.4){${\red e'}$}
\psline[linewidth=.1mm]{<-}(4.5,2)(5.2,2)
\rput(4.85,2.2){$f$}\rput(9,3){$X\kern-2pt\uparrow\kern-2pt H$}
\end{pspicture}
  \caption{Alternative construction of the $1$-skeleton of $X\kern-2pt\uparrow\kern-2pt H$.}
  \label{chapter3:actionsuniversal:figure450}
\end{figure}
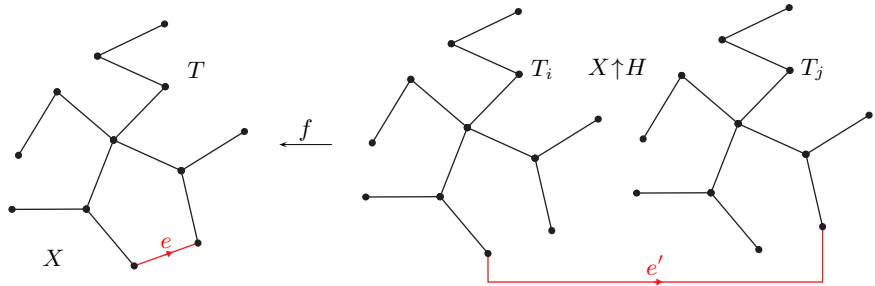

\begin{definition}[``bottom-up'' cover: version $\mathbf{2}$]
\label{chapter3:actionsuniversal:definition85}
Let $X$ be a $2$-complex and $H\subset\pi_1(X,v)$ a subgroup. Let 
$X\kern-2pt\uparrow\kern-2pt H$ be the graph 
defined as follows:
\begin{description}
\item[(1).] Let $T\subset X$ be a spanning tree and $\{g_i\,|\,i\in I\}$
be a set of (right) coset representatives for the subgroup $H$ in $\pi_1(X,v)$.
For each $i\in I$ let $T_i$ be an isomorphic copy of $T$.
\item[(2).] We now add edges to $\bigcup T_i$: let $e$ be an edge of $X$
\emph{not\/} in $T$
with start vertex $u_1$ and terminal vertex $u_2$. Let $\gamma_e$ be the
(reduced) path that travels through $T$
from $v$ to $v_1$, traverses $e$ and then travels through $T$ from $v_2$ to $v$. 
For $i\in I$ let $u_{i_1},u_{i_2}$ be the vertices of $T_i$ corresponding to 
$v_1,v_2$ under the isomorphism $T_i\cong T$. Then there is an edge $e'$
with start $u_{i_1}$ and terminal vertex $u_{j_2}$ 
if and only if $Hg_i\gamma_e=Hg_j$. See Figure \ref{chapter3:actionsuniversal:figure450}
\end{description}
\end{definition}

Define $f:X\kern-2pt\uparrow\kern-2pt H\rightarrow X^{(1)}$ by 
$f(u_i)=u$, where $u_i\in T_i$ corresponds to $u$ via the isomorphism
$T_i\cong T$, and $f(e')=e$ where $e'\in T_i$ arises from $e\in T$ as
above.

\begin{exercise}
Show that $X\kern-2pt\uparrow\kern-2pt H$ is connected and 
$f:X\kern-2pt\uparrow\kern-2pt H\rightarrow X^{(1)}$ is a covering for which the 
conclusions of Exercise \ref{chapter3:actionsuniversal:exercise300} hold.
\end{exercise}

The remainder of the construction is as in Definition \ref{chapter3:actionsuniversal:definition75}.
We leave it 
as an Exercise to show that the claims in Proposition \ref{chapter3:actionsuniversal:result500}
hold for $X\kern-2pt\uparrow\kern-2pt H$.
Once again there is choice in the construction, and these ambiguities will be ironed out
in \S\ref{chapter4:correspondences:universal}.

Because the $1$-skeletons are given by the cosets of the subgroup $H$, these bottom-up
covers are called \emph{Schreier coset diagrams\/}.

\subsection{Universal covers}\label{chapter3:actionsuniversal:universal}

A complex
is always a cover of itself (Exercise \ref{chapter3:basics:exercise160}).
In this section we show that a complex always
has another cover at the other extreme, in that it is as ``big'' as possible.

\begin{definition}[universal covers]
\label{chapter3:actionsuniversal:definition100}
A covering $f:Y\rightarrow X$ is \emph{universal\/} if and only if for any covering
$h:Z\rightarrow X$ there is a covering $g:Y\rightarrow Z$ making the diagram,
$$
\begin{pspicture}(4,1.5)
\rput(0,-.25){
\rput(1,1.5){$Y$}
\rput(3,1.5){$Z$}
\rput(2,0.4){$X$}
\rput(2,1.7){$g$}
\rput(1.125,.825){$f$}
\rput(2.78,.825){$h$}
\psline[linewidth=.1mm]{->}(1.3,1.5)(2.7,1.5)
\psline[linewidth=.1mm]{->}(1.2,1.25)(1.8,.55)
\psline[linewidth=.1mm]{->}(2.8,1.25)(2.25,.55)
}
\end{pspicture}
$$
commute.
\end{definition}

Equivalently, $Y\rightarrow X$ is universal when any other covering $Z\rightarrow X$ 
of $X$ is intermediate to it.

The construction of a universal cover is a special case of the techniques of
the previous section: write $\wtl{X}$ for the complex $X\kern-2pt\uparrow\kern-2pt H$
obtained when $H$ is the identity subgroup of $\pi_1(X,v)$.

\begin{proposition}
\label{chapter3:actionsuniversal:result700}
The $2$-complex $\wtl{X}$ is connected, simply connected and
the covering $f:\wtl{X}\rightarrow X$ of Proposition
\ref{chapter3:actionsuniversal:result500} is universal.
\end{proposition}

\begin{proof}
That $\wtl{X}$ is connected and simply connected is immediate. Let $h:Y\rightarrow X$
be a cover. Map lifting (Theorem \ref{chapter3:basics:result600}) 
gives a map $g:\wtl{X}\rightarrow Y$ with $f=hg$, as the fundamental group of 
$\wtl{X}$ is trivial. Proposition \ref{chapter3:actionsuniversal:result150} gives $g$ is 
a cover.
\qed
\end{proof}

\begin{example}
\begin{figure}
  \centering
\begin{pspicture}(12.5,4)
\rput(2,0.25){
\rput(1,1){\BoxedEPSF{chapter3.fig4100a.eps scaled 900}}
\rput(1,1){$X^\ss$}
\rput(-.2,1){$e_1$}
\rput(2.25,1){$e_1^{-1}$}
\rput(1,2.15){$e_2$}
\rput(1,-.15){$e_2^{-1}$}
}
\psline[linewidth=.3mm]{->}(3,2.6)(3,3.2)
\rput(0,0){
\rput(3,3.75){$u$}
\rput(1.6,3.5){$e_1$}
\rput(4.4,3.5){$e_2$}
\rput(3,3.5){\BoxedEPSF{chapter3.fig2100.eps scaled 750}}
}
\rput(0,0){
\rput(9,2){\BoxedEPSF{chapter3.fig4100.eps scaled 1000}}
\rput(8.9,2.2){$v$}
}
\end{pspicture}
  \caption{universal cover of the torus: fix $\gamma$ for the face $\ss$ of $X$ as 
shown. The lifts of $\aa_\ss(\gamma)$ to the vertices in the fiber of $u$ are
shown in $\wtl{X}$ on the right.}
  \label{chapter3:actionsuniversal:figure600}
\end{figure}
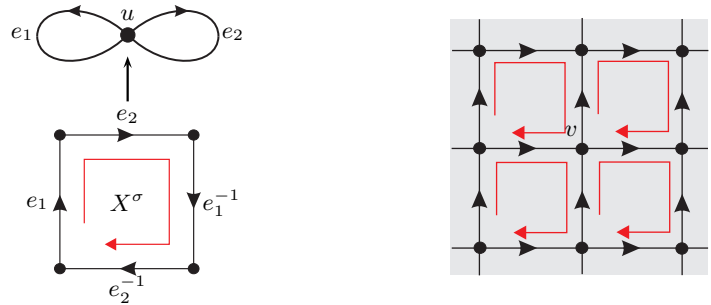

Figure \ref{chapter3:actionsuniversal:figure600} shows the result of performing this process
with the complex $X_3$ of Figure \ref{chapter3:actionsuniversal:figure300}, sewing 
faces onto the $1$-skeleton of Figure \ref{chapter3:actionsuniversal:figure400}.
Figure \ref{chapter3:actionsuniversal:figure650} shows the effect on $\wtl{X}$ of an 
extra face in $X$.  
\end{example}

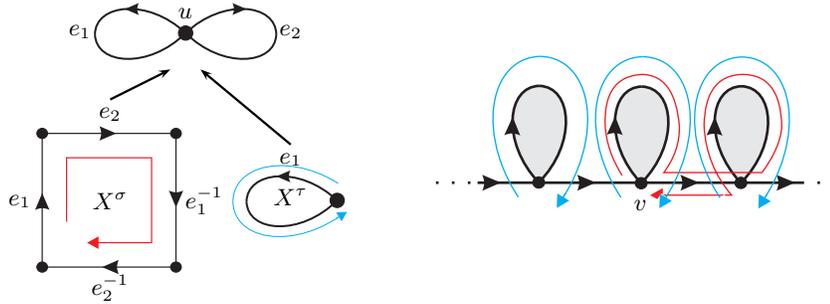
\begin{figure}
  \centering
\begin{pspicture}(12.5,4)
\rput(1,0.25){
\rput(1,1){\BoxedEPSF{chapter3.fig4100a.eps scaled 900}}
\rput(1,1){$X^\ss$}
\rput(-.2,1){$e_1$}
\rput(2.25,1){$e_1^{-1}$}
\rput(1,2.15){$e_2$}
\rput(1,-.15){$e_2^{-1}$}
}
\rput(.5,0){
\rput(4,1.25){\BoxedEPSF{chapter3.fig4200a.eps scaled 750}}
\rput(3.9,1.8){$e_1$}\rput(3.9,1.3){$X^\tau$}
}
\psline[linewidth=.3mm]{->}(2,2.6)(2.8,3)
\psline[linewidth=.3mm]{->}(4.4,2)(3.2,3)
\rput(0,0){
\rput(3,3.75){$u$}
\rput(1.6,3.5){$e_1$}
\rput(4.4,3.5){$e_2$}
\rput(3,3.5){\BoxedEPSF{chapter3.fig2100.eps scaled 750}}
}
\rput(0,0){
\rput(9,2){\BoxedEPSF{chapter3.fig4200.eps scaled 1000}}
\rput(9.05,1.2){$v$}
}
\end{pspicture}
  \caption{Complex $X$ (left) obtained by sewing another face 
$\tau$ onto the torus and its universal cover $\wtl{X}$ (right).}
  \label{chapter3:actionsuniversal:figure650}
\end{figure}

\begin{figure}
  \centering
\begin{pspicture}(12.5,4)
\rput(.5,0.25){
\rput(1,1){\BoxedEPSF{chapter3.fig4100a.eps scaled 900}}
\rput(-.2,0){$u$}
\rput(1,1){$X^\ss$}
\rput(-.2,1){$e_1$}
\rput(2.25,1){$e_1^{-1}$}
\rput(1,2.15){$e_2$}
\rput(1,-.15){$e_2^{-1}$}
}
\rput(4,0.25){
\rput(1,1){\BoxedEPSF{chapter3.fig4100b.eps scaled 900}}
\rput(-.2,0){$v$}
\rput(1,1){$X^\tau$}
\rput(-.2,1){$e_1$}
\rput(2.25,1){$e_1^{-1}$}
\rput(1,2.15){$e_2$}
\rput(1,-.15){$e_2^{-1}$}
}
\psline[linewidth=.3mm]{->}(2,2.6)(2.8,3)
\psline[linewidth=.3mm]{->}(4,2.6)(3.2,3)
\rput(0,0){
\rput(1.9,3.7){$u$}\rput(4.1,3.7){$v$}
\rput(3,3.75){$e_1$}\rput(3,3.25){$e_1$}
\rput(.5,3.5){$e_2$}\rput(5.5,3.5){$e_2$}
\rput(3,3.5){\BoxedEPSF{chapter3.fig4300.eps scaled 750}}
}
\rput(0,0){
\rput(9,2){\BoxedEPSF{chapter3.fig4100c.eps scaled 1000}}
\rput(9.05,1.2){$v$}
}
\end{pspicture}
  \caption{universal cover of an $X$ that is a degree two cover of the 
torus. The red boundary path for $\ss$ starts at $u$ whereas the blue boundary path
for $\tau$ starts at $v$.}
  \label{chapter3:actionsuniversal:figure625}
\end{figure}
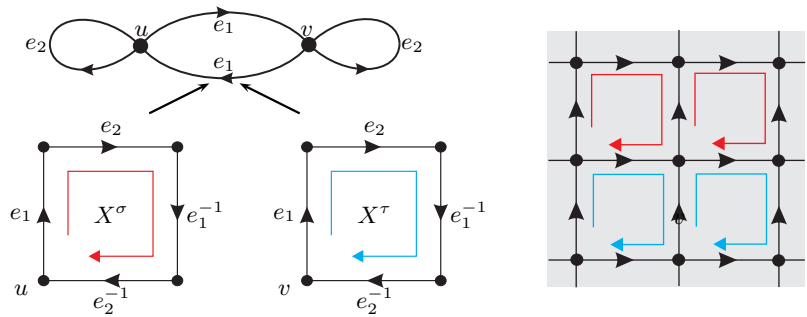

\begin{example}
Figure \ref{chapter3:actionsuniversal:figure625} shows the universal cover of 
an $X$ that is itself a degree two cover of the torus of Figure 
\ref{chapter3:actionsuniversal:figure600}.
\end{example}

\begin{exercise}
\label{chapter3:actionsuniversal:exercise400}
Show, using universal coverings, 
that if $X$ is a graph and $\gamma_1,\gamma_2$ are reduced homotopic paths 
in $X$ with the 
same start vertex, then $\gamma_1=\gamma_2$.
[\emph{hint\/}: lift the paths to $\wtl{X}$ and use properties of reduced paths in trees
to deduce that these lifts are identical.]
\end{exercise}

Here is one final result, which we will save up for later 
(\S\ref{chapter4:correspondences:universal}):

\begin{corollary}
\label{chapter3:actionsuniversal:exercise500}
Let $H$ be a subgroup of $\pi_1(X,v)$ and
let $x=H$.
Then we have an
intermediate covering $\wtl{X}_u\rightarrow
(X\kern-2pt\uparrow\kern-2pt H)_x\stackrel{f}{\rightarrow}X_v$ 
with $f_*\pi_1(X\kern-2pt\uparrow\kern-2pt H,x)=H$.
\end{corollary}

\subsection{Monodromy}\label{chapter3:actionsuniversal:monodromy}

If $f:Y\rightarrow X$ is a covering we will eventually get an action of two 
different groups
on $Y$ or parts of $Y$. The more important of these is the Galois group
of the covering, which forms the principle subject
of Chapter \ref{chapter4}.

The less important of these two actions is that of the 
fundamental group $\pi_1(X,v)$ on the fiber
$f^{-1}(v)$ of a vertex $v\in X$. 
Thus the fundamental group acts as a permutation group on the
set of vertices covering $v$. Such permutation representations of the
fundamental groups of $2$-complexes will play a key role in the proof of results
like Miller's theorem in a later Chapter.

To define this action, see that it makes sense, and is indeed a homomorphism
$$
\pi_1(X,v)\rightarrow\sym(f^{-1}(v)),
$$
we require no more than the path and homotopy lifting of \S\ref{chapter3:basics:lifting}.
So, the ``path-lifting action'' would probably be a sensible name: 
the action would then do exactly what
it says on the box! However, it is traditional in topology to 
call this action \emph{monodromy\/}, and so we 
will too.

\begin{figure}
  \centering
\begin{pspicture}(0,0)(12.5,2)
\rput(3.75,1.55){\BoxedEPSF{chapter3.fig900.eps scaled 400}}
\rput(1.15,.85){$\bullet$}\rput(2.5,.85){$\bullet$}\rput(6.35,.85){$\bullet$}
\rput(8.65,.85){$\bullet$}
\rput(4.5,.85){$\ldots$}\rput(4.05,.85){$\ldots$}\rput(4.95,.85){$\ldots$}
\rput(10.25,1){\BoxedEPSF{chapter3.fig1000.eps scaled 400}}
\psline[linewidth=.2mm]{->}(6.7,.85)(8.3,.85)
\rput(12.2,1){$\gamma$}\rput(5.5,1.6){${\red\mu}$}
\rput(8.6,.65){$v$}\rput(1.1,.65){$u$}\rput(6.35,.65){$x$}
\rput(7.5,1.05){$f$}
\rput{270}(3.7,.5){$\left.\begin{array}{c}
\vrule width 0 mm height 53 mm depth 0 pt\end{array}\right\}$}
\rput(3.7,.1){$f^{-1}(v)$}
\end{pspicture}  
  \caption{Defining the monodromy using path lifting.}
  \label{chapter3:actionsuniversal:figure1000}
\end{figure}
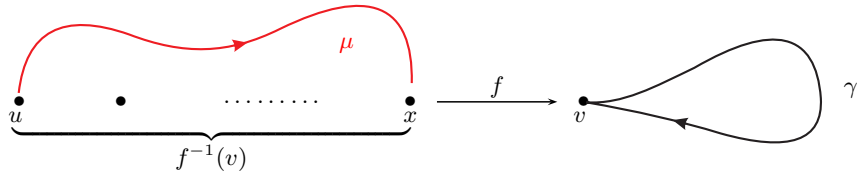

The definition is illustrated in Figure \ref{chapter3:actionsuniversal:figure1000}:
let $\gamma$ be a closed path at $v$ representing the element
$g_\gamma\in\pi_1(X,v)$. Let $\mu$ be the 
lift of $\gamma$ at a vertex $u\in f^{-1}(v)$, and let this lift have end vertex $x$. 
Let $\ss_\gamma\in\sym f^{-1}(v)$ be the permutation with $\ss_\gamma(u)=x$.

We obviously have a well-defined issue to deal with, so that the permutation $\ss_\gamma$ does
not depend on our choice of representative $\gamma$. If
$\gamma'$ is another closed path at $v$ with $g_{\gamma'}=g_\gamma$, then the paths
$\gamma,\gamma'$ are homotopic. This homotopy can be lifted, via homotopy
lifting, to a homotopy between the lifts $\mu$ and $\mu'$ at $u$,
and so the two lifts are homotopic in $Y$. 
But homotopic paths have the same endpoints! Thus $\mu'$ ends at $x$ as well,
and we get $\ss_\gamma(u)=x=\ss_{\gamma'}(u)$.

By Exercise \ref{chapter3:basics:exercise150} the lift at $u$ of the path
$\gamma_1\gamma_2$ is the path $\mu_1\mu_2$ obtained by lifting 
$\gamma_1$ and then lifting $\gamma_2$ to the
terminal vertex of $\mu_1$. In particular $\ss_{\gamma_1\gamma_2}=\ss_{\gamma_2}\ss_{\gamma_1}$,
recalling that the product is read from right to left. 

\begin{proposition}
\label{chapter3:actionsuniversal:result1000}
If $f:Y\rightarrow X$ is a covering then monodromy gives a homomorphism 
$$
\pi_1(X,v)\rightarrow\sym(f^{-1}(v)),
$$
defined by $g_\gamma\mapsto\ss_\gamma^{-1}$.
In particular, a covering of finite degree gives a homomorphism from 
$\pi_1(X,v)$ to a finite
group.
\end{proposition}

\begin{exercise}
\label{chapter3:actionsuniversal:exericse950}
Let $X$ be a $2$-complex with a single vertex  and $H\subset\pi_1(X,v)$ a subgroup with
$X\kern-2pt\uparrow\kern-2pt H$ and $f:X\kern-2pt\uparrow\kern-2pt H\rightarrow X$
the covering of \S\ref{chapter3:actionsuniversal:bottom_up}. Show that the monodromy
action corresponds to the action on the cosets given by
$Hg\mapsto H(gh^{-1})$ for $h\in\pi_1(X,v)$.
Compare with Exercise \ref{chapter4:automorphisms:exercise125}, and notice that
the action is always well defined.
\end{exercise}

Suppose that $X$ has just one vertex $v$, so that the fiber of $v$ consists
of all the vertices of $Y$. Monodromy then gives an action of $\pi_1(X,v)$
on the whole $0$-skeleton
of $Y$. The next exercise shows that in general this action cannot be extended any further
than this.

\begin{exercise}
\label{chapter3:actionsuniversal:exericse1000}
\begin{enumerate}
\item Let $X$ be the complex of Figure \ref{chapter3:actionsuniversal:figure1100}.
Describe the universal cover $\wtl{X}\rightarrow X$, showing that in particular that it is a
covering of degree $6$.
  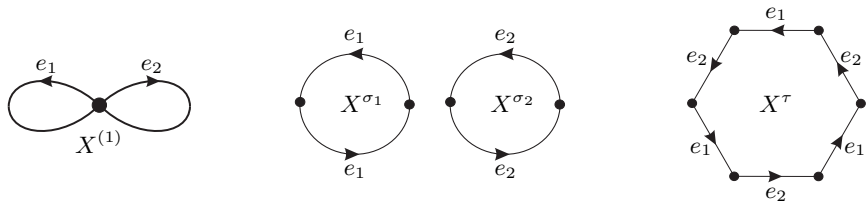
\begin{figure}
    \centering
\begin{pspicture}(0,0)(12.5,2)
\rput(-2,0){
\rput(4,1){\BoxedEPSF{chapter3.fig2100.eps scaled 750}}
\rput(4,.5){$X^{(1)}$}\rput(3.3,1.5){$e_1$}\rput(4.7,1.5){$e_2$}
}
\rput(-1,0){
\rput(6.4,1){\BoxedEPSF{chapter3.fig4400a.eps scaled 800}}
\rput(6.5,1){$X^{\ss_1}$}\rput(6.4,1.9){$e_1$}\rput(6.4,.1){$e_1$}
}
\rput(-1,0){
\rput(8.4,1){\BoxedEPSF{chapter3.fig4400a.eps scaled 800}}
\rput(8.5,1){$X^{\ss_2}$}\rput(8.4,1.9){$e_2$}\rput(8.4,.1){$e_2$}
}
\rput(11,1){\BoxedEPSF{chapter3.fig4400b.eps scaled 750}}
\rput(11,1){$X^\tau$}
\rput(11,2.2){$e_1$}\rput(9.95,1.6){$e_2$}
\rput(10,.4){$e_1$}\rput(11,-.2){$e_2$}
\rput(12.05,.4){$e_1$}\rput(12,1.6){$e_2$}
\end{pspicture}  
    \caption{complex $X$ for Exercise \ref{chapter3:actionsuniversal:exericse1000}
with three faces $\ss_1,\ss_2,\tau$.}
    \label{chapter3:actionsuniversal:figure1100}
  \end{figure}
\item Show that in $\wtl{X}$ there exists a pair of vertices joined by an edge but there is
no edge joining the images of these two vertices under the monodromy action of $\pi_1(X,v)$. 
Thus it is not possible to define an automorphism of $\wtl{X}$ at this edge. Deduce that there 
can therefore be no homomorphism from $\pi_1(X,v)$ to the automorphism group of the $1$-skeleton
extending the monodromy action on the $0$-skeleton.
\end{enumerate}
\end{exercise}

\section{Operations on coverings}\label{chapter3:operations}

In Chapter \ref{chapter1}
we had three constructions arising from a collection of complexes and
maps between them: the pushout, pullback and Higman composition. In this section
we show that all three are useful ways of creating new coverings from old.

The set-up is as follows: we have a fixed covering $f:Y\rightarrow X$ together with two
coverings intermediate to $f$ as in \S\ref{chapter3:actionsuniversal:intermediate}:
$$
Y\stackrel{g_1}{\rightarrow}Z_1\stackrel{h_1}{\rightarrow}X
\text{ and }
Y\stackrel{g_2}{\rightarrow}Z_2\stackrel{h_2}{\rightarrow}X.
$$
We then pushout the covers $Y\rightarrow Z_i$ and pullback the covers $Z_i\rightarrow X$.
\emph{Throughout this section all complexes are connected\/}.

\subsection{Pushouts of covers}\label{chapter3:operations:pushouts}

\begin{figure}
  \centering
\begin{pspicture}(12.5,3)
\rput(2.5,0.5){
\rput(0,2){$Y$}
\rput(0,0){$Z_1$}
\rput(2,2){$Z_2$}
\rput(2,0){$X$}
\psline[linewidth=.2mm]{->}(0,1.8)(0,.2)
\psline[linewidth=.2mm]{->}(.2,2)(1.8,2)
\psline[linewidth=.2mm]{->}(.2,0)(1.8,0)
\psline[linewidth=.2mm]{->}(2,1.8)(2,.2)
\psline[linewidth=.2mm]{->}(.1,1.9)(1.8,.2)
\rput(-0.2,1){$g_1$}
\rput(1,2.2){$g_2$}
\rput(1,.2){$h_1$}
\rput(2.2,1){$h_2$}
\rput*(1,1){$f$}
}
\rput(2,0){
\rput(5,1){
\rput(0,2){$Y$}
\rput(0,0){$Z_1$}
\rput(2,2){$Z_2$}
\rput(2.1,0){$Z_1\coprod_Y Z_2$}
\psline[linewidth=.2mm]{->}(0,1.7)(0,.3)
\psline[linewidth=.2mm]{->}(.3,2)(1.7,2)
\psline[linewidth=.2mm]{->}(.3,0)(1.2,0)
\psline[linewidth=.2mm]{->}(2,1.7)(2,.3)
\rput(.25,1.05){$g_1$}
\rput(1,1.8){$g_2$}
\rput(.75,.2){$t_1$}
\rput(1.8,1){$t_2$}
}
\rput(8,0){$X$}
\psline[linewidth=.2mm,linearc=.3]{->}(5,.8)(5,0)(7.8,0)
\psline[linewidth=.2mm,linearc=.3]{->}(7.25,3)(8,3)(8,.2)
\psline[linewidth=.3mm,linestyle=dotted]{->}(7.2,.8)(7.85,.15)
\rput(6,.25){$h_1$}\rput(8.2,2){$h_2$}
}
\end{pspicture}  
\caption{two intermediate covers (left) and their pushout (right).}
\label{chapter3:operations:pushouts:figure100}
\end{figure}
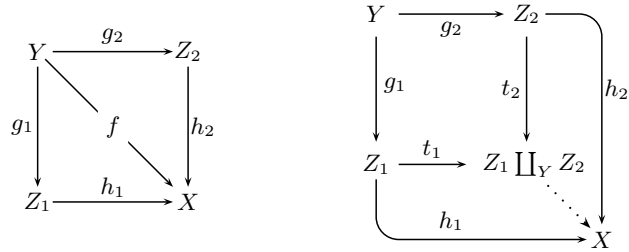

Let $f:Y\rightarrow X$ be a fixed covering of connected $2$-complexes, and
$Y\stackrel{g_i}{\rightarrow}Z_i\stackrel{h_i}{\rightarrow}X\,(i=1,2)$ be coverings 
intermediate to $f$ with the $Z_i$ distinct and connected. 
Thus we have the commuting diagram on the left 
of Figure \ref{chapter3:operations:pushouts:figure100},
with all the maps in sight coverings. 
As the $g_i$ are dimension preserving,
we can by Theorem \ref{chapter1:2complexes:result400}, form the pushout $Z_1\coprod_Y Z_2$, 
obtaining in the process maps $t_i:Z_i\rightarrow Z_1\coprod_Y Z_2\,(i=1,2)$ as the composition
$Z_i\hookrightarrow Z_1\bigcup Z_2\rightarrow Z_1\bigcup Z_2\quo$ of the 
inclusion of $Z_i$ in the disjoint union and the quotient defined in 
\S\ref{chapter1:quotients:pushouts}.
The universality of the pushout,
applied to the maps $h_i:Z_i\rightarrow X$,
gives the commuting diagram on the right of 
Figure \ref{chapter3:operations:pushouts:figure100}.
If $[z]$ is a cell of $Z_1\coprod_Y Z_2$, then $z\in Z_i$ for some $i$,
so that the map $h:Z_1\coprod_Y Z_2\rightarrow X$ sends $[z]$ to $h_i(z)\in X$.

\begin{proposition}[pushouts of covers]
\label{chapter3:operations:result100}
The maps 
$$
t_i:Z_i\rightarrow Z_1\coprod_Y Z_2,\,(i=1,2)\text{ and }h:Z_1\coprod_Y Z_2\rightarrow X,
$$
are coverings. Thus, the pushout of two intermediate coverings 
$Y\rightarrow Z_i\rightarrow X\,(i=1,2)$
is a connected intermediate covering $Y\rightarrow Z_1\coprod_YZ_2\rightarrow X$.
\end{proposition}

\begin{proof}
We show that the map $t_1g_1:Y\rightarrow Z_1\coprod_Y Z_2$ is a covering, and then two
applications of Proposition \ref{chapter3:actionsuniversal:result150} give
$t_1$ and $h$ are coverings. That $t_2$ is a covering is completely analogous.
The map $t_1g_1$
is dimension preserving, as $g_1$ and $t_1$ are, leaving us to show that the local continuity
maps are bijections. Suppose then that $u$ is a vertex of $Y$ mapping via $t_1g_1$ to
the vertex
$[v]$ of the pushout, so that there is a vertex $v$ of $Z_1$ with 
$v=g_1(u)$ and $[v]=t_1(v)$.

Starting with the surjectivity of the local continuity on edges,
suppose we have an edge $[e']$ of the pushout with start vertex $[v]$.
Thus, $e'$ 
is an edge in the disjoint union $Z_1\bigcup Z_2$ with start vertex $v'$
equivalent to $v$.
The equivalence between $v$ and $v'$ is 
realized by a sequence of lifts and covers (of vertices) through the coverings
$g_i:Y\rightarrow Z_i$. The same sequence applied to the edge $e'$, and using path
lifting, yields an edge $e\in Z_1$ equivalent
to $e'$, and so $t_1(e)=[e']$, and with $e$ having start vertex $v$. Lifting $e$ to 
$u\in Y$ gives an edge mapping via local continuity to $[e']$. Faces work the same:
start with an occurrence of $[v]$ in $[\ss']$; use face lifting to get an 
occurrence of $v$ in a face $\ss\in Z_1$ with $t_1(\ss)=[\ss']$. 

Now to the injectivity of the local continuity on edges, for which we suppose there
are edges $e_1,e_2$ starting at $u\in Y$ and mapping via $t_1g_1$ to an edge $[e]$ 
starting at $[v]$ in the pushout. Thus the edges $g_1(e_1),g_1(e_2)$ map via $t_1$ to $[e]$, hence by 
$h_1$ to $h[e]$ (starting at $h[u]$). Uniqueness of path lifting, applied first to the 
covering $h_1$ and then to $g_1$ gives $e_1=e_2$. Again, faces work the same.
\qed
\end{proof}

\subsection{Pullbacks of covers}\label{chapter3:operations:pullbacks}

\begin{figure}
  \centering
\begin{pspicture}(12.5,3)
\rput(2.5,0.5){
\rput(0,2){$Y$}
\rput(0,0){$Z_1$}
\rput(2,2){$Z_2$}
\rput(2,0){$X$}
\psline[linewidth=.2mm]{->}(0,1.8)(0,.2)
\psline[linewidth=.2mm]{->}(.2,2)(1.8,2)
\psline[linewidth=.2mm]{->}(.2,0)(1.8,0)
\psline[linewidth=.2mm]{->}(2,1.8)(2,.2)
\psline[linewidth=.2mm]{->}(.1,1.9)(1.8,.2)
\rput(-0.2,1){$g_1$}
\rput(1,2.2){$g_2$}
\rput(1,.2){$h_1$}
\rput(2.2,1){$h_2$}
\rput*(1,1){$f$}
}
\rput(3.5,-1){
\rput(5,1){
\rput(0,2){$Z_1\prod_X Z_2$}
\rput(0,0){$Z_1$}
\rput(2,2){$Z_2$}
\rput(2,0){$X$}
\psline[linewidth=.2mm]{->}(0,1.7)(0,.3)
\psline[linewidth=.2mm]{->}(.9,2)(1.7,2)
\psline[linewidth=.2mm]{->}(.3,0)(1.7,0)
\psline[linewidth=.2mm]{->}(2,1.7)(2,.3)
\rput(.25,1){$t_1$}
\rput(1.3,1.8){$t_2$}
\rput(1,.2){$h_1$}
\rput(1.75,1){$h_2$}
}
\rput(4,4){$Y$}
\psline[linewidth=.2mm,linearc=.3]{->}(4,3.8)(4,1)(4.7,1)
\psline[linewidth=.2mm,linearc=.3]{->}(4.2,4)(7,4)(7,3.2)
\rput(3.8,2){$g_1$}
\rput(6,3.8){$g_2$}
\psline[linewidth=.3mm,linestyle=dotted]{->}(4.1,3.9)(4.7,3.3)
}
\end{pspicture}
  \caption{two intermediate covers (left) and their pullback (right).}
  \label{chapter3:operations:pullbacks:figure100}
\end{figure}
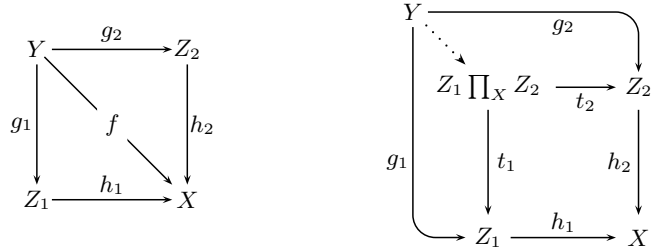

As in \S\ref{chapter3:operations:pushouts}, 
let $f:Y\rightarrow X$ be a fixed covering of connected $2$-complexes, and 
$Y\stackrel{g_i}{\rightarrow}Z_i\stackrel{h_i}{\rightarrow}X\,(i=1,2)$ be coverings 
intermediate to $f$ with the $Z_i$ connected. 
Thus we have the commuting diagram on the left of Figure \ref{chapter3:operations:pullbacks:figure100}
with all the maps in sight coverings. As the coverings $h_i$ are dimension preserving,
we may, via \S\ref{chapter1:pullhig:pullbacks}, form the pullback $Z_1\prod_X Z_2$, 
obtaining in the process maps $t_i:Z_1\prod_X Z_2\rightarrow Z_i$ given by
$t_i:z_1\times z_2\mapsto z_i$. The universality of the pullback, Theorem
\ref{chapter1:2complexes:result600}, applied to the maps $g_i:Y\rightarrow Z_i$,
gives the commuting diagram on the right of 
Figure \ref{chapter3:operations:pullbacks:figure100}. The new map
$h:Y\rightarrow Z_1\prod_X Z_2$ 
sends a cell $y\in Y$ to the cell $g_1(y)\times g_2(y)\in Z_1\prod_X Z_2$.

\begin{proposition}[pullbacks of covers]
\label{chapter3:operations:result500}
The maps 
$$
t_i:Z_1\prod_X Z_2\rightarrow Z_i\,(i=1,2)\text{ and }h:Y\rightarrow Z_1\prod_X Z_2,
$$
are coverings. Thus, the pullback of two intermediate coverings 
$Y\rightarrow Z_i\rightarrow X\,(i=1,2)$
is an intermediate covering $Y\rightarrow Z_1\prod_XZ_2\rightarrow X$.
\end{proposition}

\begin{proof}
It suffices, by Proposition \ref{chapter3:actionsuniversal:result150}
to show that $h$ is a covering. It is dimension preserving as the $g_i$ are, and so it remains to
show that the various local continuity maps are bijections. This is similar for both edges and 
faces: suppose that $v_1\times v_2$ is a vertex of the pullback and $u$ a vertex of $Y$
with $h(u)=v_1\times v_2$.
If two objects at $u$ (edges starting at $u$ or appearances of $u$ in faces)
map under $h$ to a single object at $v_1\times v_2$, then these two map via the $g_i$ to
single objects at $v_i\in Z_i$. The $g_i$ are coverings, ensuring that the original
two objects coincide, hence injectivity of the local continuity maps.

Surjectivity requires a couple more steps: start with an object at the vertex 
$v_1\times v_2$ of the pullback. It maps via the $t_i$ to objects at the $v_i\in Z_i$, and they
in turn map via the $h_i$ to the \emph{same\/} object at $v=h_i(v_i)$. The path and face lifting
provided by the covers $g_i$ give two objects at $u$ mapping to this single object at $u$,
one via $h_1g_1$ and
the other via $h_2g_2$. But then these two objects 
map via the covering $f$ to this single object at $u$, and so must be the same object.
By definition, the image via $h$ of this single object at $u$ must be the original
object at $v_1\times v_2$ that we started with.
\qed
\end{proof}

We saw in \S\ref{chapter1:pullhig:pullbacks} that the pullback is not necessarily 
connected. To get a connected covering we use the pointed version of
Exercise \ref{chapter1:2complexes:exercise650}: the $h_i:(Z_i)_{z_i}\rightarrow X_x$
are pointed coverings, giving
$$
t_ih_i:(Z_1\prod_X Z_2)_z\rightarrow X_x,
$$
a pointed covering of connected complexes for $z=z_1\times z_2$.

\begin{exercise}
Let $Y\rightarrow X\leftarrow Z$ be coverings with $Z$ a forest. Show that the 
pullback $X\prod_XZ$ is also a forest.  
\end{exercise}

\subsection{Higman compositions of covers}\label{chapter3:operations:higman}

Let $f_i:Y_i\rightarrow X\, (i=1,\ldots,n)$ be dimension preserving maps
and $\{e_{j1},e_{j2}\}$ ($j=1,\ldots,m$) a handle configuration in $\bigcup_i
Y_i$, with $Y=\hcl Y_1,\ldots,Y_n\hcr$ the Higman composition and $f=\hcl
f_1,\ldots,f_n\hcr:Y\rightarrow X$ the map of \S\ref{chapter1:pullhig:higman}.

\begin{proposition}[Higman compositions of covers]
\label{chapter3:operations:result550}
If the $f_i$ are covering maps then $f$ is a covering map. Moreover, 
if the Higman composition is connected then
$$
\deg(\hcl Y_1,\ldots,Y_n\hcr\rightarrow X)=\sum_i\deg(Y_i\rightarrow X).
$$
\end{proposition}

\begin{proof}
Let $u\in Y$ and $v\in X$ be vertices with $f(u)=v$. If $u$ is not the start or terminal
vertex of one of the edges in the handle configuration, then the edges starting at $u$ are
completely unaffected by the composition. Otherwise, the edges starting at $u$ are unchanged
in number, as Figure \ref{chapter1:2complexes:figure3100} shows. For local continuity of faces,
we have the desired bijection before composition courtesy of one of the covering maps maps $f_i$,
with the composition replacing certain occurrences of $u$ in the $\tau_\ell$ by occurrences in the
$\tau_\ell'$, while still maintaining the bijection. Thus, $f$ is a covering map. 
The degree assertion follows from 
the fact that if $e$ is the edge of Definition \ref{chapter1:2complexes:definition800}
giving rise to the handle configuration, then 
the fiber $f^{-1}(e)$ is in bijective correspondence with the disjoint union
of fibers $\bigcup_i f^{-1}_i(e)$.
\qed
\end{proof}

\begin{figure}
  \centering
\begin{pspicture}(12.5,4.5)
\rput(2,2){\BoxedEPSF{chapter3.fig5000.eps scaled 750}}
\rput(6.5,2){\BoxedEPSF{chapter3.fig2100.eps scaled 750}}
\rput(11,2){\BoxedEPSF{chapter3.fig5000a.eps scaled 750}}
\psline[linewidth=.2mm]{->}(4.5,2)(5,2)
\psline[linewidth=.2mm]{->}(8.4,2)(7.9,2)
\rput(1.5,2.3){${\red a}$}\rput(2.5,2.3){${\red a}$}
\rput(9.9,2.4){${\red a}$}\rput(12.1,2.4){${\red a}$}
\rput(5.8,2.5){$a$}
\rput(4.75,2.2){$f_1$}\rput(8.15,2.2){$f_2$}
\end{pspicture}
  \caption{coverings $f_1,f_2$ and a handle configuration (in red).}
  \label{chapter3:operations:hc:figure100}
\end{figure}
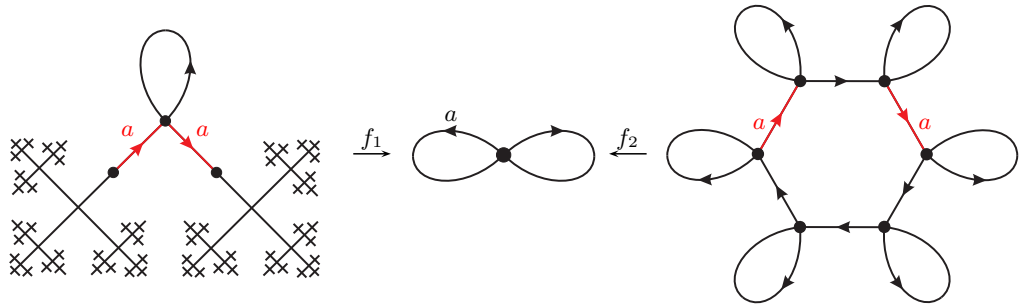

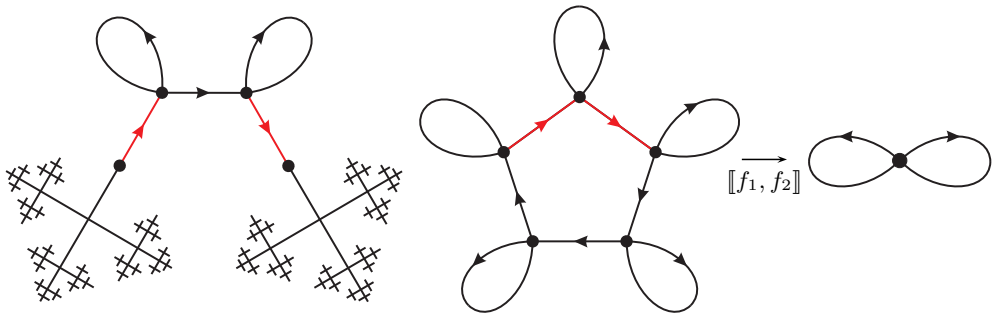
\begin{figure}
  \centering
\begin{pspicture}(12.5,4)
\rput(2.25,2){\BoxedEPSF{chapter3.fig5000b.eps scaled 750}}
\rput(7.25,2){\BoxedEPSF{chapter3.fig5000c.eps scaled 750}}
\rput(11.5,2){\BoxedEPSF{chapter3.fig2100.eps scaled 750}}
\psline[linewidth=.2mm]{->}(9.4,2)(10,2)
\rput(9.7,1.7){$\hcl f_1,f_2\hcr$}
\end{pspicture}
  \caption{the (disconnected) Higman composition resulting from the set-up in Figure
\ref{chapter3:operations:hc:figure100}.}
  \label{chapter3:operations:hc:figure200}
\end{figure}

\begin{example}
Figure \ref{chapter3:operations:hc:figure100}
shows graphs $X,Y_1,Y_2$ and two coverings $f_i:Y_i\rightarrow X$,
a handle configuration in $Y_1\bigcup Y_2$ and the resulting (disconnected) Higman composition
in Figure \ref{chapter3:operations:hc:figure200}.
\end{example}

\section{Lattices of covers}\label{chapter3:lattices}

When we introduced pullbacks in Chapter \ref{chapter1} we said that they
would act as a kind of union of $2$-complexes, with the pushout acting as a kind 
of intersection. This section makes this precise: we give the set of coverings
intermediate to a fixed covering $f:Y\rightarrow X$ the structure of a poset
in which the pullback and pushout give a join $\vee$ and a meet $\wedge$. Thus the 
intermediate coverings form a lattice (Theorem \ref{chapter3:lattices:lattice:result300} 
below). The whole business is complicated by the fact that the resulting lattice 
is \emph{slightly too big\/} for what we want it for in Chapter \ref{chapter4}.
This forces us in \S\ref{chapter3:lattices:poset}
below to consider instead intermediate coverings upto
a certain equivalence.

\subsection{\emph{Aside}: posets and lattices}\label{chapter3:lattices:aside}

We pause and take a brief look at the theory of posets and lattices and some important 
examples.
There are many books on this subject: we have followed
\cite{Stanley97}*{Chapter 3}.

Partially ordered sets (or \emph{posets\/}) formalise the idea of ordering: 
a poset is a set $P$ and a binary relation $\leq$
that is  
\emph{reflexive\/}: $x\leq x$ for all $x\in P$;
\emph{antisymmetric\/}: if $x\leq y$ and $y\leq x$ then $x=y$; and
\emph{transitive\/}: if $x\leq y$ and $y\leq z$ then $x\leq z$.
The motivating example is meant to be the integers $\Z$ with their usual ordering $\leq$,
and the usual notational conventions from there are used in general:
we write $x<y$ to mean $x\leq y$ but $x\not= y$. Elements $x,y$ with $x\leq y$ 
or $y\leq x$ are 
\emph{comparable\/}, otherwise they are \emph{incomparable\/} (a possibility that
obviously doesn't arise with the primordial example $\Z$). We say that $y$ \emph{covers\/}
$x$, written $x\prec y$, when $x<y$ and if $x\leq z\leq y$ then either $z=x$ or $z=y$. 

A \emph{morphism\/} (or just \emph{map\/}) of posets $f:P\rightarrow Q$ is an order-preserving
map of the underlying sets: if $x\leq y$ in $P$ then $f(x)\leq f(y)$ in $Q$. 
Notice that this is a one way business:
comparable elements are sent to comparable elements, but
incomparable elements are allowed to map to comparable ones. 
An \emph{anti-morphism\/} is an order-reversing map: if $x\leq y$ in $P$ then $f(y)\leq f(x)$ in $Q$.

Bijective morphisms have inverse set maps, although they may not be morphisms. A bijective
morphism with order-preserving inverse is an \emph{isomorphism\/}: $x\leq y$ in $P$ if and only if
$f(x)\leq f(y)$ in $Q$. Similarly a bijective anti-morphism with order-reversing inverse is
an \emph{anti-isomorphism\/}.

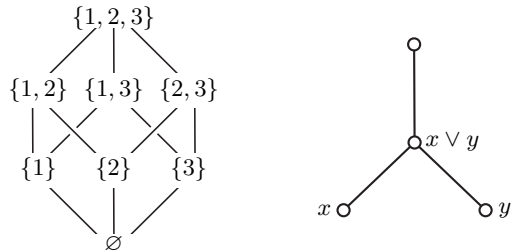
\begin{figure}
  \centering
\begin{pspicture}(0,0)(12.5,3.25)
\rput(-3.5,-.5){
\rput(7.5,2){\BoxedEPSF{cube4.eps scaled 300}}
\rput(7.5,0.5){${\varnothing}$}
\rput(6.5,1.5){${\{1\}}$}
\rput(7.5,1.5){${\{2\}}$}
\rput(8.5,1.5){${\{3\}}$}
\rput(6.5,2.55){${\{1,2\}}$}
\rput(7.5,2.55){${\{1,3\}}$}
\rput(8.5,2.55){${\{2,3\}}$}
\rput(7.5,3.5){${\{1,2,3\}}$}
}
\rput(8,1.55){\BoxedEPSF{chapter3.fig2800.eps scaled 400}}
\rput(6.8,.45){$x$}\rput(9.2,.45){$y$}\rput(8.5,1.35){$x\vee y$}
\end{pspicture}  
  \caption{Hasse diagram for the poset of subsets of the set $\{1,2,3\}$ ordered
by inclusion (left) and for a poset with four elements (right) that is not a lattice.}
  \label{chapter3:lattices:figure100}
\end{figure}

Posets are often illustrated using their \emph{Hasse diagram\/}: a graph whose vertices
are the elements of $P$ and whose edges give the covering relations. Thus, if $x\prec y$
then the vertex $y$ is drawn above the vertex $x$ with an edge connecting them. 
Two examples
of Hasse diagrams (and posets) are given in Figure 
\ref{chapter3:lattices:figure100}.

A special place is reserved for those posets which have supremums and infimums. 
If $x,y\in P$ then $z$ is an \emph{upper bound\/} for $x$ and $y$ when both $x\leq z$ and 
$y\leq z$. It is a \emph{least upper bound\/} or \emph{supremum\/} or \emph{join\/}
when it is an upper bound such that for any other upper bound $w$ we have
$z\leq w$. Similarly, $z$ is a \emph{lower bound\/} for $x$ and $y$ when both $z\leq x$ and 
$z\leq y$. It is a \emph{greatest lower bound\/} or \emph{infimum\/} or \emph{meet\/}
when it is an lower bound such that for any other lower bound $w$ we have
$w\leq z$.

It is easy to show that if $x$ and $y$ have a join then it is unique (\emph{hint\/}: any
two joins must be $\leq$ each other) and similarly for the meet. 
Write $x\vee y$ for the join and $x\wedge y$ for
the meet of $x$ and $y$.

A poset is a \emph{lattice\/} if for every pair of elements $x$ and $y$,
the join $x\vee y$ and meet $x\wedge y$ exist.

The poset on the left of Figure \ref{chapter3:lattices:figure100}
is a lattice, as can be checked directly from the Hasse diagram, but the example on the right is
not: if $x$ and $y$ are the two elements shown, then they have a 
join, but no meet.

\begin{exercise}
\label{chapter3:lattices:exercise50}
A $\hat{1}$ in a poset $P$ is a unique maximal element: for all $x\in P$ we have
$x\leq\hat{1}$. Similarly
a $\hat{0}$ in $P$ is a unique minimal element: for all $x\in P$ we have
$\hat{0}\leq x$.
Show that a finite lattice has a $\hat{0}$ and a $\hat{1}$.
\end{exercise}

\begin{exercise}
\label{chapter3:lattices:exercise100}
A poset is a \emph{meet-semilattice\/} if any two elements have a meet. 
Dually we have the notion of a \emph{join-semilattice}. Show that
if $P$ is a \emph{finite\/} meet-semilattice with a $\hat{1}$ then $P$ is a lattice
(dually, if $P$ is a finite join-semilattice with a $\hat{0}$ then $P$ is a lattice).
\end{exercise}

\begin{exercise}
\label{chapter3:lattices:exercise200}
Let $P$ and $Q$ be lattices and $f:P\rightarrow Q$ a lattice isomorphism 
(respectively anti-isomorphism). Show that $f$ sends joins to joins and meets
to meets (resp. joins to meets and meets to joins), ie: $f(x\vee y)=f(x)\vee f(y)$
and $f(x\wedge y)=f(x)\wedge f(y)$.
\end{exercise}

The most commonly occurring lattice ``in nature'' is the \emph{Boolean lattice\/} on 
a set $X$: its elements are the subsets of $X$ with $A\leq B$ if
and only if $A\subset B$. Meets and 
joins are just intersections and unions: $A\wedge B=A\cap B$ and $A\vee B=A\cup B$.

\begin{exercise}
\label{chapter3:lattices:exercise250}
Let $X$ be a finite set, $P$ the Boolean lattice on $X$ and $V$ the real vector space
with basis $X$. If $v=\sum_X\lambda_x x\in V$ define $|v|^2=\sum_X\lambda_x^2$, and let
$\Box^n:=\{v\in V :|v|^2\leq 1\}$, the \emph{$n$-dimensional cube\/}. Embed the
underlying set of $P$ in $V$ via the map sending $A\subset X$ to $\sum_{x\in A}x$,
and show that the image of $P$ is the set of vertices of $\Box^n$, while the 
vertices and edges of $\Box^n$ give the Hasse diagram for $P$. 
\end{exercise}

Another example is the lattice $L_n(\F)$ of all subpaces 
of the $n$-dimensional vector space over the field $\F$, with the 
ordering given by inclusion of one subspace in another. The meet of two subspaces
is again their intersection, but this time the union is too small to be their
join: the union of two subspaces is not a subspace! Instead we define
$U\vee V=U+V$, their sum, consisting of all vectors of the form $u+v$ for 
$u\in U$ and $v\in V$. We leave it to the reader to verify that these are indeed
infimums and supremums.

Here is one we are particularly interested in,

\begin{definition}[lattice of subgroups]
\label{chapter3:lattices:definition300}
Let $G$ be a group. The lattice of subgroups $\SSS(G)$ has as elements the subgroups
of $G$ ordered by inclusion, and with $H\wedge K=H\cap K$, $H\vee K=\langle H, K\rangle$,
the subgroup generated by $H$ and $K$.
\end{definition}

\begin{figure}
  \centering
\begin{pspicture}(0,0)(12.5,3)
\rput(6.25,3){$\SS_3$}
\rput(6.25,2){$\{\id,\ss,\ss^2\}$}
\rput(4,1){$\{\id,\tau\}$}
\rput(8,1){$\{\id,\ss\tau\}$}
\rput(10,1){$\{\id,\ss^2\tau\}$}
\rput(6.25,0){$\id$}
\psline[linewidth=.2mm]{-}(6.25,2.2)(6.25,2.8)
\psline[linewidth=.2mm]{-}(6.25,.2)(6.25,1.8)
\psline[linewidth=.2mm]{-}(6.1,2.8)(4,1.2)
\psline[linewidth=.2mm]{-}(6.4,2.8)(8,1.2)
\psline[linewidth=.2mm]{-}(6.55,2.8)(10,1.2)
\psline[linewidth=.2mm]{-}(6.1,.2)(4,.8)
\psline[linewidth=.2mm]{-}(6.4,.2)(8,.8)
\psline[linewidth=.2mm]{-}(6.55,.2)(10,.8)
\end{pspicture}  
  \caption{subgroup lattice for the symmetric group $\SS_3$ with $\ss=(1,2,3)$ and
$\tau=(2,3)$.}
  \label{chapter3:lattices:figure200}
\end{figure}
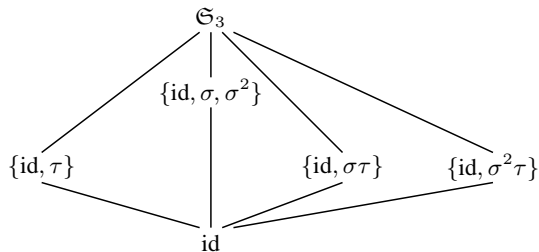

\begin{exercise}
\label{chapter3:lattices:exercise300}
Show that the set of \emph{finite index\/} subgroups of a group $G$ also
forms a lattice, with the same meet and join as $\SSS(G)$. Show the same
with finite index replaced by \emph{finitely generated\/}.
\end{exercise}

\begin{exercise}
\label{chapter3:lattices:exercise300}
Show that an isomorphism $G_1\rightarrow G_2$ of groups induces an
isomorphism of lattices $\SSS(G_1)\rightarrow\SSS(G_2)$.
\end{exercise}

\subsection{The poset of intermediate covers}
\label{chapter3:lattices:poset}

In this section and the
next, we construct a lattice whose elements are, more or less, the coverings intermediate
to a particular fixed covering $f$. In the next chapter we'll see that if we look at 
this lattice sideways and squint our eyes a little, then it looks the same as the lattice
of subgroups of the ``group of automorphisms'' of the covering $f$. 

Throughout this section $f:Y\rightarrow X$ is a fixed 
covering of connected $2$-complexes and
$Y\stackrel{g_i}{\rightarrow}(Z_i)_{x_i}\stackrel{h_i}{\rightarrow}X\,(i=1,2)$ 
are coverings 
intermediate to $f$ with the $Z_i$ connected.

We call these two intermediate coverings
\emph{equivalent\/} 
if and only if there is a isomorphism $Z_1\rightarrow Z_2$ making the diagram
on the right of Figure \ref{chapter3:lattices:figure250}
commute.

\begin{figure}
\centering
\begin{pspicture}(12.5,2.5)
\rput(-1,0){
\rput(2.5,0){
\rput(0,2){$Y$}
\rput(0,0){$Z_1$}
\rput(2,2){$Z_2$}
\rput(2,0){$X$}
\psline[linewidth=.2mm]{->}(0,1.8)(0,.2)
\psline[linewidth=.2mm]{->}(.2,2)(1.8,2)
\psline[linewidth=.2mm]{->}(.2,0)(1.8,0)
\psline[linewidth=.2mm]{->}(2,1.8)(2,.2)
\psline[linewidth=.2mm]{->}(.1,1.9)(1.8,.2)
\rput(-0.2,1){$g_1$}
\rput(1,2.2){$g_2$}
\rput(1,.2){$h_1$}
\rput(2.2,1){$h_2$}
\rput*(1,1){$f$}
}
\rput(8,0){
\rput(0,2){$Y$}
\rput(0,0){$Z_1$}
\rput(2,2){$Z_2$}
\rput(2,0){$X$}
\psline[linewidth=.2mm]{->}(0,1.8)(0,.2)
\psline[linewidth=.2mm]{->}(.2,2)(1.8,2)
\psline[linewidth=.2mm]{->}(.2,0)(1.8,0)
\psline[linewidth=.2mm]{->}(2,1.8)(2,.2)
\psline[linewidth=.3mm,linestyle=dotted]{->}(.1,.2)(1.75,1.85)
\rput(-0.2,1){$g_1$}
\rput(1,2.2){$g_2$}
\rput(1,.2){$h_1$}
\rput(2.2,1){$h_2$}
\rput*(1,1){$\cong$}
}
\rput(6,1){equivalent}\rput(7.2,1){$\Leftrightarrow$}\rput(11.5,1){commutes}
}
\end{pspicture}
\caption{}
  \label{chapter3:lattices:figure250}
\end{figure}

This is an equivalence relation on the set of coverings intermediate to
$f$, and we write $\mathcal{L}(Y\stackrel{f}{\rightarrow}X)$ or just $\mathcal{L}(Y,X)$ for the
set of equivalence classes. The notation here can become 
very cumbersome, so where possible we will write $Z\in\mathcal{L}(Y,X)$ to mean 
the equivalence class represented by the
coverings $Y\rightarrow Z\rightarrow X$ intermediate to $f$.

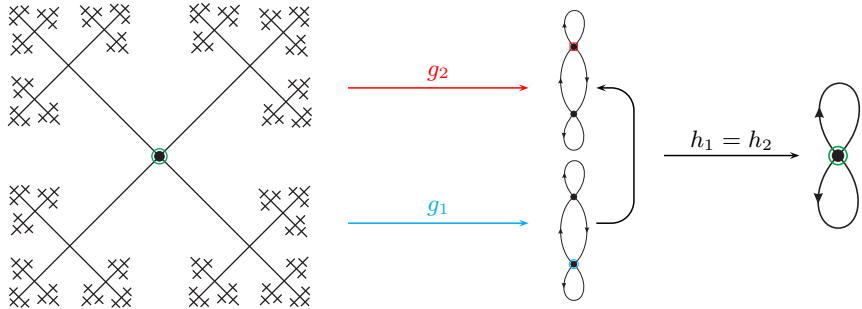
\begin{figure}
  \centering
\begin{pspicture}(0,0)(12.5,4)
\rput(2.5,2){\BoxedEPSF{chapter3.fig2400.eps scaled 500}}
\rput(8,1){\BoxedEPSF{chapter3.fig2300a.eps scaled 300}}
\rput(8,3){\BoxedEPSF{chapter3.fig2300b.eps scaled 300}}
\rput(11.5,2){\BoxedEPSF{chapter3.fig2200.eps scaled 600}}
\psline[linewidth=.2mm,linecolor=red]{->}(5,2.9)(7.4,2.9)
\psline[linewidth=.2mm,linecolor=cyan]{->}(5,1.1)(7.4,1.1)
\psline[linewidth=.2mm,linearc=.3]{<-}(8.3,2.9)(8.8,2.9)(8.8,1.1)(8.3,1.1)
\psline[linewidth=.2mm]{->}(9.2,2)(11,2)
\rput(6.2,1.3){${\cyan g_1}$}\rput(6.2,3.1){${\red g_2}$}
\rput(10.1,2.2){$h_1=h_2$}
\end{pspicture}  
  \caption{equivalent intermediate graph coverings: 
$X$ is a bouquet of two loops and $Y$ is its universal
cover, an infinite $4$-valent tree. The intermediate $Z_1,Z_2$, cover
$X$ by the same coverings $h_1=h_2$. The coverings differ in how $Y$ covers
the $Z_i$ via the $g_i$: the first sends the green vertex of $Y$ to the red vertex of $Z_1$
and the second sends it to the blue vertex of $Z_2$. There is an isomorphism 
$Z_1\rightarrow Z_2$ interchanging the two vertices.}
  \label{chapter3:lattices:figure300}
\end{figure}

It is possible to do everything in this section and the next in terms of intermediate 
coverings themselves, and not worry about equivalence at all. Nevertheless, 
the notion of equivalence will become essential later for
the accounting to come out in the wash. 
Figure \ref{chapter3:lattices:figure300} shows a pair of equivalent graph coverings. 

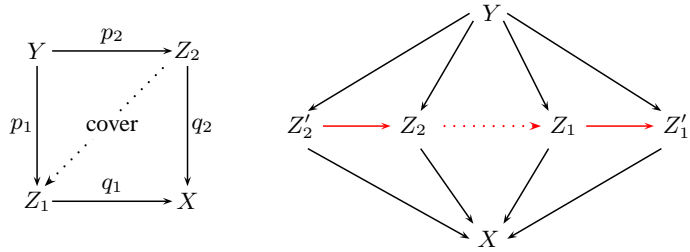
\begin{figure}
  \centering
\begin{pspicture}(12.5,3)
\rput(2,0.5){
\rput(0,2){$Y$}
\rput(0,0){$Z_1$}
\rput(2,2){$Z_2$}
\rput(2,0){$X$}
\psline[linewidth=.2mm]{->}(0,1.8)(0,.2)
\psline[linewidth=.2mm]{->}(.2,2)(1.8,2)
\psline[linewidth=.2mm]{->}(.2,0)(1.8,0)
\psline[linewidth=.2mm]{->}(2,1.8)(2,.2)
\psline[linewidth=.3mm,linestyle=dotted]{<-}(.1,.2)(1.75,1.85)
\rput(-0.2,1){$p_1$}
\rput(1,2.2){$p_2$}
\rput(1,.2){$q_1$}
\rput(2.2,1){$q_2$}
\rput*(1,1){cover}
}
\rput(2,0){
\rput(6.05,3){$Y$}
\psline[linewidth=.2mm]{->}(5.8,2.9)(5.1,1.7)
\psline[linewidth=.2mm]{->}(6.2,2.9)(6.8,1.7)
\psline[linewidth=.2mm]{->}(5.75,3)(3.6,1.7)
\psline[linewidth=.2mm]{->}(6.25,3)(8.3,1.7)
\rput(5,1.5){$Z_2$}\rput(7,1.5){$Z_1$}
\rput(3.5,1.5){$Z_2'$}\rput(8.5,1.5){$Z_1'$}
\psline[linewidth=.2mm]{->}(5.1,1.2)(5.8,.2)
\psline[linewidth=.2mm]{->}(6.8,1.2)(6.2,.2)
\psline[linewidth=.2mm]{->}(3.6,1.2)(5.75,0)
\psline[linewidth=.2mm]{->}(8.3,1.2)(6.25,0)
\psline[linewidth=.2mm,linecolor=red]{->}(3.8,1.5)(4.7,1.5)
\psline[linewidth=.2mm,linecolor=red]{->}(7.3,1.5)(8.2,1.5)
\psline[linewidth=.3mm,linecolor=red,linestyle=dotted]{->}(5.3,1.5)(6.7,1.5)
\rput(6,0){$X$}}
\end{pspicture}  
  \caption{The partial order on $\mathcal{L}(Y,X)$ is well defined.}
  \label{chapter3:lattices:poset:figure100}
\end{figure}

We now turn $\mathcal{L}(Y,X)$ into a poset: one intermediate covering is
``bigger'' than another if the first covers the second. Specifically, if
$Z_1,Z_2\in\mathcal{L}(Y,X)$ then define $Z_1\leq Z_2$ precisely when there is a covering
$Z_2\rightarrow Z_1$ making the diagram on the left of
Figure \ref{chapter3:lattices:poset:figure100} commute. 

Presupposing for a minute that this definition makes sense and gives
a partial order, we have:

\begin{definition}[poset of intermediate covers]
For a fixed covering $f:Y\rightarrow X$ of connected complexes, the set
$\mathcal{L}(Y,X)$ of equivalence classes 
of connected intermediate coverings, together with the partial order $\leq$ defined above 
is called the \emph{poset of intermediate coverings\/} (to $f$).
\end{definition}

It is not hard to check that the order $\leq$ is well defined: suppose that for $i=1,2$, 
we have intermediate coverings $Z_i'$, equivalent to the $Z_i$ via 
isomorphisms $Z_i\leftrightarrow Z_i'$ and $Z_1\leq Z_2$. 
As isomorphisms are nothing other than 
degree one coverings, the red map across the middle of the diagram on the right
of Figure \ref{chapter3:lattices:poset:figure100}
is a covering making the big outside square commute. Thus $Z_1'\leq Z_2'$,
and the order doesn't depend on which representative for the equivalence class
we choose.

\begin{lemma}
The set $\mathcal{L}(Y,X)$ is a poset.
\end{lemma}

\begin{proof}
Reflexivity and transitivity are immediate, as the identity map is a covering and the
composition of coverings is a covering. Only anti-symmetry requires a moments thought: 
suppose we have $Z_1,Z_2\in\mathcal{L}(Y,X)$ with $Z_1\leq Z_2$ and
$Z_2\leq Z_1$, so that there are coverings $Z_1\leftrightarrows Z_2$ 
making the appropriate diagrams
commute. Let $t$ be the composition $Z_1\rightarrow Z_2\rightarrow Z_1$ of 
these two. Then consideration of these commuting diagrams gives $g_1=t h_1$,
were $g_1:Y\rightarrow Z_1$ is the covering, and so by the surjectivity of $g_1$,
$t$ is the identity on $Z_1$. But then the covering $Z_1\rightarrow Z_2$ 
must be injective, ie: of degree $1$, and so an isomorphism. Thus
$Z_1=Z_2$ in $\mathcal{L}(Y,X)$.
\qed
\end{proof}

\subsection{The lattice of intermediate covers}\label{chapter3:lattices:lattice}

In the last section we introduced the poset $\mathcal{L}(Y,X)$ of equivalence classes
of connected coverings intermediate to a fixed covering $Y\rightarrow X$. In this section we
show that we have in fact a \emph{lattice\/}, with join a pullback and meet a pushout.
Because we want all our complexes to be connected and the pullback isn't 
necessarily so, everything in sight has to be pointed, and we use the pointed versions
of the pushout and pullback in Exercises \ref{chapter1:2complexes:exercise550}
and \ref{chapter1:2complexes:exercise650}.

Throughout then, $f:Y_u\rightarrow X_v$ is a fixed pointed covering of 
connected $2$-complexes. All intermediate coverings 
$Y\stackrel{g}{\rightarrow}Z_x\stackrel{h}{\rightarrow}X$ are connected and pointed,
and $\mathcal{L}(Y_u,X_v)$ is the poset of equivalence classes of pointed connected intermediate
coverings. 

We start by showing that we have a meet. Let $Y_u{\rightarrow}(Z_1)_{z_1}{\rightarrow}X_v$ and 
$Y_u{\rightarrow}(Z_2)_{z_2}{\rightarrow}X_v$
be intermediate to $f$, and $Z_1\coprod_YZ_2$ the pushout 
of the coverings $g_1:Y_u{\rightarrow}(Z_1)_{z_1}$
and $g_2:Y_u{\rightarrow}(Z_2)_{z_2}$.
Let $z=[z_1]=[z_2]$, where $q(x)=[x]$ is the quotient map 
arising from the construction of the pushout,
and $(Z_1\coprod_YZ_2)_z$ the resulting
pointed pushout.

By Proposition \ref{chapter3:operations:result100} we have a new element of 
$\mathcal{L}(Y_u,X_v)$ given by the equivalence class of the intermediate 
covering,
$$
Y_u{\rightarrow}(Z_1\coprod_YZ_2)_z{\rightarrow}X_v.
$$
We now need a technical result to ensure that the whole process is well defined:
if the $X_i$ are replaced by equivalent coverings, then the new pushout that
results is equivalent to the old one:

\begin{proposition}\label{chapter3:lattices:lattice:result100}
Let
$(V_1)_{y_1},(V_2)_{y_2}\in\mathcal{L}(Y_u,X_v)$ 
be equivalent to $(Z_1)_{z_1},(Z_2)_{z_2}$ via the isomorphisms,
$$s_1:(Z_1)_{z_1}\rightarrow(V_1)_{y_1}\text{ and }
s_2:(Z_2)_{z_2}\rightarrow(V_2)_{y_2}.
$$ 
Define a map
$s_1\amalg s_2:Z_1\bigcup Z_2\rightarrow V_1\bigcup V_2$ 
between the disjoint unions by
$s_1\amalg s_2|_{Z_1}=s_1$ and
$s_1\amalg s_2|_{Z_2}=s_2$. Then the map 
$$
s:(Z_1\coprod_Y Z_2)_z\rightarrow(V_1\coprod_Y V_2)_{y}, (y=[y_1]=[y_2]),
$$
defined by $sq=q'(s_1\amalg s_2)$,
is an isomorphism making these pointed
pushouts equivalent, where
$q:Z_1\bigcup Z_2\rightarrow Z_1\coprod_Y Z_2$ and 
$q':V_1\bigcup V_2\rightarrow V_1\coprod_Y V_2$ 
are the quotient maps arising in the pushouts.
\end{proposition}

Thus, the pushout can be extended in a well defined way to equivalence classes
of intermediate coverings, so for $(Z_1)_{z_1},(X_2)_{z_2}\in\mathcal{L}(Y_u,X_v)$
we write
$$
(Z_1\coprod_Y Z_2)_z\in\mathcal{L}(Y_u,X_v),
$$
for the pushout of these two equivalence classes.

\begin{proof}[of Proposition \ref{chapter3:lattices:lattice:result100}]
is a tedious but routine diagram
chase.
\qed
\end{proof}

Now, Proposition \ref{chapter3:operations:result100} gives coverings 
$t_i:Z_i\rightarrow Z_1\coprod_Y Z_2$ so that 
$(Z_1\coprod_Y Z_2)_z\leq (Z_i)_{z_i}\,(i=1,2)$
is a lower bound in $\mathcal{L}(Y_u,X_v)$. If $V_y$ is any other lower bound 
then we get coverings $Y_u\rightarrow(Z_i)_{z_i}\rightarrow V_y$, and by the universality
of the pushout, Proposition \ref{chapter1:2complexes:result400},
we have a map $(Z_1\coprod_Y Z_2)_z\rightarrow (V)_y$, which by Proposition 
\ref{chapter3:actionsuniversal:result150} is a covering. Thus
$V_y\leq(Z_1\coprod_Y Z_2)_z$, and the pushout is the meet of the two 
equivalence classes $(Z_i)_{z_i}$. We write 
$$
(Z_1)_{z_1}\wedge (Z_2)_{z_2}=(Z_1\coprod_Y Z_2)_z.
$$

Now to joins, which are similar. Let $Y_u{\rightarrow}(Z_1)_{z_1}{\rightarrow}X_v$ and 
$Y_u{\rightarrow}(Z_2)_{z_2}{\rightarrow}X_v$
be intermediate to $f$, and $Z_1\prod_XZ_2$ the pullback
of the coverings $h_1:(Z_1)_{z_1}\rightarrow X_v$
and $h_2:(Z_2)_{z_2}\rightarrow X_v$.
Let $z=z_1\times z_2$, a vertex of the pullback, 
and $(Z_1\prod_XZ_2)_z$ the
pointed pullback consisting of the connected component
containing $z$.

We have a well-definedness result analogous to Proposition 
\ref{chapter3:lattices:lattice:result100}:

\begin{proposition}\label{chapter3:lattices:lattice:result200}
Let
$(V_1)_{y_1},(V_2)_{y_2}\in\mathcal{L}(Y_u,X_v)$ 
be equivalent to $(Z_1)_{z_1},(Z_2)_{z_2}$ via the isomorphisms,
$$s_1:(Z_1)_{z_1}\rightarrow(V_1)_{y_1}\text{ and }
s_2:(Z_2)_{z_2}\rightarrow(V_2)_{y_2}.
$$ 
Then the map 
$$
s:(Z_1\prod_X Z_2)_{z_1\times z_2}\rightarrow(V_1\prod_Z V_2)_{y_1\times y_2},
$$
defined by $i(z_1\times z_2)=i_1(z_1)\times i_2(z_2)$ is an isomorphism making the pointed
pullbacks equivalent.
\end{proposition}

Proposition \ref{chapter3:operations:result500} gives coverings 
$t_i:(Z_1\prod_X Z_2)_z\rightarrow (Z_i)_{z_i}$ 
for $z=z_1\times z_2$,
so that the 
$(Z_i)_{z_i}\leq (Z_1\coprod_Y Z_2)_z\,(i=1,2)$ and the pullback is 
an upper bound in $\mathcal{L}(Y_u,X_v)$. If $V_y$ is any other upper bound 
we get coverings $V_y\rightarrow (Z_i)_{z_i}\rightarrow X_v$, and by the universality
of the pullback, Proposition \ref{chapter1:2complexes:result600},
we have a map $(V)_y\rightarrow(Z_1\prod_X Z_2)_z$. 
Proposition 
\ref{chapter3:actionsuniversal:result150} again 
gives the map $V_y\rightarrow(Z_1\prod_X Z_2)_z$
is a covering. Thus
$(Z_1\prod_X Z_2)_z\leq V_y$, and the pullback is the join of the two 
equivalence classes $(Z_i)_{z_i}$. We write 
$$
(Z_1)_{z_1}\vee (Z_2)_{z_2}=(Z_1\prod_X Z_2)_z.
$$

Finally, recall from Exercise \ref{chapter3:lattices:exercise50}
that a $\hat{1}$ in a poset is a unique maximal element and a $\hat{0}$
is a unique minimal element. 

\begin{theorem}[lattice of intermediate coverings]
\label{chapter3:lattices:lattice:result300}
The poset $\mathcal{L}(Y_u,X_v)$ of pointed connected covers
intermediate to a fixed covering $f:Y_u\rightarrow X_v$
is a lattice 
with join $(Z_1)_{z_1}\vee(Z_2)_{z_2}$ the pullback 
$(Z_1\prod_X Z_2)_{z_1\times z_2}$, meet $(Z_1)_{z_1}\wedge(Z_2)_{z_2}$
the pushout $(Z_1\coprod_Y Z_2)_{[z_i]}$, unique minimal element $\hat{0}=X_v$
and unique maximal element $\hat{1}=Y_u$.
\end{theorem}

\section{Notes on Chapter \thechapter}\label{chapter3:notes}

\chapter{Galois Theory}\label{chapter4}

Galois theory arises whenever we have the following situation:
$\mathscr{A}$ is some object and
$\gal(\mathscr{A})$ is its group of ``symmetries'', or Galois
group. If $\mathscr{B}\subset\mathscr{A}$ is a sub-object, then
there is a subgroup $H\subset\gal(\mathscr{A})$ consisting of those
symmetries of $\mathscr{A}$ that act trivially on $\mathscr{B}$. On the otherhand, if
$H\subset\gal(\mathscr{A})$ is a subgroup, there is a sub-object
$\mathscr{B}\subset\mathscr{A}$
on which the action of $H$ has been cancelled out.

The first serious theorem in any Galois theory then says that this
correspondence between sub-objects of $\mathscr{A}$ and subgroups 
of $\gal(\mathscr{A})$
is perfect: the sub-objects of
$\mathscr{A}$ form a lattice, as do the subgroups of
$\gal(\mathscr{A})$, and these two lattices are \emph{anti}-isomorphic.

In the classical Galois theory, $\mathscr{A}$ is
an extension $E\subset F$ of fields and $\gal(\mathscr{A})$ the
field automorphisms of $E$ fixing $F$ pointwise. For us,
$\mathscr{A}$ will be a covering $Y\rightarrow X$ of $2$-complexes,
and $\gal(\mathscr{A})$ the automorphisms of $Y$ that permute the
fibers of the covering.

Throughout this chapter we make the running assumption that
\emph{all complexes are connected.}

\section{Galois groups}\label{chapter4:galoisgroups}

\subsection{Automorphisms and Galois
  groups}\label{chapter4:galoisgroups:automorphisms}

\begin{definition}[automorphism of a covering]
\label{chapter4:automorphisms:definition100}
Let $f:Y\rightarrow X$ be a covering
with $Y,X$ connected. A {\em covering automorphism\/} 
(or \emph{Galois automorphism\/}) of $f$ is an 
isomorphism $a:Y\rightarrow Y$ making the diagram,
$$
\begin{pspicture}(0,0)(4,1.5)
\rput(0,-.25){
\rput(1,1.5){$Y$}
\rput(3,1.5){$Y$}
\rput(2,0.4){$X$}
\rput(2,1.7){$a$}
\rput(1.125,.825){$f$}
\rput(2.78,.825){$f$}
\psline[linewidth=.1mm]{->}(1.3,1.5)(2.7,1.5)
\psline[linewidth=.1mm]{->}(1.2,1.25)(1.8,.55)
\psline[linewidth=.1mm]{->}(2.8,1.25)(2.25,.55)
}
\end{pspicture}
$$
commute. If $f:Y_u\rightarrow X_v$ is a pointed covering
then $a$ is a pointed isomorphism $a:Y_u\rightarrow Y_{u'}$ making the
diagram commute.
\end{definition}

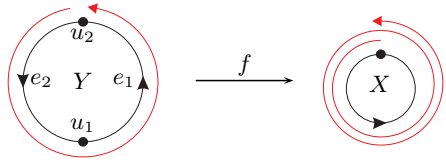
\begin{figure}
  \centering
\begin{pspicture}(0,0)(12.5,2)
\rput(4,1){\BoxedEPSF{chapter3.fig200.eps scaled 750}}
\rput(8,1){\BoxedEPSF{chapter3.fig300.eps scaled 750}}
\rput(4,1){$Y$}\rput(7.95,.95){$X$}\rput(6.15,1.2){$f$}
\psline[linewidth=.2mm]{->}(5.5,1)(6.8,1)
\rput(4,.4){$u_1$}\rput(4,1.6){$u_2$}\rput(4.55,1){$e_1$}\rput(3.45,1){$e_2$}
\end{pspicture}
  \caption{a simple graph covering: the two vertices of $Y$ cover the 
single vertex of $X$ and the two arcs of $Y$ similarly. There is a
single non-trivial covering automorphism interchanging the vertices
$u_1$ and $u_2$ and the edges $e_1$ and $e_2$.}
  \label{chapter4:automorphisms:figure50}
\end{figure}

A covering automorphism is thus an automorphism of $Y$ that permutes
the fibers of the covering. 
In topology covering automorphisms are often called {\em deck transformations\/}.
Figure \ref{chapter4:automorphisms:figure50} shows a graph
covering with exactly two covering automorphisms. One is the identity
$\id:Y\rightarrow Y$, and the other is the automorphism of $Y$ that
interchanges the vertices $u_1$ and $u_2$ and interchanges the edges
$e_1$ and $e_2$.

Notice that $Y$ has other automorphisms, but they are not covering
automorphisms. For example, the map $Y\rightarrow Y$ that fixes the
vertices $u_1$ and $u_2$, and sends the edge $e_1$ to $e_2^{-1}$ and
$e_2$ to $e_1^{-1}$ is an automorphism, but it is not a covering
automorphism as $e_1$ and $e_2^{-1}$ lie in different fibers of the covering.

\begin{figure}
  \centering
\begin{pspicture}(0,0)(12.5,3)
\rput(-9.5,0){
\rput(12.5,1.5){\BoxedEPSF{chapter1.fig4200.eps scaled 500}}
\rput(2.5,-4.5){
\rput(11.6,7.4){$\ss_1$}
\rput(11.6,4.6){$\ss_2$}
\rput(12.45,6.4){${\red\partial\ss_1}$}
\rput(12.45,5.8){${\red\partial\ss_2}$}
\rput(9.7,5.45 ){${v}_1$}
\rput(10.2,6.55){${v}_2$}
\rput(10.6,5.45){${e}_2$}
\rput(9.4,6.6){${e}_1$}
}}
\rput(8.7,0){
\rput(1.25,1.25){\BoxedEPSF{chapter1.fig4300.eps scaled 500}}
\rput(2.2,.2){$\ss$}\rput(.4,1.55){$v$}\rput(2.2,2.4){$v$}\rput(1.2,1.05){$e$}
\rput(1.2,2.5){$e$}
}
\rput(-.25,0){
\psline[linewidth=.2mm]{->}(6,1.5)(8,1.5)
\rput(7,1.75){$f$}
}
\rput(0,1.5){$Y$}\rput(12.5,1.5){$X$}
\psline[linewidth=.2mm,linearc=.3,linecolor=blue]{<->}(1,.5)(.5,.5)(.5,2.5)(1,2.5)
\rput(.35,1.5){${\blue a}$}
\end{pspicture}
  \caption{The non-trivial automorphism of the covering of Figure
    \ref{chapter4:automorphisms:figure50} extends to a unique covering
    automorphism of the covering of 
    Figure \ref{chapter3:basics:figure300}. The faces $\ss_1$ and
    $\ss_2$ are interchanged by $\aa$ and the face isomorphism
    $Y^{\ss_1}\rightarrow Y^{\ss_2}$ is a $1/2$-turn.}
\label{chapter4:automorphisms:figure100}  
\end{figure}
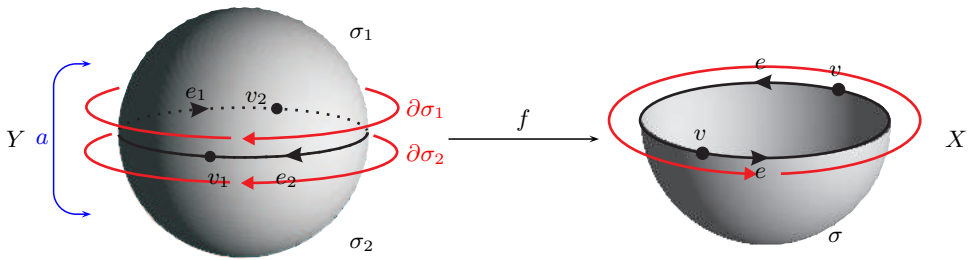

Figure \ref{chapter4:automorphisms:figure100} 
extends this example to
a covering of $2$-complexes
as in
\S\ref{chapter3:basics:coverings}.
The identity automorphism $\id: Y\rightarrow Y$ is a covering
automorphism, as is the automorphism
that interchanges the vertices $u_i$, the edges $e_i$ and the faces
$\ss_i$. These are the only two.

One can check that 
there is no automorphism of the $Y$ of Figure \ref{chapter4:automorphisms:figure100}  
that
fixes the vertices and edges but interchanges the faces. 
In particular notice that distinct covering automorphisms of $Y$ restrict to
distinct covering automorphisms of the $1$-skeleton of $Y$.

\begin{exercise}
\label{chapter4:automorphisms:exercise100}  
If $f:Y\rightarrow X$ is a covering and $a_1,a_2:Y\rightarrow Y$
covering automorphisms, show that their composition $a_2a_1$ and
$a_1^{-1}$ are covering automorphisms. Show that the identity map
$\id:Y\rightarrow Y$ is a covering automorphism. Deduce that the set
of covering automorphisms forms a group.
\end{exercise}

\begin{definition}[Galois group of a covering]
\label{chapter4:automorphisms:definition200}  
The covering automorphisms of the covering $f:Y\rightarrow X$ form a
group called the \emph{Galois group\/} of the covering, denoted
$$
\gal(Y\stackrel{f}{\rightarrow}X)\text{ or just }\gal(Y,X).
$$
If the covering is
pointed we write $\gal(Y_u,X_v)$.
\end{definition}

We now come to some basic properties of the action of the Galois group:
Recall from \S\ref{chapter1:2complexes:maps} that a group acts freely on a $2$-complex $Y$
precisely when it acts freely on the vertices of $Y$.

\begin{lemma}\label{chapter4:automorphisms:result100}  
(i). The action of $\gal(Y,X)$ on $Y$ is orientation preserving.
(ii). The effect of a covering automorphism $a\in\gal(Y,X)$
is completely determined by the image of a single vertex. In particular, 
the Galois group acts freely on $Y$.
\end{lemma}

\begin{proof}
(i). Let $a\in\gal(Y,X)$. If $x$ is an edge or face of $Y$ then both 
$x$ and $a(x)$ lie in the same 
fiber of the covering $f$, so that
if $a(x)=x^{-1}$ then $f(x)=f(x)^{-1}$, a contradiction. Thus, the Galois group
acts without inversions.
(ii). Let $a\in\gal(Y,X)$ with $a(u)=u'$ 
for $u,u'\in Y$ vertices. 
If $\gamma$ is a path in $Y$ starting at $u$ (and covering $f(\gamma)$), then $a(\gamma)$ 
is a path starting at $u'$ and also covering $f(\gamma)$, as
$a$ is a covering automorphism. By uniqueness of path lifting,
$a(\gamma)$ must be the lift of $f(\gamma)$ to $u'$. Thus the effect
of $a$ on the $1$-skeleton is completely determined by $a(u)=u'$.
Faces are similar, using the uniqueness of face lifting.
\qed
\end{proof}

\begin{figure}
  \centering
\begin{pspicture}(0,0)(12.5,4)
\rput(3.75,3.15){$u'$}\rput(3.75,1.8){$u$}
\rput(9,3.5){$x'$}\rput(9,2.2){$x$}
\rput(6,2.3){$\gamma$}\rput(6,1.1){$f(\gamma)$}
\rput(2.5,.1){$X$}\rput(2.5,1.4){$Y$}\rput(2.5,2.7){$Y$}
\psbezier[linewidth=.2mm,linecolor=blue]{->}(3.55,1.8)(2.95,2.4)(2.95,2.5)(3.55,3.1)
\psbezier[linewidth=.2mm,linecolor=blue]{->}(9.2,2.2)(9.8,2.8)(9.8,2.85)(9.2,3.45)
\rput(3.4,2.45){${\blue a}$}\rput(9.4,2.825){${\blue a}$}
\rput(0,0){\psline[linewidth=.4mm,linestyle=dotted,linecolor=red]{->}(4,3)(4,.7)}
\rput(1.25,.4){\psline[linewidth=.4mm,linestyle=dotted,linecolor=red]{->}(4,3)(4,.7)}
\rput(2.5,0){\psline[linewidth=.4mm,linestyle=dotted,linecolor=red]{->}(4,3)(4,.7)}
\rput(3.4,.35){\psline[linewidth=.4mm,linestyle=dotted,linecolor=red]{->}(4,3)(4,.7)}
\rput(4.7,.4){\psline[linewidth=.4mm,linestyle=dotted,linecolor=red]{->}(4,3)(4,.7)}
\rput(6.25,2){\BoxedEPSF{chapter4.fig100.eps scaled 1000}}
\end{pspicture}
  \caption{cover and lift}
\label{chapter4:automorphisms:figure200}  
\end{figure}
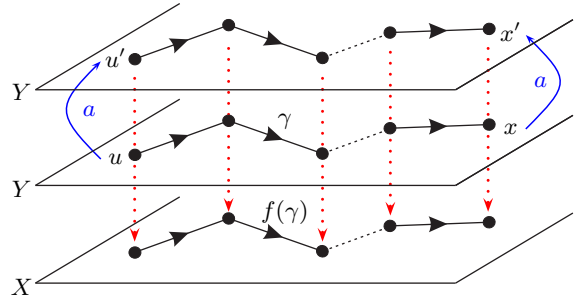

The technique used in the proof of Lemma
\ref{chapter4:automorphisms:result100} is called ``cover and lift''
(see Figure \ref{chapter4:automorphisms:figure200}): start with the 
path $\gamma$; it covers the path $f(\gamma)$ and this in turn lifts
to the image path $a(\gamma)$. 

\begin{lemma}\label{chapter4:automorphisms:result200}  
Let $f:Y\rightarrow X$ be a covering and $f:Y^{(1)}\rightarrow
X^{(1)}$ its restriction to the $1$-skeletons. The restriction of any
$a\in\gal(Y,X)$ to the one skeletons is a covering automorphism,
and this induces an injective homomorphism $\gal(Y,X)\rightarrow\gal(Y^{(1)},X^{(1)})$.
\end{lemma}

\begin{proof}
As $ab$ restricted to the $1$-skeleton is just $b$ restricted to
the $1$-skeleton followed by $a$ restricted to the $1$-skeleton, we
get a homomorphism. Lemma \ref{chapter4:automorphisms:result100} shows that a covering
automorphism $a$ is completely determined by its efect on the
$1$-skeleton, thus giving injectivity.
\qed
\end{proof}

\subsection{Constructing automorphisms}
\label{chapter4:galoisgroups:construction}

The explicit construction of automorphisms is achieved by the following
result:

\begin{proposition}\label{topological:galois:result200}
Let $f:Y\rightarrow X$ be a covering with $f(u)=v$, 
and $u'$ another vertex in the fiber $f^{-1}(v)$.
Then there is a covering automorphism $a\in\gal(Y,X)$ with $a(u)=u'$
if and only if
for any
closed path $\gamma$ at $v$ with lifts $\mu_1,\mu_2$ 
at $u,u'$, we have
\begin{equation*}
  \label{eq:1}
\text{$\mu_1$ is closed $\Leftrightarrow$
$\mu_2$ is closed.
}
\tag{$\dag$}
\end{equation*}
\end{proposition}

Indeed, the covering automorphism $a$ is unique by Lemma
\ref{chapter4:automorphisms:result100} and
comes about as follows: use ``cover and
lift'' as in Figure \ref{chapter4:automorphisms:figure200} to get the effect of $a$ on the
vertices. There is a well-defined issue, and condition (\dag) is exactly
what is needed to resolve it. The effect of $a$ on the edges is given by
path-lifting 
and on the faces by face lifting. 

\begin{proof}
As automorphisms send closed paths to closed paths and non-closed paths
to non-closed paths, the only if direction is clear. On the other hand, 
if $x$ is a vertex of $Y$ and $\mu$ a path from $u$ to $x$, then
define $a(x)$ to be the terminal vertex of the lift of $f(\mu)$ to
$u'$.
If $e$ is an edge let $a(e)$ be the lift of $f(e)$ to the vertex
$a s(e)$. If $\ss$ is a face and $x$ some vertex appearing in its
boundary, let $a(\ss)$
be the lift of $f(\ss)$ to $a(x)$.

If $\mu'$ is another path from $u$ to $x$
(so that $\mu_1=\mu'\mu^{-1}$ is closed at $u$) then
$f(\mu')f(\mu)^{-1}$ is a closed path at $v$ lifting to
$\mu_1$, hence by ($\dag$), it lifts to a closed path at $u'$. Thus
the lifts to $u'$ of $f(\mu)$ and $f(\mu')$ end at the same vertex and so
$a$ is well defined on the vertices.
There is also a choice of boundary vertex involved in the definition of
$a$ on the edges and faces: we chose the vertex $s(e)$ in the
boundary of the edge $e$ and the vertex $x$ in the boundary of the
face $\ss$. We show first that an arbitrary choice extends $a$ to a
covering automorphism, and then appeal to Lemma
\ref{chapter4:automorphisms:result100}(ii) to see that a
different choice gives the same $a$.

We have that $sa(e)$ is by definition equal to $a s(e)$, and the
uniqueness of the lift of $f(e)^{-1}$ to $ta(e)$ gives
$a(e^{-1})=a(e)^{-1}$. If $\ss$ is a face, then splicing together the two diagrams provided
by the cover $\ss\mapsto f(\ss)$ and the lift $a(\ss)\mapsto f(\ss)$
gives the required commuting diagram for $\ss\mapsto a(\ss)$, and so
$a$ is a map. Interchanging the roles of the vertices $u$ and $u'$
gives a map $b:Y\rightarrow Y$ with the uniqueness of path and face
lifting giving $ab=ba=\id$ on $Y$. It is easy to see that
$a$ is dimension-preserving, so that we have an
automorphism. Finally, $x$ and 
$a(x)$ lie in the same fiber of the covering, for any cell
$x$, whence $fa=f$. 
\qed
\end{proof}

\begin{example}
Figure \ref{chapter4:galoisgroups:figure250} shows a covering with a
pair of vertices (ringed) satisfying (\ref{eq:1}) of the
Proposition with the obvious $1/2$-turn rotation the resulting
automorphism. 

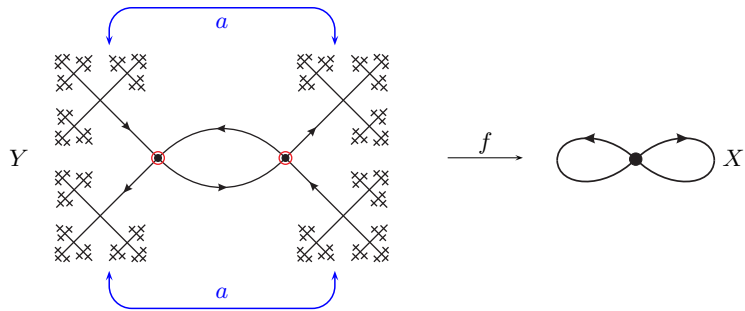
\begin{figure}[h]
  \centering
\begin{pspicture}(0,0)(12,4)
\rput(.5,0.5){
\rput(3,1.5){\BoxedEPSF{chapter4.fig4100.eps scaled 500}}
\rput(8.5,1.5){\BoxedEPSF{chapter3.fig2100.eps scaled 650}}
\psline[linewidth=.1mm]{->}(6,1.5)(7,1.5)
}
\psline[linewidth=.2mm,linearc=.3,linecolor=blue]{<->}(2,3.5)(2,4)(5,4)(5,3.5)
\psline[linewidth=.2mm,linearc=.3,linecolor=blue]{<->}(2,.5)(2,0)(5,0)(5,.5)
\rput(3.5,3.8){${\blue a}$}\rput(3.5,.2){${\blue a}$}
\rput(.8,2){$Y$}\rput(10.3,2){$X$}\rput(7,2.2){$f$}
\end{pspicture}
  \caption{vertices satisfying condition (\ref{eq:1}) of 
Proposition \ref{topological:galois:result200},
yielding a covering automorphism $a\in\gal(Y,X)$}
\label{chapter4:galoisgroups:figure250}
\end{figure}
\end{example}

\subsection{From covers to subgroups} 
\label{chapter4:galoisgroups:coverstosubgroups}

If $f:Y\rightarrow X$ is a covering and $\gal(Y,X)$ is its Galois
group, then in this section and the next we show how an intermediate covering
$Y\rightarrow Z\rightarrow X$ 
gives rise to a subgroup of $\gal(Y,X)$ and vice-versa. Indeed, recalling the
definition of equivalent intermediate coverings from 
\S\ref{chapter3:lattices:poset}, 
equivalent coverings give the same subgroup. Thus, an equivalence
class of intermediate coverings gives rise to a subgroup of the Galois
group. 

\begin{figure}[h]
  \centering
\begin{pspicture}(0,0)(12,4.75)
\rput(4.5,2.5){
\rput(0.8,1.5){$Y$}
\rput(3.2,1.5){$Y$}
\rput(2,0.4){$Z_2$}
\rput(1.95,-1.1){$Z_1$}
\rput(1.95,-2.6){$X$}
\rput(2,1.7){$a$}
\psline[linewidth=.1mm]{->}(1.1,1.5)(2.9,1.5)
\psline[linewidth=.1mm]{->}(1,1.25)(1.75,.55)
\psline[linewidth=.1mm]{->}(3,1.25)(2.25,.55)
\psline[linewidth=.1mm]{->}(1.95,.2)(1.95,-.85)
\psline[linewidth=.1mm]{->}(1.95,-1.3)(1.95,-2.35)
\psbezier[linewidth=.1mm,showpoints=false]{->}(.8,1.25)(.8,.55)(.8,.05)(1.75,-.85)
\psbezier[linewidth=.1mm,showpoints=false]{->}(3.2,1.25)(3.2,.55)(3.2,.05)(2.15,-.85)
\psbezier[linewidth=.1mm,showpoints=false]{->}(.7,1.25)(.7,.05)(.7,-1.1)(1.75,-2.35)
\psbezier[linewidth=.1mm,showpoints=false]{->}(3.3,1.25)(3.3,.05)(3.3,-1.1)(2.15,-2.35)
\pspolygon[linearc=.3,linecolor=red,
showpoints=false](.3,1.9)(3.7,1.9)(3.7,1)(1.95,0)(.3,1)
\pspolygon[linearc=.3,linecolor=blue,
showpoints=false](.1,2.1)(3.9,2.1)(3.9,0)(1.95,-1.5)(.1,0)
}
\end{pspicture}
  \caption{from intermediate covers to subgroups of the Galois group.}
\label{chapter4:galoisgroups:figure260}
\end{figure}
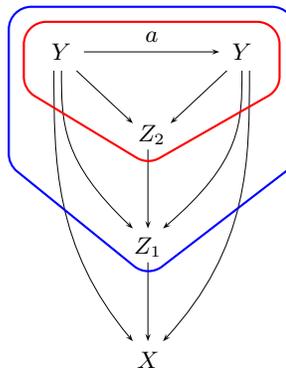

Let $Y\rightarrow Z_2\rightarrow Z_1\rightarrow X$ be a sequence of
coverings intermediate to $f$,
and let $a$ be an element of the Galois group $\gal(Y,Z_2)$, so that the
triangle ringed in red in Figure \ref{chapter4:galoisgroups:figure260} commutes.

The unlabeled maps are the coverings. It is easy to see that the
triangle ringed in blue must then also commute, so that $a$ can be identified
with an element of $\gal(Y,Z_1)$, and we get a map
$\gal(Y,Z_2)\rightarrow\gal(Y,Z_1)$.
The following is then immediate,

\begin{lemma}\label{chapter4:automorphisms:result1000}  
The map $\gal(Y,Z_2)\rightarrow\gal(Y,Z_1)$ is an injective homomorphism.
\end{lemma}

From now on we will just identify $\gal(Y,Z_2)$ with the subgroup of
$\gal(Y,Z_1)$ consisting of those $a$ that make the red triangle commute
in the diagram above. In particular, when $Y\rightarrow Z\rightarrow
X$ is intermediate, we can identify $\gal(Y,Z)$ with a subgroup of $\gal(Y,X)$.

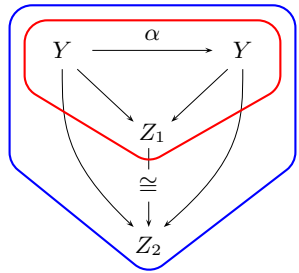
\begin{figure}[h]
  \centering
\begin{pspicture}(0,0)(12,3.5)
\rput(4.5,1.45){
\rput(0.8,1.5){$Y$}
\rput(3.2,1.5){$Y$}
\rput(2,0.4){$Z_1$}
\rput(1.95,-1.1){$Z_2$}
\rput(2,1.7){$\aa$}
\psline[linewidth=.1mm]{->}(1.2,1.5)(2.8,1.5)
\psline[linewidth=.1mm]{->}(1,1.25)(1.75,.55)
\psline[linewidth=.1mm]{->}(3,1.25)(2.25,.55)
\psline[linewidth=.1mm]{->}(1.95,.2)(1.95,-.85)
\rput*(1.95,-.3){$\cong$}
\psbezier[linewidth=.1mm,showpoints=false]{->}(.8,1.25)(.8,.55)(.8,.05)(1.75,-.85)
\psbezier[linewidth=.1mm,showpoints=false]{->}(3.2,1.25)(3.2,.55)(3.2,.05)(2.15,-.85)
\pspolygon[linearc=.3,linecolor=red,
showpoints=false](.3,1.9)(3.7,1.9)(3.7,1)(1.95,0)(.3,1)
\pspolygon[linearc=.3,linecolor=blue,
showpoints=false](.1,2.1)(3.9,2.1)(3.9,0)(1.95,-1.5)(.1,0)
}
\end{pspicture}  
  \caption{the subgroup of the Galois group does not depend on the representative of the
equivalence class.}
  \label{chapter4:galoisgroups:figure265}
\end{figure}

Now suppose we have a pair $Y\rightarrow Z_i\rightarrow X$, $(i=1,2)$
of coverings intermediate to $f$ that are equivalent
via an isomorphism $Z_1\rightarrow Z_2$. 
Then in the diagram of Figure \ref{chapter4:galoisgroups:figure265}
the red triangle commutes if and only if the blue triangle commutes
(just interchange the roles of $Z_1$ and $Z_2$). In particular, an
$a\in\gal(Y,X)$ lies in the subgroup $\gal(Y,Z_1)$ if and only if it
lies in the subgroup $\gal(Y,Z_2)$. 

These two subgroups thus coincide, and we can associate in a
well-defined manner a subgroup of $\gal(Y,X)$ with an equivalence
class of coverings intermediate to $f$.

\begin{example}
Figure \ref{chapter4:automorphisms:figure270} is the pair of equivalent
coverings from 
Figure \ref{chapter3:lattices:figure300},
with $X$ a bouquet of two loops, $Y$
its universal cover (the infinite $4$-valent tree) and $Z_1,Z_2$ 
shown in the middle of Figure \ref{chapter4:automorphisms:figure270}.

\begin{figure}[h]
  \centering
\begin{pspicture}(0,0)(12,4)
\rput(2,2){\BoxedEPSF{chapter4.fig3000.eps scaled 500}}
\psline[linewidth=.2mm,linecolor=cyan]{->}(4.2,2)(5.2,2)
\rput(6,2){\BoxedEPSF{chapter3.fig2300.eps scaled 600}}
\psline[linewidth=.2mm,linecolor=red]{->}(7.8,2)(6.8,2)
\rput(10,2){\BoxedEPSF{chapter4.fig3100.eps scaled 500}}
\rput(5.6,1.1){$x_1$}\rput(6.4,2.9){$x_2$}
\rput(2.25,2){$u$}\rput(9.7,2){$u$}
\end{pspicture}
  \caption{}
\label{chapter4:automorphisms:figure270}  
\end{figure}

The coverings $Y \rightarrow Z_i$ are shown using the little 
circles and squares.
The Galois group $\gal(Y_u,Z_{x_i})$ of the covering $Y_u\rightarrow
Z_{x_i}$ acts regularly on the vertices in each fiber. In
particular, an automorphism $a\in\gal(Y_u,Z_{x_1})$ of the lefthand version of $Y$ gives an
automorphism $a\in\gal(Y_u,Z_{x_2})$ of the righthand version and vice-versa.
\end{example}

\subsection{Subgroups to covers from the ``top down''} 
\label{chapter4:galoisgroups:subgroupstocovers}

We constructed covers from the bottom-up in
\S\ref{chapter3:actionsuniversal:bottom_up}. Here
is a reverse process:
if $f:Y\rightarrow X$ is a covering and $\gal(Y,X)$ its Galois group,
we show how a subgroup of $\gal(Y,X)$ gives rise to an intermediate covering
$Y\rightarrow Z\rightarrow X$. Passing to the
equivalance class of this covering, we get a subgroup giving
rise to an equivalence class of coverings. 

\begin{lemma}\label{chapter4:automorphisms:result1100}  
If $H_1\subset H_2\subset\gal(Y,X)$ are subgroups then there are coverings
$$
Y\rightarrow Y/H_1
\rightarrow Y/H_2
\rightarrow X
$$
intermediate to $f$.
\end{lemma}

\begin{proof}
By Proposition \ref{chapter1:2complexes:result200} and 
Lemma \ref{chapter4:automorphisms:result100}(i)   
we can form the quotient complexes $Y/H_i\,(i=1,2)$, as the subgroups
are acting orientation preservingly. The action is also free by
Lemma \ref{chapter4:automorphisms:result100}(ii), 
so the quotient maps $Y\rightarrow Y/H_i$ are coverings by 
Proposition \ref{chapter3:actionsuniversal:result100}.
Two applications of Proposition \ref{chapter3:actionsuniversal:result150}
give that $Y/H_1\rightarrow Y/H_2$ and
$Y/H_2\rightarrow X$ are coverings.
\qed
\end{proof}

In particular, letting $H_1=H_2=H$, we can associate to
$H\subset\gal(Y,X)$ the intermediate covering 
$Y\rightarrow Y/H\rightarrow X$,
and by passing to its equivalence class we get, associated to $H$, an element of the
lattice $\mathcal{L}(Y,X)$ of intermediate covers.

\subsection{The Galois group and the fundamental group} 
\label{chapter4:galoisgroups:fundamentalgroup}

Let $f:Y\rightarrow X$ be a covering with $f(u)=v$ and let $G=\pi_1(X,v)$
and $H=f_*\pi_1(Y,u)\subset G$. 
We saw in \S\ref{chapter4:galoisgroups:construction} that elements of $G$ can be 
associated in a well defined way with the vertices of $f^{-1}(v)$: let
$\gamma,\gamma'$ be closed paths at $v$ such that $g_\gamma=g_{\gamma'}$
in $\pi_1(X,v)$. Lifting the homotopic $\gamma,\gamma'$ to $u$ gives 
$\mu,\mu'$ homotopic. These must therefore end at the same vertex $u'$.

When is there a covering automorphism sending $u$ to $u'$?

\begin{lemma}\label{chapter4:automorphisms:result1200}  
The vertices $u,u'$ have property ($\dag$) of Proposition 
\ref{topological:galois:result200} if and only if $g\in G$
normalizes the subgroup $H$, ie: $g^{-1}Hg=H$.
\end{lemma}

\begin{proof}
Write $g=g_\gamma$ and recall from Corollary 
\ref{chapter3:basics:result500}  that a closed path $\lambda$ at $v$
represents an element $h\in H$ if and only if the lift $\nu$ of
$\lambda$ to $u$ is closed. 
The result follows by observing that if $\mu$ is the lift of
$\gamma$ to $u$ and $\nu'$ the lift of
$\lambda$ to $u'=t(\mu)$, then
$\mu^{-1}\nu'\mu$ is the the lift of $\gamma^{-1}\lambda\gamma$ to $u$.
\qed
\end{proof}

\begin{figure}
  \centering
\begin{pspicture}(0,0)(12.5,4.5)
\rput(3.75,3.6){$u'$}\rput(3.75,2.3){$u$}
\rput(9,2.7){${\blue x}$}
\rput(6,2.3){$\nu$}
\rput(8.1,1.1){$\gamma$}
\rput(2.5,.2){$X$}\rput(2.5,1.9){$Y$}\rput(2.5,3.2){$Y$}
\rput(0,0){\psline[linewidth=.4mm,linestyle=dotted,linecolor=red]{->}(4,3.4)(4,.9)}
\rput(6.25,2.25){\BoxedEPSF{chapter4.fig200.eps scaled 1000}}
\rput*(9.35,4){${\blue a_g(x)}$}
\rput*(4.45,2.9){$\mu$}
\end{pspicture}
  \caption{finding the image of a vertex $x$ under the covering
    automorphism $a_g$.}
\label{chapter4:automorphisms:figure300}  
\end{figure}
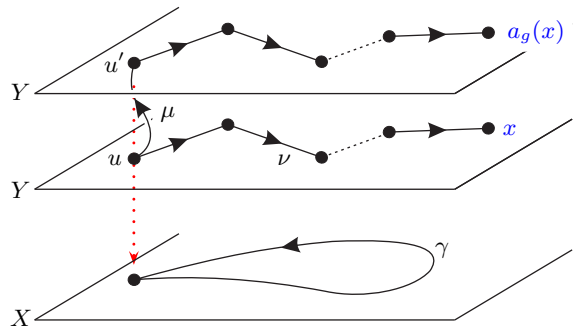

Writing $N_G(H)$ for the 
elements of $G$ normalizing $H$, the lemma provides for $g\in N_G(H)$
a (unique) covering
automorphism $a_g\in\gal(Y,X)$ sending $u$ to $u'$. 
Figure \ref{chapter4:automorphisms:figure300} reminds us 
how this covering automorphism comes about:
let $g=g_\gamma$ and lift to
$\mu$ at $u\in Y$ with $u'=t(\mu)$.
To get the image $a_g(x)$ of a vertex $x\in Y$ let $\nu$ be a
path from $u$ to $x$ and lift $f(\nu)$ to $u'$. Then $a_g(x)$ is the terminal vertex
of this lift.

\begin{exercise}
\label{chapter4:automorphisms:exercise125}    
Let $X$ be a $2$-complex with a single vertex and $H\subset\pi_1(X,v)=G$
a subgroup with $f:X\kern-2pt\uparrow\kern-2pt H\rightarrow X$ 
the covering of \S\ref{chapter3:actionsuniversal:bottom_up}.
Show that for $k\in\pi_1(X,v)$, and $Hg$ a vertex of 
$X\kern-2pt\uparrow\kern-2pt H$ as in Definition 
\ref{chapter3:actionsuniversal:definition85}, we have $a_k(Hg)=H(kg)$.
Compare with Exercise \ref{chapter3:actionsuniversal:exericse950}, and notice that
the action is well defined, ie: $Hg_1=Hg_2\Rightarrow H(kg_1)=H(kg_2)$
if and only if $k\in N_G(H)$.
\end{exercise}

Here is the
result that ties together the fundamental groups of $X$ and $Y$ and
the Galois group of the covering:

\begin{proposition}\label{chapter4:automorphisms:result1300}  
The map $g\mapsto a_{g^{-1}}$ is a surjective homomorphism
$N_G(H)\rightarrow\gal(Y,X)$ with kernel $H$, inducing an isomorphism,
$$
N_G(H)/H\stackrel{\cong}{\rightarrow}\gal(Y,X).
$$
\end{proposition}

\begin{proof}
Successively lifting representatives for $g$ and $h$, and using the
scheme in Figure \ref{chapter4:automorphisms:figure300} gives $a_{gh}=a_ga_h$.
Writing $\theta(g)=a_{g^{-1}}$ thus gives
$\theta(gh)=\theta(h)\theta(g)$ and so $\theta$ is a homomorphism
(remember: we read $gh$ from left to right in $\pi_1(X,v)$ but
$a_ga_h$ from right to left in the Galois group). The elements of $H$
are represented by the $\gamma$ lifting to closed paths at $u$, hence
giving an $a_h$ that fixes $u$, and so must be the identity covering
automorphism as automorphisms are completely determined by their
effect on a single vertex. Thus the kernel is $H$. If
$a\in\gal(Y,X)$ sends $u$ to $u'$ and $\mu$ is a path from $u$ to
$u'$, then $g=g_{f(\mu)}$
normalizes $f_*\pi_1(Y,u)$, and so $a_g=a$.
\qed
\end{proof}

\subsection{Excising simply-connected subcomplexes}
\label{chapter4:galoisgroups:excision}

We saw in \S\ref{chapter3:basics:liftingsimplyconnected} 
that if $Y\rightarrow X$ is a covering and $Z\subset X$ 
a connected, simply connected subcomplex, then $Z$ 
can be lifted to $Y$ to give a
new covering, with $Z$ and its lift ``excised''. In this section we show
that this process has no effect on the Galois group.

Let $f:Y_u\rightarrow X_v$ be a covering and $Z\subset X$
connected and simply-connected with $v\in Z$. Let 
$f^{-1}(Z)=\bigcup_I Z_i\subset Y$ with the
$Z_i$ the connected components, simply connected by Proposition 
\ref{chapter3:basics:result650}.
Let $f':Y/Z_i\rightarrow X/Z$ be the induced covering of Theorem
\ref{chapter3:basics:result660} with
$$
\begin{pspicture}(0,0)(4,2)
\rput(1.2,1.8){$Y$}\rput(1.2,0.2){$X$}
\rput(2.9,1.8){$Y/Z_i$}\rput(2.9,0.2){$X/Z$}
\psline[linewidth=.1mm]{->}(1.2,1.5)(1.2,.5)
\psline[linewidth=.1mm]{->}(2.8,1.5)(2.8,.5)
\psline[linewidth=.1mm]{->}(1.5,1.8)(2.5,1.8)
\psline[linewidth=.1mm]{->}(1.5,.2)(2.5,.2)
\rput(2,2.05){$q'$}\rput(2,.4){$q$}\rput(1,1){$f$}\rput(3.05,1){$f'$}
\end{pspicture}
$$
commuting and $q,q'$ the respective quotient maps. From Chapter
\ref{chapter2} the induced homomorphism
$$
q_*:\pi_1(X,v)\rightarrow\pi_1(X/Z,q(v))
$$
is an isomorphism.

\begin{proposition}
\label{chapter4:automorphisms:result1400}  
The isomorphism $q_*$ induces an isomorphism
$$
\gal(Y,X)\stackrel{\cong}{\rightarrow}\gal(Y/Z_i,X/Z).
$$
\end{proposition}

\begin{proof}
``Hit'' the commuting diagram above with the
$\pi_1$-functor of Chapter \ref{chapter2} to get a commuting diagram
of fundamental groups and group homomorphisms with 
$$
q'_*:\pi_1(Y,u)\rightarrow\pi_1(Y/Z_i,q'(u)),
$$
surjective. Apply Exercse \ref{chapter4:galoisgroups:exercise100} and
Proposition \ref{chapter4:automorphisms:result1300}. 
\qed
\end{proof}

\begin{exercise}
\label{chapter4:galoisgroups:exercise100}  
Suppose we have the following commutative diagram of groups and homomorphisms,
$$
\begin{pspicture}(0,0)(4,2)
\rput(1.2,1.8){$H$}\rput(1.2,0.2){$G$}
\rput(2.8,1.8){$H'$}\rput(2.8,0.2){$G'$}
\psline[linewidth=.1mm]{->}(1.2,1.5)(1.2,.5)
\psline[linewidth=.1mm]{->}(2.8,1.5)(2.8,.5)
\psline[linewidth=.1mm]{->}(1.5,1.8)(2.5,1.8)
\psline[linewidth=.1mm]{->}(1.5,.2)(2.5,.2)
\rput(2,2){$\theta$}\rput(2,.4){$\varphi$}\rput(1,1){$\psi$}\rput(3.05,1){$\psi'$}
\end{pspicture}
$$
with $\theta:H\rightarrow H'$ surjective and $\varphi:G\rightarrow G'$ an
isomorphism. Let $N:=N_G(\psi(H))$ be the normalizer in $G$ of the image
of $H$, and $N':=N_{G'}(\psi'(H'))$ similarly. Show that $\varphi$
induces an isomormphism
$$
N/\psi(H)\stackrel{\cong}{\rightarrow}N'/\psi'(H').
$$
\end{exercise}

\begin{example}
Figure \ref{chapter4:automorphisms:figure400} shows a covering
$Y\rightarrow X$ (left), a spanning tree $Z$ for $X$ consisting of the
single red edge, and the lifts $\bigcup Z_i$ of $Z$ to $Y$ (in
red). On the right we have the excised versions $X/Z$ and $Y/Z_i$. 
The Galois group $\gal(Y,X)\cong\Z$ is generated by the 
covering automorphism $a$
shown. The group is even more transparently $\cong\Z$ in the righthand version.
\begin{figure}[h]
  \centering
\begin{pspicture}(0,0)(12,4)
\rput(3.5,3){\BoxedEPSF{chapter4.fig4400.eps scaled 550}}
\rput(10,3.5){\BoxedEPSF{chapter4.fig4600.eps scaled 700}}
\rput(3.5,1){\BoxedEPSF{chapter4.fig4500.eps scaled 550}}
\rput(10,1){\BoxedEPSF{chapter3.fig2000.eps scaled 550}}
\psline[linewidth=.1mm]{->}(7,3)(8,3)
\psline[linewidth=.1mm]{->}(5,1)(9,1)
\psline[linewidth=.1mm]{->}(3.5,2.2)(3.5,1.6)
\psline[linewidth=.1mm]{->}(10,2.5)(10,1.5)
\rput(1.5,2.75){$Y$}\rput(2,1){$X$}\rput(11,1){$X/Z$}\rput(10,3.5){$Y/Z_i$}
\rput(3.5,2.75){${\blue a}$}
\end{pspicture}
  \caption{lattice excision with 
$\Z\cong\gal(Y,X)\stackrel{\cong}{\rightarrow}\gal(Y/Z_i,X/Z)\cong\Z$.}
\label{chapter4:automorphisms:figure400}  
\end{figure}
\end{example}

\begin{example}
Figure \ref{chapter4:automorphisms:figure410} 
revisits Example \ref{chapter3:basics:example410}:
a covering
$Y\rightarrow X$ (right), a spanning tree $Z$ for $X$
and the lifts $\bigcup Z_i$ of $Z$ to $Y$.
On the left we have the excised versions $X/Z$ and $Y/Z_i$. 
The two Galois groups are $\cong\Z/2$, generated by the $a$'s
shown.
\begin{figure}[h]
  \centering
\begin{pspicture}(0,0)(12,5)
\rput(0.4,0){
\rput(11,4){\BoxedEPSF{chapter3.fig6100b.eps scaled 750}}
\rput(8.9,4){\BoxedEPSF{chapter3.fig6100.eps scaled 750}}
\rput(6.8,4){\BoxedEPSF{chapter3.fig6100a.eps scaled 750}}
\rput(8.9,5){$Y$}\rput(6.8,4){$Y^{\tau_1}$}\rput(11,4){$Y^{\tau_2}$}
\psline[linewidth=.2mm,linecolor=blue]{<->}(8.35,3.45)(9.45,4.55)
\rput(8.9,4.2){${\blue a}$}
}
\rput(-.4,0){
\rput(4.8,4){\BoxedEPSF{chapter3.fig6400b.eps scaled 750}}
\rput(2.9,4){\BoxedEPSF{chapter3.fig6400.eps scaled 750}}
\rput(1,4){\BoxedEPSF{chapter3.fig6400a.eps scaled 750}}
\rput(2.9,5){$Y/Z_i$}\rput(1,3){$Y/Z_i^{q(\tau_1)}$}\rput(4.85,3){$Y/Z_i^{q(\tau_2)}$}
\psline[linewidth=.2mm,linecolor=blue]{<->}(2.5,4)(3.3,4)
\rput(2.9,4.2){${\blue a}$}
}
\rput(1.4,-3){
\rput(9,4){\BoxedEPSF{chapter3.fig6200.eps scaled 750}}
\rput(7,4){\BoxedEPSF{chapter3.fig6100c.eps scaled 750}}
\rput(9,4){$X$}\rput(7,4){$X^\ss$}
}
\rput(0.4,-3){
\rput(2.8,4){\BoxedEPSF{chapter3.fig6300.eps scaled 750}}
\rput(1.1,4){\BoxedEPSF{chapter3.fig6400c.eps scaled 750}}
\rput(2.8,4){$X/Z$}\rput(1.1,5){$X/Z^{q(\ss)}$}
}
\psline[linewidth=.1mm]{->}(7,1)(4.3,1)
\psline[linewidth=.1mm]{<-}(5.4,4)(6,4)
\psline[linewidth=.1mm]{->}(2.5,3)(2.5,2)
\psline[linewidth=.1mm]{->}(9.3,3)(9.3,2)
\end{pspicture}
  \caption{lattice excision with 
$\Z/2\cong\gal(Y,X)\stackrel{\cong}{\rightarrow}\gal(Y/Z_i,X/Z)\cong\Z/2$.}
\label{chapter4:automorphisms:figure410}  
\end{figure}
\end{example}

\section{Galois covers}\label{chapter4:galoiscovers}

To get some nice theorems about the Galois group of a cover we must 
restrict our attention to nice covers. We will see in this section that 
a nice cover is one that is highly symmetric.

\subsection{Galois covers}\label{chapter4:galoiscovers:covers}

\begin{proposition}\label{chapter4:galoiscovers:result100}
Let $f:Y\rightarrow X$ be a covering and $v$ a vertex of $X$. Then the
following are equivalent:
\begin{description}
\item[(i).] For any closed path $\gamma$ at
$v$, the lifts of $\gamma$ to each vertex of 
the fiber $f^{-1}(v)$
are either all closed or all non-closed.
\item[(ii).] The Galois group $\gal(Y,X)$ acts regularly on 
the fiber $f^{-1}(v)$.
\item[(iii).] if $u\in f^{-1}(v)$, then $f_*\pi_1(Y,u)$ is a normal subgroup
  of $\pi_1(X,v)$.
\end{description}
\end{proposition}

\begin{proof}
The equivalences follow immediately from Proposition
\ref{topological:galois:result200} and Lemma \ref{chapter4:automorphisms:result1200}.
\qed
\end{proof}

\begin{figure}
  \centering
\begin{pspicture}(0,0)(12,4)
\rput(2,3.5){\BoxedEPSF{chapter3.fig2100.eps scaled 650}}
\rput(2,1.5){\BoxedEPSF{chapter4.fig2000.eps scaled 1000}}
\rput(6.75,2){\BoxedEPSF{fig18a.eps scaled 500}}
\rput(11,2){\BoxedEPSF{chapter3.fig2500.eps scaled 1000}}
\rput(.5,3.5){$X=$}
\end{pspicture}
  \caption{some Galois coverings $Y\rightarrow X$ with Galois groups
    $\Z$, the free group of rank $2$ and $\Z\times\Z$.}  
\label{chapter4:galoiscovers:figure100}
\end{figure}
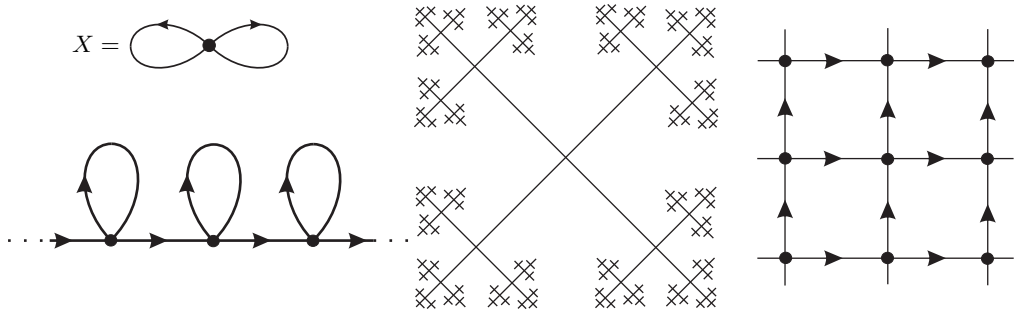

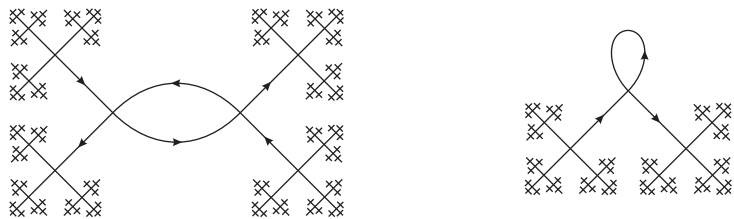
\begin{figure}
  \centering
\begin{pspicture}(0,0)(12,3)
\rput(9.5,1.5){\BoxedEPSF{fig18e.eps scaled 500}}
\rput(3.5,1.5){\BoxedEPSF{fig18g.eps scaled 500}}
\end{pspicture}
  \caption{some non-Galois coverings of the complex $X$ of Figure
    \ref{chapter4:galoiscovers:figure100}.}  
\label{chapter4:galoiscovers:figure200}
\end{figure}

Loosely, a covering is Galois when the
vertices of $Y$ are completely interchangeable with each other, or put
another way,
$Y$ looks the same from the viewpoint of any of its vertices. Figure
\ref{chapter4:galoiscovers:figure100} shows some examples and Figure
\ref{chapter4:galoiscovers:figure200} some non-examples.

\begin{definition}[Galois coverings]
\label{chapter4:galoiscovers:definition100}
A covering $Y\rightarrow X$ is \emph{Galois\/} (via $v$),
if it satisfies any (hence all) of the conditions of Proposition
\ref{chapter4:galoiscovers:result100}.
\end{definition}

\emph{Regular\/} is common alternative terminology to Galois, for reasons that
the second part of Proposition \ref{chapter4:galoiscovers:result100} makes clear.
As the effect of a covering automorphism on the fiber $f^{-1}(v)$ is 
completely determined by the image of a single vertex, a regular action is equivalent 
to a transitive action when talking about Galois groups.

\begin{lemma}\label{chapter4:galoiscovers:result200}
Let $f:Y\rightarrow X$ be a Galois covering via some vertex $v$ of
$X$. Then,
\begin{description}
\item[(i).] If $v'$ is another vertex of $X$, then the covering is Galois
  via $v'$.
\item[(ii).] If $Y\stackrel{g}{\rightarrow}Z\stackrel{h}{\rightarrow}X$ are
  intermediate coverings and $x\in h^{-1}(v)$, then $g:Y\rightarrow
  Z$ is Galois via $x$.
\item[(iii).] the induced covering $f:Y^{(1)}\rightarrow X^{(1)}$ on the
  $1$-skeletons is Galois, 
and there is an isomorphism,
$$
\gal(Y,X)\stackrel{\cong}{\rightarrow}\gal(Y^{(1)},X^{(1)}).
$$
\end{description}
\end{lemma}

In light of the first part of the Lemma, we will call a covering
$f:Y\rightarrow X$
\emph{Galois\/} without reference to a vertex in $X$.

\begin{proof}[of Lemma \ref{chapter4:galoiscovers:result200}]
For part (i), let $\lambda$ be a path in $X$ from $v'$ to $v$ and
$\gamma$ a closed path at $v'$. If $u_1',u_2'$ are vertices in the
fiber $f^{-1}(v')$ of $v'$, let $\nu_1,\nu_2$ be the lifts
of $\lambda$ to the $u_i'$. Suppose the $\nu_i$ finish at vertices
$u_i$, necessarily in the fiber $f^{-1}(v)$ of $v$. If $\mu_1,\mu_2$ are the
lifts of $\gamma$ to the $u_i'$, then the
$\nu_i^{-1}\mu_i\nu_i$ are the lifts of
$\lambda^{-1}\gamma\lambda$ to the $u_i$. The question of whether the
$\mu_i$ are closed or not becomes the question of whether the
$\nu_i^{-1}\mu_i\nu_i$ are closed or not. In
particular, if the covering is Galois via $v$ it is Galois via $v'$.

In the second part suppose that $\gamma$ is a closed path at $x$, so
that $h(\gamma)$ is a closed path at $v$, and let $u_1,u_2$ be vertices
of $Y$ in the fiber of $x$. Let $\mu_1,\mu_2$ be the lifts
(through $f$) of
$h(\gamma)$ to $u_1,u_2$. Then they are also the lifts of $\gamma$
(through $g$) to $u_1,u_2$, as $g(\mu_1),g(\mu_2)$ are closed
paths at $x$ covering $h(\gamma)$, hence must be $\gamma$. The lifts
$\mu_1,\mu_2$ are either both closed or both not closed as $f$
is Galois, and thus $g$ is Galois also.

For part (iii) we have an injective homomorphism
$\gal(Y,X)\rightarrow\gal(Y^{(1)},X^{(1)})$ from Lemma
\ref{chapter4:automorphisms:result200}, with the 
image of $\gal(Y,X)$
acting
transitively on the fiber of some vertex, so 
$\gal(Y^{(1)},X^{(1)})$ also
acting transitively, hence regularly. Thus the induced covering of graphs
is Galois. Suppose $a\in\gal(Y^{(1)},X^{(1)})$ is a graph covering
automorphism and that $a(u)=u'$ for some vertices $u,u'\in
Y^{(1)}$. Then there is a covering automorphism $\widehat{a}\in\gal(Y,X)$
that sends $u$ to $u'$ by the regularity of the action of $\gal(Y,X)$.
It must restrict on the $1$-skeleton to
$a$ and so the homomorphism is surjective. 
\qed
\end{proof}

\begin{proposition}\label{chapter4:galoiscovers:result300}
Let $f:Y\rightarrow X$ be a Galois covering and
$H\subset\gal(Y,X)$ a subgroup of index $[\gal(Y,X):H]$.
If
$Y\rightarrow Y/H\rightarrow X$ is the intermediate covering
of \S\ref{chapter4:galoisgroups:coverstosubgroups},
then,
$$
\deg(Y/H\rightarrow X)
=[\gal(Y,X):H].
$$
\end{proposition}

In particular, taking $H$ to be the trivial subgroup we get,
when $Y\rightarrow X$ is Galois,
that the
order $|\gal(Y,X)|$ of the Galois group is equal to the degree
$\deg(Y\rightarrow X)$ of the covering.

\begin{proof}
If a group $G$ acts regularly on a set and $H$ is a subgroup, 
then the number of $H$-orbits is the index $[G:H]$.
The result now follows as the $H$-orbits on the fiber (via $f$) of a
vertex $v\in X$
are precisely the vertices of $Y/H$ covering $v$ via
$Y/H\rightarrow X$.
\qed
\end{proof}

\begin{proposition}\label{chapter4:galoiscovers:result400}
If $f:Y\rightarrow X$ is Galois with $f(u)=v$ and $g\in\pi_1(X,v)$, then
the map $g\mapsto a_{g^{-1}}$ of Proposition
\ref{chapter4:automorphisms:result1300} induces an isomorphism
$$
\pi_1(X,v)/f_*\pi_1(Y,u)\stackrel{\cong}{\rightarrow}\gal(Y,X).
$$
If
$Y_u\stackrel{g}{\rightarrow}Z_x\stackrel{h}{\rightarrow}X_v$
is an intermediate covering, we get a subgroup
$h_*\pi_1(Z,x)\subset\pi_1(X,v)$, with 
$h^*\pi_1(Z,x)/f^*\pi_1(Y,u)$ mapping via this isomorphism to $\gal(Y,Z)$.
\end{proposition}

\begin{proof}
The isomorphism follows immediately from Proposition
\ref{chapter4:automorphisms:result1300}, as the normalizer
is $\pi_1(X,v)$ and $H=f_*\pi_1(Y,u)$. 
The image $h_*\pi_1(Z,x)$ consists of those $g_\gamma\in\pi_1(X,v)$
with $\gamma$
lifting to a closed path $\mu$
at $x$. Thus, it is those $g$ such that performing the process of Figure
\ref{chapter4:automorphisms:figure300},
we get for any $x\in Y$,
the cells $x$ and $a_g(x)$ lie in the same fiber of the
covering $g:Y\rightarrow Z$. But $\gal(Y,Z)$ consists precisely of
those $a\in\gal(Y,X)$ that permute the fibers of the covering
$q:Y\rightarrow Z$. 
\qed
\end{proof}

\begin{figure}
  \centering
\begin{pspicture}(0,0)(12,4)
\rput(4,3){\BoxedEPSF{fig118.eps scaled 300}}
\rput(4,.6){\BoxedEPSF{fig117.eps scaled 300}}
\rput(8,3){\BoxedEPSF{chapter4.fig1000.eps scaled 750}}
\rput(8,.6){\BoxedEPSF{chapter4.fig1100.eps scaled 750}}
\psline[linewidth=.1mm]{->}(4,1.9)(4,1.4)
\psline[linewidth=.1mm]{->}(8,1.9)(8,1.4)
\end{pspicture}  
  \caption{a Galois covering of $2$-complexes (left) and the induced
    Galois covering of graphs (right). The left one illustrates the
    inclusion of groups $1\hookrightarrow\Z/2$ and the right one the
    inclusion $2\Z\hookrightarrow\Z$. The Galois groups are both
  $\Z/2$.}
  \label{chapter4:galoiscovers:figure300}
\end{figure}
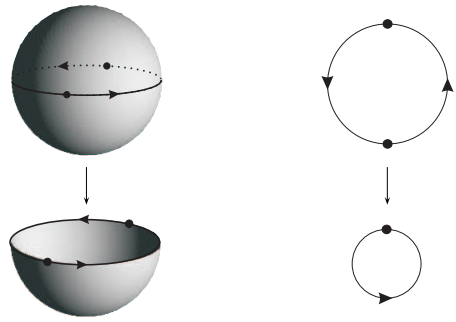

\begin{example}\label{chapter4:galoiscovers:example100}
On the left of Figure \ref{chapter4:galoiscovers:figure300}
we have the covering of \S\ref{chapter4:galoisgroups:automorphisms}. It is
easy to check that this covering is Galois, either by looking at the
lifts of closed paths in $X$, or by observing that $Y$ is simply
connected, hence $f_*\pi_1(Y,u)=1$, a normal subgroup of
$\pi_1(X,v)$. Lemma \ref{chapter4:galoiscovers:result200} tells
us that the Galois group of this covering is the same as for the
induced covering of the $1$-skeletons, shown on the right. Indeed,
these coverings have degree $2$, and so the two Galois groups are
isomorphic to $\Z/2$ by Proposition \ref{chapter4:galoiscovers:result300}.

The covering of $2$-complexes gives us no more than the covering of
the underlying graphs, at least when it comes to their Galois groups. They do
illustrate different things though as Proposition
\ref{chapter4:galoiscovers:result400} shows. The $2$-complex covering
depicts the inclusion of groups $1\hookrightarrow\Z/2$, and the graph
covering the inclusion $2\Z\hookrightarrow\Z$.  
\end{example}

\subsection{Intermediate covers}
\label{chapter4:galoiscovers:intemediatecovers}

\begin{figure}
  \centering
\begin{pspicture}(0,0)(12,4)
\rput(6,2){\BoxedEPSF{chapter4.fig2100.eps scaled 1000}}
\rput(6.2,.55){$v$}\rput(4.9,1.9){$x_1$}\rput(7,2.1){$x_2$}
\rput(4.1,3.3){$u_1$}\rput(7.65,3.7){$u_2$}\rput(5.8,3.6){$a$}
\rput(5.8,2.2){$a'$}
\rput(2.2,2.7){$Y$}\rput(2.2,1.45){$Z$}\rput(2.2,.2){$X$}
\psline[linewidth=.1mm]{->}(2.2,2.4)(2.2,1.7)
\psline[linewidth=.1mm]{->}(2.2,1.2)(2.2,.5)
\end{pspicture}  
  \caption{defining a map $\theta:\gal(Y,X)\rightarrow\gal(Z,X)$ when
    $Z\rightarrow X$ is Galois.}
  \label{chapter4:galoiscovers:figure400}
\end{figure}
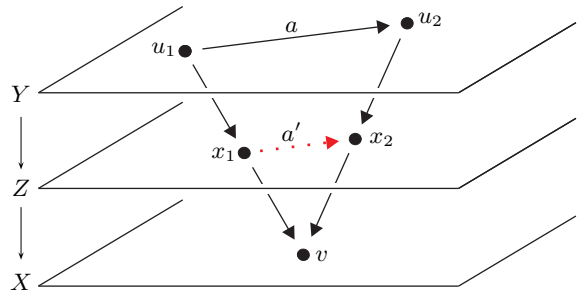

Let $Y\rightarrow X$ be a covering and $Y\rightarrow
Z\rightarrow X$ a covering intermediate to it. When $Z\rightarrow X$
is Galois we can define a map $\theta:\gal(Y,X)\rightarrow\gal(Z,X)$
using the scheme in Figure \ref{chapter4:galoiscovers:figure400}. Let $a$ be
a covering automorphism in $\gal(Y,X)$, and $u_1,u_2$ vertices of $Y$
with $a(u_1)=u_2$. In particular, the $u_i$ lie in the fiber
$f^{-1}(v)$ of the covering $Y\rightarrow X$. They also cover vertices
$x_1,x_2$ of $Z$ via the covering $Y\rightarrow Z$, with these two
lying in the fiber of $v$ via the covering $Z\rightarrow X$. As
$Z\rightarrow X$ is Galois, there is a covering automorphism
$a'\in\gal(Z,X)$ sending $x_1$ to $x_2$. 

\begin{proposition}\label{chapter4:galoiscovers:result500}
Let $f:Y\rightarrow X$ be a covering with $Y\rightarrow
Z\rightarrow X$ intermediate to it. 
\begin{description}
\item[(i).] If $Z\rightarrow X$ is Galois, then the map 
  $\theta:\gal(Y,X)\rightarrow\gal(Z,X)$defined by $\theta(a)=a'$ is a homomorphism with kernel
  $\gal(Y,Z)$.
\item[(ii).] If $Y\rightarrow X$ is Galois and $\gal(Y,Z)$ is a normal
  subgroup of $\gal(Y,X)$, then $Z\rightarrow X$ is Galois,
  and $\theta$ is surjective. 
\end{description}
In particular, if $Y\rightarrow X$ is Galois then $Z\rightarrow X$ is
Galois if and only if $\gal(Y,Z)$ is a normal subgroup of $\gal(Y,X)$,
in which case we get,
$$
\gal(Y,X)/\gal(Y,Z)\cong\gal(Z,X).
$$
\end{proposition}

\begin{proof}
The $Z\rightarrow X$ Galois condition is used in the definition of
$\theta$ and has done its job. That $\theta$ is well defined is illustrated in Figure
\ref{chapter4:galoiscovers:figure500}: suppose we choose different vertices
$u_1',u_2'$ in $Y$ covering $x_1',x_2'\in Z$
and with $\aa(u_1')=u_2'$. 
If $\gamma$ is a path from $u_1$ to $u_1'$, then lifting $f(\gamma)\in X$
to $u_2$ gives a path to $u_2'$ (as $a(u_1')=u_2'$).
The middle paths in $Z$
are the images of these two paths in $Y$ 
under the covering $Y\rightarrow Z$, and these must be the lifts
to the $x_i'$ of $f(\gamma)$. Cover and lift gives the
same covering automorphism in $\gal(Z,X)$ sending $x_1$ to $x_2$ and
$x_1'$ to $x_2'$.

\begin{figure}
  \centering
\begin{pspicture}(0,0)(12,5.25)
\rput(3,2.5){\BoxedEPSF{chapter4.fig2200.eps scaled 750}}
\rput(2.5,4.4){$a$}\rput(3.5,5.05){$a$}
\rput(.6,4.25){$u_1$}\rput(4.35,4.25){$u_2$}\rput(1.7,5){$u_1'$}\rput(5.4,5){$u_2'$}
\rput(1.4,2.15){$x_1$}\rput(3.55,2.15){$x_2$}\rput(2.45,2.85){$x_1'$}\rput(4.65,2.8){$x_2'$}
\rput(9,2.5){\BoxedEPSF{chapter4.fig2300.eps scaled 750}}
\rput(6,3.5){$u_1$}\rput(7.2,4.5){$u_2$}
\rput(9.7,3.5){$u_1'$}\rput(10.9,4.5){$u_2'$}
\rput(6.65,1.9){$x_1$}\rput(9.8,1.9){$x_2$}
\rput(10.4,.4){$\gamma$}\rput(6.2,4.7){$\mu$}
\rput(6.7,4.1){$a$}\rput(8,3.7){$b$}\rput(9,4.6){$b$}
\end{pspicture}  
  \caption{$\theta$ is well defined (left) and $Y\rightarrow X$ Galois
    implies $Z\rightarrow X$ is Galois (right).}
  \label{chapter4:galoiscovers:figure500}
\end{figure}
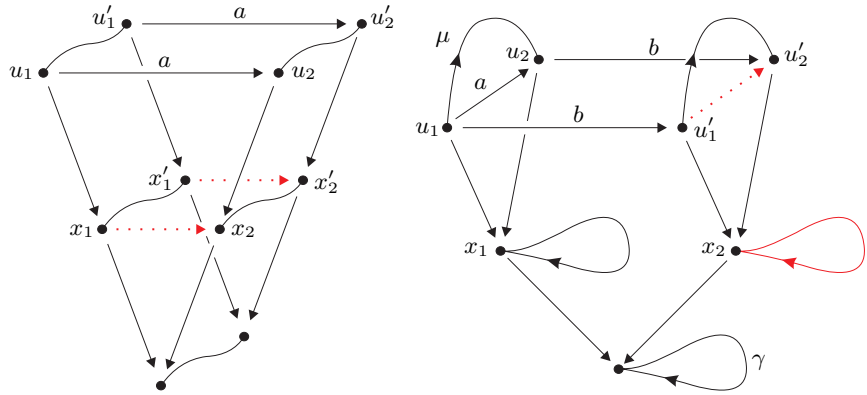

Splice together two copies of Figure
\ref{chapter4:galoiscovers:figure400} to see that $\theta$ is
a homomorphism and $a$ is in the kernel if and only if $x$ and
$a(x)$ cover the same vertex of $Z$ for all $x$; such $a$ are
precisely the subgroup $\gal(Y,Z)$ as in
\S\ref{chapter4:galoisgroups:coverstosubgroups}.

The second part is illustrated on the right of Figure
\ref{chapter4:galoiscovers:figure500}. Start at the bottom of the
picture and move
clockwise around it. Suppose the closed path $\gamma$ in $X$
lifts to a closed path at $x_1$ in $Z$, and then lift this to a
path $\mu\in Y$ starting at $u_1$ and
finishing at $u_2$. Let $x_2$ be some other
vertex of $Z$ covering $v$ and $u_1'$ a vertex of $Y$ covering
$x_2$. As $Y\rightarrow X$ is Galois there is a covering automorphism
$b\in\gal(Y,X)$ with $b(u_1)=u_1'$, and as $Y\rightarrow Z$ is
Galois (Lemma \ref{chapter4:galoiscovers:result200}(ii)), 
there is an $a\in\gal(Y,Z)$ with $a(u_1)=u_2$. If the
image path $b(\mu)$ ends at the vertex $u_2'$, then the
conjugate $bab^{-1}$ maps $u_1'$ to $u_2'$. By the normality
condition, this automorphism is also in $\gal(Y,Z)$, and
so $u_2'$ must cover $x_2$.
Finally, the image under
the covering $Y\rightarrow Z$ of $b(\mu)$, a closed path at
$x_2$, is the lift of $\gamma$ to $x_2$. Repeating the arguement
with the roles of $x_1$ and $x_2$ interchanged, we get the lift of
$\gamma$ to $x_1$ is closed if and only if the lift to $x_2$ is
closed. Thus $Z\rightarrow X$ is Galois.

If $a'\in\gal(Z,X)$ sends $x_1$ to $x_2$, and $u_1,u_1'$ are
vertices of $Y$ covering them, then there is an $a\in\gal(Y,X)$ mapping
$u_1$ to $u_1'$, as $Y\rightarrow X$ is Galois. Thus,
$\theta$ is surjective.
\qed
\end{proof}

\subsection{Universal covers and excision}
\label{chapter4:galoiscovers:universalcovers}

In Chapter \ref{chapter3} we had said quite a bit about coverings
before we gave our first non-trivial example of a covering of an arbitrary complex
$X$. This example was the universal cover $\wtl{X}\rightarrow X$ of
\S\ref{chapter3:actionsuniversal:universal}, and it will now give
our first non-trivial example of a Galois covering of an arbitrary
complex $X$:

\begin{proposition}\label{chapter4:galoiscovers:result600}
If $X$ is a $2$-complex, then the universal cover $\wtl{X}_u\rightarrow
X_v$ is a Galois covering with
$$
\gal(\wtl{X},X)\cong\pi_1(X,v).
$$
Moreover, if $\wtl{X}_u\rightarrow Z_x\stackrel{g}{\rightarrow}X_v$ is
intermediate, then the subgroup $g_*\pi_1(Z,x)\subset\pi_1(X,v)$ maps
via this homomorphism to $\gal(\wtl{X},Z)$.
\end{proposition}

\begin{proof}
By Proposition \ref{chapter3:actionsuniversal:result700} the
universal cover $\wtl{X}$ is simply connected and so 
$f_*\pi_1(\wtl{X},u)$ is the trivial subgroup, hence
is normal, so by Proposition \ref{chapter4:galoiscovers:result100} the covering is
Galois. The rest is then an immediate application of Proposition
\ref{chapter4:galoiscovers:result400}. 
\qed
\end{proof}

\begin{exercise}\label{chapter4:galoiscovers:exercise100}
Show that if $G$ is a group, then there is a covering $Y\rightarrow X$
of $2$-complexes with $\gal(Y,X)\cong G$.
\end{exercise}


Returning now to the excision of \S\ref{chapter4:galoisgroups:excision}, we show that this process sends
Galois covers to Galois covers. 
Recall that $f:Y_u\rightarrow X_v$ is a covering and $Z\subset X$
connected and simply-connected with $v\in Z$. Also $f^{-1}(Z)=\bigcup_I Z_i\subset Y$ with the
$Z_i$ connected and simply connected and
$f':Y/Z_i\rightarrow X/Z$ the induced covering of \S\ref{chapter3:basics:lifting}.

\begin{proposition}\label{chapter4:galoiscovers:result700}
The induced covering $f':Y/Z_i\rightarrow X/Z$ is Galois if and only
if $f:X\rightarrow Z$ is Galois.
\end{proposition}

\begin{proof}
In the proof of Proposition \ref{chapter4:automorphisms:result1400} we
``hit'' the commuting diagram of $2$-complexes with
$\pi_1$, although we didn't show the result. Here it is:
$$
\begin{pspicture}(0,0)(4,2)
\rput(.5,1.8){$\pi_1(Y,u)$}\rput(.5,0.2){$\pi_1(X,v)$}
\rput(3.5,1.8){$\pi_1(Y/Z_i,q'(u))$}\rput(3.5,0.2){$\pi_1(X/Z,q(v))$}
\psline[linewidth=.1mm]{->}(.5,1.5)(.5,.5)
\psline[linewidth=.1mm]{->}(3.5,1.5)(3.5,.5)
\psline[linewidth=.1mm]{->}(1.2,1.8)(2.3,1.8)
\psline[linewidth=.1mm]{->}(1.2,.2)(2.3,.2)
\rput(1.9,2.05){$q'_*$}\rput(1.9,.4){$q_*$}\rput(.3,1){$f_*$}\rput(3.75,1){$f'_*$}
\end{pspicture}
$$
with the top homomorphism surjective and the bottom an
isomorphism. This gives the image $f'_*\pi_1(Y/Z_i,q'(u))\subset\pi_1(X/Z,q(v))$
equaling the image $q_*f_*\pi_1(Y,u)$. In particular,
$\pi_1(Y,u)\subset\pi_1(X,v)$ is normal if and only if
$\pi_1(Y/Z_i,q'(u))\subset\pi_1(X/Z,q(v))$ is normal.
\qed
\end{proof}

\section{Galois correspondences}\label{chapter4:correspondences}

In its purest form, the Galois correspondence says that one lattice is
the same as another turned upside down. The two
lattices concerned are the lattice of covers intermediate to a fixed
cover $Y\rightarrow X$, and the lattice of subgroups of the Galois group of
$Y\rightarrow X$.

\subsection{The Galois
  correspondence}\label{chapter4:correspondences:correspondence}

We return to the set-up of the last section of Chapter \ref{chapter3}, where
$f:Y_u\rightarrow X_v$ is a fixed pointed covering of 
$2$-complexes, connected as
always, with $\mathcal{L}=\mathcal{L}(Y_u,X_v)$ the lattice of
equivalence classes of pointed connected intermediate covers and
$\gal(Y_u,X_v)$ the Galois group. Let $\SSS=\SSS(Y_u,X_v)$, be the lattice
of subgroups of $\gal(Y_u,X_v)$ as in Definition
\ref{chapter3:lattices:definition300}. 

\begin{theorem}[Galois correspondence]
\label{galois:correspondence:result300}
If $Y_u{\rightarrow}X_v$ is a Galois covering of $2$-complexes, then the map 
$$
\Phi:\mathcal{L}\rightarrow\SSS,
$$
that associates to the equivalence class of 
$Y_u{\rightarrow}Z_x{\rightarrow}X_v$ the subgroup 
$\gal(Y_u,Z_x)$ is a lattice 
anti-isomorphism.
Its inverse $\Psi$ is the map that associates to the subgroup
$H\subset\gal(Y_u,X_v)$ the equivalence class of 
$Y_u{\rightarrow}Y/H_{q(u)}
{\rightarrow}X_v$.
In particular, 
\begin{description}
\item[(i).] $\gal(Y_u,Y/H_{q(u)})=H$ and $Y_u\rightarrow Z_x\rightarrow X_v$ is
  equivalent to $Y_u\rightarrow Y/\gal(Y,Z)_{q(u)}\rightarrow X_v$;
\item[(ii).] equivalence classes of covers
$Y_u{\rightarrow}Z_x{\rightarrow}X_v$ with $Z_x\rightarrow X_v$ Galois
correspond to normal subgroups of the Galois group $\gal(Y_u,X_v)$;
\item[(iii).] given the intermediate covering
  $Y_u{\rightarrow}Z_x{\rightarrow}X_v$, we have,
$$
[\gal(Y_u,X_v):\gal(Y_u,Z_x)]=\deg(Z_x\rightarrow X_v);
$$
\item[(iv).] given the subgroup $H\subset\gal(Y_u,X_v)$, we have
$$
\deg(Y/H_{q(u)}{\rightarrow}X_v)=[\gal(Y_u,X_v),H].
$$
\end{description}
\end{theorem}

\begin{proof}
We show that both $\Phi$ and $\Psi$ are order-reversing
bijections, recalling the order on $\mathcal{L}$
from \S\ref{chapter3:lattices:poset}, and that the order on
$\SSS$ is inclusion of subgroups.
Suppose we have two equivalence classes of intermediate coverings with
representatives $Z_x,Z_y'\in\mathcal{L}$, and $Z_x\leq
Z_y'$.
By Lemma \ref{chapter4:automorphisms:result1000}  we get
$\gal(Y_u,Z_y')\subset\gal(Y_u,Z_x)$, ie: $\Phi(Z_y')\leq\Phi(Z_x)$.
On the otherhand,
if $H_1\subset H_2$ in $\SSS$,
then Lemma \ref{chapter4:automorphisms:result1100}
gives $Y_u\rightarrow Y/H_1\rightarrow Y/H_2\rightarrow X_v$, ie: 
$\Psi(H_2)\leq \Psi(H_1)$. Thus $\Phi$ and $\Psi$ are lattice
anti-morphisms.

We now show that the composition $\Psi\Phi$ is the identity map on the
lattice $\mathcal{L}$.
Let $Y_u\rightarrow Z_x\rightarrow X_v$ be intermediate
with $Y_u\rightarrow Z_x$ Galois by Lemma
\ref{chapter4:galoiscovers:result200}.
We also have the intermediate covering
$Y_u\rightarrow Y_u/\gal(Y_u,Z_x)\rightarrow Z_x$
with 
$Y_u/\gal(Y_u,Z_x)\rightarrow Z_x$ degree $1$ by Proposition
\ref{chapter4:galoiscovers:result300}, hence an isomorphism by
Corollary \ref{chapter3:basics:result800}. 

\begin{figure}
  \centering
\begin{pspicture}(0,0)(4,2)
\rput(-1,0.2){
\rput(3,1.5){$Y/\gal(Y_u,Z_x)$}\rput(3,0){$Z_{x}$}
\psline[linewidth=.1mm]{->}(3,1.2)(3,.4)
\rput(1.3,.75){$Y_u$}\rput(4.7,.75){$X_v$}
\psline[linewidth=.1mm]{->}(1.6,.5)(2.7,.14)
\rput(-.5,0){\psline[linewidth=.1mm]{->}(2.1,.9)(2.7,1.26)}
\psline[linewidth=.1mm]{->}(3.3,.14)(4.4,.5)
\rput(.5,0){\psline[linewidth=.1mm]{->}(3.3,1.26)(3.9,.9)}
}
\end{pspicture}  
  \caption{The composition $\Psi\Phi=\id_\mathcal{L}$.}
 \label{chapter4:correspondences:figure100}
\end{figure}
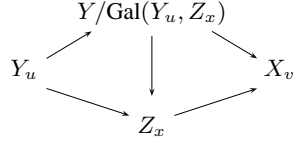

What results is the diagram
of Figure \ref{chapter4:correspondences:figure100},
with the whole square commuting as $Z$ and $Y/\gal(Y,Z)$ are intermediate
to $Y\rightarrow X$,
and the left triangle commuting
as $Y/\gal(Y,Z)$ is intermediate to $Y\rightarrow Z$.
The
right triangle thus commutes as well,
and the map across the middle is an isomorphism. 
Thus, the intermediate coverings
$Y_u\rightarrow Z_x\rightarrow X_v$ and 
$Y_u\rightarrow Y_u/\gal(Y_u,Z_x)\rightarrow X_v$
are equivalent, and so we have $\Psi\Phi=\id_{\mathcal{L}}$.

Now for $\Phi\Psi$: if $H\subset\gal(Y_u,X_v)$ with $q:Y\rightarrow Y/H$ the
quotient map,
then $qa=q$ for any $a\in H$, 
so $H\subset\gal(Y_u,(Y/H)_{q(u)})$ with the covering
$Y_u\rightarrow Y/H_{q(u)}$ intermediate, hence Galois. Proposition 
\ref{chapter4:galoiscovers:result300} gives the index of $H$ in 
$\gal(Y_u,(Y/H)_{q(u)})$ as the degree of the covering
$Y/H_{q(u)}\rightarrow Y/H_{q(u)}$, ie: $H=
\gal(Y_u,(Y/H)_{q(u)})$,
and we have $\Phi\Psi=\id_{\SSS}$. The anti-morphisms $\Phi$ and $\Psi$
are thus anti-isomorphisms.

The claims of part (i) are direct applications of $\Phi\Psi=\id_\mathcal{L}$
and $\Psi\Phi=\id_\SSS$. 
The correspondence in (ii) between intermediate Galois coverings and normal
subgroups follows from Proposition
\ref{chapter4:galoiscovers:result500}. Part (iii) is just part (i) applied to 
Proposition \ref{chapter4:galoiscovers:result300} and part (iv) \emph{is\/} just 
Proposition \ref{chapter4:galoiscovers:result300}.
\qed
\end{proof}

\begin{example}\label{chapter4:correspondence:example100}
Figure \ref{chapter4:correspondences:figure150} shows a covering
$Y_u\rightarrow X_v$ with the Galois group $\gal(Y_u,X_v)$ of order
four, generated by the $1/4$-turn rotation $a$ shown. Thus
$\gal(Y_u,X_v)\cong\Z_4$. The fiber of the vertex $v$ consists of
the the four vertices at the center of $Y$. In particular the
Galois group acts transitively (hence regularly) on this fiber and
the covering is Galois. 

\begin{figure}[h]
  \centering
\begin{pspicture}(0,0)(12,5)
\rput(3.5,2.5){\BoxedEPSF{chapter4.fig4200.eps scaled 500}}
\rput(9.5,2.5){\BoxedEPSF{chapter4.fig4000.eps scaled 500}}
\psline[linewidth=.1mm]{->}(6.5,2.5)(7.5,2.5)
\rput(3.5,2.5){${\blue a}$}\rput(3.5,3.65){$u$}\rput(9.5,2.55){$v$}
\end{pspicture}
  \caption{Galois covering $Y\rightarrow X$ of degree four and covering 
automorphism $a$.}
\label{chapter4:correspondences:figure150}
\end{figure}
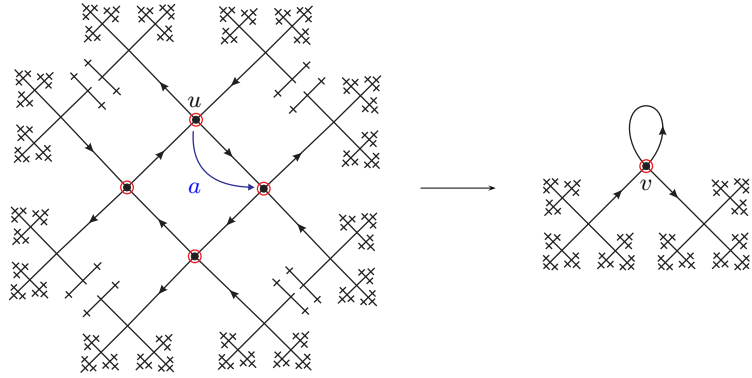

The subgroup lattice is very simple, just $1\subset\Z_2\subset\Z_4$,
and so by the Galois correspondence there is a single equivalence
class of intermediate coverings $Y_u\rightarrow Z_x\rightarrow X_v$
with the degree of the covering $Z\rightarrow X$ the index
$[\Z_4:\Z_2]=2$, and $Z$ the quotient of $Y$ by the automorphism
$a^2$ (see Figure \ref{chapter4:correspondences:figure175}).

\begin{figure}[h]
  \centering
\begin{pspicture}(0,0)(12,4)
\rput(6,2){\BoxedEPSF{chapter4.fig4100.eps scaled 500}}
\rput(2.25,2){$Z=Y/\langle a^2\rangle=$}
\end{pspicture}
  \caption{quotient $Y/\Z_2$ corresponding to the single equivalence class
of intermediate coverings.}
\label{chapter4:correspondences:figure175}
\end{figure}
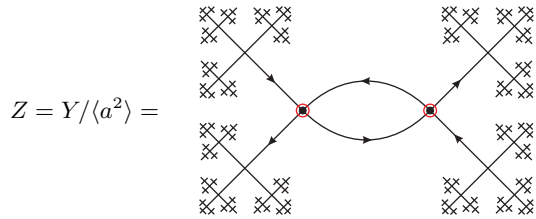
\end{example}

Returning to generalities, when two lattices are anti-isomorphic, 
then the join of two elements in one
corresponds to the meet in the other, and vice-versa. 
Following this through with $\Phi(-)=\gal(Y_u,-)$ we
get the following,

\begin{corollary}\label{galois:correspondence:result400}
Let $Y_u{\rightarrow}X_v$ be Galois
with $Y_u\rightarrow Z_{x}\rightarrow X_v$ and
$Y_u\rightarrow Z'_{y}\rightarrow X_v$ representatives of elements in the lattice
$\mathcal{L}(Y_u,X_v)$.
Then,
\begin{align*}
\gal(Y_u,(Z\prod_XZ')_z)
&=\gal(Y_u,Z_{x})\cap\gal(Y_u,Z'_{y}),\\
\gal(Y_u,(Z\coprod_YZ')_z)&=
\langle\gal(Y_u,Z_{x}),\gal(Y_u,Z'_{y})\rangle,
\end{align*}
where in the first case $z=x\times y$ in the formation of the pullback $(Z\prod_XZ')_z$,
and in the second case $z=[x]=[y]$ in the formation of the pushout $(Z\coprod_YZ')_z$.
\end{corollary}

Similarly the anti-isomorphism $\Psi(-)=(Y/-)_{q(u)}$ sends joins to
meets and meets to joins, giving

\begin{corollary}\label{galois:correspondence:result500}
Let $Y_u{\rightarrow}X_v$ be Galois and $H_1,H_2\subset \gal(Y_u,X_v)$.
Then the intermediate coverings,
\begin{align*}
Y_u\rightarrow Y/\langle H_1,H_2\rangle_{q(u)}\rightarrow X_v
\text{ and }
Y_u\rightarrow(Y/H_1\coprod_Y Y/H_2)_{q(u)}
\rightarrow X_v,\\
Y_u\rightarrow Y/(H_1\cap H_2)_{q(u)}\rightarrow X_v\text{ and }
Y_u\rightarrow(Y/H_1\prod_XY/H_2)_w\rightarrow X_v.
\end{align*}
are equivalent, where the $q$ are the appropriate 
quotient maps and $w=q_1(u)\times q_2(u)$ with $q_i:Y\rightarrow Y/H_i\,(i=1,2)$ the
the quotient map.
\end{corollary}


Here is another simple application of the Galois correspondence.
In Exercise \ref{chapter4:galoiscovers:exercise100} we showed that any
group could be realized as the 
Galois group of some covering $Y\rightarrow X$ of $2$-complexes.
If the reader will allow us the liberty of ``borrowing'' a theorem
from a later chapter, then we can realize the group as the Galois
group of a covering of \emph{graphs}:

\begin{corollary}[Inverse Galois theorem]\label{galois:correspondence:result700}
Let $G$ be a group. Then there is a Galois covering of graphs $Y\rightarrow X$
with $\gal(Y,X)\cong G$.
\end{corollary}

\begin{proof}
The theorem to be borrowed is that there is a free group $F$ and a
normal subgroup $H\lhd F$ with $G\cong F/H$. As the fundamental groups
of graphs are free, choose a graph $X$ with $\pi_1(X,v)\cong F$.
Consider the universal
cover $f:\wtl{X}_u\rightarrow X_v$ with  
$$
\gal(\wtl{X}_u,X_v)\cong\pi_1(X,v),
$$
by
Proposition \ref{chapter4:galoiscovers:result600}.
Identifying $H$ with
the corresponding subgroup of $\gal(\wtl{X},X)$ under the isomorphism
above, we get an intermediate covering 
$$
\wtl{X}\rightarrow\wtl{X}/H\rightarrow X,
$$
with $\gal(\wtl{X},\wtl{X}/H)=H$.  In particular, Proposition
\ref{chapter4:galoiscovers:result500} gives $\wtl{X}/H\rightarrow
X$ is a Galois covering with Galois group isomorphic to $G$.
\qed
\end{proof}

\begin{exercise}[completely irregular covers $\mathbf{=}$ malnormal
  subgroups]
\label{chapter4:correspondences:exercise200}
A subgroup $H\subset G$ is {\em malnormal\/} whenever
$g\in G\setminus H$ implies $gHg^{-1}\cap H=\{1\}$. 
Let $f:Y\rightarrow X$ be a covering with $f(u)=v$, and let $\gamma$ be a homotopically
non-trivial closed path at $v$ whose lift to $u$ is also closed. Call $f$ 
\emph{completely irregular\/} if the lift of $\gamma$ to every other vertex of $f^{-1}(v)$
is non-closed. 

Show that equivalence classes of covers $Y_u\rightarrow Z_x\rightarrow X_v$ with 
$Z_x\rightarrow X_v$ completely irregular correspond to malnormal subgroups
of the Galois group $\gal(Y_u,X_v)$.
\end{exercise}

\subsection{Galois correspondence for the universal cover}
\label{chapter4:correspondences:universal}

There is a more familiar version of the Galois correspondence using the
fact that the universal cover $\wtl{X}\rightarrow X$ is Galois, and
translating the various ingredients of Theorem
\ref{galois:correspondence:result300} into this
special setting. The aim is to remove all mention of the
top cover, and just focus on coverings of $X$.

Let $X$ be a $2$-complex and $v\in X$ a vertex. Two coverings 
$f_i:Y_i\rightarrow X\,(i=1,2)$ with $f_i(y_i)=v_i$ are \emph{equivalent\/} 
when there is an isomorphism $(Y_1)_{y_1}\rightarrow (Y_2)_{y_2}$ making
$$
\begin{pspicture}(0,0)(4,1.5)
\rput(0,-.25){
\rput(1,1.5){$Y_1$}
\rput(3,1.5){$Y_2$}
\rput(2,0.4){$X$}
\rput(2,1.7){$\cong$}
\rput(1.125,.825){$f_1$}
\rput(2.78,.825){$f_2$}
\psline[linewidth=.1mm]{->}(1.3,1.5)(2.7,1.5)
\psline[linewidth=.1mm]{->}(1.2,1.25)(1.8,.55)
\psline[linewidth=.1mm]{->}(2.8,1.25)(2.25,.55)
}
\end{pspicture}
$$
commute. Let $\mathcal{L}=\mathcal{L}(X_v)$ be the set of equivalence classes of 
coverings.

\begin{exercise}
\label{chapter4:correspondences:exercise100}
Show that two coverings are equivalent in the sense just defined if and only if
the intermediate coverings
$$
\wtl{X}_u\rightarrow(Y_i)_{y_i}\stackrel{f_i}{\rightarrow}X_v
$$
are equivalent, and so there is an isomorphism of lattices 
$\mathcal{L}(X_v)\rightarrow\mathcal{L}(\wtl{X}_u,X_v)$.
\end{exercise}

Identifying $\mathcal{L}(\wtl{X}_u,X_v)$ with $\mathcal{L}(X_v)$,
the Galois correspondence sends (the equivalence class of) the
cover $f:Y_y\rightarrow X_v$ to the subgroup $\gal(\wtl{X}_u,Y_y)\subset\gal(\wtl{X}_u,X_v)$.
Proposition \ref{chapter4:galoiscovers:result600} in turn gives an
isomorphism
$$
\varphi:\gal(\wtl{X}_u,X_v)\rightarrow\pi_1(X,v),
$$
which sends the subgroup $\gal(\wtl{X}_u,Y_y)$ to
$f_*\pi_1(Y,y)\subset\pi_1(X,v)$. 

Let $\SSS=\SSS(X_v)$ be the lattice of
subgroups of the fundamental group $\pi_1(X,v)$ ordered, as usual, by
inclusion. A translation of 
Theorem \ref{galois:correspondence:result300} then gives:

\begin{corollary}[Galois correspondence for $\mathbf{\wtl{X}}$: first go]
\label{galois:correspondence:result800}
If $X$ is a $2$-complex, then the map 
$$
\Phi:\mathcal{L}\rightarrow\SSS,
$$
that associates to the equivalence class of the covering 
$f:Y_y \rightarrow X_v$ the subgroup 
$f_*\pi_1(Y,y)\subset\pi_1(X,v)$, is a lattice 
anti-isomorphism.
Its inverse $\Psi$ is the map that associates to the subgroup
$H\subset\pi_1(X,v)$ the equivalence class of
$\wtl{X}/\varphi^{-1}(H)_{q(u)}{\rightarrow}X_v$ where 
$q:\wtl{X}\rightarrow\wtl{X}/\varphi^{-1}(H)$ is the quotient map.

Moreover, equivalence classes of Galois covers
$Y_y\rightarrow X_v$ 
correspond to normal subgroups of $\pi_1(X,v)$.
\end{corollary}

The mantra is thus, ``covers of $X$ correspond to subgroups of the
fundamental group of $X$''.

Corollary \ref{galois:correspondence:result800} is not quite perfect:
our aim was to remove all reference to the top covering ($\wtl{X}$
in this case) and focus entirely on $X$, but our top-down approach
to associating a cover to a subgroup has hard
wired $\wtl{X}$ into the picture. The solution is to return 
to the bottom-up construction of
\S\ref{chapter3:actionsuniversal:bottom_up}.

\begin{proposition}[bottom-up$\mathbf{=}$top-down]
\label{galois:correspondence:result900}
Let $H\subset\pi_1(X,v)$ be a subgroup and
$(X\kern-2pt\uparrow\kern-2pt H)_x\rightarrow X_v$ the bottom-up cover of $X$
constructed in \S\ref{chapter3:actionsuniversal:bottom_up}.
Then the coverings
$$
(X\kern-2pt\uparrow\kern-2pt H)_x\stackrel{}{\rightarrow}X_v\text{ and }
\wtl{X}/\varphi^{-1}(H)_{q(u)}{\rightarrow}X_v
$$
are equivalent.
\end{proposition}

\begin{proof}
By Corollary
\ref{chapter3:actionsuniversal:exercise500} the cover 
$(X\kern-2pt\uparrow\kern-2pt H)_x\rightarrow X_v$ corresponds to $H\subset\pi_1(X,v)$,
hence is equivalent to $\wtl{X}/\varphi^{-1}(H)_{q(u)}{\rightarrow}X_v$
by Corollary \ref{galois:correspondence:result800}.
\qed
\end{proof}

\begin{corollary}[Galois correspondence for $\mathbf{\wtl{X}}$: second go]
\label{galois:correspondence:result1000}
If $X$ is a $2$-complex, then the map 
$$
\Phi:\mathcal{L}\rightarrow\SSS,
$$
that associates to the equivalence class of the covering 
$f:Y_y \rightarrow X_v$ the subgroup 
$f_*\pi_1(Y,y)\subset\pi_1(X,v)$, is a lattice 
anti-isomorphism.
Its inverse $\Psi$ is the map that associates to the subgroup
$H\subset\pi_1(X,v)$ the equivalence class of
$(X\kern-2pt\uparrow\kern-2pt H)_x\rightarrow X_v$.

Moreover, equivalence classes of Galois covers
$Y_y\rightarrow X_v$ 
correspond to normal subgroups of $\pi_1(X,v)$.
\end{corollary}

\subsection{Lattice excision}
\label{chapter4:correspondences:lattice}

We have seen (\S\ref{chapter4:galoisgroups:excision}) that the
excision of simply-connected subcomplexes has no 
effect on Galois groups. Another consequence of the Galois
correspondence is that these excisions have no effect on lattices
of intermediate coverings.

Recalling the set-up, let $f:Y_u\rightarrow X_v$ be a covering and $Z\subset X$
connected and simply-connected with $v\in Z$. Let $f^{-1}(Z)=\bigcup_I Z_i\subset Y$ 
a disjoint union with the
$Z_i$ connected and simply connected. Finally, $q:X\rightarrow X/Z$,
$q':Y\rightarrow Y/Z_i$ are the quotient maps and 
$f':Y/Z_i\rightarrow X/Z$ the induced covering.

\begin{corollary}[lattice excision]
\label{galois:correspondence:result800}
There is an 
isomorphism of lattices 
$$
\LL(Y_u,X_v)\stackrel{\cong}{\rightarrow}
\LL((Y/Z_i)_{q'(u)},X/Z_{q(v)}),
$$
that sends the equivalence class of 
$Y_u\rightarrow W_{x}\stackrel{g}{\rightarrow}X_v$ to the 
equivalence class of 
$(Y/Z_i)_{q'(u)}\rightarrow (W/Z'_j)_{q''(x)}\rightarrow(X/Z)_{q(v)}$,
with $g^{-1}(Z)=\bigcup_J Z'_j$,
and Galois coverings $W_{x}\stackrel{g}{\rightarrow}X_v$ to Galois
coverings $(W/Z'_j)_{q''(x)}\rightarrow(X/Z)_{q(v)}$.
\end{corollary}

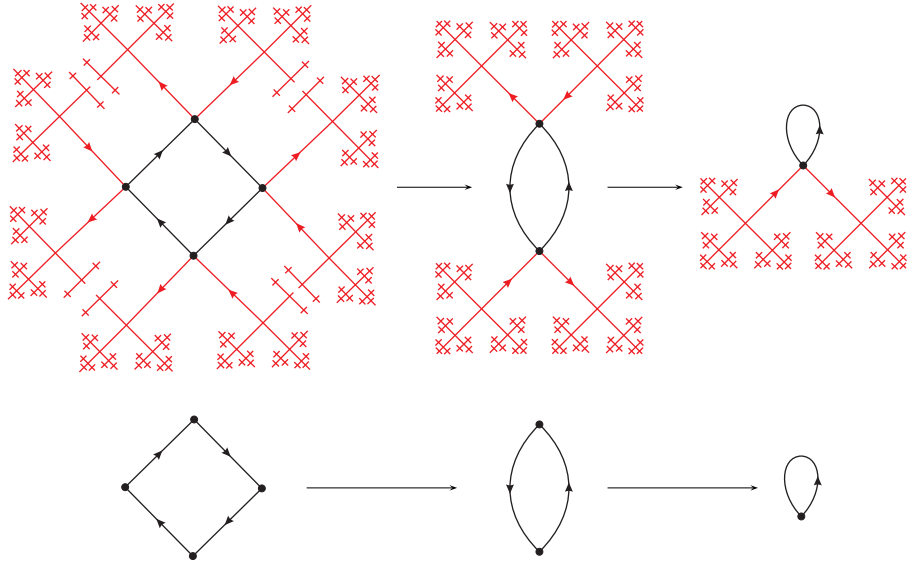
\begin{figure}[h]
  \centering
\begin{pspicture}(0,0)(12,8)
\rput(2.5,5.5){\BoxedEPSF{chapter4.fig4200a.eps scaled 500}}
\rput(7.1,5.5){\BoxedEPSF{chapter4.fig4100a.eps scaled 500}}
\rput(10.6,5.5){\BoxedEPSF{chapter4.fig4000a.eps scaled 500}}
\rput(2.5,1.5){\BoxedEPSF{chapter4.fig5000b.eps scaled 500}}
\rput(7.1,1.5){\BoxedEPSF{chapter4.fig5000.eps scaled 500}}
\rput(10.6,1.5){\BoxedEPSF{chapter4.fig5000a.eps scaled 500}}
\psline[linewidth=.1mm]{->}(5.2,5.5)(6.2,5.5)
\psline[linewidth=.1mm]{->}(8,5.5)(9,5.5)
\psline[linewidth=.1mm]{->}(4,1.5)(6,1.5)
\psline[linewidth=.1mm]{->}(8,1.5)(10,1.5)
\end{pspicture}
  \caption{lattice excision applied to Example \ref{chapter4:correspondence:example100}.}
\label{chapter4:correspondences:figure500}
\end{figure}

\begin{proof}
We show the result first in the case that
$Y\rightarrow X$ is Galois, leaving off the pointings
for clarity. By Proposition \ref{chapter4:automorphisms:result1400}  
the Galois groups $\gal(Y,X)$ and $\gal(Y/Z_i,X/Z)$ are isomorphic
and by Exercise \ref{chapter3:lattices:exercise300} this induces a lattice
isomorphism $\SSS(Y,X)\rightarrow\SSS(Y/Z_i,X/Z)$ between the subgroup lattices
of the two Galois groups. Proposition \ref{chapter4:galoiscovers:result700}
gives the induced covering $Y/Z_i\rightarrow X/Z$ is Galois.
Now apply the Galois correspondence \emph{twice\/},
\begin{equation}
  \label{eq:3}
\LL(Y,X)\rightarrow\SSS(Y,X)\rightarrow\SSS(Y/Z_i,X/Z)\rightarrow\LL(Y/Z_i,X/Z),  
\tag{\ddag}
\end{equation}
to give 
a composition of an isomorphism and two anti-isomorphisms, with net
effect the isomorphism we seek.

It is possible to find the image of $Y\rightarrow W\rightarrow X$ by
brute force. Alternatively, commuting diagrams like that in the proof
of Proposition \ref{chapter4:galoiscovers:result700} give the
images under the isomorphism
$q_*:\pi_1(X,v)\rightarrow\pi_1(X/Z,q(v))$ of the subgroups
$f_*\pi_1(Y,u)$ and $g_*\pi_1(W,x)$ to be 
$f_*'\pi_1(Y/Z_i,q'(u))$ and $g'_*\pi_1(W/Z'_j,q'(x))$. 
In particular, and by Proposition
\ref{chapter4:galoiscovers:result400}, the isomorphism
$\gal(Y,X)\rightarrow\gal(Y/Z_i,X/Z)$, and hence the isomorphism
$\SSS(Y,X)\rightarrow\SSS(Y/Z_i,X/Z)$, sends the subgroup $\gal(Y,W)$ 
to $\gal(Y/Z_i,W/Z'_j)$. Following this through
in (\ref{eq:3}) gives
$$
(Y\rightarrow W\rightarrow X)\mapsto\gal(Y,W)\mapsto\gal(Y/Z_i,W/Z_j')
\mapsto(Y/Z_i\rightarrow W/Z_j'\rightarrow X/Z).
$$
We've already observed that 
$q_*g_*\pi_1(W,x)=g'_*\pi_1(W/Z'_j,q'(x))$, hence 
Galois covers correspond to Galois covers.

Suppose now that $Y\rightarrow X$ is an arbitrary cover and let
$h:\wtl{X}\rightarrow X$ be the universal cover of
\S\ref{chapter3:actionsuniversal:universal}. Thus, we have
intermediate coverings
$$
\wtl{X}\rightarrow Y\rightarrow X,
$$
hence intermediate coverings $\wtl{X}/Z'_j\rightarrow Y/Z_i\rightarrow
X/Z$, with $\bigcup Z_j'=h^{-1}(Z)\subset\wtl{X}$. 
It is clear that $\LL(Y,X)$ is a sublattice of
$\LL(\wtl{X},X)$ and $\LL(Y/Z_i,X/Z)$ is a sublattice of
$\LL(\wtl{X}/Z'_j,X/Z)$. Applying lattice excision to the Galois cover
$\wtl{X}\rightarrow X$ gives an isomorphism
$\LL(\wtl{X},X)\rightarrow\LL(\wtl{X}/Z'_j,X/Z)$ and it is easy to
check $\LL(Y,X)$ is sent to $\LL(Y/Z_i,X/Z)$.
\qed
\end{proof}

In particular
there are homeomorphisms
\begin{align*}
(W_1\prod_X W_2)/Z_k&\rightarrow 
(W_1/Z_{1j})\prod_{X/Z}(W_2/Z_{2j}),\\
(W_1\coprod_Y W_2)/Z_k&\rightarrow 
(W_1/Z_{1j})\coprod_{Y/Z_i}(W_2/Z_{2j}).
\end{align*}
where $W_1,W_2$ are intermediate to $Y\rightarrow X$,
and the $Z_{ij},Z_{2j},Z_k$ are the components of the 
preimages of $Z$ via the various coverings.

\begin{example}
We revisit Example \ref{chapter4:correspondence:example100}, with the
covering $Y\rightarrow X$ shown in Figure
\ref{chapter4:correspondences:figure500}. A spanning tree $Z$ for $X$
is given in red and the lifts $\bigcup Z_i$ in $Y$.
We saw in Example \ref{chapter4:correspondence:example100} that the
trees were just a distraction when it came to intermediate covers:
with the Galois group $\cong\Z_4$ there was just a single equivalence
class of intermediate cover. Indeed the same effect could have been
achieved by ignoring all the red trees and focusing on the covering of
a single loop by the central square.
\end{example}

\section{Notes on Chapter \thechapter}\label{chapter4:notes}


\chapter{Generators and Relations}\label{chapter5}

\chapter{The Topological Dictionary}\label{chapter}

\chapter{Amalgams}\label{chapter8}

\chapter{The Arboreal Dictionary}\label{chapter9}

\chapter{Ends}\label{chapter10}

\chapter{Appendix}\label{appendix}

\section{Comparision with ``proper''  topology}\label{}

\subsection{The category of CW complexes}

\subsection{A functor}

\chapter{Hints for the Exercises}
\label{chapter.hints}

\backmatter




\end{document}